\documentclass{amsart}

\usepackage{amsmath,amssymb}

\begin{document}

\newcommand{\nc}{\newcommand}

\nc{\pr}{\noindent{\em Proof. }} 
\nc{\g}{\mathfrak g}
\nc{\n}{{\bf n}} 
\nc{\opn}{\overline{\n}}
\nc{\h}{\mathfrak h}
\renewcommand{\b}{\mathfrak b}
\nc{\Ug}{U(\g)} 
\nc{\Uh}{U(\h)} 
\nc{\Un}{U(\n)}
\nc{\Uopn}{U(\opn)}
\nc{\Ub}{U(\b)} 
\nc{\p}{\mathfrak p}
\renewcommand{\l}{\mathfrak l}
\nc{\z}{{\bf z}} 
\renewcommand{\h}{\mathfrak h}
\nc{\m}{{\bf m}}
\renewcommand{\k}{{\mathfrak k}}
\nc{\opk}{\overline{\k}}
\nc{\opb}{\overline{\b}}
\nc{\e}{{\varepsilon}}
\nc{\gr}{\bullet}
\nc{\ra}{\rightarrow}
\nc{\Ho}{{\rm Hom}}
\nc{\mK}{{\mathfrak{K}}}

\nc{\bg}{{\bf g}}
\nc{\W}{\mathcal{W}}
\nc{\D}{\mathcal{D}}
\renewcommand{\P}{\mathcal{P}}
\renewcommand{\O}{\mathcal{O}}
\nc{\Sw}{\mathcal{S}}
\nc{\R}{\mathbb{R}}
\nc{\Q}{\mathbb{Q}}
\nc{\N}{\mathbb{N}}
\nc{\E}{\mathbb{E}}
\renewcommand{\L}{\mathbb{L}}
\renewcommand{\P}{\mathbb{P}}
\nc{\G}{\Gamma}
\nc{\K}{\mathcal{K}}
\nc{\A}{\mathcal{A}}
\nc{\mB}{\mathcal{B}}
\nc{\B}{\mathbb{B}}
\nc{\M}{\mathcal{M}}
\nc{\Hm}{{\rm Hom}}
\nc{\Sm}{{\rm Sym}}
\nc{\dm}{{\rm deg_m}}
\nc{\Ml}{{\bf Mult}}
\nc{\bt}{\begin{theorem}}
\nc{\et}{\end{theorem}}
\nc{\be}{\begin{equation}}
\nc{\ee}{\end{equation}}
\nc{\bp}{\begin{proposition}}
\nc{\ep}{\end{proposition}}
\nc{\bd}{\begin{definition}}
\nc{\ed}{\end{definition}}
\nc{\bl}{\begin{lemma}}
\nc{\el}{\end{lemma}}
\nc{\bpr}{\begin{proof}}
\nc{\epr}{\end{proof}}
\nc{\bc}{\begin{corollary}}
\nc{\ec}{\end{corollary}}
\nc{\lab}{\label}
\nc{\br}{\begin{remark}}
\nc{\er}{\end{remark}}

\newtheorem{theorem}{Theorem}{}
\newtheorem{lemma}[theorem]{Lemma}{}
\newtheorem{corollary}[theorem]{Corollary}{}
\newtheorem{conjecture}[theorem]{Conjecture}{}
\newtheorem{proposition}[theorem]{Proposition}{}
\newtheorem{axiom}{Axiom}{}
\newtheorem{remark}[theorem]{Remark}{}
\newtheorem{example}{Example}{}
\newtheorem{exercise}{Exercise}{}
\newtheorem{definition}{Definition}{}

\renewcommand{\thetheorem}{\thesubsection.\arabic{theorem}}

\renewcommand{\thelemma}{\thesubsection.\arabic{lemma}}

\renewcommand{\theproposition}{\thesubsection.\arabic{proposition}}

\renewcommand{\thecorollary}{\thesubsection.\arabic{corollary}}

\renewcommand{\theremark}{\thesubsection.\arabic{remark}}

\renewcommand{\thedefinition}{\thesubsection.\arabic{definition}}

\title{A tree--free approach to 3D Yang--Mills Langevin dynamic. \\ Analytic estimates and the existence of a model for a regularity structure}

\author{A. Sevostyanov}

\address{Institute of Pure and Applied Mathematics,
University of Aberdeen \\ Aberdeen AB24 3UE, United Kingdom \\ e-mail: a.sevastyanov@abdn.ac.uk }

\begin{abstract}
Using the multi--index approach to regularity structures due to F. Otto et al., we construct a regularity structure and a model for it associated to the stochastic Langevin equation for the 3D Euclidean Yang--Mills functional. For the model we also obtain global stochastic and global pointwise weighted Besov type estimates which hold almost surely. The model is defined as a limit of a sequence of smooth models introduced with the help of a mollified noise. When the mollification is removed the sequence converges in a certain topology defined with the help of the stochastic estimates. To obtain these results we develop the multi--index approach for systems of equations with vector--valued white noises. This project is motivated by the problem for constructing 3D Euclidean Yang--Mills measure and by the earlier results of the author on the related problem of canonical quantization of the Yang--Mills field on the Minkowski space.
\end{abstract}

\keywords{Stochastic differential equation, Yang-Mills field, regularity structures}

\maketitle

\pagestyle{myheadings}
\markboth{A. SEVOSTYANOV}{A TREE--FREE APPROACH TO 3D YANG--MILLS LANGEVIN DYNAMIC}




\renewcommand{\theequation}{\thesubsection.\arabic{equation}}


\tableofcontents


\section{Introduction}\label{Int}

\subsection{Canonical quantization of the Yang--Mills field on the Minkowski space and 3D Euclidean Yang--Mills measure}

\setcounter{equation}{0}
\setcounter{theorem}{0}

It is well known that the Yang--Mills theory on the Minkowski space is the main ingredient of the most complete contemporary model in quantum field theory, called the Standard Model, which describes weak, strong and electromagnetic interactions of elementary particles. The problem for defining a mathematically rigorous quantum Yang--Mills theory remains open for many decades, and became central in quantum field theory at least since the Standard Model was introduced in the end of 1960s. One of the aspects of this problem, a rigorous definition of the quantized Yang--Mills Hamiltonian and the proof of the mass gap property for it, is stated as a Clay Mathematical Institute Millennium problem (see \cite{JW}).
 
In paper \cite{S} it was shown that a natural candidate for the canonically quantized  Hamiltonian of the Yang--Mills field on the Minkowski space is the Ornstein–Uhlenbeck operator associated to the yet to be defined 3D Euclidean Yang--Mills measure (see \cite{FS}, Chapter 3 for the Hamiltonian formulation for the Yang--Mills field on the Minkowski space). 

Actually, in \cite{S} it was observed that the classical Hamiltonian of the Yang--Mills field on the trivial $\mathtt K$--bundle over the Minkowski space associated to the adjoint representation of a compact Lie group $\mathtt K$ with Lie algebra $\k$ belongs to a class of Hamiltonians which can be described in the toy finite--dimensional case as follows.

Let $M$ be a Riemannian manifold with a metric $\langle \cdot,\cdot\rangle$. Consider a Hamiltonian on the symplectic manifold $T^*M$ of the form 
\begin{equation}\label{classh}
\frac12(\langle p,p\rangle +\langle v(x),v(x)\rangle ),	
\end{equation}
where $p\in T_xM\simeq T^*_xM$ is the momentum and $v={\rm grad}~ \phi$ is a potential vector field. According to the results of \cite{S}, this Hamiltonian has a family of canonical quantizations of the form 
$$
\frac 12\sum_{a=1}^n\xi_a^*(x)\xi_a(x):L_2(M,d\mu)\rightarrow L_2(M,d\mu),		
$$
where $\xi_a(x)$, $a=1,\ldots, {\rm dim}~M$ is an orthonormal basis of $T_xM$, and $\xi_a^*(x)$ is the operator formally adjoint to $\xi_a(x)$ with respect to the canonical scalar product in the space $L_2(M,d\mu)$ of square integrable functions on $M$ with respect to the measure $d\mu=\psi e^{-2\phi}dx$, where $dx$ is the Lebesgue measure on $M$ associated to the Riemannian metric, and $\psi$ is an arbitrary smooth strictly positive function on $M$.

The function $\psi$ is a quantization parameter which can be chosen arbitrarily if we do not require that the quantized components of the momentum in the Hamiltonian formalism should be self--adjoint in $L_2(M,d\mu)$. This restriction can be often lifted in quantum field theory when $M$ is infinite--dimensional and the Hamiltonian formalism momentum $p$ has no physical meaning, hence its quantized components do not need to be self--adjoint (do not confuse $p$ with the relativistic momentum which is, of course, a physical observable). 

On the other hand, the possibility of introducing an extra factor $\psi$ in the definition of the measure $d\mu$ becomes primarily important in quantum field theory because the usual Lebesgue measure does not exist on infinite--dimensional spaces, but choosing $\psi$ appropriately one can expect to find proper probabilistic counterparts of measures of the type $d\mu=\psi e^{-2\phi}dx$ which can be defined on such spaces.

The key observation of \cite{S} is that the Hamiltonian of the Yang--Mills field on the Minkowski space is of type (\ref{classh}), where $M$ is the space of gauge equivalence classes of connections $A$ on the trivial $\mathtt K$--bundle over $\R^3$ associated to the adjoint representation of $\mathtt K$, and $\phi$ is the Chern--Simons functional which we denote by $CS$.\footnote{Strictly speaking, the Chern--Simons functional defined on the space of connections is only locally gauge invariant. Under some gauge transformations it may shift by a topological constant. In the informal discussion in this introduction we do not consider this problem.} A natural metric on $M$ is induced by the natural metric $\langle \cdot,\cdot\rangle$ on the space of $\k$-valued differential one--forms on $\mathbb{R}^3$, see formula (\ref{prod}) below. 

The discussion in this introduction is rather informal, and we do not rigorously define the space of gauge equivalence classes of connections on $\mathbb{R}^3$ required for quantization. This is a subtle matter as probabilistic measures on infinite--dimensional function spaces are typically supported on spaces which consist of very singular generalized functions (see \cite{CCHS,CCHS1} for a rigorous definition of such spaces in the case of gauge equivalence classes of connections on a two--dimensional or a three--dimensional torus).   

A natural choice of the function $\psi$ for which the measure $d\mu$ is expected to exist in the case of the Yang--Mills field is such that
\be\lab{YMmden}
\psi e^{-2\phi}=\exp(-YM(A)-\frac12 m^2\langle (A-Z),(A-Z)\rangle ),	
\ee 
and hence 
$$
\psi=\exp(-YM(A)-\frac12 m^2\langle (A-Z),(A-Z)\rangle +2CS(A) ),
$$
where $YM(A)$ is the Euclidean Yang--Mills functional on $\R^3$ (see formula (\ref{YM}) below), $m>0$ is a  parameter, and we identify the affine space of connections with the space of $\k$--valued one--forms on $\R^3$ by fixing an origin $Z$, so that the expression $\langle (A-Z),(A-Z)\rangle$ is gauge invariant\footnote{In \cite{S} a slightly different choice of $\psi$ is used. The choice of $\psi$ as in (\ref{YMmden}) is more natural and symmetric, especially in the case of non--abelian $\mathtt K$.}.

As explained in \cite{S}, the necessity for introducing a non--zero parameter $m$ in the  right hand side of formula (\ref{YMmden}) for the ``density function'' of the corresponding conjectural measure $d\mu$, called the 3D Euclidean Yang--Mills measure, is already clear in the abelian case when ${\mathtt K}=U(1)$, and in this case the corresponding measure can be constructed by classical methods of functional analysis. In fact, for ${\mathtt K}=U(1)$ the expression in the exponent in the right hand side of (\ref{YMmden}) is quadratic in $A$, and the gauge action is reduced to adding differentials of functions, so that one can fix a ``gauge'' (e.g the Coulomb gauge as in \cite{S}) which provides a model for the space of gauge equivalence classes of connections, and make the 3D Euclidean Yang--Mills measure manifestly Gaussian on this model. Note that the corresponding quantization of the Yang--Mills Hamiltonian in this case, given by the Ornstein–Uhlenbeck operator associated to this Gaussian measure, is different from that used in quantum electrodynamics. 

As shown in \cite{S}, this operator exhibits one important property which is also expected from the quantized self--adjoint Yang--Mills Hamiltonian in the non--abelian case and which is primarily important for applications in quantum chromodynamics: this quantized Hamiltonian has a gap in the spectrum which separates the zero eigenvalue corresponding to the ground state from the rest of the spectrum.

As it is clear from the preceding discussion the key step in the problem of quantizing the classical Hamiltonian of the Yang--Mills field defined on the Minkowski space is now reduced to the problem for defining the 3D Euclidean Yang--Mills measure.  
In the next part of the Introduction we shall discuss this problem in more detail.


\subsection{3D Euclidean Yang--Mills measure and the Langevin equation associated to the 3D Euclidean Yang--Mills functional}\lab{stat}

\setcounter{equation}{0}
\setcounter{theorem}{0}

Firstly, we recall some basic facts from gauge theory.
Let $\mathtt K$ be a compact simple Lie group, $\k$ its Lie algebra and $\g$ the
complexification of $\k$. We denote by $(\cdot ,\cdot )$ the Killing form of $\g$.
Recall that the restriction of this form to $\k$ is non-degenerate and negative definite. For the Lie bracket $[\cdot,\cdot]$ on $\k$ it will be important to consider a family of Lie brackets $g[\cdot,\cdot]$ on the vector space $\k$, and the corresponding family of Lie brackets on $\g$, parametrized by $g\in \R$ which is called the coupling constant. We shall denote the vector space $\k$ equipped with the  Lie bracket corresponding to $g$ by $\k_g$, and the corresponding compact simple Lie group by ${\mathtt K}_g$.

We shall consider the Yang--Mills functional on the affine space of smooth
connections on the trivial ${\mathtt K}_g$-bundle over the standard Euclidean space
${\mathbb R}^{3}$ and associated to the adjoint representation of ${\mathtt K}_g$. Fixing a trivialization of this bundle and the
trivial connection as an origin in the affine space of connections we can identify this space
with the space $\Omega^1(\mathbb{R}^{3},\k_g)$ of $\k_g$-valued 1-forms on ${\mathbb R}^{3}$.
Let $A$ be such a connection.

Let $\Omega^*(\mathbb{R}^{3},\k_g)$ be the space of
$\k$-valued differential forms on $\mathbb{R}^{3}$. We
define a scalar product on this space, whenever it is finite, by
\be\lab{prod}
\langle \omega_1,\omega_2\rangle =-\int_{\mathbb{R}^{3}}(\omega_1\wedge,*\omega_2)=
-\int_{\mathbb{R}^{3}}*(\omega_1\wedge,*\omega_2)d^{3}x,~
\omega_{1,2} \in \Omega^*(\mathbb{R}^{3},\k_g)
\ee
where $*$ stands for the Hodge star operation associated to a standard
Euclidean metric on $\mathbb{R}^{3}$, and we evaluate the Killing form $(\cdot,\cdot)$ on the
values of $\omega_1$ and $*\omega_2$ and also take their exterior product. Fix orthonormal Cartesian coordinates $(x_1,x_2,x_3)$ on $\R^3$ corresponding to the Euclidean metric.

We recall that the covariant derivative
$d_A:\Omega^n(\mathbb{R}^{3},\k_g)\rightarrow \Omega^{n+1}(\mathbb{R}^{3},\k_g)$ associated to $A$
is defined by $d_A\omega =d\omega + g[A\wedge \omega]$, and the operator $d_A^*$
formally adjoint to $d_A$ with respect to scalar product (\ref{prod}) is
equal to $-*d_A*$.

Let $\mathcal{K}_g$ be the corresponding gauge group, i.e. the group
of ${\mathtt K}_g$-valued smooth maps ${\hat{g}}:\mathbb{R}^{3}\rightarrow {\mathtt K}_g$. The Lie algebra of $\mathcal{K}_g$
is isomorphic to the space $\Omega^0(\mathbb{R}^{3},\k_g)$ of $\k_g$-valued smooth functions on ${\mathbb R}^{3}$.

The gauge group $\mathcal{K}_g$ acts on the space of connections 
by
\be\lab{Gaugeact}
{\hat{g}}\circ A = -d{\hat{g}}{\hat{g}}^{-1}+{\hat{g}}A{\hat{g}}^{-1},~{\hat{g}}\in \mathcal{K}_g,~A\in \Omega^1(\mathbb{R}^{3},\k_g),
\ee
where we denote $d{\hat{g}}{\hat{g}}^{-1}={\hat{g}}^*\theta_R$, ${\hat{g}}A{\hat{g}}^{-1}=\mathrm{Ad}{\hat{g}}(A)$, and
$\theta_R$ is the right-invariant Maurer--Cartan form on ${\mathtt K}_g$. The corresponding vector fields for the action of the Lie algebra $\Omega^0(\mathbb{R}^{3},\k_g)$ of $\mathcal{K}_g$ have the form
\be\lab{Gaugeactalg}
X_A=-dX+g[X,A]=-d_AX,~X\in \Omega^0(\mathbb{R}^{3},\k_g),~A\in \Omega^1(\mathbb{R}^{3},\k_g),
\ee
where $X_A=-d_AX\in \Omega^1(\mathbb{R}^{3},\k_g)$ should be regarded as a tangent vector to the space of connections at point $A$, and the tangent space is naturally identified with $\Omega^1(\mathbb{R}^{3},\k_g)$.

Denote by $F_A$ the curvature 2-form of this connection, $F_A=dA + \frac{g}{2} [A\wedge A]$.
Here as usual we denote by $[A\wedge A]$ the
operation which takes the
exterior product of $\k$-valued 1-forms and the commutator of their values in
$\k$. The Yang--Mills functional $YM$ evaluated at $A$ is defined by
the formula
\be\lab{YM}
YM(A) = -\frac 12 \int_{{\mathbb R}^{3}} (F_A \wedge , * F_A),
\ee
where we evaluate the Killing form on the
values of $F$ and $*F$ and also take their exterior product. 

The Yang--Mills functional can be compactly rewritten in terms of scalar product (\ref{prod}),
$$
YM(A)=\frac 12 \langle F_A,F_A\rangle.
$$

According to \cite{PW} a probabilistic measure with the ``density function'' given by (\ref{YMmden}) can be sought as an invariant measure for the Markov process constructed with the help of the Langevin equation associated to the functional 
$$
-YM(A)-\frac12 m^2\langle (A-Z),(A-Z)\rangle 	
$$ 
which appears in the exponent in formula (\ref{YMmden}).

The Langevin equation associated to this functional is the gradient flow equation with an extra white noise term. Formally it has the form
\be\lab{Lang1}
\partial_tA=-d_A^*F_A-m^2(A-Z)+\xi,
\ee
where $\xi=\sum_{i=1}^3 \xi_i dx_i$, $\xi_i$, $i=1,2,3$ are independent white noises on $\R\times \R^3$ with values in the adjoint representation of ${\mathtt K}_g$, with the covariance induced by the Killing form $(\cdot,\cdot)$, $t$, called ``time'', is a coordinate on $\R$ in the product $\R\times \R^3$, and $(x_1,x_2,x_3)$ is the fixed Cartesian coordinate system on $\R^3$ in it.

Recall that we fixed a trivialization of the trivial ${\mathtt K}_g$--bundle over $\R^3$ associated to the adjoint representation and an origin in the affine space of connections. We can assume that it is done in such a way that $Z=0$ after identification of the affine space of connections with $\Omega^1(\mathbb{R}^{3},\k_g)$. Then equation (\ref{Lang1}) takes the form
\be\lab{Lang}
\partial_tA=-d_A^*F_A-m^2A+\xi.
\ee

One of the problems which arise in the course of the study of equation (\ref{Lang}) is that the partial differential operator applied to $A$ in the right hand side is not elliptic. To overcome this difficulty one adds to the right hand side of (\ref{Lang}) an extra term of the form $d_AH(A)$, where $H(A)$ is a sufficiently regular $\k_g$--valued functional of $A$, 
\be\lab{LangH}
\partial_tA=-d_A^*F_A+d_AH(A)-m^2A+\xi,
\ee
Since $d_AH(A)$ has the form of an infinitesimal gauge transformation (\ref{Gaugeactalg}) with $X=H(A)$, at least formally solutions of equations (\ref{Lang1}) and (\ref{LangH}) are gauge equivalent in law under time-dependent gauge transformations (see e.g. \cite{CCHS1,C1}).

A common choice of $H(A)$, that we are going to use, is $H(A)=-d^*A$ (Coulomb gauge) which yields the so-called DeTurck-Zwanziger term $-d_Ad^*A$ (see \cite{DeT,Zw}) in the right hand side of (\ref{LangH}). With this choice of $H(A)$ equation (\ref{LangH}) takes the form
\be\lab{LangdT}
\partial_tA=-d_A^*F_A-d_Ad^*A-m^2A+\xi,
\ee
or, in terms of the coordinates $(t,x_1,x_2,x_3)$ on $\R\times \R^3$ and of the corresponding components of $A=\sum_{i=1}^3A_i dx_i$,
\be\lab{LangHc}
\partial_tA_i=\Delta A_i-m^2A_i+g[A_j , 2\partial_j A_i - \partial_iA_j + g[A_j , A_i]]+\xi_i,~i=1,2,3,
\ee
where $\Delta=\partial_i\partial_i$ is the Laplace operator, and here and below we take sums over all possible values of repeated indices in products.

Equation (\ref{LangHc}) is highly singular due to presence of the non-linear and of the white noise terms in the right hand side. In its original formal form it has no solutions of any reasonable class.

Now we discuss this difficulty in more detail. Firstly, it is well known that the white noise on $\R^4$ belongs almost surely to a weighted Besov space. Therefore, it suffices to consider equation (\ref{LangdT}) with $\xi$ scalar components of which belong to this space, and hence we can assume that these components are less singular than generic elements of the white noise space.

More precisely, let $\D=C^\infty_0(\R^4)$ be the space of smooth compactly supported real--valued functions on $\R^4$, $\D'$ the corresponding space of distributions, i.e. the strong dual of $\D$, $\Sw\subset C^\infty(\R^4)$ the Schwartz space of smooth real--valued functions on $\R^4$ decreasing at infinity faster than any polynomial, $\Sw'$ the space of tempered distributions, i.e. its topological dual. For a compact subset $\mK\subset \R^4$, let $C^\infty(\mK)$ be the space of smooth real--valued functions on $\mK$, $\D(\mK)=C^\infty_0(\mK)$ the space of smooth real--valued functions supported on $\mK$. 

We fix a centered Gaussian measure on $\Sw'$ which equips it with the structure of a white noise space. Then for any $s<-\frac 52$ scalar components of the white noise term in equation (\ref{LangdT}) belong almost surely to the weighted Besov spaces $B^s_{\infty,\infty}(w)$, $w(x)=(1+\left\|x\right\|^2)^{\frac 52}$, $x=(t,x_1,x_2,x_3)\in \R^4$, $\left\|x\right\|:=\sqrt{t^2+x_1^2+x_2^2+x_3^2}$, with exponent $s$ associated to a non-Euclidean scaling of $\R^4=\R\times\R^3$ (see \cite{HS}, Lemma 2.23). In Section \ref{ptwise} we shall prove a more general result for recentered maps of a model for a regularity structure.

For completeness we recall the definition of weighted Besov spaces.
For $x=(t,x_1,x_2,x_3)\in \R^4$, let $|x|:=\sqrt[4]{t^2+x_1^4+x_2^4+x_3^4}$, $B_R=\{x\in \R^4:|x|\leq R\}$, $R>0$. 
For any real--valued function $\varphi(x)$, $x=(t,x_1,x_2,x_3)\in \R^4=\R\times \R^3$ we define $\varphi_x^{\lambda}(t',y_1,y_2,y_3)=\lambda^{-d}\varphi(\lambda^{-2}(t-t'),\lambda^{-1}(x_1-y_1),\lambda^{-1}(x_2-y_2),\lambda^{-1}(x_3-y_3))$, where $d=5$ is the parabolic dimension of $\R^4=\R\times \R^3$.

Let $s<0$ and fix $r\in \N:=\{0,1,2,\ldots\}$, $r>|s|$. Then the Besov space $B^s_{\infty,\infty}(w)$ can be defined as the space of tempered distributions $\zeta\in \Sw'$ satisfying the following condition
\be\lab{bes}
|\zeta(\varphi_x^{\lambda})|\lesssim w(x)\lambda^{-s}
\ee
which holds uniformly over $x\in \R^4$, $\varphi \in \mB^r:=\{\varphi\in \D(B_1):\left\|\varphi\right\|_{C^r}\leq 1\}$, $\lambda \in \R$, $0<\lambda\leq 1$. 

Here and below the symbol $C\lesssim D$ for any two real--valued quantities (i.e. real numbers or real--valued functions) means that $C\leq a D$, where $a\geq 0$. When using the symbol $C\lesssim D$ we often specify if the constant $a$ in the inequality depends on the variables on which $C$ and $D$ depend. For instance, in (\ref{bes}) the constant in the inequality only depends on $\zeta$.

The first natural idea in the course of the study of equation (\ref{LangdT}) is to replace the singular white noise term with a mollification $\xi^\rho(x):=\xi(\eta_x^\rho)$, $\rho>0$, where $\eta\in \Sw$ satisfies $\eta(y)\geq 0$, $y\in \R^4$, $\int_{\R^4}\eta(y)dy=\int_{\R^4}|\eta(y)|dy=1$, and by $\xi(\eta_x^\rho)$ we mean that the scalar components of $\xi$ are evaluated at $\eta_x^\rho$. To preserve the symmetries of equation (\ref{LangdT}) with respect to spatial reflections and coordinate permutations we shall also require that $\eta(t,x_1,x_2,x_3)$ is even in the spatial coordinates $(x_1,x_2,x_3)$ and is invariant under all permutations of these coordinates. 

However, it turns out that in the limit $\rho\to 0$, when the components of $\xi^\rho$ tend to the corresponding components of $\xi$ in $\Sw'$, solutions to the regularized equation (\ref{LangdT}) with $\xi^\rho$ instead of $\xi$ in the right hand side have no limits in any relevant topology, even for $\xi$ with components in the spaces $B^s_{\infty,\infty}(w)$, $s<-\frac 52$. 

To circumvent this obstruction a more deep modification of equation (\ref{LangdT}) is required which is achieved by adding to the right hand side extra terms, called counter terms, depending on $\rho$. A proper realization of this idea initiated by M. Hairer in pioneering paper \cite{H1} for a class of non-linear stochastic differential equations and developed in \cite{BCCH,BHZ,CH} is very complex. 

Firstly, equations of this class, called subcritical, with mollified white noises, can be lifted to equations in spaces of the so-called modelled distributions which are functions with values in other spaces called model spaces of regularity structures. The lifted equations can be solved and pushed down to equations in Besov spaces using a reconstruction operator defined in terms of smooth models for regularity structures. The smooth models depend on the mollification parameter $\rho$ and one shows that there is an appropriate choice of a smooth model such that when $\rho\to 0$ the pushed down solutions converge in a reasonable sense. However, the pushed down ``renormalized'' equations, which are the original equations with mollified noises and extra counter terms, often have no meaning in the limit $\rho\to 0$ since coefficients in the counter terms may diverge when $\rho\to 0$.

Note that the method suggested in \cite{BCCH,BHZ,CH,H1} yields a local existence theorem for the Cauchy problem, where $t$ plays the role of time, similar to Picard's theorem for the local existence and uniqueness of solutions to ordinary differential equations.

This method was applied in \cite{CCHS} to equation (\ref{LangdT}) with the spatial part defined on the three--dimensional torus (see also \cite{CCHS1} for the case of the two--dimensional torus). Note that $m=0$ in this case. The corresponding counter term obtained in \cite{CCHS,CCHS1}, that should be added to the right hand side of equation (\ref{LangdT}) with a mollified white noise, has the form $cA$, where $c$ is a constant which depends on $\rho$. 

Note that the mollified noise $\xi^\rho$ is not manifestly gauge covariant. But the gauge covariance can be recovered in the limit $\rho\to 0$, at least in law. 

The renormalized equation (\ref{LangdT}) with a mollified white noise has the form
\be\lab{LangdTrm}
\partial_tA+m^2A=-d_A^*F_A-d_Ad^*A+cA+\xi^\rho,
\ee
or, in terms of the coordinates, 
\be\lab{LangHcrm}
(\partial_t-\Delta+m^2)A_i=g[A_j , 2\partial_j A_i - \partial_iA_j + g[A_j , A_i]]+cA_i+\xi_i^\rho,~i=1,2,3.
\ee

In papers \cite{BOS,BOT,LOTT,LOT,LO,OSSW,OST,T} an alternative multi--index approach to regularity structures and to their applications for solving non-linear parabolic stochastic partial differential equations was developed. This approach has several benefits.    

Firstly, the multi--index combinatorics is much simpler than the combinatorics of colored trees used in \cite{BCCH,BHZ,CH}. Secondly, in the multi--index framework one can obtain  solutions to the corresponding stochastic differential equations using a priori estimates and the continuity method (see \cite{BOS}). As in the theory of parabolic partial differential equations this method is more likely to yield Cauchy problem solutions defined at least for any finite time interval, while only local in time solutions can be obtained in the original Hairer's framework using Picard's iterations and the fixed point theorem.    

In the case of Yang-Mills measures Hairer's approach indeed yields local in time solutions to the corresponding renormalized equation (\ref{LangdTrm}), with the spatial part defined on the two--dimensional or the three--dimensional torus, and $m=0$, as described in \cite{CCHS} and \cite{CCHS1}. In the two--dimensional case the Yang-Mills measure obtained long time ago by other means (see \cite{Lev,Sen} or \cite{Dri,GKS} for earlier works) is an invariant measure for the corresponding Markov process defined on a space of gauge equivalence classes of connections (see \cite{CS}). But the construction of global solutions to equation (\ref{LangdTrm}), and of the corresponding Yang-Mills measure in the three--dimensional case are still central open problems.

The purpose of the project initiated in this paper is to develop the multi--index approach for solving equation (\ref{LangdTrm}). In particular, in this paper we construct a regularity structure associated to equation (\ref{LangdTrm}), prove the existence of a model for it, and obtain stochastic and pointwise estimates for the model. For this purpose we develop the multi--index approach for systems of equations and vector--valued noises. 

From now on we only consider equation (\ref{LangdTrm}) associated to the Yang-Mills functional on $\R^3$. We mainly follow the multi--index approach as it is presented in \cite{BOT}, and refer the reader to \cite{BOT} for motivations of the definitions given in this paper for equation (\ref{LangdTrm}). 

According to \cite{BOT}, the renormalization parameter $c$ in equation (\ref{LangdTrm}) should be a polynomial in the coupling constant $g$ without constant term. The constant $m$ is a free parameter of the theory. As we mentioned in the first part of the Introduction, for the existence of the Yang-Mills measure in the abelian case on $\R^3$ when $\k=\R$, $g=0$, and the polynomial $c$ reduces to zero, $m^2$ must be strictly positive (see also \cite{H}, Example 1.3 and Remark 1.16 in \cite{C1}). Therefore we naturally assume that in the non--Abelian case $m>0$ is a constant independent of $\rho$, and $c$ is a polynomial in $g$ without constant term the coefficients of which depend on $\rho$.

Equation (\ref{LangHcrm}) is our main object of study. Note that the operator $L:=\partial_t-\Delta+m^2$, which appears in its left hand side, has much better properties than the heat operator commonly used in both Hairer's and the multi--index approaches. In particular, at infinity its Green function behaves like a function from the Schwartz space. 

Note that $m$ plays the role of a mass parameter in the theory, and it is expected that the condition $m>0$ ensures convergence of the corresponding Markov process on a space of gauge equivalence classes of connections, existence of an invariant Yang--Mills measure for it, and the spectral gap for its generator which should play the role of the Hamiltonian of the canonically quantized Yang--Mills theory on the Minkowski space as shown in \cite{S}. The last fact is expected to explain the mass gap phenomenon in quantum chromodynamics.

In conclusion we briefly describe the structure of the paper. More comments on the content can be found in the beginning of each section.

In Section \ref{regm}, after preliminary discussion of scaling properties of functions from Besov spaces and of some properties of the operator $L$, we introduce the model space and the smooth model for a regularity structure associated to equation (\ref{LangHcrm}). We follow the philosophy of the multi--index approach in the form described in \cite{BOT} for scalar equations. But our definition of the model space adapted for systems of equations with vector--valued noises requires vector--valued analogues of multi--indices which appear in Section \ref{mspac}. 

The definition of the recentered maps and of the canonical lift for a smooth model of the regularity structure in Section \ref{recm} is given by applying a combination of the original Hairer's method from \cite{H1} and of the multi--index approach, the combinatorics being close to the multi--index approach, and the analytic formulas resemble Hairer's definition based on the notion of regularizing kernels. The reason for this combined approach originates in the fact that we prefer the simplified multi--index combinatorics but consider the case when $\xi$ is a white noise, which is important for preserving gauge covariance, while in \cite{BOT} the case of less rough noises is considered. The last fact allows to introduce recentered maps using the heat kernel itself, not its cut--off version, called a regularizing kernel, as in \cite{H3}. 

In turn, this possibility is based on the scaling properties of noises. White noises scale (in law) with rational exponents under dilations on the underlying space $\R^4$ which a priori forces to use estimates for the regularizing kernel in H\"{o}lder spaces with integer exponents in the theory of regularity structures. But actually such estimates can only be obtained in Zygmund spaces which are different from H\"{o}lder spaces at integer exponents. 

One way to circumvent this obstruction is to manually distort the exponents in the definition of regularity structures which is done in the original Hairer's approach. In the setting of this paper the need for such distortion is already clear from the fact that the components of the white noise $\xi$ belong to the spaces $B^s_{\infty,\infty}(w)$ with $s<-\frac 52$. This method is applicable in the case of white noises, while the other way is to consider less rough noises which already scale in law with irrational exponents, as it is done in \cite{BOT}. 

Since gauge covariance, which is only preserved in the case of the white noise, is primarily important for our purposes we follow the former path. Its disadvantage is that in \cite{BOT} uniform estimates at all scales are obtained for the model, while in the original Hairer's approach estimates are local and the scaling parameter takes values in a bounded interval, the last restriction being impossible to overcome in the case of white noises. Since we stick to Hairer's definition of the recentered maps in its analytic part, we still use a cut--off version of the kernel of the Green function of the operator $L$ in this definition. We already mentioned that this Green function behaves like a Schwartz function at infinity unlike the heat kernel itself which is not even integrable. Nevertheless, the Green function of $L$ is not compactly supported which makes it unsuitable for Hairer's definition of recentered maps.  

The structure group maps are defined in Section \ref{StrG} with the help of automorphisms relating the recentered maps and the canonical lift which appear in Section \ref{Frec}. The definition of these automorphisms is simple and explicit. Note that the multi--index description of the regularity structure is dual in a sense to the one in Hairer's approach, and the definition of automorphisms relating the recentered maps and the canonical lift is obtained by ``dualizing'' a similar definition in terms of regularity structures where it becomes inductive and rather inexplicit (see \cite{H1}, \S 8.2 or \cite{C2}, Section 2). Formulas of the same type also appeared in \cite{LOT,LO}.

In Section \ref{Moddeq} we also show, following Hairer's approach, how to construct a natural lift of equation (\ref{LangHcrm}) to an equation in a space of modelled distributions.  

Section \ref{Anal} contains analytic tools required for obtaining stochastic estimates for models. After introducing some function spaces and recalling some classical Sobolev type inequalities we state Schauder estimates for regularizing kernels and the reconstruction theorem for germs in Section \ref{Shaud}. We need stochastic pointed and local versions of the Schauder estimates which are close to those obtained in \cite{BCZ} in the deterministic case (see also \cite{FH}, Chapter 14). A stochastic version of Hairer's reconstruction theorem for germs is stated in a form close to Theorem 3.2 in \cite{BL} (see also \cite{CZ}). The stochastic versions of these two statements are required in the multi--index approach to stochastic estimates for models as suggested in \cite{BOT,LOTT} where their global stochastic versions are used. Section \ref{Diffrec} also contains necessary results on Malliavin derivatives of the recentered maps used in the proof of the stochastic estimates in the multi--index approach.  

Section \ref{Stohe} is central in this paper. After recalling in Section \ref{BPHZsubs} the so--called generalized BPHZ condition used to fix the renormalization parameter $c(\rho)$ we confirm the translational invariance in law for the recentered and for the structure group maps in Section \ref{Trinvsubs} and proceed with the statement of the stochastic estimates for them in Section \ref{Mainstat}. 

The proof of these estimates presented in Sections \ref{Baseind}, \ref{Indstep}, \ref{indstGGS}, \ref{indstdP<2}, \ref{indstepPb} and \ref{indstGGSpp} is quite lengthy. We proceed by induction over the indices from the index set of the regularity structure following the idea of \cite{BOT,LOTT}. However, in contrast to \cite{BOT,LOTT}, our proof is arranged as a linear chain of propositions. The other difference is that, since we use local and pointed Schauder estimates for a regularizing kernel on a bounded scale, which naturally appears in a definition of the relevant H\"{o}lder spaces, the range of the scaling parameter shrinks after each iteration of the induction step, and to run the induction step again this parameter must be greater than one. As a consequence, starting from a large enough range for the rescaling parameter in the base case, we can obtain stochastic estimates for recentered maps and for structure group maps labeled by indices from the index set of the regularity structure uniformly bounded by a given constant which can be chosen arbitrarily large. The BPHZ conditions in our approach should also be fixed at each iteration of the induction step separately for recentered maps labeled by the corresponding indices from the index set. 

In Section \ref{conv} we show that the smooth models converge in a certain space when the mollification of the white noise is removed, i.e. when $\rho\to 0$. The proof is based, as in \cite{HS,T}, on a Cauchy property for smooth models which is established in the same manner as the stochastic estimates.

Finally, in Section \ref{ptwise} using Kolmogorov continuity arguments we obtain weighted Besov type global pointwise estimates for models which hold almost surely. This is done in Theorem \ref{GestMNprop} for the structure group maps of the model and in Theorem \ref{PestNprop} for the recentered maps.

\vskip 0.5cm
\noindent
{\bf Acknowledgements.} The author is grateful to Thierry L\'{e}vy and to Ilya Chevyrev for stimulating discussions at the initial stage of this project.



\section{A regularity structure and a model for the Langevin equation with mollified noise}\lab{regm}

In this section we introduce a regularity structure and a smooth model for it associated to equation (\ref{LangHcrm}). We start by defining the underlying combinatorics relying on the notion of multi--indices which appear in Section \ref{scaleh}. In the definition of the smooth model for the regularity structure we explicitly use the kernel of the heat operator with an extra mass term the properties of which are recalled in Section \ref{heatsubs}. Then we proceed with an algebraic description of the model space in Section \ref{mspac}. Using this model space recentered maps of a smooth model for the corresponding regularity structure and a canonical lift for it associated to equation (\ref{LangHcrm}) are defined in Section \ref{recm}. In the algebraic part of their definition in Lemma \ref{lF} we follow the strategy of \cite{BOT}. But, as we mentioned in the introduction, analytic formulas (\ref{defP}) for the recentered maps have some similarity with the original Hairer's definition given in \cite{H} as well as provisional local analytic estimates for them in the end of Section \ref{recm}. The structure group maps of the smooth model are defined with the help of model space automorphisms relating the recentered maps and the canonical lift as in the original Hairer's approach. However, we introduce these automorphisms in Section \ref{Frec} using the multi--index combinatorics in terms of which their definitions are much more simple and explicit. The definition of the structure group automorphisms resembles Hairer's definition and their properties derived in Lemma \ref{Gprop} are similar to those obtained in \cite{BOT} in the case of a single scalar equation. In the end of this section we explicitly relate the constructed smooth model to equation (\ref{LangHcrm}). To this end we define in Section \ref{Moddeq} a lift of equation (\ref{LangHcrm}) to an equation for modelled distributions and show that after applying the so-called reconstruction operator this equation is reduced to equation (\ref{LangHcrm}).  


\subsection{Scale invariance and homogeneity}\lab{scaleh}

\setcounter{equation}{0}
\setcounter{theorem}{0}

Recall that our final goal is to obtain the $\rho\to 0$ limit of solutions to equation (\ref{LangdTrm}) with $\xi\in B^s_{\infty,\infty}(w)$, $s<-\frac 52$. The structure of the left hand side of equation (\ref{LangHcrm}) implies that such solutions should belong at least to the local version of the space $B^\alpha_{\infty,\infty}$, $\alpha:=s+2$.
Therefore, following \cite{BOT}, in view of definition (\ref{bes}) in terms of a scaling, it is natural to consider scaling properties of equation (\ref{LangdTrm}) with respect to the following transformation
\be\lab{scal}
\xi(x)\mapsto \lambda^{-(\alpha-2)}\xi(R_\lambda x),\rho\mapsto \lambda^{-1}\rho, A\mapsto \lambda^{-\alpha}A(R_\lambda x),~\lambda\in \R,~\lambda> 0,
\ee
where for $x=(t,x_1,x_2,x_3)$,  we define $R_\lambda x=(\lambda^2t,\lambda x_1,\lambda x_2,\lambda x_3)$. From now on we also assume that $\alpha=-\frac 12-\varepsilon$, where $0<\varepsilon<\frac 12$, $\varepsilon \not\in \mathbb{Q}$ can be chosen arbitrarily small.

Denoting $x_0:=t$ and introducing multi--indices $\n=(n_0,n_1,n_2,n_3)\in \N^4$ corresponding to the coordinates $(x_0,x_1,x_2,x_3)$ we arrive at the following induced transformation for partial derivatives 
\be\lab{scaln}
\partial^{\n}A_i(x)\mapsto  \lambda^{|{\n}|-\alpha}\partial^{\n}A_i(R_\lambda x),~i=1,2,3,
\ee
where for ${\n}=(n_0,n_1,n_2,n_3)$ one defines $\partial^{\n}=\prod_{i=0}^3\partial_i^{n_i}$, and
\be\lab{nhom}
|{\n}|=2n_0+n_1+n_2+n_3.
\ee

Equation (\ref{LangdTrm}) becomes invariant under transformations (\ref{scal}), (\ref{scaln}) if one postulates that these transformations act on $g$, $m^2$ and $c$ as follows
\be\lab{scalc}
g\mapsto \lambda^{\alpha+1} g,~m^2 \mapsto \lambda^2 m^2,~c\mapsto \lambda^2 c.
\ee

The relevant combinatorics can be presented in the most convenient way using the set of multi--indices $\Ml$ over the set $\{{\bg}\}\coprod \N^4$, where $\bg$ is a formal variable. The multi--indices are functions $\beta:\{{\bg}\}\coprod \N^4\to \N$ taking non--zero values only at finitely--many elements of $\{{\bg}\}\coprod \N^4$. Let ${\bf 0}=(0,0,0,0)\in \N^4$, and denote by $\delta_{\bg}$ (resp. $\delta_{\n}$) the multi--index equal to 1 at ${\bg}$ (resp. ${\n}\in \N^4$) and zero otherwise.

For any multi--index $\beta$ we introduce its homogeneity  $|\beta|$ by
\be\lab{hom}
|\beta|=(\alpha+1)\beta({\bg})+\sum_{{\n}\in \N^4}\beta({\n})(|{\n}|-\alpha)+\alpha.
\ee

Then writing 
\be\lab{c'}
c=\sum_{\tiny \begin{array}{c}\beta=k\delta_{\bg}+\delta_{\bf 0} \\ k\in \N\setminus \{0\}\end{array}}c_\beta g^{\beta({\bg})}
\ee
for some $c_\beta=c_\beta(\rho)\in \R$ we arrive at the following transformation law for $c_\beta$ induced by (\ref{scalc})
\be
c_\beta\mapsto \lambda^{2-|\beta|}c_\beta.
\ee

This transformation rule and (\ref{scal}) imply that if we fix $\lambda^{-1}\rho=1$ and consider the limit $\rho\to 0$ removing the regularization of $\xi$ then one can expect $c_\beta(\rho) \to \infty$ for $|\beta|<2$ and $c_\beta(\rho) \to 0$ for $|\beta|>2$ (note that for non-zero $c_\beta$ one has $|\beta|\not\in \mathbb{Q}$). Therefore it is natural to assume that $c_\beta=0$ for $|\beta|\geq 2$. In Sections \ref{Stohe} and \ref{conv} we shall prove the main results of this paper under this assumption.

As by (\ref{hom}) $|k\delta_{\bg}+\delta_{\bf 0}|=k(\alpha+1)$, one has from (\ref{c'})
$$
c=\sum_{0<k<\frac{2}{1+\alpha}}c_k g^k,~c_k=c_k(\rho)\in \R,
$$
and under the assumption $\alpha=-\frac 12-\varepsilon$, where $0<\varepsilon\leq\frac{1}{10}$, $\varepsilon \not\in \mathbb{Q}$, we infer 
$$
\frac{2}{1+\alpha}=\frac{4}{1-2\e}\leq 5,
$$
and hence
$$
c=\sum_{k=1}^4c_k g^k.
$$

From now on we shall always assume that $\varepsilon \not\in \mathbb{Q}$ and satisfies 
\be\lab{erest}
0<\varepsilon\leq\frac 1{100}.
\ee


\subsection{The heat operator with a mass term}\lab{heatsubs}

\setcounter{equation}{0}
\setcounter{theorem}{0}

To introduce the smooth model for a regularity structure associated to equation (\ref{LangHcrm}) we shall need some properties of the operator $L=\partial_t-\Delta+m^2$ which appears in its left hand side. These properties are well known. 

For $x=(t,x_1,x_2,x_3)\in \R^4$, denote $\bar{x}=(x_1,x_2,x_3)\in \R^3$. Let $$\|\bar{x}\|=\sqrt{x_1^2+x_2^2+x_3^2}$$ be the standard norm on $\R^3$. The Green function $K_m$ of the operator $L$ can be obtained using the Fourier transform. It has the form
\be\lab{GrL}
K_m(x)=\frac{\theta(t)}{(4\pi t)^{\frac 32}}e^{-\frac{\|\bar{x}\|^2}{4t}-m^2t}=K_H(x)e^{-m^2t},
\ee
where 
$$
\theta(t)=\left\{\begin{array}{ll} 1 & {\rm if}~t>0 \\ 0 & {\rm if}~t\leq 0 \end{array}, \right.
$$
and $K_H(x)=\frac{\theta(t)}{(4\pi t)^{\frac 32}}e^{-\frac{\|\bar{x}\|^2}{4t}}$ is the heat kernel. 

For $x=(x_0,x_1,x_2,x_3)\in \R^4$ let $\|x\|:=\sqrt{x_0^2+x_1^2+x_2^2+x_3^2}$ be the standard Euclidean norm on $\R^4$. We also denote by $|x|$ the parabolic distance from $x\in \R^4$ to the origin defined by $|x|=\sqrt[4]{x_0^2+x_1^4+x_2^4+x_3^4}$. The parabolic metric $d(x,y):=|x-y|$, $x,y\in \R^4$ is a genuine metric. 

For any $R>0$, $x\in \R^4$ we denote by $B_R(x):=\{y\in \R^4:|x-y|\leq R\}$ the parabolic ball of radius $R$ centered at $x$. We shall also write for brevity $B_R:=B_R(0)$.

Fix a function $\varsigma \in \D$ supported in $B_1$ and invariant under each spatial reflection $R_i:\R\times \R^3\to \R\times \R^3$  in the plane defined by the equation $x_i=0$, $i=1,2,3$ and under all permutations of the coordinates $x_i$, $i=1,2,3$. 

Let $K(x):=K_m(x)\varsigma(x)$, $x\in \R^4\setminus \{0\}$. Then $K$ is a 2-regularizing kernel according to Definition 14.9 in \cite{FH}.

Denote the corresponding integral operator by $\K$,
\be\lab{Kdef}
\K f(x)=\int_{\R^4}K(x-y)f(y)dy=\int_{\R^4}K(y)f(x-y)dy.
\ee

Let also $\K_m$ be the integral operator with the kernel $K_m$,
$$
\K_m f(x)=\int_{\R^4}K_m(x-y)f(y)dy=\int_{\R^4}K_m(y)f(x-y)dy.
$$

Many functions that we shall introduce in the next section will belong to the space $\O$ of smooth functions on $\R^4$ of polynomial growth which is defined by
$$
\O=\{f\in C^\infty(\R^4):~\forall~{\n}\in\N^4~\exists~N({\n}),C({\n})\geq 0:~|\partial^{\n}f(x)|\leq C({\n})(1+\|x\|^2)^{N({\n})}\}.
$$

Tempered distributions have the following useful property (see Theorem IX.4(a) in \cite{RS}).
\bl\lab{SS'conv}
For any $\xi\in \Sw'$ and $\eta\in \Sw$ one has $(\xi*\eta)(x):=\xi(\eta^1_x)\in \O$. In particular, for any $\rho>0$, one has $\xi^\rho\in \O$.
\el

Now we show that the operator $\K$ acts in the space $\O$. 

For any ${\n}=(n_0,n_1,n_2,n_3)\in \N^4$ and $x=(x_0,x_1,x_2,x_3)\in \R^4$ denote $x^{{\n}}=x_0^{n_0}x_1^{n_1}x_2^{n_2}x_3^{n_3}$. We call $|{\n}|$ the parabolic degree of the monomial $x^{\n}$. For $n\in \N$ let $\P_{n}\subset \O$ be the subspace generated by monomials of parabolic degrees up to $n$. We call $\P_{n}$ the space of polynomials of parabolic degree at most $n$. For convenience we also define $\P_{n}=\{0\}$ for $n\in {\mathbb{Z}}$, $n<0$.

We proceed by observing that, since $K$ is compactly supported and absolutely integrable,  for any $n\in \N$ the space $\P_{n}$ and the space of polynomials of degree at most $n$ are invariant under the action of the operator $\K$. We state this obvious fact as a lemma for future references.
\bl\lab{KPol}
(i) Let $P$ be a polynomial on $\R^4$ of parabolic degree at most $n$. Then $\K P$ is well defined and is a polynomial of parabolic degree at most $n$.

(ii) Let $P$ be a polynomial on $\R^4$ of degree at most $n$. Then $\K P$ is well defined and is a polynomial of degree at most $n$.
\el

As a corollary of this lemma we immediately obtain the following result.
\bp\lab{KO}
$\K$ is well defined as an operator $\K:\O\to \O$.
\ep

\bpr
First recall that $K$ is compactly supported and absolutely integrable. Hence for any $f\in C^\infty(\R^4)$ one has $\K f\in C^\infty(\R^4)$, and for any ${\n}\in \N^4$
$$
\partial^{\n}(\K f)(x)=(\K \partial^{\n}f)(x).
$$

In particular, for any $f\in \O$, one has $\K f\in C^\infty(\R^4)$, and for any ${\n}\in \N^4$ we have using Lemma \ref{KPol}(ii)
$$
|\partial^{\n}(\K f)(x)|=|(\K \partial^{\n}f)(x)|\leq \int_{\R^4}|K(y)||\partial^{\n}f(x-y)|\leq
$$
$$
\leq \int_{\R^4}|K(y)|C({\n})(1+\|x-y\|^2)^{N({\n})}\leq C'({\n})(1+\|x\|^2)^{N({\n})}
$$
for some $C({\n}),C'({\n}),N({\n})\geq 0$. Therefore by the definition $\K f\in \O$.

\epr


\subsection{The index set and the model space for a regularity structure}\lab{mspac}

\setcounter{equation}{0}
\setcounter{theorem}{0}

In this section we introduce the model space for a regularity structure for equation (\ref{LangHcrm}). The motivation for this definition will be clarified in the next two sections.

Fix a total normal order on the set of multi--indices $\N^4$ (e.g. one can use the lexicographic order), and for ${\n}\in \N^4$ let $W_{\n}:=\k\otimes \R^3$, $\Sm_{\n}$ the symmetric algebra of $W_{\n}$, and $\Sm_{\n}^*$ the symmetric algebra of $W_{\n}^*$. The tensor products of algebras 
$$
\Sm_{{\n}_0}\otimes \Sm_{{\n}_1}\otimes\ldots\otimes \Sm_{{\n}_p},
$$
$$
\Sm_{{\n}_0}^*\otimes \Sm_{{\n}_1}^*\otimes\ldots\otimes \Sm_{{\n}_p}^*,
$$
where ${\n}_0<{\n}_1<\ldots <{\n}_p<\ldots $ is the fixed normal order on $\N^4$, form direct systems of algebras with respect to natural inclusions. Denote their direct limits by $\Sm$ and $\Sm^*$, respectively. $\Sm$ and $\Sm^*$ are associative algebras elements of which can be regarded as linear combinations of formal infinite tensor products of the form
$$
v_0\otimes v_1 \otimes \ldots \otimes v_p\otimes 1 \otimes 1 \otimes \ldots ,
$$
$$
v_0^*\otimes v_1^* \otimes \ldots \otimes v_p^* \otimes 1 \otimes 1 \otimes \ldots,
$$
respectively, where $v_i\in \Sm_{{\n}_i}$, $v_i^*\in \Sm_{{\n}_i}^*$, $i=0,\ldots, p$, $p\in \N$, and only finitely many terms in the tensor products are different from 1. 

For any multi--index $\beta=l\delta_\bg+\sum_{i=0}^p k_i\delta_{{\n}_i}$ denote 
\be\lab{defWb}
W_\beta:=W_{{\n}_0}^{\odot k_0}\otimes W_{{\n}_1}^{\odot k_1}\otimes \ldots \otimes W_{{\n}_p}^{\odot k_p}\subset \Sm,
\ee 
\be\lab{defWb*}
W_\beta^*:={W_{{\n}_0}^*}^{\odot k_0}\otimes {W_{{\n}_1}^*}^{\odot k_1}\otimes \ldots \otimes {W_{{\n}_p}^*}^{\odot k_p}\subset \Sm^*,
\ee
where $\odot$ stands for the symmetrized tensor power. Note that according to this definition for $l\in \N$ one has $W_{l\delta_\bg}=\R$, and
$W_\beta^*$ is naturally dual to $W_\beta$. This gives rise to a natural pairing
\be\lab{par}
\Sm^*\otimes \Sm\to \R,
\ee
and for any multi--indices $\beta$ and $\gamma$ the multiplications in $\Sm$ and in $\Sm^*$ give rise to linear maps
\be\lab{Wprod}
W_\beta \otimes W_\gamma \to W_{\beta+\gamma}, 
\ee
\be\lab{W*prod}
W_\beta^* \otimes W_\gamma^* \to W_{\beta+\gamma}^*
\ee
which are isomorphisms by a dimension count.

Let $V=\k\otimes \R^3$, and
$$
V_\beta=V_\beta^{-}:=W_\beta^*\otimes V\simeq \Hm(W_\beta, V).
$$
We have  natural evaluation maps
$$
V_\beta \otimes W_\beta \to V,~V_\beta^{-} \otimes W_\beta \to V.
$$ 

For any multi--index $\beta$ we introduce $[\beta] \in \mathbb{Z}$ by 
\be\lab{pop}
[\beta]=\beta(\bg)-\sum_{{\n}\in \N^4}\beta({\n}).
\ee
Let $\M'$ be the subset of multi--indices $\beta$ which either satisfy the condition $[\beta]\geq 0$ or have one of the following two forms:
\be\lab{b2}
\beta=\delta_\bg+\delta_{\n}+\delta_{{\n}'},~{\n},{\n}'\in \N^4,
\ee
\be\lab{b3}
\beta=2\delta_\bg+\delta_{\n}+\delta_{{\n}'}+\delta_{{\n}''},~{\n},{\n}',{\n}''\in \N^4.
\ee 

Following the terminology of \cite{BOT} we call multi--indices of the form $\delta_{\n}$, ${\n}\in \N^4$ purely polynomial. Denote by $\M_{pp}:=\{\delta_{\n}:{\n}\in \N^4\}$ the subset of purely polynomial multi--indices. Let $\M_{\geq 0}:=\{\beta\in \Ml:[\beta]\geq 0\}\subset \M'$, and $\M:=\M'\cup \M_{pp}$. As in \cite{BOT}, multi--indices from the set $\M$ are called populated.

In the following lemma we obtain the local finiteness property for the index set corresponding to $\M$ in the sense of Hairer.
\bl\lab{locf}
The subset $H$ of $\R$ of the homogeneities $|\beta|$ of all possible multi--indices $\beta$ and the set $I:=\{|\beta|:\beta \in \M\}\subset H\subset \R$ are subsets of the set $(\alpha+1)\N+\N-\alpha\N+\alpha$. In particular, $H$ and $I$ are locally finite, i.e. any bounded open subinterval of $\R$ contains finite numbers of elements of both $H$ and of $I$. $\alpha$ is the minimal element for both $H$ and $I$. In particular, for any $r\in \R$ the set $\{\beta\in \Ml:|\beta|\leq r\}$ is finite. Moreover, $\beta=0$ is the only multi--index for which $|\beta|=\alpha$.
\el

\bpr
Formula (\ref{hom}) implies that the set $H$ is a subset of the set $(\alpha+1)\N+\N-\alpha\N+\alpha$. Since $\alpha+1$ and $-\alpha$ are strictly positive, we deduce from formula (\ref{hom}) that $\alpha$ is the minimal element of $H$ and of $I$, and $\beta=0$ is the only multi--index for which $|\beta|=\alpha$. 

For the same reason for any $r\in \R$ the set $((\alpha+1)\N+\N-\alpha\N+\alpha)\cap (-\infty,r)$ is finite. Since $H\subset (\alpha+1)\N+\N-\alpha\N+\alpha$, we infer that $H$ is locally finite, and hence $I\subset H$ is locally finite as well. 

\epr

Let 
$$
T=\bigoplus_{\beta\in \M}V_\beta^*\oplus \bigoplus_{\beta\in \M}{V_\beta^{-}}^*.
$$
$T$ is called the model space for a regularity structure for equation (\ref{LangHcrm}). This choice of $T$ is not unique. Another possible choice would be 
$$
T':=\bigoplus_{\beta\in \M}W_\beta \oplus \bigoplus_{\beta\in \M}W_\beta
$$ 
(compare with a similar discussion in \cite{BCCH,BHZ}), so that $T\simeq T'\otimes V^*$. Our choice is convenient for the description of the recentered maps of a smooth model in the next section. However, the structure group maps of the smooth model defined in Section \ref{StrG} are induced by endomorphisms of $T'$.

In the notation of \cite{H3} the index set for $T$ is $I\cup (I-2)$, so $T=\bigoplus_{r\in  I\cup (I-2)}T_r$, and for $r\in I\cup (I-2)$ the corresponding graded component of $T$ has the form 
$$
T_r=\bigoplus_{\beta\in \M, |\beta|=r}V_\beta^*\oplus \bigoplus_{\beta\in \M,|\beta|-2=r}{V_\beta^{-}}^*,
$$
and
$$
\bigoplus_{{\n}\in \N^4}V_{\delta_{\n}}^* 
$$
is the polynomial sector of $T$.

Since the set $I$ is locally finite and bounded from below with the minimal element $\alpha$, the index set $I\cup (I-2)$ is locally finite and bounded from below with the minimal element $s=\alpha-2$. Hence $I\cup (I-2)$ satisfies the conditions in the definition of index sets in \cite{H3}. Moreover, $V_{\delta_{\n}}^*\simeq \Hm(V,V)$ for all ${\n}\in \N^4$ which agrees with the fact that for systems of equations with vector--valued noises the polynomial sector of the regularity structure should be extended accordingly.


\subsection{Recentered maps of the smooth model and the canonical lift in the case of the mollified noise}\lab{recm}

\setcounter{equation}{0}
\setcounter{theorem}{0}

Following the strategy of \cite{BOT} we proceed with defining a smooth model for the regularity structure $T$. Firstly, we introduce an equation for defining recentered maps of the model using induction over the homogeneity. For this purpose we rewrite equation (\ref{LangdTrm}) in a convenient form.

For smooth $V$--valued functions $U,U',U''\in C^\infty(V):=C^\infty(\R^4)\otimes V$ define $\A(U,U')$, $\mB(U,U',U'')\in C^\infty(V)$ by
\be\lab{AB}
\A(U,U')=A^{ijk}_l[U_j,\partial_i U_k']e^l,~\mB(U,U',U'')=B_l^{ijk}[U_j[U_k',U_i'']]e^l,
\ee
where 
$$
A^{ijk}_l=2\delta_{ij}\delta_{kl}-\delta_{il}\delta_{jk},~B_l^{ijk}=\delta_{il}\delta_{jk},
$$
$\delta_{ij}$ is the Kronecker delta, $e^l$, $l=1,2,3$ is the orthonormal basis of the Cartesian coordinate system in $\R^3$ fixed in Section \ref{stat}, $U_i, U_i',U_i''\in C^\infty(\R^4)\otimes \k$, $i=1,2,3$ are the components of $U,U',U''\in C^\infty(V)=C^\infty(\R^4)\otimes \k\otimes \R^3$ with respect to this basis, so $U=(U_1,U_2,U_3)$, etc., the commutators of the components $U_i, U_i',U_i''\in C^\infty(\R^4)\otimes \k$, $i=1,2,3$ in (\ref{AB}) are taken in $\k$, and the usual product of functions is taken in $C^\infty(\R^4)$.
The partial derivatives $\partial_i$, $i=1,2,3$ in (\ref{AB}) act on the factor $C^\infty(\R^4)$ in the tensor product $C^\infty(\R^4)\otimes \k$. 

Note that $A\in C^\infty(V)$ in (\ref{LangdTrm}), and the definitions of $\A$ and $\mB$ imply that equation (\ref{LangdTrm}) can be rewritten in the form
\be\lab{LangdTrmAB}
LA=g\A(A,A) + g^2 \mB(A,A,A)+cA+\xi^\rho,
\ee
where the operator $L$ acts on the factor $C^\infty(\R^4)$ in the tensor product $C^\infty(V)=C^\infty(\R^4)\otimes V$.

Let $\O(V):=\O\otimes V$, and for any multi--index $\beta$, let $\O_\beta:=\O\otimes V_\beta\simeq \O\otimes W_\beta^*\otimes V\simeq \O(V)\otimes W_\beta^*$,
$\O_\beta^-:=\O\otimes V_\beta^-\simeq \O\otimes W_\beta^*\otimes V\simeq \O(V)\otimes W_\beta^*$. By this definition $\O_\beta=\O_\beta^-$, but it will be convenient to use different symbols for the same space in different contexts.

Since $\O\subset C^\infty(\R^4)$ one has $\O(V)\subset C^\infty(V)$, and for $U,U',U''\in \O(V)$ one has $\A(U,U')$, $\mB(U,U',U'')\in \O(V)$ as $\O$ is closed under differentiation and multiplication.

Now assume that for all purely polynomial multi--indices $\delta_{\n}$, ${\n}\in \N^4$ elements
\be
\Phi_{\delta_{\n}}\in \P_{\delta_{\n}}:=\P_{|{\n}|}\otimes V_{\delta_{\n}}\subset \O_{\delta_{\n}}
\ee
are fixed.  

Let
\be\lab{Ph0}
\Phi_{0}^-:=\xi^\rho\in \O(V)\simeq \O_0^-,~ 
\Phi_{0}:=\K\Phi_{0}^-=\K\xi^\rho\in \O_0^-=\O_0.
\ee 
Here and below the operator $\K$ acts on the $\O$--factor in any tensor product $\O\otimes V_\beta^-= \O_\beta^-$.

Suppose that for a multi--index $\beta\neq 0$ elements $\Phi_{\gamma}\in \O_\gamma$ for $|\gamma|<|\beta|$ are given. Then introduce $\Phi_{\beta}^-\in \O_\beta^-$ by
\be\lab{defP-}
\Phi_{\beta}^-=\sum_{\tiny \begin{array}{c}\beta_1,\beta_2\in \Ml: \\ \beta_1+\beta_2+\delta_\bg=\beta \end{array}}\A(\Phi_{\beta_1},\Phi_{\beta_2})+
\ee
$$
+\sum_{\tiny \begin{array}{c}\beta_1,\beta_2,\beta_3\in \Ml: \\ \beta_1+\beta_2+\beta_3+2\delta_\bg=\beta\end{array}}\mB(\Phi_{\beta_1},\Phi_{\beta_2},\Phi_{\beta_3})+\sum_{\tiny \begin{array}{c} k\in  \{1,2,3,4\},\beta_1\in \Ml: \\ k\delta_\bg+\beta_1=\beta\end{array}}c_k\Phi_{\beta_1},
$$
where we assume that $\Phi_{\beta}^-=0$ if all sums in the right hand side are void, $\A$ and $\mB$ are evaluated on the $\O(V)$--factors of the elements $\Phi_{\beta_i}\in \O_{\beta_i}\simeq \O(V)\otimes W_{\beta_i}^*$, $i=1,2,3$, and product (\ref{W*prod}) is applied to the $W_{\beta_i}^*$--factors of these elements. 

Correctness of this inductive definition is confirmed in the following lemma along with some other important properties of $\Phi_{\beta}^-$. 
\bl\lab{lF}
Assume that for a multi--index $\beta\neq 0$ elements $\Phi_{\gamma}\in \O_\gamma$ for $|\gamma|<|\beta|$ are given. Then the following statements are true.

(i) Formula (\ref{defP-}) defines an element $\Phi_{\beta}^-\in \O_\beta^-$, i.e. all sums in the right hand side of (\ref{defP-}) are finite, and only $\Phi_{\gamma}$ with $|\gamma|<|\beta|$ appear in the right hand side of (\ref{defP-}). Moreover, $\Phi_{\beta}^-=0$ if $\beta\in \M_{pp}$.

(ii) Assume that $\Phi_{\gamma}=0$ for $\gamma\not \in \M_{\geq 0}\cup \M_{pp}$ with $|\gamma|<|\beta|$. Then $\Phi_{\beta}^-=0$ if $\beta \not\in \M'$.

(iii) Under the assumption in part (ii), if $[\beta]<0$ then $\Phi_{\beta}^-=0$ if $\beta\not\in \M'\setminus \M_{\geq 0}$, and if $\beta\in \M'\setminus \M_{\geq 0}$ then $|\beta|\in \N$ and $\Phi_{\beta}^-\in \P_{|\beta|-2}\otimes V_{\beta}^-\subset \O_\beta^-$.

(iv) Assume that $\Phi_{\gamma}=0$ for $\gamma$ with $[\gamma]<-1$ and $|\gamma|<|\beta|$. Then $\Phi_{\beta}^-=0$ if $[\beta]<-1$.
\el 

\bpr
(i) All sums in the right hand side of (\ref{defP-}) are finite since the set of multi--indices is the free commutative monoid generated by the elements $\delta_\bg$ and $\delta_{\n}$, ${\n}\in \N^4$, so every its element is represented as a sum of other elements in finitely many ways. 

By (\ref{hom}) for $\beta=\beta_1+\beta_2+\delta_\bg$ (resp. $\beta=\beta_1+\beta_2+\beta_3+2\delta_\bg$, and $\beta=k\delta_\bg+\beta_1$, $k\in \{1,2,3,4\}$) one has
$|\beta|=|\beta_1|+|\beta_2|+1$ (resp. $|\beta|=|\beta_1|+|\beta_2|+|\beta_3|+2$, and $|\beta|=k(\alpha+1)+|\beta_1|$, $k\in \{1,2,3,4\}$, $\alpha+1>0$). The conditions $|\beta_i|\geq \alpha=-\frac 12 -\e$ for all $i$ and $\alpha+1=\frac 12 -\e$ ensure that in all three cases only $\Phi_{\gamma}$ with $|\gamma|<|\beta|$ appear in the right hand side of (\ref{defP-}), and $\Phi_{\beta}^-=0$ if $\beta\in \M_{pp}$ as $\beta=\beta_1+\beta_2+\delta_\bg$, $\beta=\beta_1+\beta_2+\beta_3+2\delta_\bg$, and $\beta=k\delta_\bg+\beta_1$, $k\in \{1,2,3,4\}$ which appear in the sums in the right hand side of (\ref{defP-}) are not purely polynomial.

(ii) For all $\beta_i$, $i=1,2,3$ in the non--zero terms in the right hand side of (\ref{defP-}) we have either $[\beta_i]\geq 0$ or $\beta_i=\delta_{\n}$ for some ${\n}\in \N^4$, and in the latter case $[\beta_i]=-1$.

For a sum $\beta=\beta_1+\beta_2+\delta_\bg$ (resp. $\beta=\beta_1+\beta_2+\beta_3+2\delta_\bg$, and $\beta=k\delta_\bg+\beta_1$, $k\in \{1,2,3,4\}$) one has by definition (\ref{pop}) that $[\beta]=[\beta_1]+[\beta_2]+1$ (resp. $[\beta]=[\beta_1]+[\beta_2]+[\beta_3]+2$, and $[\beta]=k+[\beta_1]$, $k\in \{1,2,3,4\}$). Therefore if for at least one $\beta_i$ in the first two sums $[\beta_i]\geq 0$ then $[\beta]\geq 0$, and in the last sum we always have $[\beta]\geq 0$.

Otherwise, $\beta=\delta_{\n}+\delta_{{\n}'}+\delta_\bg$ for some ${\n},{\n}'\in \N^4$ (resp. $\beta=\delta_{\n}+\delta_{{\n}'}+\delta_{{\n}''}+2\delta_\bg$ for some ${\n},{\n}',{\n}''\in \N^4$), and in all cases $\beta \in \M'$.

(iii) If $[\beta]<0$ then $\Phi_{\beta}^-=0$ if $\beta\not\in \M'\setminus \M_{\geq 0}$ by part (ii).

For $\beta \in \M'\setminus \M_{\geq 0}$, from the last paragraph of the proof of part (ii) we deduce that $\beta=\delta_{\n}+\delta_{{\n}'}+\delta_\bg$ for some ${\n},{\n}'\in \N^4$ or $\beta=\delta_{\n}+\delta_{{\n}'}+\delta_{{\n}''}+2\delta_\bg$ for some ${\n},{\n}',{\n}''\in \N^4$. From formula (\ref{hom}) we obtain that for $\beta=\delta_{\n}+\delta_{{\n}'}+\delta_\bg$ one has $|\beta|=|{\n}|+|{\n}'|+1$, and for $\beta=\delta_{\n}+\delta_{{\n}'}+\delta_{{\n}''}+2\delta_\bg$ one has $|\beta|=|{\n}|+|{\n}'|+|{\n}''|+2$. In particular, 
\be\lab{n1}
|{\n}|+|{\n}'|=|\delta_{\n}+\delta_{{\n}'}+\delta_\bg|-1=|\beta|-1
\ee 
in the former case, and 
\be\lab{n2}
|{\n}|+|{\n}'|+|{\n}''|=|\delta_{\n}+\delta_{{\n}'}+\delta_{{\n}''}+2\delta_\bg|-2=|\beta|-2
\ee 
in the latter case.

Since for any ${\n}\in \N^4$ one has $\Phi_{\delta_{\n}}\in \P_{\delta_{\n}}=\P_{|{\n}|}\otimes V_{\delta_{\n}}$, for $\beta=\delta_{\n}+\delta_{{\n}'}+\delta_\bg$ the terms in the first sum in the right hand side of (\ref{defP-}) belong to $\P_{|{\n}|+|{\n}'|-1}\otimes V_{\beta}^-$, and for $\beta=\delta_{\n}+\delta_{{\n}'}+\delta_{{\n}''}+2\delta_\bg$ the terms in the second sum in the right hand side of (\ref{defP-}) belong to $\P_{|{\n}|+|{\n}'|+|{\n}''|}\otimes V_{\beta}^-$. By (\ref{n1}) and (\ref{n2}) this implies $\Phi_{\beta}^-\in \P_{|\beta|-2}\otimes V_{\beta}^-\subset \O_\beta^-$.

(iv) For all $\beta_i$, $i=1,2,3$ in the non--zero terms in the right hand side of (\ref{defP-}) we have $[\beta_i]\geq -1$.

For a sum $\beta=\beta_1+\beta_2+\delta_\bg$ (resp. $\beta=\beta_1+\beta_2+\beta_3+2\delta_\bg$, and $\beta=k\delta_\bg+\beta_1$, $k\in \{1,2,3,4\}$) one has in this case by definition (\ref{pop}) that $[\beta]=[\beta_1]+[\beta_2]+1\geq -1$ (resp. $[\beta]=[\beta_1]+[\beta_2]+[\beta_3]+2\geq -1$, and $[\beta]=k+[\beta_1]\geq 0$, $k\in \{1,2,3,4\}$). This implies $\Phi_{\beta}^-=0$ if $[\beta]<-1$. 

\epr

\br\lab{covrem}
Let $O(\R^3)$ be the group of orthogonal transformations of $\R^3$ equipped with the scalar product fixed in the Introduction, ${\rm Aut}(\k)$ the group of automorphisms of $\k$.

Note that equation (\ref{LangdTrm}) is manifestly covariant under the natural action of $O(\R^3)$ on the space of $\k$--valued differential one--forms on $\R^3$ depending also on $t\in \R$, which is induced by the action of $O(\R^3)$ on the $\R^3$--factor with coordinates $(x_1,x_2,x_3)$ in $\R^4=\R\times \R^3$ with coordinates $(t,x_1,x_2,x_3)$, and under the natural action of ${\rm Aut}(\k)$ on this space induced by its action on $\k$. Since the operator $L$ is obviously invariant under both actions, the right hand side of the equivalent equation (\ref{LangHcrm}) is also covariant under the resulting natural action of ${\rm Aut}(\k)\times O(\R^3)$. As a consequence, the equivalent equation (\ref{LangdTrmAB}) and its right hand side enjoy the same covariance property under the natural action of ${\rm Aut}(\k)\times O(\R^3)$. 

Note that ${\rm Aut}(\k)\times O(\R^3)$ also naturally acts on $V=\k\otimes \R^3$, $W_\n=\k\otimes \R^3$, $\n\in \N^4$, and hence on $W_\beta$, $W_\beta^*$, $\beta\in \Ml$ defined by (\ref{defWb}) and (\ref{defWb*}). These actions induce natural actions of ${\rm Aut}(\k)\times O(\R^3)$ on $\O_\beta$ and $\O_\beta^-$, $\beta\in \Ml$.

This fact, together with the covariance property of the right hand side of (\ref{LangdTrmAB}) under the action of ${\rm Aut}(\k)\times O(\R^3)$ justified in the end of the previous paragraph, implies that equations (\ref{Ph0}), (\ref{defP-}) are also covariant with respect to the natural action of ${\rm Aut}(\k)\times O(\R^3)$. We shall partially use this fact in the proof of Proposition \ref{BPHZc}.
\er

Now using (\ref{defP-}) and Lemma \ref{lF} we apply an inductive procedure to construct elements $\Phi_{\beta}$ with two different choices of $\Phi_{\delta_{\n}}$, ${\n}\in \N^4$ which yield the recentered maps and the canonical lift.

In the first case for $x\in \R^4$ we denote $\Phi_\beta=\Pi_{x\beta}$, $\Phi_\beta^-=\Pi_{x\beta}^-$ and for all purely polynomial multi--indices $\delta_{\n}$, ${\n}\in \N^4$ we define
\be\lab{Pipol}
\Pi_{x\delta_{\n}}(y):=(y-x)^{\n}{\rm Id}_{W_{\delta_\n}}^V\in \P_{\delta_{\n}}= \P_{|{\n}|}\otimes V_{\delta_{\n}}\subset \O\otimes V_{\delta_{\n}}=\O_{\delta_{\n}},~{\n}\in \N^4,
\ee
where ${\rm Id}_{W_{\delta_\n}}^V\in W_{\delta_{\n}}^*\otimes V= V_{\delta_{\n}}$ is the natural operator identifying $W_{\delta_\n}$ and $V$.


Assuming that for a non--purely polynomial  multi--index $\beta$ elements $\Pi_{x\gamma}\in \O_\gamma$ for $|\gamma|<|\beta|$ are given, and $\Pi_{x\gamma}=0$ for $\gamma\not\in \M_{\geq 0}\cup \M_{pp}$ with $|\gamma|<|\beta|$, we introduce $\Pi_{x\beta}^-=\Phi_\beta^-\in \O_\beta^-$ by (\ref{defP-}).

Then we define 
\be\lab{defP}
\Pi_{x\beta}(y)=(\K\Pi_{x\beta}^-)(y)-\sum_{{\n}\in \N^4: |{\n}|<|\beta|}\partial^{\n}(\K\Pi_{x\beta}^-)(x)\frac{(y-x)^{\n}}{{\n}!},
\ee
where for ${\n}=(n_0,n_1,n_2,n_3)$ ${\n}!=\prod_{i=0}^3n_i!$ and the sum in the last term is void if $|\beta|\leq 0$.

By Lemma \ref{lF} (ii) and (iii) $\Pi_{x\beta}^-=0$ if $\beta \not\in \M'$, and
 if $[\beta]<0$ then $\Pi_{x\beta}^-\in \P_{|\beta|-2}\otimes V_{\beta}\subset \O_\beta^-$. Therefore, since any polynomial of parabolic degree at most $n$ is equal to its Taylor polynomial of parabolic degree $n$, and by Lemma \ref{KPol} (i) the subspace $\P_n\otimes V_{\beta}\subset \O_\beta^-$ is invariant under the action of the operator $\K$, definition (\ref{defP}) implies that $\Pi_{x\beta}=0$ if $[\beta]<0$, i.e. $\Pi_{x\beta}$ satisfies the same property as $\Pi_{x\gamma}$ with $|\gamma|<|\beta|$ in the induction assumption.

This assumption obviously holds for $\Pi_{x0}$ as $[0]=0$, and $\alpha=|0|$ is the minimal possible value of the homogeneity which is only achieved at $\beta=0$. Therefore the inductive procedure outlined above uniquely defines $\Pi_{x\beta}$ and $\Pi_{x\beta}^-$ for all $\beta$. Moreover, $\Pi_{x\beta}^-=0$ if $\beta \not\in \M'$, and $\Pi_{x\beta}=0$ for $\beta\not\in \M_{\geq 0}\cup \M_{pp}$,
\be\lab{Pi-0}
\Pi_{x\beta}^-=0~{\rm for }~\beta \not\in \M',
\ee
\be\lab{Pi0}
\Pi_{x\beta}=0~{\rm for }~[\beta]<0,~\beta\neq \delta_{\n},~{\n}\in \N^4.
\ee

Also, by definition (\ref{Pipol}) and  (\ref{defP}) one has 
\be\lab{Pder}
\partial^\n(\Pi_{x\beta})(y)=\partial^\n(\K\Pi_{x\beta}^-)(y) ~{\rm for ~any }~ x,y\in \R^4,~\beta\in \Ml,~\n\in \N^4 
\ee
$$
~{\rm such~ that }~|\n|>|\beta|.
$$

In the second case we denote $\Phi_\beta={\bf \Pi}_{\beta}$, $\Phi_\beta^-={\bf \Pi}_{\beta}^-$ and for all purely polynomial multi--indices $\delta_{\n}$, ${\n}\in \N^4$ we define
\be\lab{Pibpol}
{\bf \Pi}_{\delta_{\n}}(y):=y^{\n}{\rm Id}_{W_{\delta_\n}}^V\in \P_{\delta_{\n}}= \P_{|{\n}|}\otimes V_{\delta_{\n}}\subset \O\otimes V_{\delta_{\n}}=\O_{\delta_{\n}},~{\n}\in \N^4.
\ee

Assuming that for a non--purely polynomial  multi--index $\beta$ elements ${\bf \Pi}_{\gamma}\in \O_\gamma$ for $|\gamma|<|\beta|$ are given, and ${\bf \Pi}_{\gamma}=0$ if $[\gamma]<-1$, we introduce ${\bf \Pi}_{\beta}^-=\Phi_\beta^-\in \O_\beta^-$ by (\ref{defP-}).

Then we define 
\be\lab{defPb}
{\bf \Pi}_{\beta}=\K{\bf \Pi}_{\beta}^-.
\ee

By Lemma \ref{lF} (i) ${\bf \Pi}_{\beta}^-=0$ if $\beta\in \M_{pp}$, and by Lemma \ref{lF} (iv) ${\bf \Pi}_{\beta}^-=0$ if $[\beta]<-1$, and hence by (\ref{defPb}) ${\bf \Pi}_{\beta}=0$ if $[\beta]<-1$. Obviously for ${\bf \Pi}_{0}$ one has $[0]=0\geq -1$, and $\alpha=|0|$ is the minimal possible value of the homogeneity which is only achieved at $\beta=0$. Also, for ${\bf \Pi}_{\delta_{\n}}$, ${\n}\in \N^4$ we have $[\delta_{\n}]=-1$.

Therefore the inductive procedure described above uniquely defines ${\bf \Pi}_{\beta}$ and ${\bf \Pi}_{\beta}^-$ for all $\beta$. Moreover, ${\bf \Pi}_{\beta}^-=0$ if $[\beta]<-1$ or $\beta\in \M_{pp}$, and ${\bf \Pi}_{\beta}=0$ if $[\beta]<-1$,
\be\lab{BP0}
{\bf \Pi}_{\beta}^-={\bf \Pi}_{\beta}=0~{\rm for }~[\beta]<-1,~{\bf \Pi}_{\delta_\n}^-=0,~\n\in \N^4.
\ee

\br
Note that by the definition the kernel $K_m$ is invariant with respect to that action of $O(\R^3)$ on $\R^4=\R\times \R^3$ considered in Remark \ref{covrem}. Therefore if $\varsigma \in \D$ is invariant with respect to this action of $O(\R^3)$ then the kernel $K(x)=K_m(x)\varsigma(x)$ is also $O(\R^3)$--invariant. Thus the operator $\K$ is $O(\R^3)$--invariant. Together with covariance of Taylor polynomials of a fixed parabolic degree with respect to spatial orthogonal transformations, which preserve the parabolic degree, this implies that formula (\ref{defP}) is covariant with respect to the natural action of ${\rm Aut}(\k)\times O(\R^3)$ on $\O_\beta$ and $\O_\beta^-$, $\beta\in \Ml$ in the sense that if $\theta \in {\rm Aut}(\k)\times O(\R^3)$ then for $\beta \not\in \M_{pp}$ one has
\be\lab{Pcov}
(\theta(\Pi_{\theta(x)\beta}))(y)=
\ee
$$
=(\K(\theta(\Pi_{\theta(x)\beta}^-)))(y)-\sum_{{\n}\in \N^4: |{\n}|<|\beta|}\partial^{\n}(\K(\theta(\Pi_{\theta(x)\beta}^-))(x)\frac{(y-x)^{\n}}{{\n}!}.
$$

Clearly, a similar property holds for formula (\ref{defPb}),
\be\lab{Pbcov}
\theta({\bf \Pi}_{\beta})=\K(\theta({\bf \Pi}_{\beta}^-)),~\beta \not\in \M_{pp}.
\ee
\er

Now we define recentered maps in the sense of Hairer as follows. 
For any vector space $W$, denote by ${\rm Id}_W$ its identity linear automorphism.
Let $\beta\in \Ml$, and define
\be\lab{tau}
\tau_\beta:={\rm Id}_{V_\beta^*}\in {\rm End}(V_\beta^*)\simeq V_\beta^*\otimes V_\beta,
\ee
and 
\be\lab{tau-}
\tau_\beta^-:={\rm Id}_{{V_\beta^-}^*}\in {\rm End}({V_\beta^-}^*)\simeq {V_\beta^-}^*\otimes V_\beta^-.
\ee

For any $\Phi_\beta\in \O_\beta$, $\Phi_\beta^-\in \O_\beta^-$, $\beta \in \Ml$ we define a linear map 
$$
\Phi:{\bf T}\to \O, ~{\bf T}:=\bigoplus_{\beta\in \Ml}V_\beta^*\oplus \bigoplus_{\beta\in \Ml}{V_\beta^{-}}^*,
$$
by
\be\lab{defPP}
(\Phi\otimes {\rm Id}_{V_\beta})\tau_\beta=\Phi_{\beta},~ (\Phi\otimes {\rm Id}_{V_\beta^-})\tau_\beta^-=\Phi_{\beta}^-.
\ee

We can apply this definition to $\Pi_{x\beta}\in \O_\beta$, $\Pi_{x\beta}^-\in \O_\beta^-$, $\beta \in \Ml$ to get a linear map
\be\lab{PixdefbfT}
\Pi_x:{\bf T}\to \O.
\ee
Due to properties (\ref{Pi-0}) and (\ref{Pi0}) $\Pi_x|_{V_\beta^*}=0$ for $\beta\not\in \M_{\geq 0}\cup \M_{pp}$, and $\Pi_x|_{{V_\beta^-}^*}=0$ for $\beta \not\in \M'$. Therefore $\Pi_x$ is completely determined by its restriction to $T\subset {\bf T}$, 
\be\lab{Pixdef}
\Pi_x:T\to \O.
\ee
The maps $\Pi_x$, $x\in \R^4$ are called recentered maps (we follow the terminology of \cite{C2}). 

Similarly, applying definition (\ref{defPP}) to  ${\bf \Pi}_\beta$ and ${\bf \Pi}_\beta^-$, $\beta \in \Ml$ we can define a linear map
\be\lab{Pidef}
{\bf \Pi}:{\bf T}\to \O.
\ee
Due to property (\ref{BP0}) ${\bf \Pi}|_{V_\beta^*}=0$ and ${\bf \Pi}|_{{V_\beta^-}^*}=0$ for $[\beta ]<-1$. Therefore ${\bf \Pi}$ is completely determined by its restriction to 
$$
\bigoplus_{\beta\in \Ml:[\beta]\geq -1}V_\beta^*\oplus \bigoplus_{\beta\in \Ml:[\beta]\geq -1}{V_\beta^{-}}^*\subset {\bf T},
$$ 
$$
{\bf \Pi}:\bigoplus_{\beta\in \Ml:[\beta]\geq -1}V_\beta^*\oplus \bigoplus_{\beta\in \Ml:[\beta]\geq -1}{V_\beta^{-}}^*\to \O.
$$
The map ${\bf \Pi}$ is called the canonical lift.

In conclusion, we obtain certain local analytic bounds for $\Pi_{x\beta}$ and $\Pi_{x\beta}^-$ which are similar to those obtained in \cite{H1}, Proposition 8.27. These bounds are provisional in our approach. In fact, in Section \ref{recestptwsgl} we shall derive global  uniform estimates for $\Pi_{x\beta}$ and $\Pi_{x\beta}^-$ which hold almost surely.
\bp\lab{Pibond}
Fix any norm on $V$, and for every $\beta\in \Ml$ fix any norm on $W_\beta$. Denote by $|\cdot |$ the induced operator norm on $V_\beta=V_\beta^-={\rm Hom}(W_\beta,V)$. Then for every compact set $\mK\subset \R^4$, every $\rho>0$, and every $\beta\in \Ml$ there exist constants $C_{\mK,\beta,\n}(\rho)$, $C_{\mK,\beta,\n}^-(\rho)\geq 0$, $\n\in \N^4$ such that
\be\lab{Pibnd1}
|\partial^{\n}(\Pi_{x\beta})(y)|\leq C_{\mK,\beta,\n}(\rho)|x-y|^{|\beta|-|\n|},~x,y\in \mK,~x\neq y,~\n\in \N^4,
\ee
\be\lab{Pibnd2}
|\partial^{\n}(\Pi_{x\beta}^-)(y)|\leq C_{\mK,\beta,\n}^-(\rho)|x-y|^{|\beta|-2-|\n|},~x,y\in \mK,~x\neq y,~\n\in \N^4.
\ee
\ep

\bpr
For $\beta \in \M_{pp}$ the bound in (\ref{Pibnd1}) follows from definition (\ref{Pipol}) and the bound in (\ref{Pibnd2}) is trivial by (\ref{Pi-0}).

If $\beta\not \in \M_{pp}$ and $|\n|> |\beta|$ the bound in (\ref{Pibnd1}) is obvious as $\Pi_{x\beta}\in \O_\beta=\O\otimes V_\beta\subset C^\infty(\R^4)\otimes V_\beta$, and if $\beta\not \in \M_{pp}$ and $|\n|> |\beta|-2$ the bound in (\ref{Pibnd2}) is obvious as $\Pi_{x\beta}^-\in \O_\beta^-=\O\otimes V_\beta^-\subset C^\infty(\R^4)\otimes V_\beta^-$.

If $|\n|\leq |\beta|$ the bound in (\ref{Pibnd1}) follows from the formula obtained by applying $\partial^{\n}$ to (\ref{defP}) and from the expression for the remainder in the generalized Taylor formula (see \cite{H1}, Appendix A and the proof of Proposition 8.27 in \cite {H1}). If $|\n|\leq |\beta|-2$ the bound in (\ref{Pibnd2}) follows from (\ref{Pibnd1}) and from the formula obtained by applying $\partial^{\n}$ to (\ref{defP-}) in the case of $\Pi_{x\beta}^-$.

\epr

Taking the limit  $y\to x$ in (\ref{Pibnd1}) and (\ref{Pibnd2}) when $|\n|< |\beta|$ and $|\n|<|\beta|-2$, respectively, we immediately obtain the following corollary.
\bc
\be\lab{Pvan}
\partial^\n(\Pi_{x\beta})(x)=0 ~{\rm for~ any }~ x\in \R^4,~\beta\in \Ml,~\n\in \N^4 ~{ \rm such ~that }~|\n|<|\beta|,
\ee
and
\be\lab{P-van}
\partial^\n(\Pi_{x\beta}^-)(x)=0 ~{\rm for~ any }~ x\in \R^4,~\beta\in \Ml,~\n\in \N^4 
\ee
$$
{\rm such~ that }~|\n|<|\beta|-2.
$$
\ec


\subsection{Automorphisms relating the recentered maps and the canonical lift}\lab{Frec}

\setcounter{equation}{0}
\setcounter{theorem}{0}

In this section, as a technical tool for defining structure group maps we introduce certain automorphisms relating the recentered maps and the canonical lift. For this purpose we consider a formal power series algebra in terms of which these automorphisms are defined. Firstly, we introduce the relevant notions and terminology.
  
Let $f_i$, $i=1,\ldots, {\rm dim}(V)$ be a linear basis of $V=\k\otimes \R^3$, and $\z_\bg$, ${\z}^i_{\n}$, $i=1,\ldots, {\rm dim}(V)$, ${\n}\in \N^4$ formal variables. Consider the algebra 
$$
\Sm[[{\z}]]:=\Sm[[{\z}_{\bg},\{{\z}^i_{\n}:i=1,\ldots, {\rm dim}(V), {\n}\in \N^4\}]]
$$ 
of formal power series in the variables $\z_\bg$, $\z^i_{\n}$, $i=1,\ldots, {\rm dim}(V)$, ${\n}\in \N^4$ with coefficients in the algebra $\Sm$.

For any monomial in $\Sm[[\z]]$ in the variables $\z_\bg$, $\z^i_{\n}$, $i=1,\ldots, {\rm dim}(V)$, ${\n}\in \N^4$ we introduce its $\Ml$--valued multi--index degree $\dm$ by requiring that $\dm(\z_\bg)=\delta_\bg$, $\dm(\z^i_{\n})=\delta_\n$, $i=1,\ldots, {\rm dim}(V)$, ${\n}\in \N^4$ and that the multi--index degree of the product of two monomials is equal to the sum of their multi--index degrees. 

For any element $Z\in \Sm[[{\z}]]$ and any multi--index $\beta$ we denote by $(Z)_\beta$ the sum of all monomials of multi--index degree $\beta$ which appear in $Z$. This sum is always finite since the set of multi--indices is the free commutative monoid generated by the elements $\delta_\bg$ and $\delta_{\n}$, ${\n}\in \N^4$, and hence each its element is represented as a sum of other elements in finitely many ways.

Identifying $V$ with $W_{\delta_{\n}}$ we define elements 
$$
\z^{\delta_{\n}}:=\sum_{i=1}^{{\rm dim}(V)}f_i\z^i_{\n}\in W_{\delta_{\n}}[\{\z^i_{\n}:i=1,\ldots, {\rm dim}(V)\}]\subset \Sm[[\z]].
$$
For any multi--index $\beta=l\delta_\bg+\sum_{i=1}^pk_i\delta_{{\n}_i}$, let 
$$
\z^\beta=\z_\bg^l(\z^{\delta_{{\n}_1}})^{k_1}\ldots (\z^{\delta_{{\n}_p}})^{k_p}\in W_{\beta}[\z_\bg, \{\z^j_{{\n}_i}:j=1,\ldots, {\rm dim}(V),i=1,\ldots, p\}]\subset \Sm[[\z]],
$$
where the product in $\Sm[[\z]]$ is used in the right hand side.

\br\lab{indepbas}
For each $\n\in \N^4$ one can identify the real vector space spanned by the formal variables $\z^i_{\n}$, $i=1,\ldots, {\rm dim}(V)$ with $V^*=W_{\delta_{\n}}^*$ in such a way that $\z^i_{\n}$, $i=1,\ldots, {\rm dim}(V)$ is the linear basis of $W_{\delta_{\n}}^*$ dual to the basis $f_i$, $i=1,\ldots, {\rm dim}(V)$ of $W_{\delta_{\n}}=V$. Under this identification and using the natural isomorphism $W_{\delta_{\n}}\otimes W_{\delta_{\n}}^*\simeq {\rm End}(W_{\delta_{\n}})$ we have $\z^{\delta_{\n}}={\rm Id}_{W_{\delta_{\n}}}\in {\rm End}(W_{\delta_{\n}})\simeq W_{\delta_{\n}}\otimes W_{\delta_{\n}}^*$, and hence the definitions of $\z^{\delta_{\n}}$, $\n\in \N^4$, and, as a consequence, of $\z^\beta$, $\beta\in \Ml$ do not depend on the choice of the basis $f_i$, $i=1,\ldots, {\rm dim}(V)$. 

From now on we shall always assume that the real vector space spanned by the formal variables $\z^i_{\n}$, $i=1,\ldots, {\rm dim}(V)$ is identified with $W_{\delta_{\n}}^*$ in such a way that $\z^i_{\n}$, $i=1,\ldots, {\rm dim}(V)$ is the linear basis of $W_{\delta_{\n}}^*$ dual to the basis $f_i$, $i=1,\ldots, {\rm dim}(V)$ of $W_{\delta_{\n}}$. 
\er

By the definition the coefficients of the polynomial $\z^\beta$ in the variables $\z_\bg$ and $\z^j_{{\n}_i}$, $j=1,\ldots, {\rm dim}(V)$, $i=1,\ldots, p$ form a linear basis of $W_\beta$, and hence any $\Phi_{\beta}\in \O(V)\otimes W_\beta^*\simeq \O_\beta$ is completely determined by its values on these coefficients, i.e. by the polynomial 
$$
\Phi_{\beta}(\z^\beta)\in \O(V)[\z_\bg,\{\z^j_{{\n}_i}:j=1,\ldots, {\rm dim}(V),i=1,\ldots, p\}],
$$
in $\z_\bg$, $\z^j_{{\n}_i}$, $j=1,\ldots, {\rm dim}(V)$, $i=1,\ldots, p$ the coefficients of which are obtained by evaluating $\Phi_{\beta}\in \O_\beta\simeq \O(V)\otimes W_\beta^*$ on the coefficients of the polynomial $\z^\beta\in W_{\beta}[\z_\bg,\{\z^j_{{\n}_i}:j=1,\ldots, {\rm dim}(V),i=1,\ldots, p\}]$ with the help of the inclusion 
$$
\O_\beta \otimes W_\beta \simeq \O(V)\otimes W_\beta^*\otimes W_\beta\subset \O\otimes \Sm^*\otimes \Sm
$$ 
and of the natural map 
\be\lab{parO}
\O\otimes \Sm^*\otimes \Sm\to \O
\ee
induced by pairing (\ref{par}).


If $\Phi_{\beta}\in \O_\beta$ are defined for all multi--indices $\beta$ it is also convenient to introduce the formal power series
\be\lab{Pb}
\Phi(\z):=\sum_{\beta\in \Ml}\Phi_{\beta}(\z^\beta)\in \O(V)[[\z_\bg,\{\z^i_{\n}, i=1,\ldots, {\rm dim}(V), {\n}\in \N^4\}]],
\ee
and any element of $\O(V)[[\z_\bg,\{\z^i_{\n}, i=1,\ldots, {\rm dim}(V), {\n}\in \N^4\}]]$ uniquely determines elements of $\O_\beta$ for all multi--indices $\beta$.

Now for all $x\in \R^4$ we define endomorphisms $F_x^*$ of the algebra $\R[[\z]]:=\R[[\z_\bg,\{\z^i_{\n}, i=1,\ldots, {\rm dim}(V), {\n}\in \N^4\}]]\subset \Sm[[\z]]$ by requiring that
\be\lab{Fdef}
F_x^*1=1,~F_x^*\z_\bg=\z_\bg,
\ee
\be\lab{Fdef1}
({\rm Id}_{W_{\delta_{\n}}}^V\otimes F_x^*)\z^{\delta_{\n}}=\sum_{\tiny\begin{array}{c} \m\in \N^4:\\ |\n|\leq |\m|\end{array}}{\m \choose \n}(-x)^{\m-\n}({\rm Id}_{W_{\delta_{\m}}}^V\otimes {\rm Id}_{\R[[\z]]})\z^{\delta_\m}-
\ee
$$
-\frac{1}{\n!}\sum_{\beta\in \M':|{\bf n}|\leq |\beta|}\sum_{\tiny\begin{array}{c} \m\in \N^4:\\ |\m|\leq |\beta|-|{\bf n}|\end{array}}\frac{(-x)^\m}{\m!}\partial^{\n+\m}(\K\Pi_{x\beta}^-)(x)(\z^\beta),
$$
where  
$$
{\m \choose \n}=\frac{\m!}{\n!(\m-\n)!},
$$ 
and as usual for $\n=(n_0,n_1,n_2,n_3)$, $\m=(m_0,m_1,m_2,m_3)$ we define ${\m \choose \n}=0$ unless $m_i\geq n_i$ for $i=0,1,2,3$. In the first line in (\ref{Fdef1}) each $\z^{\delta_\m}$ is regarded as an element of $W_{\delta_{\m}}\otimes \R[[\z]]\supset W_{\delta_{\n}}[\{\z^i_{\n}:i=1,\ldots, {\rm dim}(V)\}]$, so that (\ref{Fdef1}) must hold in $V[[\z_\bg,\{\z^i_{\n}, i=1,\ldots, {\rm dim}(V), {\n}\in \N^4\}]]\simeq V\otimes \R[[\z]]$.

Since $V$ is finite--dimensional this definition completely determines the action of $F_x^*$ on the formal variables $\z_\n^i$ and $\z_\bg$. Below we shall see that the natural extension of $F_x^*$ defined by (\ref{Fdef}), (\ref{Fdef1}) to endomorphisms of $\R[[\z]]$ is well defined.

Definition (\ref{Fdef}), (\ref{Fdef1}) is obtained by dualizing a similar definition of a dual map for regularity structures (see \cite{H1}, \S 8.2 or \cite{C2}, Section 2). The maps $F_x^*$, $x\in \R^4$ actually induce automorphisms of the model space $T$. We shall not need this fact. A similar property will be proved below for the structure group maps.

To confirm that $F_x^*$ can be extended to automorphisms of the algebra $\R[[\z]]$ we firstly note that formula (\ref{hom}) for the homogeneity can be rewritten in the form
\be
|\beta|=\alpha(1+[\beta])+\beta(g)+\sum_{{\n}\in \N^4}\beta({\n})|{\n}|.
\ee
This implies that for $\beta\in \M_{\geq 0}$  one has $|\beta|\in \alpha(\N\setminus \{0\})+\N$. Since $\alpha \not\in \Q$, we deduce that for $\beta\in \M_{\geq 0}$ one has $|\beta|\not\in \Q$,
\be\lab{poph}
\{|\beta|:\beta\in \M_{\geq 0}\}\subset \R\setminus \Q.
\ee

Therefore, if in the last sum in (\ref{Fdef1}) $|\beta|=|\n|\in \N$ then $[\beta]<0$. In this case by Lemma \ref{lF} (iii) $\Pi_{\beta}^-\in \P_{|\beta|-2}\otimes V_{\beta}^-\subset \O_\beta^-$. This implies $\partial^{\n+\m}(\K\Pi_{x\beta}^-)(x)=0$, and hence (\ref{Fdef1}) takes the form
\be\lab{Fdef2}
({\rm Id}_{W_{\delta_{\n}}}^V\otimes F_x^*)\z^{\delta_{\n}}=\sum_{\tiny\begin{array}{c} \m\in \N^4:\\ |\n|\leq |\m|\end{array}}{\m \choose \n}(-x)^{\m-\n}({\rm Id}_{W_{\delta_{\m}}}^V\otimes {\rm Id}_{\R[[\z]]})\z^{\delta_\m}-
\ee
$$
-\frac{1}{\n!}\sum_{\beta\in \M_{\geq 0}:|{\bf n}|< |\beta|}\sum_{\tiny\begin{array}{c} \m\in \N^4:\\ |\m|< |\beta|-|{\bf n}|\end{array}}\frac{(-x)^\m}{\m!}\partial^{\n+\m}(\K\Pi_{x\beta}^-)(x)(\z^\beta),
$$
as if $[\beta]\geq 0$ one has $|\beta|\not\in \Q$ by (\ref{poph}), 
so $|\beta|\leq |\n|\in \N$ implies $|\beta|<|\n|$, and $|\m|\leq |\beta|-|{\bf n}|$ implies $|\m|< |\beta|-|{\bf n}|$.

Since the right hand sides of the second expression in (\ref{Fdef}) and of (\ref{Fdef2}) do  not contain zero powers of the formal variables $\z_\bg$, $\z^i_{\n}, i=1,\ldots, {\rm dim}(V), {\n}\in \N^4$ all $F_x^*$ can be naturally extended by multiplicativity to continuous automorphisms of the algebra $\R[[\z]]$. 



Next, we study properties of the endomorphisms $F_x^*$. For this purpose we introduce some relevant notions and notation.
First note that any linear map $F^*:\R[\z]\to \R[[\z]]$, where $\R[\z]$  is the polynomial algebra in the variables $\z_\bg$, $\z^i_{\n}, i=1,\ldots, {\rm dim}(V), {\n}\in \N^4$ naturally extends to a linear map $\Sm[\z]\to \Sm[[\z]]$, and any continuous linear endomorphism $F^*$ of the vector space $\R[[\z]]$ equipped with the formal power series algebra topology naturally induces a continuous linear endomorphism of the vector space $\Sm[[\z]]$. We denote these extensions by the same symbols. 

For any $\beta,\gamma\in \Ml$ and any  linear map $F^*:\R[\z]\to \R[[\z]]$ we can define linear maps $(F)_{\beta}^\gamma: W_\beta\to W_\gamma$ by
\be\lab{Fmatr}
(F^*\z^\gamma)_\beta=(F)_{\beta}^\gamma \z^\beta.
\ee

Now we obtain a triangularity property of the automorphisms $F_x^*$ with respect to the homogeneity of multi--indices.
For this purpose we note that from (\ref{Fdef2}) one can infer
\be\lab{Fdg}
(F_x)_{\delta_\m}^\gamma=\left\{\begin{array}{cc} {\m \choose \n}(-x)^{\m-\n}{\rm Id}_{W_{\delta_{\m}}}^{W_{\delta_{\n}}} & {\rm if }~ \gamma =\delta_\n,~\n\in \N^4 \\ 0 & {\rm else} \end{array}\right. ,
\ee
where ${\rm Id}_{W_{\delta_{\m}}}^{W_{\delta_{\n}}}: W_{\delta_{\m}}\to W_{\delta_{\n}}$ is the natural isomorphism.

From formula (\ref{Fdef2}) we also obtain that 
$$
F_x^*\z^{\delta_\n}=\z^{\delta_\n}+\sum_{\tiny\begin{array}{c}\beta\in \M_{\geq 0}\cup \M_{pp} : \\ |{\bf n}|< |\beta|\end{array}}(F_x)_{\beta}^{\delta_\n} \z^\beta
$$
Together with (\ref{Fdef}) this implies by multiplicativity that for $\gamma=l\delta_\bg+\sum_{i=1}^p\delta_{\m_i}$, $l\in \N$, $\m_i\in \N^4$ one has
\be\lab{F*zg}
F_x^*\z^\gamma=\z^\gamma+\z_\bg^l \hspace{-0.5cm} \sum_{\tiny\begin{array}{c}\beta_i\in \M_{\geq 0}\cup \M_{pp}: \\|{\bf m}_i|\leq |\beta_i|,\\ i=1,\ldots,p \end{array}}\hspace{-0.5cm} (F_x)_{\beta_1}^{\delta_{\m_1}} \z^{\beta_1}\ldots (F_x)_{\beta_p}^{\delta_{\m_p}} \z^{\beta_p},
\ee
where in each term of the sum at least one $\beta_i$, $i=1,\ldots,p$ satisfies $|{\bf m}_i|< |\beta_i|$. Therefore by linearity of the function $|\cdot|-\alpha$ the homogeneity of each such term equal to $|l\delta_\bg+\beta_1+\ldots+\beta_p|$ satisfies
$$
|l\delta_\bg+\beta_1+\ldots+\beta_p|-\alpha=|l\delta_\bg|-\alpha+\sum_{i=1}^p(|\beta_i|-\alpha)>|l\delta_\bg|-\alpha+\sum_{i=1}^p(|\m_i|-\alpha)=
$$
\be\lab{bgsum}
=|l\delta_\bg|-\alpha+\sum_{i=1}^p(|\delta_{\m_i}|-\alpha)=|\gamma|-\alpha,
\ee
i.e. $|l\delta_\bg+\beta_1+\ldots+\beta_p|>|\gamma|$.

Thus for any monomial $Z\in \R[\z]$ of multi--index degree $\gamma$ one has 
\be\lab{triang}
(F_x^*Z-Z)_\beta = 0 ~{\rm for }~ |\beta|\leq|\gamma|.
\ee

Later we shall also need a similar property of the automorphisms $F_x^*$ with respect to the modified homogeneity  $|\cdot|_{\prec}$ defined by
\be\lab{homm}
|\beta|_\prec:=|\beta|+\frac{d}{2}([\beta]+1),~\beta\in \Ml.
\ee

Note that for $\beta\in \M_{\geq 0}\cup \M_{pp}$ one has $[\beta]\geq -1$, and hence by the definition for such $\beta$ one has $|\beta|_\prec\geq |\beta|$. Moreover, for any $\n\in\N^4$ $[\delta_\n]=-1$, and hence by the definition $|\delta_\n|_\prec=|\delta_\n|=|\n|$.

Using these observations we can immediately rewrite (\ref{F*zg}) in the form
$$
F_x^*\z^\gamma=\z^\gamma+\z_\bg^l \hspace{-0.5cm} \sum_{\tiny\begin{array}{c}\beta_i\in \M_{\geq 0}\cup \M_{pp}: \\|{\bf m}_i|_\prec\leq |\beta_i|_\prec,\\ i=1,\ldots,p \end{array}}\hspace{-0.5cm} (F_x)_{\beta_1}^{\delta_{\m_1}} \z^{\beta_1}\ldots (F_x)_{\beta_p}^{\delta_{\m_p}} \z^{\beta_p},
$$
where in each term of the sum at least one $\beta_i$, $i=1,\ldots,p$ satisfies $|{\bf m}_i|_\prec< |\beta_i|_\prec$. Therefore, similarly to (\ref{bgsum}), by linearity of the function $|\cdot|_\prec-\alpha-\frac d2$ the modified homogeneity of each such term equal to $|l\delta_\bg+\beta_1+\ldots+\beta_p|_\prec$ satisfies
$$
|l\delta_\bg+\beta_1+\ldots+\beta_p|_\prec-\alpha-\frac d2 =|l\delta_\bg|_\prec-\alpha-\frac d2 +\sum_{i=1}^p(|\beta_i|_\prec-\alpha-\frac d2 )>
$$
$$
>|l\delta_\bg|_\prec-\alpha-\frac d2 +\sum_{i=1}^p(|\m_i|_\prec-\alpha-\frac d2 )=|l\delta_\bg|_\prec-\alpha-\frac d2 +\sum_{i=1}^p(|\delta_{\m_i}|_\prec-\alpha-\frac d2 )=|\gamma|_\prec-\alpha-\frac d2 ,
$$
i.e. $|l\delta_\bg+\beta_1+\ldots+\beta_p|_\prec>|\gamma|_\prec$.

Thus for any monomial $Z\in \R[\z]$ of multi--index degree $\gamma$ one has 
\be\lab{triangm}
(F_x^*Z-Z)_\beta = 0 ~{\rm for }~ |\beta|_\prec\leq|\gamma|_\prec.
\ee

In the next lemma we study triangularity properties of endomorphisms of $\R[[\z]]$, similar to (\ref{triang}), (\ref{triangm}), in a more general setting.
\bl\lab{Flemc}
(i) Let $\R[\z]$ be the polynomial algebra in the formal variables $\z_\bg$, $\z^i_{\n}, i=1,\ldots, {\rm dim}(V), {\n}\in \N^4$, and $F^*:\R[\z]\to \R[[\z]]$ a linear map. Assume that for any monomial $Z\in \R[\z]$ of multi--index degree $\gamma$  
\be\lab{triangF'}
(F^*Z)_\beta = 0 ~{\rm for }~ |\beta|+N_\beta<|\gamma|,
\ee
where $N_\beta\in \R$ may depend on $\beta$. Then $F^*$ induces a continuous linear map $F^*:\R[[\z]]\to \R[[\z]]$ and a linear map
\be\lab{FdefT}
F:\bigoplus_{\beta\in \Ml}W_\beta\to \bigoplus_{\beta\in \Ml}W_\beta
\ee
defined by 
\be\lab{FdefT1}
\sum_{\beta\in \Ml}F(\z^\beta)=\sum_{\gamma\in \Ml}F^*(\z^\gamma),
\ee
where in the left hand side $F$ acts on the coefficients of the polynomial $\z^\beta\in W_{\beta}[\z_\bg,\{\z^j_{{\n}_i}:j=1,\ldots, {\rm dim}(V),i=1,\ldots, p\}]$ which are elements of $W_\beta$, $\beta\in \Ml$. 

Moreover, for any $\beta\in \Ml$
\be\lab{triang1'}
F:W_\beta\to \bigoplus_{\tiny\begin{array}{c} \gamma \in \Ml: \\ |\gamma|\leq |\beta|+N_\beta\end{array} }W_\gamma,
\ee
and if for any $\gamma\in \Ml$ we denote by $P_\gamma$ the natural projection operator $P_\gamma: \bigoplus_{\beta\in \Ml}W_\beta\to W_\gamma$ then $(F)_{\beta}^\gamma=P_\gamma F P_\beta$.

(ii) If $F^*:\R[\z]\to \R[[\z]]$ is a linear map such that for any monomial $Z\in \R[\z]$ of multi--index degree $\gamma$
\be\lab{triangF}
(F^*Z-Z)_\beta = 0 ~{\rm for }~ |\beta|\leq|\gamma|
\ee
then $F^*$ induces a continuous linear map $F^*:\R[[\z]]\to \R[[\z]]$ and a linear map (\ref{FdefT}). 

Moreover, one has
\be\lab{triang1}
F:W_\beta\to W_\beta\oplus\bigoplus_{\tiny\begin{array}{c} \gamma \in \Ml: \\ |\gamma|< |\beta|\end{array}}W_\gamma,
\ee
and 
\be\lab{triang2}
(F)_{\beta}^\beta={\rm Id}_{W_\beta},
\ee
i.e.
\be\lab{triang3}
(F-{\rm Id}_{W_\beta}):W_\beta\to \bigoplus_{\tiny\begin{array}{c}\gamma \in \Ml: \\ |\gamma|< |\beta|\end{array}}W_\gamma.
\ee
Also, if
\be\lab{triangm1}
(F^*Z-Z)_\beta = 0 ~{\rm for }~ |\beta|_\prec\leq|\gamma|_\prec
\ee
then
\be\lab{triang3m}
(F-{\rm Id}_{W_\beta}):W_\beta\to \bigoplus_{\tiny\begin{array}{c} \gamma \in \Ml: \\ |\gamma|_\prec< |\beta|_\prec\end{array}}W_\gamma.
\ee

The maps $F,F^*:\R[[\z]]\to \R[[\z]]$ are invertible, and ${F^*}^{-1}$ is continuous.
Denote by $\rm Id$ the identity linear endomorphism of $\bigoplus_{\gamma \in \Ml}W_\gamma$.
Then $F-{\rm Id}$ is locally nilpotent, and the inverse to $F$ can be defined by $F^{-1}=\sum_{n=0}^\infty ({\rm Id}-F)^n$, where the series in the right hand side converges in the strong operator sense, so that $F^{-1}$ satisfies (\ref{triang3}), and, if (\ref{triangm1}) holds, it satisfies (\ref{triang3m}) as well.

(iii) Continuous linear endomorphisms $F^*$ of $\R[[\z]]$ satisfying (\ref{triangF}) form a group, and the corresponding linear maps $F$ form a group as well.

(iv) The algebra automorphisms $F_x^*$, $x\in \R^4$ induce linear maps $F_x$ defined in (\ref{FdefT}), (\ref{FdefT1}) with $F$ replaced by $F_x$ and $F^*$ replaced by $F_x^*$. The maps $F_x$ are independent of the choice of the basis $f_i$, $i=1,\ldots, {\rm dim}(V)$ and satisfy (\ref{triang1}), (\ref{triang2}), (\ref{triang3}), and (\ref{triang3m}).
\el

\bpr
(i) Note that the expression in the right hand side of (\ref{FdefT1}) is well defined as
\be\lab{defF*}
\sum_{\gamma\in \Ml}(F^*(\z^\gamma))_\beta=\sum_{\tiny\begin{array}{c} \gamma\in \Ml: \\ |\gamma|\leq |\beta|+N_\beta\end{array}}F^\gamma_\beta \z^\beta,
\ee
and by Lemma \ref{locf} the sum in the right hand side of the last expression is finite. This observation together with (\ref{triangF'}) immediately imply all statements in part (i).

The first claim in part (ii) and (\ref{triang1}) follow from (i) with $N_\beta=0$, and (\ref{triangF}) also determines $(F)_{\beta}^\beta$ as described in (\ref{triang2}). Property (\ref{triang3m}) follows from (\ref{triangm1}).

Identity (\ref{triang3}) implies that $F-{\rm Id}$ is locally nilpotent, and the inverse to $F$ can be defined by $F^{-1}=\sum_{n=0}^\infty ({\rm Id}-F)^n$, where the series in the right hand side converges in the strong operator sense, so that $F^{-1}$ satisfies (\ref{triang3}), and, if (\ref{triangm1}) holds, $F^{-1}$ satisfies (\ref{triang3m}) as well. This property also implies that $F^*$ is invertible and that ${F^*}^{-1}$ is continuous by the first claim in part (ii) which is already proved.

Part (iii) immediately follows from (ii).

(iv) Since the algebra automorphisms $F_x^*$ satisfy (\ref{triangF}) and (\ref{triangm1}) (see (\ref{triang}) and (\ref{triangm})), by part (ii) they induce linear maps $F_x$ defined in (\ref{FdefT}), (\ref{FdefT1}), with $F$ replaced by $F_x$ and $F^*$ replaced by $F_x^*$. The maps $F_x$ satisfy (\ref{triang1}), (\ref{triang2}), (\ref{triang3}), and (\ref{triang3m}). By Remark \ref{indepbas} and by explicit formulas (\ref{Fdef}), (\ref{Fdef1}) the maps $F_x$ are also independent of the choice of the basis $f_i$, $i=1,\ldots, {\rm dim}(V)$.

\epr

Now we define an action of the automorphisms $F_x^*$ on the canonical lift ${\bf \Pi}$ to relate it to the recentered maps $\Pi_x$. We start by introducing some notation and by deriving some preliminary formulas.

For any $\Phi(\z)$ as in (\ref{Pb}) and any continuous linear endomorphism $F^*$ of the vector space $\R[[\z]]$ we define
\be\lab{FPdef'}
(F^*\Phi)(\z):= \sum_{\gamma\in \Ml}(F^*\Phi)_\gamma(\z^\gamma):=F^*(\sum_{\gamma\in \Ml}\Phi_\gamma(\z^\gamma)), 
\ee
where the coefficients of the formal power series in the last expression in the right hand side are obtained by applying $F^*$ to the formal variables $\z_\bg$, and $\z^i_{\n}, i=1,\ldots, {\rm dim}(V), {\n}\in \N^4$ in the expression $\sum_{\gamma\in \Ml}\Phi_\gamma(\z^\gamma)\in \O(V)[[\z_\bg,\{\z^i_{\n}, i=1,\ldots, {\rm dim}(V), {\n}\in \N^4\}]]$. By the definition of $F^*$ the expression in the right hand side of (\ref{FPdef'}) is well defined as a formal power series in the formal variables $\z_\bg$, $\z^i_{\n}, i=1,\ldots, {\rm dim}(V), {\n}\in \N^4$ with coefficients in the vector space of $V$--valued functions on $\R^4$.

If a continuous linear map $F^*:\R[[\z]]\to\R[[\z]]$ satisfies (\ref{triangF'}), which is the case for the automorphisms $F_x^*$ (see \ref{triang}), then by (\ref{FdefT1}) and by Lemma \ref{locf} definition (\ref{FPdef'}) can be rewritten in the form
\be\lab{FPdef}
(F^*\Phi)(\z)= \sum_{\gamma\in \Ml}(F^*\Phi)_\gamma(\z^\gamma)= \hskip 5cm
\ee
$$
\hskip 2cm =\sum_{\gamma\in \Ml}\Phi_\gamma(F^*\z^\gamma)\in \O(V)[[\z_\bg,\{\z^i_{\n}, i=1,\ldots, {\rm dim}(V), {\n}\in \N^4\}]], 
$$
where the coefficients of the formal power series in the last expression in the right hand side are obtained by evaluating $\Phi_{\gamma}\in \O_\gamma\simeq \O(V)\otimes W_\gamma^*\subset \O(V)\otimes \Sm^*$ on the coefficients of the series $F^*\z^\gamma\in \Sm[[\z]]$ using natural map (\ref{parO}).

Now assume that $F^*$ is an algebra endomorphism of $\R[[\z]]$.
Let $\gamma,\gamma_1,\gamma_2\in \Ml$ be such that $\gamma=\gamma_1+\gamma_2$. Using the multiplicative property $\z^\gamma=\z^{\gamma_1}\z^{\gamma_2}$ which follows from the definition of $\z^\beta$, $\beta\in \Ml$, definition  (\ref{Fmatr}), and the fact that $F^*$ is an algebra endomorphism we have for any $\beta\in \Ml$
\be\lab{Fmult}
(F)_{\beta}^\gamma \z^\beta=(F^*\z^\gamma)_\beta=(F^*(\z^{\gamma_1}\z^{\gamma_2}))_\beta=(F^*(\z^{\gamma_1})F^*(\z^{\gamma_2}))_\beta=
\ee
$$
=\sum_{\tiny\begin{array}{c}\beta_1,\beta_2\in \Ml: \\ \beta_1+\beta_2=\beta\end{array}}(F^*\z^{\gamma_1})_{\beta_1}(F^*\z^{\gamma_2})_{\beta_2}=\sum_{\tiny\begin{array}{c}\beta_1,\beta_2\in \Ml: \\ \beta_1+\beta_2=\beta\end{array}} (F)_{\beta_1}^{\gamma_1} (\z^{\beta_1})(F)_{\beta_2}^{\gamma_2} (\z^{\beta_2}),
$$
where in the second line we also used the multiplication rule in $\Sm[[\z]]$ in the form
$$
(ZW)_\beta=\sum_{\tiny\begin{array}{c}\beta_1,\beta_2\in \Ml: \\ \beta_1+\beta_2=\beta\end{array}}(Z)_{\beta_1}(W)_{\beta_2},~Z,W\in \Sm[[\z]],~\beta\in \Ml.
$$

Recalling linear isomorphisms (\ref{Wprod}) we can rewrite (\ref{Fmult}) as follows
\be\lab{Fsum}
(F)_{\beta}^\gamma=\sum_{\tiny\begin{array}{c}\beta_1,\beta_2\in \Ml: \\ \beta_1+\beta_2=\beta\end{array}} (F)_{\beta_1}^{\gamma_1} (F)_{\beta_2}^{\gamma_2},
\ee
where $(F)_{\beta_1}^{\gamma_1} (F)_{\beta_2}^{\gamma_2}$ is defined as the composition
$$
W_\beta\xrightarrow{\sim} W_{\beta_1}\otimes W_{\beta_2}\stackrel{F_{\beta_1}^{\gamma_1}\otimes F_{\beta_2}^{\gamma_2}}{\rightarrow} W_{\gamma_1} \otimes W_{\gamma_2}\xrightarrow{\sim} W_\gamma,
$$
and the first and the last maps in the composition are an isomorphism (\ref{Wprod}) and an inverse isomorphism, respectively.

 
Now we are in a position to state the main result of this subsection.
The relation (\ref{PPb}) in the proposition below is an analogue in the formal power series approach of a similar relation between the canonical lift and the recentered maps in the original Hairer's framework (see formula (8.26) in \cite{H1} or Proposition 2.1 in \cite{C2}).
\bp
(i) For all $x\in \R^4$, let 
\be\lab{xP}
\Pi_x(\z)=\sum_{\beta\in \Ml}\Pi_{x\beta}(\z^\beta),~\Pi_x^-(\z)=\sum_{\beta\in \Ml}\Pi_{x\beta}^-(\z^\beta),
\ee
be the formal power series defined as in (\ref{Pb}), and 
\be\lab{FxP}
(F_x^*{\bf \Pi})(\z)=\sum_{\beta\in \Ml}{\bf \Pi}_{\beta}(F_x^*\z^\beta),~(F_x^*{\bf \Pi}^-)(\z)=\sum_{\beta\in \Ml}{\bf \Pi}_{\beta}^-(F_x^*\z^\beta)
\ee
the formal power series defined as in (\ref{FPdef}). Then
\be\lab{Px}
\Pi_x(\z)=(F_x^*{\bf \Pi})(\z),~{\rm and } ~\Pi_x^-(\z)=(F_x^*{\bf \Pi}^-)(\z).
\ee

(ii) For all $x\in \R^4$ relations (\ref{Px}) are equivalent to
\be\lab{PPb}
\Pi_x={\bf \Pi}F_x,
\ee
where $\Pi_x,{\bf \Pi}:{\bf T}\to \O$ are defined by (\ref{PixdefbfT}) and (\ref{Pidef}), respectively, and we denote by the same symbol $F_x$ the natural extension of $F_x$ defined by (\ref{FdefT}),
\be\lab{Fxdef}
F_x:\bigoplus_{\beta\in \Ml}W_\beta\to \bigoplus_{\beta\in \Ml}W_\beta,
\ee
to an automorphism of 
\be\lab{TT}
{\bf T}=\bigoplus_{\beta\in \Ml}V_\beta^*\oplus \bigoplus_{\beta\in \Ml}{V_\beta^{-}}^*\simeq
(\bigoplus_{\beta\in \Ml}W_\beta\oplus \bigoplus_{\beta\in \Ml}W_\beta)\otimes V^*
\ee
acting as in (\ref{Fxdef}) on both summands $\bigoplus_{\beta\in \Ml}W_\beta$ in the last expression in (\ref{TT}), and as the identity transformation on the factor $V^*$ of the tensor product in that expression. 
\ep

\bpr
(i) Firstly, we rewrite (\ref{Px}) in terms of components using definitions (\ref{FxP}) and recalling (\ref{triang}),
\be\lab{Px1}
\Pi_{x\beta}(\z^\beta)=\sum_{\gamma\in \Ml:|\gamma|\leq |\beta|}{\bf \Pi}_{\gamma}((F_x^*\z^\gamma)_\beta)=(F_x^*{\bf \Pi})_{\beta}(\z^\beta),~\beta\in \Ml,
\ee
\be\lab{Px2}
\Pi_{x\beta}^-(\z^\beta)=\sum_{\gamma\in \Ml:|\gamma|\leq |\beta|}{\bf \Pi}_{\gamma}^-((F_x^*\z^\gamma)_\beta)=(F_x^*{\bf \Pi}^-)_{\beta}(\z^\beta),~\beta\in \Ml.
\ee

If $\beta=\delta_\n,~{\n}\in \N^4$ is purely polynomial then we have by (\ref{Pipol}), (\ref{Pibpol})
$$ 
\Pi_{x\delta_{\n}}(y)=(y-x)^{\n}{\rm Id}_{W_{\delta_\n}}^V,~{\bf \Pi}_{\delta_{\n}}(y)=y^{\n}{\rm Id}_{W_{\delta_\n}}^V.
$$

Thus from (\ref{Fdef2}) we obtain for any $y\in \R^4$ using the binomial formula that
$$
(F_x^*{\bf \Pi})_{\delta_{\n}}(y)(\z^{\delta_{\n}})=\sum_{\gamma\in \Ml:|\gamma|\leq |\n|}{\bf \Pi}_{\gamma}(y)((F_x^*\z^\gamma)_{\delta_{\n}})=
$$
$$
=\sum_{\m\in \N^4 :|\m|\leq |\n|}{\bf \Pi}_{\delta_\m}(y)((F_x^*\z^{\delta_{\m}})_{\delta_{\n}})=
$$
$$
=\sum_{\tiny\begin{array}{c} \m\in \N^4:\\ |\m|\leq |\n|\end{array}}{\n \choose \m}y^{\m}(-x)^{\n-\m}{\rm Id}_{W_{\delta_\n}}^V(\z^{\delta_\n})=(y-x)^{\n}{\rm Id}_{W_{\delta_\n}}^V(\z^{\delta_\n})=\Pi_{x\delta_{\n}}(y)(\z^{\delta_\n})
$$
which proves (\ref{Px1}) for all polynomial $\beta$.

Similarly, by (\ref{Pi-0}) and (\ref{BP0})
$$
\Pi_{x\delta_{\n}}^-={\bf \Pi}_{\delta_{\n}}^-=0.
$$
Thus from (\ref{Fdef2}) we obtain 
$$
(F_x^*{\bf \Pi}^-)_{\delta_{\n}}(\z^{\delta_{\n}})=\sum_{\gamma\in \Ml:|\gamma|\leq |\n|}{\bf \Pi}_{\gamma}^-((F_x^*\z^\gamma)_{\delta_{\n}})=
$$
$$
=\sum_{\m\in \N^4 :|\m|\leq |\n|}{\bf \Pi}_{\delta_\m}^-((F_x^*\z^{\delta_{\m}})_{\delta_{\n}})=0=\Pi_{x\delta_{\n}}^-(\z^{\delta_\n})
$$
which proves (\ref{Px2}) for all polynomial $\beta$.

For non--polynomial $\beta\in \Ml$ we shall prove (\ref{Px1}) and (\ref{Px2}) by induction over $|\beta|$.

For $\beta=0$, which is the only multi--index with $|\beta|=\alpha$, (\ref{Px1}) and (\ref{Px2}) trivially hold as in this case in (\ref{Px1}) and (\ref{Px2}) $(F_x^*\z^\gamma)_0=0$ for $\gamma\neq 0$ by (\ref{Fdef1}) and $(F_x^*1)_0=1$ by (\ref{Fdef}), so in the right hand sides of (\ref{Px1}) and (\ref{Px2}) only the terms with $\gamma=0$ are non--trivial, and by (\ref{Ph0}) these terms are equal to $\Pi_{x0}={\bf \Pi}_0=\K\xi^\rho$ and $\Pi_{x0}^-={\bf \Pi}_0^-=\xi^\rho$, respectively.

Now assume that for some $\eta\not\in \M_{pp}$, $\eta\neq 0$ identities (\ref{Px1}) and (\ref{Px2}) hold for all $\beta\in \Ml$ such that $|\beta|<|\eta|$. Then by (\ref{Pibpol}) and (\ref{Fdef2})
$$
(F_x^*{\bf \Pi})_{\eta}(y)(\z^\eta)=\sum_{\gamma\in \Ml:|\gamma|\leq |\eta|}{\bf \Pi}_{\gamma}(y)((F_x^*\z^\gamma)_\eta)=
$$
$$
=\sum_{\tiny\begin{array}{c}\gamma\in \Ml, \gamma\not\in \M_{pp}: \\ |\gamma|\leq |\eta|\end{array}}{\bf \Pi}_{\gamma}(y)((F_x^*\z^\gamma)_\eta)+\sum_{\n\in \N^4:|\n|\leq |\eta|}{\bf \Pi}_{\delta_\n}(y)((F_x^*\z^{\delta_\n})_\eta)=
$$
$$
=\sum_{\tiny\begin{array}{c}\gamma\in \Ml, \gamma\not\in \M_{pp}: \\ |\gamma|\leq |\eta|\end{array}}{\bf \Pi}_{\gamma}(y)((F_x^*\z^\gamma)_\eta)-
$$
$$
-\sum_{\n\in \N^4:|\n|\leq |\eta|}\sum_{\tiny\begin{array}{c} \m\in \N^4:\\ |\m|< |\eta|-|\n|\end{array}}\frac{y^\n(-x)^\m}{\n!\m!}\partial^{\n+\m}(\K\Pi_{x\eta}^-)(x)(\z^\eta)=
$$
$$
=\sum_{\tiny\begin{array}{c}\gamma\in \Ml, \gamma\not\in \M_{pp}: \\ |\gamma|\leq |\eta|\end{array}}{\bf \Pi}_{\gamma}(y)((F_x^*\z^\gamma)_\eta)-\sum_{\n\in \N^4:|\n|< |\eta|}\frac{(y-x)^\n}{\n!}\partial^{\n}(\K\Pi_{x\eta}^-)(x)(\z^\eta),
$$
where at the last step we also used the binomial formula.

Now we rewrite the first sum in the previous formula using (\ref{defPb}) and the fact that ${\bf \Pi}_{\gamma}^-=0$ for $\gamma\in \M_{pp}$ by (\ref{BP0}). This yields
\be\lab{F1}
(F_x^*{\bf \Pi})_{\eta}(y)(\z^\eta)=
\ee
$$
=\sum_{\tiny\begin{array}{c}\gamma\in \Ml, \gamma\not\in \M_{pp}: \\ |\gamma|\leq |\eta|\end{array}}\K{\bf \Pi}_{\gamma}^-(y)((F_x^*\z^\gamma)_\eta)-\sum_{\n\in \N^4:|\n|< |\eta|}\frac{(y-x)^\n}{\n!}\partial^{\n}(\K\Pi_{x\eta}^-)(x)(\z^\eta)=
$$
$$
=\sum_{\gamma\in \Ml: |\gamma|\leq |\eta|}\K{\bf \Pi}_{\gamma}^-(y)((F_x^*\z^\gamma)_\eta)-\sum_{\n\in \N^4:|\n|< |\eta|}\frac{(y-x)^\n}{\n!}\partial^{\n}(\K\Pi_{x\eta}^-)(x)(\z^\eta)=
$$
$$
=\K(F_x^*{\bf \Pi}^-)_{\eta}(y)(\z^\eta)-\sum_{\n\in \N^4:|\n|< |\eta|}\frac{(y-x)^\n}{\n!}\partial^{\n}(\K\Pi_{x\eta}^-)(x)(\z^\eta),
$$
where at the last step we also used the fact that $F_x^*$ acts on the formal variables $\z_\bg$, $\z^i_{\n}$, $i=1,\ldots, {\rm dim}(V)$, ${\n}\in \N^4$, and hence it commutes with the integral operator $\K$.

Now we show that (\ref{Px2}) holds for $\beta=\eta$, i.e. $(F_x^*{\bf \Pi}^-)_{\eta}=\Pi_{x\eta}^-$, and hence the first term in the right hand side of (\ref{F1}) coincides with $\K\Pi_{x\eta}^-(y)(\z^\eta)$. Indeed, evaluating (\ref{defP-}) for $\Phi_{\beta}^-={\bf \Pi}_{\beta}^-$, $\Phi_{\beta_i}={\bf \Pi}_{\beta_i}$, $i=1,2,3$ on $\z^\beta$, using the multiplicativity $\z^\beta=\z_\bg\z^{\beta_1}\z^{\beta_2}$ in the first sum, $\z^\beta=\z_\bg^2\z^{\beta_1}\z^{\beta_2}\z^{\beta_3}$ in the second sum, and $\z^\beta=\z_\bg^k\z^{\beta_1}$ in the last sum, taking the sum over all $\beta\neq 0$, and adding identity (\ref{Ph0}) in the form ${\bf \Pi}_0^-=\xi^\rho$ we obtain
$$
{\bf \Pi}^-(\z)=\z_\bg \A({\bf \Pi}(\z),{\bf \Pi}(\z))
+\z_\bg^2 \mB({\bf \Pi}(\z),{\bf \Pi}(\z),{\bf \Pi}(\z))+c(\z_\bg){\bf \Pi}(\z)+\xi^\rho,
$$
where $c(\z_\bg)=\sum_{k=1}^\infty c_k\z_\bg^k$, and the last sum is actually finite.

Applying $F_x^*$ to this identity and recalling (\ref{Fdef}) we deduce using that $F_x^*$ is an automorphism
$$
(F_x^*{\bf \Pi}^-)(\z)=\z_\bg \A((F_x^*{\bf \Pi})(\z),(F_x^*{\bf \Pi})(\z))
+
$$
$$
+\z_\bg^2 \mB((F_x^*{\bf \Pi})(\z),(F_x^*{\bf \Pi})(\z),(F_x^*{\bf \Pi})(\z))+c(\z_\bg)(F_x^*{\bf \Pi})(\z)+\xi^\rho.
$$

Consider the $\eta$--component of this identity,
\be\lab{F2}
(F_x^*{\bf \Pi}^-)_{\eta}=\sum_{\tiny \begin{array}{c}\beta_1,\beta_2\in \Ml: \\ \beta_1+\beta_2+\delta_\bg=\eta \end{array}}\A((F_x^*{\bf \Pi})_{\beta_1},(F_x^*{\bf \Pi})_{\beta_2})+
\ee
$$
+\sum_{\tiny \begin{array}{c}\beta_1,\beta_2,\beta_3\in \Ml: \\ \beta_1+\beta_2+\beta_3+2\delta_\bg=\eta\end{array}}\mB((F_x^*{\bf \Pi})_{\beta_1},(F_x^*{\bf \Pi})_{\beta_2},(F_x^*{\bf \Pi})_{\beta_3})+$$
$$
+\sum_{\tiny \begin{array}{c} k\in \{1,2,3,4\},\beta_1\in \Ml: \\ k\delta_\bg+\beta_1=\eta\end{array}}c_k(F_x^*{\bf \Pi})_{\beta_1}.
$$
By Lemma \ref{lF} (i) for all $\beta_i$, $i=1,2,3$ in the right hand side we have $|\beta_i|<|\eta|$, and hence by the induction assumption $(F_x^*{\bf \Pi})_{\beta_i}=\Pi_{x\beta_i}$, so (\ref{F2}) takes the form
$$
(F_x^*{\bf \Pi}^-)_{\eta}=\sum_{\tiny \begin{array}{c}\beta_1,\beta_2\in \Ml: \\ \beta_1+\beta_2+\delta_\bg=\eta \end{array}}\A(\Pi_{x\beta_1},\Pi_{x\beta_2})+
$$
$$
+\sum_{\tiny \begin{array}{c}\beta_1,\beta_2,\beta_3\in \Ml: \\ \beta_1+\beta_2+\beta_3+2\delta_\bg=\eta\end{array}}\mB(\Pi_{x\beta_1},\Pi_{x\beta_2},\Pi_{x\beta_3})+\sum_{\tiny \begin{array}{c} k\in \{1,2,3,4\},\beta_1\in \Ml: \\ k\delta_\bg+\beta_1=\eta\end{array}}c_k\Pi_{x\beta_1}=\Pi_{x\eta}^-,
$$
where at the last step we used definition (\ref{defP-}) in the case of $\Pi_{x\eta}^-$. This establishes (\ref{Px2}) for $\beta=\eta$, and hence (\ref{F1}) takes the form
$$
(F_x^*{\bf \Pi})_{\eta}(y)(\z^\eta)=
$$
$$
=\K\Pi_{x\eta}^-(y)(\z^\eta)-\sum_{\n\in \N^4:|\n|< |\eta|}\frac{(y-x)^\n}{\n!}\partial^{\n}(\K\Pi_{x\eta}^-)(x)(\z^\eta)=\Pi_{x\eta}(y)(\z^\eta),
$$
where at the last step we used (\ref{defP}).
This confirms (\ref{Px1}) for $\beta=\eta$ and completes the proof of the induction step.

(ii) 
Consider, for instance, the first relation in (\ref{Px}) and expand its left hand side and right hand side using (\ref{xP}) and (\ref{FxP}), respectively,
$$
\sum_{\beta\in \Ml}\Pi_{x\beta}(\z^\beta)=\sum_{\beta\in \Ml}{\bf \Pi}_{\beta}(F_x^*\z^\beta).
$$

By the first relation in (\ref{defPP}) the last displayed identity takes the form
\be\lab{F*}
\sum_{\beta\in \Ml}(\Pi_x\otimes {\rm Id}_{V_\beta})(\tau_\beta)(\z^\beta)=\sum_{\beta\in \Ml}({\bf \Pi}\otimes {\rm Id}_{V_\beta})(\tau_\beta)(F_x^*\z^\beta),
\ee
or, as ${\rm Id}_{V_\beta}={\rm Id}_{W_\beta^*}\otimes {\rm Id}_{V}$ is an identity transformation and $\tau_\beta={\rm Id}_{V_\beta^*}={\rm Id}_{W_\beta}\otimes {\rm Id}_{V^*}$,
\be\lab{F3}
\sum_{\beta\in \Ml}(\Pi_x\otimes {\rm Id}_{V_\beta})(\tau_\beta)(\z^\beta)=({\bf \Pi}\otimes {\rm Id}_{V})\sum_{\beta\in \Ml}(F_x^*\z^\beta\otimes {\rm Id}_{V^*})=
\ee
$$
=({\bf \Pi}\otimes {\rm Id}_{V})(\sum_{\beta\in \Ml}F_x^*\z^\beta)\otimes {\rm Id}_{V^*},
$$
where 
$$
F_x^*\z^\beta\otimes {\rm Id}_{V^*}\in W_\beta\otimes V^*\otimes V[[\z_\bg,\{\z^j_{{\n}}:j=1,\ldots, {\rm dim}(V),\n\in \N^4\}]]\simeq 
$$
$$
\simeq V_\beta^*\otimes V[[\z_\bg,\{\z^j_{{\n}}:j=1,\ldots, {\rm dim}(V),\n\in \N^4\}]].
$$

Recalling (\ref{FdefT}) and (\ref{FdefT1}) and the definition of the extension $F_x:{\bf T}\to {\bf T}$ in the statement of part (ii) of this proposition we can reduce relation (\ref{F3}) to
\be\lab{F4}
\sum_{\beta\in \Ml}(\Pi_x\otimes {\rm Id}_{V_\beta})(\tau_\beta)(\z^\beta)=({\bf \Pi}F_x\otimes {\rm Id}_{V})\sum_{\beta\in \Ml}(\z^\beta\otimes {\rm Id}_{V^*}).
\ee

Applying the arguments which led from (\ref{F*}) to (\ref{F3}) in the opposite order we infer
$$
\sum_{\beta\in \Ml}(\Pi_x\otimes {\rm Id}_{V_\beta})(\tau_\beta)(\z^\beta)=\sum_{\beta\in \Ml}({\bf \Pi}F_x\otimes {\rm Id}_{V_\beta})(\tau_\beta)(\z^\beta).
$$

Similarly one can obtain
$$
\sum_{\beta\in \Ml}(\Pi_x\otimes {\rm Id}_{V_\beta^-})(\tau_\beta^-)(\z^\beta)=\sum_{\beta\in \Ml}({\bf \Pi}F_x\otimes {\rm Id}_{V_\beta^-})(\tau_\beta^-)(\z^\beta).
$$
The last two relations are obviously equivalent to (\ref{PPb}). This completes the proof of part (ii).

\epr


\subsection{Structure group automorphisms and the structure group}\lab{StrG}

\setcounter{equation}{0}
\setcounter{theorem}{0}

Now we are in a position to define the last component of our smooth model, the structure group maps. 

Following \cite{H}, for any $x,y\in \R^4$ we define structure group automorphisms $\G_{xy}$ by
\be\lab{Gdef}
\G_{xy}=F_x^{-1}F_y:\bigoplus_{\beta\in \Ml}W_\beta\to \bigoplus_{\beta\in \Ml}W_\beta,
\ee
and denote by the same symbols their extensions to automorphisms of $\bf T$
acting as in (\ref{Gdef}) on both summands $\bigoplus_{\beta\in \Ml}W_\beta$ in the last expression in (\ref{TT}), and as the identity transformation on the factor $V^*$ of the tensor product in that expression. Note that by part (iv) of Lemma \ref{Flemc} the definition of $\G_{xy}$ does not depend on the choice of the basis $f_i$, $i=1,\ldots, {\rm dim}(V)$.

Using (\ref{FdefT1}) with $F=\G_{xy}$ one can define automorphisms $\G_{xy}^*$ of the algebra $\Sm[[\z]]$ which are actually induced by automorphisms of the subalgebra $\R[[\z]]\subset \Sm[[\z]]$ by the definition of $F_x$ and of $\G_{xy}=F_x^{-1}F_y$.

The automorphisms $\G_{xy}$ have the following properties which are similar to properties of analogous automorphisms in \cite{BOT} in the case of a single equation for real--valued functions.
\bp\lab{Gprop}
(i) For any $x,y\in \R^4$ and any $\beta\in \Ml$ one has
\be\lab{triangG}
\G_{xy}:W_\beta\to W_\beta\oplus \bigoplus_{\gamma \in \Ml:|\gamma|< |\beta|}W_\gamma,
\ee
and 
\be\lab{triangG1}
(\G_{xy})_{\beta}^\beta={\rm Id}_{W_\beta},
\ee
i.e
$$
(\G_{xy}-{\rm Id}_{W_\beta}):W_\beta\to \bigoplus_{\gamma \in \Ml:|\gamma|< |\beta|}W_\gamma.
$$

Also,
\be\lab{triangGm}
(\G_{xy}-{\rm Id}_{W_\beta}):W_\beta\to \bigoplus_{\gamma \in \Ml:|\gamma|_\prec< |\beta|_\prec}W_\gamma.
\ee

(ii) For any $x,y,z\in \R^4$ one has
\be
\G_{xy}\G_{yz}=\G_{xz},
\ee
and $\G_{xx}$ is the identity transformation of $\bigoplus_{\beta\in \Ml}W_\beta$, so that the automorphisms $\G_{xy}$ form a group.

(iii) For any $x,y\in \R^4$ one has
\be\lab{PG}
\Pi_y=\Pi_x\G_{xy},~\Pi_y(\z)=(\G_{xy}^*\Pi_x)(\z),~\Pi_y^-(\z)=(\G_{xy}^*\Pi_x^-)(\z)
\ee

(iv) For any $x,y\in \R^4$
\be\lab{G1g}
\G_{xy}^*1=1,~\G_{xy}^*\z_\bg=\z_\bg,
\ee
and
\be\lab{Gn}
(\G_{xy})_{\delta_\m}^\gamma=\left\{\begin{array}{cc} {\m \choose \n}(x-y)^{\m-\n}{\rm Id}_{W_{\delta_{\m}}}^{W_{\delta_{\n}}} & {\rm if }~ \gamma =\delta_\n,~\n\in \N^4,~n_i\leq m_i,~i=1,2,3,4 \\ 0 & {\rm else} \end{array}\right. .
\ee

(v) For any $\gamma\in \M_{\geq 0}\cup \M_{pp}$ and any $\beta\not\in \M_{\geq 0}\cup \M_{pp}$ one has $(\G_{xy})_{\beta}^\gamma =0$,
\be\lab{tr1}
\gamma\in \M_{\geq 0}\cup \M_{pp},~\beta\not\in \M_{\geq 0}\cup \M_{pp}\Longrightarrow (\G_{xy})_{\beta}^\gamma =0,
\ee
and for any $\gamma\in \M'\setminus \M_{\geq 0}$ and any $\beta\not\in \M'$ one has $(\G_{xy})_{\beta}^\gamma =0$,
\be\lab{tr2}
\gamma\in \M'\setminus \M_{\geq 0},~\beta\not\in \M'\Longrightarrow(\G_{xy})_{\beta}^\gamma =0.
\ee
In particular, for any $\gamma\in \M$ and any $\beta\not\in \M$ one has $(\G_{xy})_{\beta}^\gamma =0$,
\be\lab{tr3}
\gamma\in \M,~\beta\not\in \M\Longrightarrow (\G_{xy})_{\beta}^\gamma =0.
\ee

Moreover, if $\beta, \gamma \in \M'\setminus \M_{\geq 0}$ are such that $(\G_{xy})_\beta^\gamma\neq 0$ then $\m:=\sum_{i=0}^\infty\gamma(\n_i)\n_i$, ${\bf k}:=\sum_{i=0}^\infty\beta(\n_i)\n_i\in \N^4$ satisfy $m_i\leq k_i$ for $i=0,1,2,3$, $W_\beta$ is naturally isomorphic to $W_\gamma$, and
\be\lab{m'}
(\G_{xy})_\beta^\gamma=A_\beta^\gamma (x-y)^{{\bf k}-\m}{\rm Id}_{W_\beta}^{W_\gamma},
\ee
where $A_\beta^\gamma\neq 0$ is a combinatorial numeric factor.

(vi) For every $\beta\in \Ml$ fix the norm on $W_\beta$ as in Proposition \ref{Pibond}. Denote by $|\cdot |$ the induced operator norms on ${\rm Hom}(W_\beta,W_\gamma)$, $\gamma,\beta\in \Ml$. Then for every compact set $\mK\subset \R^4$, every $\rho>0$, and every $\beta\in \Ml$ and $\gamma\in \M$ there exists a constant $D_{\mK,\beta}^\gamma(\rho)\geq 0$ such that
\be\lab{Gbond}
|(\G_{xy})_\beta^\gamma|\leq D_{\mK,\beta}^\gamma(\rho)|x-y|^{|\beta|-|\gamma|},~x,y\in \mK,~x\neq y.
\ee

(vii) Let 
$$
P_{W}:\bigoplus_{\beta\in \Ml}W_\beta\to \bigoplus_{\beta\in \M}W_\beta=:W,
$$
and
$$
P_T:{\bf T}=\bigoplus_{\beta\in \Ml}V_\beta^*\oplus \bigoplus_{\beta\in \Ml}{V_\beta^{-}}^*\to \bigoplus_{\beta\in \M}V_\beta^*\oplus \bigoplus_{\beta\in \M}{V_\beta^{-}}^*=T
$$
be the natural projection operators. Then for $x,y\in \R^4$ the linear maps $\G_{xy}^{W}:=P_{W}\G_{xy}P_{W}$ (resp. $\G_{xy}^T:=P_T\G_{xy}P_T$ for $\G_{xy}$ acting on $\bf T$) are automorphisms of $W$ (resp. $T$) which form a group called the structure group of the regularity structure $T$, and for any $x,y,z\in \R^4$ one has
\be\lab{strGG}
\G_{xy}^{W}\G_{yz}^{W}=\G_{xz}^{W},~
\G_{xy}^T\G_{yz}^T=\G_{xz}^T.
\ee

\ep

\bpr
(i) Properties (\ref{triangG}), (\ref{triangG1}) and (\ref{triangGm}) follow from the definition of $\G_{xy}$ and from (\ref{triang1}), (\ref{triang2}) and (\ref{triang3m}), respectively, for $F_x$ and $F_x^{-1}$, $x\in \R^4$.   

(ii) We have by the definition of $\G_{xy}$
$$
\G_{xy}\G_{yz}=F_x^{-1}F_yF_y^{-1}F_z=F_x^{-1}F_z=\G_{xz}.
$$

Obviously, $\G_{xx}$ is the identity transformation by the definition.

(iii) By (\ref{PPb}) we have
$$
{\bf \Pi}=\Pi_xF_x^{-1}=\Pi_yF_y^{-1},
$$
which is equivalent to the first relation in (\ref{PG}). The last two relations in (\ref{PG}) follow from (\ref{Px}) in a similar way.

(iv) (\ref{G1g}) follows immediately from (\ref{Fdef}), and by (\ref{Fdg})
\be\lab{Fdg-1}
(F_x^{-1})_{\delta_\m}^\gamma=\left\{\begin{array}{cc} {\m \choose \n}x^{\m-\n}{\rm Id}_{W_{\delta_{\m}}}^{W_{\delta_{\n}}} & {\rm if }~ \gamma =\delta_\n,~\n\in \N^4 \\ 0 & {\rm else} \end{array}\right. .
\ee

Together with (\ref{Fdg}) this yields
$$
(\G_{xy})_{\delta_\m}^\gamma=\sum_{{\bf k}\in \N^4}(F_x^{-1})_{\delta_{\bf k}}^\gamma(F_y)_{\delta_\m}^{\delta_{\bf k}}.
$$
If $\gamma\not\in \M_{pp}$ the right hand side of the last expression is zero by (\ref{Fdg-1}), and if $\gamma=\delta_\n$, $\n\in \N^4$ we get using (\ref{Fdg}) and (\ref{Fdg-1}) 
$$
(\G_{xy})_{\delta_\m}^{\delta_\n}=\sum_{{\bf k}\in \N^4}(F_x^{-1})_{\delta_{\bf k}}^{\delta_\n}(F_y)_{\delta_\m}^{\delta_{\bf k}}=
$$
$$
=\hspace{-0.5cm}\sum_{\tiny \begin{array}{c} {\bf k}\in \N^4: \\ |\n|\leq |{\bf k}|\leq |\m|\end{array}}\hspace{-0.5cm}{{\bf k} \choose \n}x^{{\bf k}-\n}{\m \choose {\bf k}}(-y)^{\m-{\bf k}}{\rm Id}_{W_{\delta_{\m}}}^{W_{\delta_{\n}}}.
$$

Observing that in the last expression
$$
{{\bf k} \choose \n}{\m \choose {\bf k}}={\m \choose \n}{\m-\n \choose {\bf k}-\n}
$$
and introducing the new variable of summation ${\bf l}:={\bf k}-\n$ we obtain
$$
(\G_{xy})_{\delta_\m}^{\delta_\n}={\m \choose \n}\hspace{-0.5cm}\sum_{\tiny \begin{array}{c} {\bf l}\in \N^4: \\ 0\leq |{\bf l}|\leq |\m-\n|\end{array}}\hspace{-0.5cm}{\m-\n \choose {\bf l}}x^{\bf l}(-y)^{\m-\n-{\bf l}}{\rm Id}_{W_{\delta_{\m}}}^{W_{\delta_{\n}}}=
$$
$$
={\m \choose \n}(x-y)^{\m-\n}{\rm Id}_{W_{\delta_{\m}}}^{W_{\delta_{\n}}},
$$
where at the last step we used the binomial formula. This completes the proof of (\ref{Gn}).

(v) and (vi) Fix a compact subset $\mK\subset \R^4$. Property (\ref{tr3}) follows from (\ref{tr1}) and (\ref{tr2}). 

If $\gamma=0$ then $[\gamma]=0$, and $\gamma\in \M_{\geq 0}$. Thus property (\ref{tr2}) is void, and (\ref{tr1}) follows from (\ref{Fmatr}) with $F=\G_{xy}$ and from the first formula in (\ref{G1g}) which imply
$$
(\G_{xy})_{\beta}^0=\left\{\begin{array}{cc} {\rm Id}_{W_0}  & {\rm if }~ \beta =0  \\ 0 & {\rm else} \end{array}\right. .
$$

The last formula also implies (\ref{Gbond}) in the case when $\gamma=0$ and $\beta\in \Ml$ is arbitrary.

We prove (\ref{tr1}), (\ref{tr2}) and (\ref{Gbond}) by induction over $|\beta|$. We follow the arguments in \S 1.18 of \cite{BOT}.

$\beta=0$ is the only possible for which $|\beta|=\alpha$ is minimal possible. In this case $\beta\in \M_{\geq 0}\subset \M'\subset \M$, so that properties (\ref{tr1}) and (\ref{tr2}) are void.

Now assume that for some $\beta\in \Ml$ the statements of (v) and (vi) are true for all $\gamma\in \Ml$ and all $\beta'\in \Ml$ with $|\beta'|<|\beta|$. As we already observed in the beginning of the proof of parts (v) and (vi), the statements are true for $\gamma=0$ and for any $\beta$. Therefore we can assume that $\gamma=l\delta_\bg+\sum_{i=1}^p\delta_{{\m}_i}\neq 0$, $\m_i\in \N^4$.

\bl\lab{Gbgmult}
For any multi--index $\gamma=l\delta_\bg+\sum_{i=1}^p\delta_{{\m}_i}\neq 0$, $\m_i\in \N^4$, $l\in \N$ and any algebra endomorphism $F^*$ of $\R[[\z]]$ such that $F^*\z_\bg=\z_\bg$ one has
\be\lab{Gb}
(F)_{\beta}^\gamma= \hspace{-1cm}  \sum_{\tiny\begin{array}{c}\beta_1,\ldots,\beta_p\in \Ml: \\ \beta_1+\ldots+\beta_p+l\delta_\bg=\beta\end{array}} \hspace{-1cm}  (F)_{\beta_1}^{\delta_{\m_1}}\ldots (F)_{\beta_p}^{\delta_{\m_p}},
\ee
where we assume that an empty sum is zero and an empty product is 1. 

Moreover, if $p>1$ or $l>0$, and $F^*$ satisfies (\ref{triang1}) then in the non--zero terms in the sum in (\ref{Gb}) $|\beta_i|<|\beta|$, $i=1,\ldots,p$.
\el

\bpr
Firstly, iterating (\ref{Fsum}) and using $F^*\z_\bg=\z_\bg$ yields (\ref{Gb}).

Since $(F)_0^{\delta_{\m_i}}=0$, $i=1,\ldots,p$ by (\ref{triang1}) as $|\delta_{\m_i}|=|\m_i|\geq 0>|0|=\alpha$, we may assume that $\beta_i\neq 0$,  $i=1,\ldots,p$ in (\ref{Gb}).

From definition (\ref{hom}) it follows that $|\cdot|-\alpha$ is a linear function on $\Ml$ achieving its minimum equal to zero at $\beta=0$, and hence for each non--zero term in (\ref{Gb}) we have
\be\lab{bs}
|\beta|-\alpha=\sum_{i=1}^p(|\beta_i|-\alpha)+|l\delta_\bg|-\alpha,~|\beta_i|-\alpha\geq 0,~i=1,\ldots, p.
\ee

Recall that $\beta_i\neq 0$, $i=1,\ldots,p$, and hence by Lemma \ref{locf} $|\beta_i|-\alpha>0$, $i=1,\ldots,p$. Observe also that $|l\delta_\bg|-\alpha=l(\alpha+1)\geq 0$. Therefore, if $p>1$ or $l>0$, by (\ref{bs}) $|\beta_i|-\alpha<|\beta|-\alpha$, $i=1,\ldots,p$, or $|\beta_i|<|\beta|$, $i=1,\ldots,p$.

\epr

Note that by part (i) of this proposition Lemma \ref{Gbgmult} holds for $F=\G_{xy}$.
Now if $\gamma=l\delta_\bg+\sum_{i=1}^p\delta_{{\m}_i}\neq 0$, $\m_i\in \N^4$, $\gamma\in \M'$ then $l\neq 0$, and hence by the previous lemma for all $\beta_i$ in the right hand side of (\ref{Gb}) with $F=\G_{xy}$ one has $|\beta_i|<|\beta|$. So by the induction assumption $(\G_{xy})_{\beta_i}^{\delta_{\m_i}}\neq 0$ only if $\beta_i\in \M_{\geq 0}\cup\M_{pp}$, and estimates (\ref{Gbond}) hold for all such $(\G_{xy})_{\beta_i}^{\delta_{\m_i}}$. In all these cases 
\be\lab{bg0}
\beta_i\in \M_{\geq 0}\cup\M_{pp},~[\beta_i]\geq -1,~ i=1,\ldots, p, 
\ee
and the triangle inequality along with estimates (\ref{Gbond}) for $(\G_{xy})_{\beta_i}^{\delta_{\m_i}}$ in each term of sum (\ref{Gb}) with $F=\G_{xy}$ yield for $x,y\in \mK$
$$
|(\G_{xy})_{\beta}^\gamma|\leq C_\beta^\gamma \hspace{-1.2cm} \sum_{\tiny\begin{array}{c}\beta_1,\ldots,\beta_p\in \Ml: \\ \beta_1+\ldots+\beta_p+l\delta_\bg=\beta\end{array}}\hspace{-1.2cm} |(\G_{xy})_{\beta_1}^{\delta_{\m_1}}|\ldots |(\G_{xy})_{\beta_2}^{\delta_{\m_p}}|\leq 
$$
$$
\leq D_{\mK,\beta}^\gamma(\rho)|x-y|^{\sum_{i=1}^p(|\beta_i|-|\m_i|)}=D_{\mK,\beta}^\gamma(\rho)|x-y|^{|\beta|-|\gamma|}
$$
for some $C_\beta^\gamma>0$ and $D_{\mK,\beta}^\gamma(\rho)\geq 0$ as by (\ref{bs}) 
$$
|\beta|=\sum_{i=1}^p|\beta_i|+|l\delta_\bg|-p\alpha,
$$
and similarly from (\ref{hom})
$$
|\gamma|=\sum_{i=1}^p|\delta_{\m_i}|+|l\delta_\bg|-p\alpha=\sum_{i=1}^p|\m_i|+|l\delta_\bg|-p\alpha,
$$
so 
$$
\sum_{i=1}^p(|\beta_i|-|\m_i|)=|\beta|-|\gamma|.
$$
This proves the induction step for (vi) in the case when $\gamma\not\in \M_{pp}$.

Note also that $[\gamma]=l-p$ by the definition and that by the linearity of the function $[\cdot]$ one has 
\be\lab{pbe}
[\beta]=l-p+\sum_{i=1}^p([\beta_i]+1),
\ee
where $\beta_i$, $i=1,\ldots, p$ are introduced in (\ref{Gb}).
Therefore if $\gamma\in \M_{\geq 0}$ then $l-p\geq 0$, and this fact together with (\ref{bg0}) imply that all terms in the right hand side of (\ref{pbe}) are non--negative, so $[\beta]\geq 0$, and $\beta\in \M_{\geq 0}$. This establishes the induction step for (v) in the case when $\gamma\in \M_{\geq 0}$.

If $\gamma \in \M'\setminus \M_{\geq 0}$ then by the definition of $\M'$
\be\lab{g1}
\gamma=\delta_\bg+\delta_{\m_1}+\delta_{\m_2},~\m_1,\m_2\in \N^4,
\ee
or
\be\lab{g2}
\gamma=2\delta_\bg+\delta_{\m_1}+\delta_{\m_2}+\delta_{\m_3},~\m_1,\m_2,\m_3\in \N^4.
\ee

In both cases $[\gamma]=l-p=-1$. Therefore by (\ref{pbe}) and by (\ref{bg0}) either $[\beta]\geq 0$, i.e. $\beta\in \M_{\geq 0}$, or $[\beta]=-1$ if $[\beta_i]=-1$ for all $i$. In the latter case by (\ref{bg0}) $\beta_i=\delta_{{\bf k}_i}$ for some ${\bf k}_i\in \N^4$ and for all $i$. So that for $\gamma$ as in (\ref{g1}) or in (\ref{g2})
$$
\beta=\delta_\bg+\delta_{{\bf k}_1}+\delta_{{\bf k}_2},~{\bf k}_1,{\bf k}_2\in \N^4,
$$
or
$$
\beta=2\delta_\bg+\delta_{{\bf k}_1}+\delta_{{\bf k}_2}+\delta_{{\bf k}_3},~{\bf k}_1,{\bf k}_2,{\bf k}_3\in \N^4,
$$
respectively. In both cases $\beta \in \M'\setminus \M_{\geq 0}$, and, as $(\G_{xy})_{\delta_{{\bf k}_i}}^{\delta_{\m_i}}={{\bf k}_i \choose \m_i}(x-y)^{{\bf k}_i-\m_i}{\rm Id}_{W_{\delta_{{\bf k}_i}}}^{W_{\delta_{\m_i}}}$, we obtain in this case that by (\ref{Gb}) with $F=\G_{xy}$ that $W_\beta\simeq W_\gamma$, and 
$$
(\G_{xy})_\beta^\gamma=A_\beta^\gamma (x-y)^{{\bf k}-\m}{\rm Id}_{W_\beta}^{W_\gamma},~\m:=\sum_{i=0}^\infty\gamma(\n_i)\n_i,~{\bf k}:=\sum_{i=0}^\infty\beta(\n_i)\n_i\in \N^4,
$$
where $\m$ and $\bf k$ satisfy $m_i\leq  k_i$ for $i=0,1,2,3$, and $A_\beta^\gamma\neq 0$ is a combinatorial numeric factor. This proves (\ref{m'}). 

We also established the induction step for (v) in the case when $\gamma \in \M'\setminus \M_{\geq 0}$.

It remains to justify the induction step for (v) and (vi) in the case when $\gamma\in \M_{pp}$. For this purpose we rewrite the second relation in (\ref{PG}) in terms of components using (\ref{FdefT1}), (\ref{Fmatr}) and (\ref{FPdef}). This yields
$$
\Pi_{y\beta}(w)=\sum_{\eta\in \Ml}\Pi_{x\eta}(w)\circ (\G_{xy})_\beta^\eta, w\in \R^4,
$$
where $\circ$ stands for the composition of the linear maps $\Pi_{x\eta}(w):W_\eta\to V$ and $(\G_{xy})_\beta^\eta:W_\beta\to W_\eta$ defined for each $x,y,w\in \R^4$. 

Since by (\ref{Pi0}) for any $\delta\in \Ml$ one has $\Pi_{x\delta}=0$ for $\delta\not\in \M_{\geq 0}\cup \M_{pp}$ and by part (i) of this proposition $(\G_{xy})_\beta^\eta=0$ for $|\eta|> |\beta|$, we obtain from the last formula
\be\lab{pin}
\Pi_{y\beta}(w)=\sum_{\tiny \begin{array}{c} \eta\in \M_{\geq 0}: \\ |\eta|\leq |\beta|\end{array}}\Pi_{x\eta}(w)\circ (\G_{xy})_\beta^\eta+\sum_{\m\in \N^4}\Pi_{x\delta_\m}(w)\circ (\G_{xy})_\beta^{\delta_\m}=
\ee
$$
=\sum_{\tiny \begin{array}{c} \eta\in \M_{\geq 0}: \\ |\eta|\leq |\beta|\end{array}}\Pi_{x\eta}(w)\circ (\G_{xy})_\beta^\eta+\sum_{\m\in \N^4}(w-x)^{\m}{\rm Id}_{W_{\delta_\m}}^V\circ (\G_{xy})_\beta^{\delta_\m},
$$
where at the last step we used (\ref{Pipol}).

By the induction step for (v) proved for $\eta\in \M_{\geq 0}$ we have $(\G_{xy})_\beta^\eta=0$ for $\beta\not\in \M_{\geq 0}$. Also, by (\ref{Pi0}) we obtain $\Pi_{x\beta}=0$ for $\beta\not\in \M_{\geq 0}\cup \M_{pp}$. Thus by (\ref{pin}) in this case
$$
\sum_{\m\in \N^4}(w-x)^{\m}{\rm Id}_{W_{\delta_\m}}^V\circ (\G_{xy})_\beta^{\delta_\m}=0
$$
if $\beta\not\in \M_{\geq 0}\cup \M_{pp}$, and, since the monomials $(w-x)^{\m}$ are linearly independent for different $\m$, we deduce that $(\G_{xy})_\beta^{\delta_\m}=0$ for all $\m\in \N^4$ if $\beta\not\in \M_{\geq 0}\cup \M_{pp}$. This establishes the induction step for (v) and for $\gamma\in \M_{pp}$ and completes the proof of (v). 

To justify the induction step for (vi) in the case when $\gamma=\delta_\n$, $\n\in \N^4$ we apply $\partial^\n$ with respect to $w$ to (\ref{pin}) and evaluate the obtained expression at $w=x$. After rearranging terms and dividing by $\n!$ we obtain 
$$
{\rm Id}_{W_{\delta_\n}}^V\circ (\G_{xy})_\beta^{\delta_\n}=\frac{1}{\n!}\partial^\n(\Pi_{y\beta})(x)
-\frac{1}{\n!}\sum_{\tiny \begin{array}{c}\eta\in \M_{\geq 0}: \\ |\eta|\leq |\beta|\end{array}}\partial^\n(\Pi_{x\eta})(x)\circ (\G_{xy})_\beta^\eta.
$$

Note that by (\ref{Pvan}) one has $\partial^\n(\Pi_{x\eta})(x)=0$ for any $\eta\in \Ml$, $\n\in \N^4$ such that $|\n|<|\eta|$, and hence the previous expression takes the form
$$
{\rm Id}_{W_{\delta_\n}}^V\circ (\G_{xy})_\beta^{\delta_\n}=\frac{1}{\n!}\partial^\n(\Pi_{y\beta})(x)
-\frac{1}{\n!}\sum_{\tiny \begin{array}{c}\eta\in \M_{\geq 0}: \\ |\eta|\leq \min{\{|\n|, |\beta|\}}\end{array}}\partial^\n(\Pi_{x\eta})(x)\circ (\G_{xy})_\beta^\eta.
$$

By Lemma (\ref{locf}) the sum in the right hand side of the last expression is finite. Hence, using the triangle inequality in the right hand side of it, estimate (\ref{Pibnd1}) for the first term in the right hand side and estimates (\ref{Gbond}) for $|(\G_{xy})_\beta^\eta|$ with $\eta \in \M'$, $|\eta|\leq |\beta|$ which hold by the induction assumption and by the induction step already proved for $\beta$ and any $\eta\in \M'\supset \M_{\geq 0}$, we obtain 
\be\lab{Gbnd}
|(\G_{xy})_\beta^{\delta_\n}|\leq \frac{1}{\n!}|\partial^\n(\Pi_{y\beta})(x)|
+\frac{1}{\n!}\sum_{\tiny \begin{array}{c}\eta\in \M_{\geq 0}: \\ |\eta|\leq \min{\{|\n|,|\beta|\}}\end{array}}|\partial^\n(\Pi_{x\eta})(x)||(\G_{xy})_\beta^\eta|\leq
\ee
$$
\leq \frac{1}{\n!}|x-y|^{|\beta|-|\n|}+\frac{1}{\n!}\sum_{\tiny \begin{array}{c}\eta\in \M_{\geq 0}: \\ |\eta|\leq \min{\{|\n|,|\beta|\}}\end{array}}|\partial^\n(\Pi_{x\eta})(x)|D_{\mK,\beta}^\eta(\rho)|x-y|^{|\beta|-|\eta|}, ~x,y\in \mK.
$$

Since $\Pi_{x\eta}\in \O_\eta=\O\otimes V_\eta\subset C^\infty(\R^4)\otimes V_\eta$, their partial derivatives $\partial^\n(\Pi_{x\eta})(x)$ are bounded on the compact subset $\mK\subset\R^4$, and for $|\eta|\leq |\n|$ one has $|x-y|^{|\beta|-|\eta|}\leq C_\mK |x-y|^{|\beta|-|\n|}$ for any $x,y\in \mK$ and some $C_\mK \geq 0$. Therefore (\ref{Gbnd}) implies (\ref{Gbond}) for $\gamma=\delta_\n$ as by (\ref{hom}) $|\n|=|\delta_\n|$ in (\ref{Gbnd}). This establishes the induction step for (vi) in the case when $\gamma\in \M_{pp}$ and completes the proof of part (vi).

(vii) By (\ref{tr3}) $P_{W}\G_{xy}P_{W}=P_{W}\G_{xy}$, and hence by part (ii) of this proposition
$$
\G_{xy}^{W}\G_{yz}^{W}=P_{W}\G_{xy}P_{W} P_{W}\G_{yz}P_{W}=P_{W}\G_{xy}\G_{yz}P_{W}=P_{W}\G_{xz}P_{W}=\G_{xz}^{W}.
$$
This proves the first relation in (\ref{strGG}) and together with part (ii) of this proposition justifies that the linear maps $\G_{xy}^{W}$, $x,y\in \R^4$ are automorphisms of $W$ which form a group.

A similar claim for the linear maps $\G_{xy}^T$, $x,y\in \R^4$ and the second relation in (\ref{strGG}) immediately follow from the definition of the maps $\G_{xy}^T$, $x,y\in \R^4$.

\epr

In \cite{BOT,LOTT} arguments based on formulas (\ref{Gb}) and (\ref{pin}) are called the algebraic and the three--point argument, respectively. They are common technical proof tools in the multi--index approach to regularity structures.


\subsection{Modelled distributions and the lift of the Langevin equation to an equation for modelled distributions}\lab{Moddeq}

\setcounter{equation}{0}
\setcounter{theorem}{0}

Estimates (\ref{Gbond}) are similar to those obtained in Proposition 8.27 in \cite{H1}. Together with parts (i), (iii), (vii) of Proposition \ref{Gprop}, (\ref{Pipol}), (\ref{Pi-0}), (\ref{Pi0}) and Proposition (\ref{Pibond}) they imply that $(\Pi_x,\G_{xy}^T)$ is a smooth model for the regularity structure $T$ in the sense of Hairer (see \cite{H1}, Definition 2.17).

Estimates (\ref{Pibnd1}), (\ref{Pibnd2}) and (\ref{Gbond}) are rather provisional as well as the definition of modelled distributions given in this section. Actually in Section \ref{Stohe} we shall obtain stochastic estimates which hold uniformly for all $\rho >0$, and in Section \ref{conv} we show that the models $(\Pi_x,\G_{xy}^T)$ converge in a certain sense when $\rho\to 0$ in such a way that the estimates still hold after taking the limit. The proper definition of modelled distributions should be then given using the limit of the sequence of the smooth models for which global pointwise estimates similar to (\ref{Pibnd1}), (\ref{Pibnd2}) and (\ref{Gbond}) are obtained in Section \ref{ptwise}. 

The purpose of this section is just to show that the lift of equation (\ref{LangHcrm}) to an equation for modelled distributions is reduced to the original equation after applying the so-called reconstruction operator (see Theorem 3.10 in \cite{H1} for its definition). 

We start by defining the space of modelled distributions. For any $\kappa\in \R$, let
$$
W_{<\kappa}=\bigoplus_{\tiny\begin{array}{c}\beta\in \M :\\|\beta|<\kappa\end{array}}W_\beta.
$$
Let $f:\R^4\to W=\bigoplus_{\beta\in \M}W_\beta$ be a function. For any $\gamma\in \Ml$ we denote $f_\gamma(x)=(f(x))_\gamma:=P_\gamma (f(x))$, where $P_\gamma:\bigoplus_{\beta\in \Ml}W_\beta\to W_\gamma$ is the natural projection operator. Note that by Proposition \ref{Gprop}(i) for any $x,y\in \R^4$ and any $f:\R^4\to W_{<\kappa}$ we have $\G_{xy}^{W}f(y)\in W_{<\kappa}$. The space of modelled distributions $\D^\kappa(W)$ is defined as the space of all functions $f:\R^4\to W_{<\kappa}$ such that for every compact set $\mK\subset \R^4$ one has
\be\lab{distn}
|||f|||_{\kappa,\mK}:=\sup_{x\in \mK}\sup_{\tiny\begin{array}{c}\beta\in \M :\\ |\beta|<\kappa\end{array}}|f_\beta(x)|+
\ee
$$
+\sup_{\tiny\begin{array}{c} x,y\in \mK, \\ x\neq y\end{array}}\sup_{\tiny\begin{array}{c}\beta\in \M :\\|\beta|<\kappa\end{array}}\frac{|(f(x)-\G_{xy}^{W}f(y))_\beta|}{|x-y|^{\kappa-|\beta|}}<+\infty,
$$
where $|\cdot|$ in the first term in the right hand side and in the numerator of the second term is the norm on $W_\beta$ fixed in Proposition \ref{Pibond}.

As in \cite{H1} one can also define $\D^\kappa(W)$ as the set of equivalence classes of $W$--valued functions on $\R^4$ such that a function $f$ is equivalent to $g$ if $f_\beta =g_\beta$ for all $\beta$ with $|\beta|<\kappa$, and in formula (\ref{distn}) the representative $f$ with $f_\beta=0$ for all $\beta$ with $|\beta|\geq \kappa$ is used. Using this definition $\D^{\kappa'}(W)$ can be regarded as a subspace in $\D^\kappa(W)$ for any $\kappa'\geq \kappa$.

We call a modelled distribution $f\in \D^\kappa(W)$ coherent if
\be\lab{bcoh}
f_{\delta_\bg}=g,~{\rm and }~ f_{\beta_1+\beta_2}=f_{\beta_1}f_{\beta_2}
\ee
for all $\beta_1$, $\beta_2\in \M$ such that $\beta_1+\beta_2\in \M$ and $|\beta_1|$, $|\beta_2|$, $|\beta_1+\beta_2|<\kappa$, where the product $f_{\beta_1}f_{\beta_2}$ is defined using multiplication (\ref{Wprod}) in the target space and pointwise product of functions.
From now on we assume that $g\neq 0$. Then conditions (\ref{bcoh}) imply that 
\be\lab{bcoh1}
f_0=1,
\ee
and for any $\beta=l\delta_\bg+\sum_{i=0}^pk_i \delta_{\n_i}\in \M$ one has 
\be\lab{bcoh2}
f_\beta=g^l\prod_{i=0}^pf_{\delta_{\n_i}}^{k_i}.
\ee

We denote by $\D^\kappa_{\rm coh}(W)$ the subset of coherent modelled distributions in $\D^\kappa(W)$.  $\D^\kappa_{\rm coh}(W)$ is not a vector subspace in $\D^\kappa(W)$. But it inherits a topology from $\D^\kappa(W)$.

Since (\ref{LangHcrm}) is a system of equations, we shall also need the notion of $V$--valued modelled distributions. For any $\kappa\in \R$, let
$$
T_{<\kappa}=\bigoplus_{\tiny\begin{array}{c}\beta\in \M :\\|\beta|<\kappa\end{array}}V_\beta^*\oplus \bigoplus_{\tiny\begin{array}{c}\beta\in \M :\\|\beta|<\kappa+2\end{array}}{V_\beta^{-}}^*.
$$
Let ${\bf f}:\R^4\to T\otimes V=(\bigoplus_{\beta\in \M}V_\beta^*\oplus \bigoplus_{\beta\in \M}{V_\beta^{-}}^*)\otimes V$ be a function. For any $\gamma\in \Ml$ we denote ${\bf f}_\gamma(x)=({\bf f}(x))_\gamma:=(P_\gamma\otimes {\rm Id}_V) ({\bf f}(x))$ (resp ${\bf f}_\gamma^-(x)=({\bf f}(x))_\gamma^-:=(P_\gamma^- \otimes {\rm Id}_V)({\bf f}(x))$), where $P_\gamma:{\bf T}\to V_\gamma$ (resp. $P_\gamma^-:{\bf T}\to V_\gamma^-$) is the natural projection operator. Note that by Proposition \ref{Gprop}(i) for any $x,y\in \R^4$ we have $(\G_{xy}^{T}\otimes {\rm Id}_V){\bf f}(y)\in T_{<\kappa}\otimes V$. 

For all $\beta\in \Ml$ fix norms on $V_\beta^*\otimes V$ and on ${V_\beta^-}^*\otimes V$. The space of modelled distributions $\D^\kappa(T)$ is defined as the space of all functions $f:\R^4\to T_{<\kappa}\otimes V$ such that for every compact set $\mK\subset \R^4$ one has
$$
||{\bf f}||_{\kappa,\mK}:=\sup_{x\in \mK}\sup_{\tiny\begin{array}{c}\beta\in \M :\\ |\beta|<\kappa\end{array}}|{\bf f}_\beta(x)|+\sup_{x\in \mK}\sup_{\tiny\begin{array}{c}\beta\in \M :\\ |\beta|<\kappa+2\end{array}}|{\bf f}_\beta^-(x)|+
$$
\be\lab{distn1}
+\sup_{\tiny\begin{array}{c} x,y\in \mK, \\ x\neq y\end{array}}\sup_{\tiny\begin{array}{c}\beta\in \M :\\|\beta|<\kappa\end{array}}\frac{|({\bf f}(x)-(\G_{xy}^{T}\otimes {\rm Id}_V){\bf f}(y))_\beta|}{|x-y|^{\kappa-|\beta|}}+
\ee
$$
+\sup_{\tiny\begin{array}{c} x,y\in \mK, \\ x\neq y\end{array}}\sup_{\tiny\begin{array}{c}\beta\in \M :\\|\beta|<\kappa+2\end{array}}\frac{|({\bf f}(x)-(\G_{xy}^{T}\otimes {\rm Id}_V){\bf f}(y))_\beta^-|}{|x-y|^{\kappa+2-|\beta|}}<+\infty,
$$
where $|\cdot|$ in the first and in the second term in the right hand side and in the numerator of the third and of the fourth terms is the fixed norm on $V_\beta^*\otimes V$ or ${V_\beta^-}^*\otimes V$. Since all such norms are equivalent the definition does not depend on the choice of the norms.

One can also define $\D^\kappa(T)$ as the set of equivalence classes of $T\otimes V$--valued functions on $\R^4$ such that a function $f$ is equivalent to $g$ if ${\bf f}_\beta={\bf g}_\beta$ for all $\beta$ with $|\beta|<\kappa$, and ${\bf f}_\beta^-={\bf g}_\beta^-$ for all $\beta$ with $|\beta|<\kappa+2$, and in formula (\ref{distn1}) the representative ${\bf f}$ with ${\bf f}_\beta=0$ for all $\beta$ with $|\beta|\geq \kappa$ and ${\bf f}_\beta^-=0$ for all $\beta$ with $|\beta|\geq \kappa+2$ is used. Using this definition $\D^{\kappa'}(T)$ can be regarded as a subspace in $\D^\kappa(T)$ for any $\kappa'\geq \kappa$,
\be\lab{embD}
\D^{\kappa'}(T)\subset \D^\kappa(T),~\kappa'\geq \kappa.
\ee

For any $f\in \D^{\kappa}(W)$, let
\be\lab{Ff}
{\bf F}(f)(x)=\sum_{\tiny\begin{array}{c}\beta\in \M :\\|\beta|<\kappa\end{array}}\tau_\beta({\rm Id}_{V_\beta}\otimes f_\beta(x))\in \D^\kappa(T),~x\in \R^4,
\ee
where by $\tau_\beta({\rm Id}_{V_\beta}\otimes f_\beta(x))$ it is meant that the second tensor factor of 
$$
\tau_\beta={\rm Id}_{V_\beta^*}\in {\rm End}(V_\beta^*)\simeq V_\beta^*\otimes V_\beta,
$$
which belongs to $V_\beta={\rm Hom}(W_\beta,V)$, is evaluated on $f_\beta(x)\in W_\beta$, and the first tensor factor is not affected, so that the result is an element of $V_\beta^*\otimes V$.

Similarly we define
\be\lab{F-f}
{\bf F}^-(f)(x)=\sum_{\tiny\begin{array}{c}\beta\in \M :\\|\beta|<\kappa\end{array}}\tau_\beta^-({\rm Id}_{V_\beta^-}\otimes f_\beta(x))\in \D^{\kappa-2}(T).
\ee

Following \cite{H1}, \S5 we also introduce the abstract integration map
$$
\mathcal{I}:\bigoplus_{\beta\in \M }{V_\beta^-}^*\to \bigoplus_{\beta\in \M}V_\beta^*
$$
by requiring that 
\be\lab{IInt}
\mathcal{I}|_{{V_\beta^-}^*}=\left\{\begin{array}{ll} {\rm Id}_{{V_\beta^-}^*}^{V_\beta^*} &{\rm if }~ \beta\in \M' \\
0 & {\rm if }~ \beta\in \M_{pp}\end{array}\right. ,
\ee
where ${\rm Id}_{{V_\beta^-}^*}^{V_\beta^*}: {V_\beta^-}^*\to V_\beta^*$ is the natural isomorphism.

For any $r> 0$, $r\not \in \N$, let
\be\lab{Br}
B^{r,loc}_{\infty,\infty}:=\{f:\R^4\to \R:|f|_{B^r_{\infty,\infty}(\mK)}:=\max_{|\n|\leq [r]}\sup_{x\in \mK}|\partial^{\n}f(x)|+ 
\ee
$$
+\max_{|\n|=[r]}\sup_{\tiny\begin{array}{c} x,y\in \mK, \\ x\neq y\end{array}}\frac{|\partial^{\n}f(x)-\partial^{\n}f(y)|}{|x-y|^{\{r\}}}<+\infty~\text{\rm for any compact subset }\mK\subset \R^4\},
$$
be the local H\"{o}lder space of exponent $r$.

For any $s<0$ define also the Besov space $B^{s,loc}_{\infty,\infty}$ as the space of tempered distributions $\zeta\in \Sw'$ such that some $r>|s|$ and for any compact subset $\mK\subset \R^4$ 
$$
\sup_{0<\lambda\leq 1}\sup_{x\in \mK}\sup_{\varphi \in \mB^r}|\zeta(\varphi_x^{\lambda})|\lambda^{-s}<\infty.
$$

Let $\mathcal{R}:\D^\kappa(T)\to B^{s,loc}_{\infty,\infty}\otimes V$ be the natural extension of the reconstruction operator acting on $T$--valued modelled distributions as defined in Theorem 3.10 in \cite{H1}. As it is mentioned in Remark 3.15 in \cite{H1}, in the case when $\Pi_x(\varpi)\in B^{r,loc}_{\infty,\infty}$ for some $0<r<1$ and for all $\varpi\in T$, which is the case in our situation as $\Pi_x(\varpi)\in \O\subset B^{r,loc}_{\infty,\infty}$ for any $r>0$, the reconstruction operator has a very simple form which implies the following expression for $\mathcal{R}{\bf f}$ if ${\bf f}\in \D^\kappa(T)$ for some $0<\kappa<1$,
\be\lab{rec}
(\mathcal{R}{\bf f})(x)=(\Pi_x\otimes {\rm Id}_V)({\bf f}(x))(x).
\ee
In fact, the right hand side of the last formula defines an element of $B^{u,loc}_{\infty,\infty}\otimes V$, where $u=\min\{r,\kappa\}$. 

Indeed, by (\ref{distn1}) and (\ref{Br}) for any compact subset $\mK\subset \R^4$ one has for all $x\in \mK$
\be\lab{R1}
|(\mathcal{R}{\bf f})(x)|\leq C_B\max_{\varpi\in B}\sup_{x\in \mK}|\Pi_x(\varpi)(x)|\sup_{x\in \mK}|{\bf f}(x)|\leq 
\ee
$$
\leq C_B\max_{\varpi\in B}|\Pi_x(\varpi)|_{B^r_{\infty,\infty}(\mK)}||{\bf f}||_{\kappa,\mK},
$$
where $B$ is a basis of $T_{<\kappa}$, $C_B > 0$, $|\cdot |$ in the left hand side is a norm on $V$, $|\cdot |$ in the right hand side evaluated at ${\bf f}(x)$ is the norm on $T$ induced by the norms on $V_\beta^*\otimes V$ and ${V_\beta^-}^*\otimes V$ as in (\ref{distn1}). For all $x,y\in \mK$ we also obtain by (\ref{distn1})  and (\ref{Br}) 
$$
|(\mathcal{R}{\bf f})(x)-(\mathcal{R}{\bf f})(y)|\leq |(\Pi_x\otimes {\rm Id}_V)({\bf f}(x))(x)-(\Pi_x\otimes {\rm Id}_V)({\bf f}(x))(y)|+
$$
$$
+|(\Pi_y\otimes {\rm Id}_V)((\G^{T}_{yx}\otimes {\rm Id}_V){\bf f}(x)-{\bf f}(y))(y)|\leq
$$
\be\lab{R2}
\leq C_B'\max_{\varpi\in B}|\Pi_x(\varpi)(x)|_{B^r_{\infty,\infty}(\mK)}||{\bf f}||_{\kappa,\mK}|x-y|^{r}+
\ee
$$
+C_{B,\mK}\max_{\varpi\in B}|\Pi_x(\varpi)(x)|_{B^r_{\infty,\infty}(\mK)}||{\bf f}||_{\kappa,\mK}|x-y|^{\kappa}\leq
$$
$$
\leq C_B^\mK\max_{\varpi\in B}|\Pi_x(\varpi)(x)|_{B^r_{\infty,\infty}(\mK)}||{\bf f}||_{\kappa,\mK}|x-y|^{\min\{r,\kappa\}},
$$
where $C_B'$, $C_{B,\mK}$, $C_B^\mK> 0$, and in the second term in the right hand side we used the fact that by (\ref{Pibnd1}) and (\ref{Pibnd2}) one has $\Pi_x(\varpi)(x)=0$ for any $\varpi\in T_v$, $v>0$, and hence when applying (\ref{distn1}) in the second term in the right hand side we can use that for some $C_\mK$, $C'_\mK> 0$ and for all $x,y\in \mK$ one also has $|x-y|^{\kappa-|\beta|}\leq C_\mK|x-y|^{\kappa}$ for $|\beta|<0$ and $|x-y|^{\kappa+2-|\beta|}\leq C'_\mK|x-y|^{\kappa}$ for $|\beta|-2<0$.

Estimates (\ref{R1}) and (\ref{R2}) imply that 
\be\lab{fBloc}
\mathcal{R}{\bf f}\in B^{u,loc}_{\infty,\infty}\otimes V,~u=\min\{r,\kappa\}.
\ee

Now for $\kappa>2$, following the ideas of Sections 5 and 7 of \cite{H1}, we consider the following equation for $f\in \D^\kappa_{\rm coh}(W)$  
$$
{\bf F}(f)(x)=(\mathcal{I}\otimes {\rm Id}_V){\bf F}^-(f)(x)+
$$
\be\lab{eqreg}
+\sum_{\n\in \N^4:|\n|<\kappa}\frac{1}{\n!}\tau_{\delta_\n}({\rm Id}_{V_{\delta_\n}}\otimes {\rm Id}_V^{W_{\delta_\n}}(\partial^\n\K_m(\mathcal{R}{\bf F}^-(f))(x)-
\ee
$$
-\partial^\n\K(\Pi_xP_{\leq |\n|-2}\otimes {\rm Id}_V)({\bf F}^-(f)(x))(x)+\partial^\n\widehat{u}_0(x))),
$$
where ${\rm Id}_V^{W_{\delta_\n}}:V\to W_{\delta_\n}$ is the natural isomorphism, for any $a\in \R$
$$
P_{\leq a}:{\bf T}=\bigoplus_{\beta\in \Ml}V_\beta^*\oplus \bigoplus_{\beta\in \Ml}{V_\beta^{-}}^*\to \bigoplus_{\tiny \begin{array}{c} \beta\in \Ml: \\|\beta|\leq a \end{array}}V_\beta^*\oplus \bigoplus_{\tiny \begin{array}{c}\beta\in \Ml: \\ |\beta|-2\leq a\end{array}}{V_\beta^{-}}^*
$$
is the natural projection operator, and 
\be\lab{inc}
\widehat{u}_0(x)=\widehat{u}_0(t,\bar{x})=\int_{\R^3}K_m(t,\bar{x}-\bar{y})u_0(\bar{y})d\bar{y}
\ee
for some $u_0:\R^3\to V$ which plays the role of an initial data.

Equation (\ref{eqreg}) resembles equation (7.18) in \cite{H1}.
We claim that solutions to equation (\ref{eqreg}) satisfy (\ref{LangHcrm}).
  
In this section we do not discuss the existence of solutions to equation (\ref{eqreg}). Actually one can show that the right hand side of (\ref{eqreg}) belongs to $\D^\kappa(T)$. We do not prove this result as we shall not need it in this paper. We only note that, for any $\n\in \N^4$ with $|\n|<2+\varepsilon$, $0<\varepsilon<1$, the $\n$-th term in the right hand side of (\ref{eqreg}) is a genuine function as by (\ref{embD}), (\ref{F-f}) and (\ref{fBloc}) one has $\mathcal{R}{\bf F}^-(f)\in B^{\varepsilon,loc}_{\infty,\infty}\otimes V$, so by the classical Schauder estimates (see e.g. \cite{Kryl}, Theorem 8.9.2) we obtain $\K_m\mathcal{R}{\bf F}^-(f)\in B^{\varepsilon+2,loc}_{\infty,\infty}\otimes V$, and $\partial^\n\K_m(\mathcal{R}{\bf F}^-(f))\in B^{\varepsilon+2-|\n|,loc}_{\infty,\infty}\otimes V$ for $|\n|<\varepsilon+2$. The third term in the right hand side of (\ref{eqreg}) is a function by the definition, and for a large class of initial data $u_0$ the last term is smooth for $t>0$, as $\widehat{u}_0(t,\bar{x})$ is smooth for $t>0$ if, e.g., $u_0$ is locally integrable. In particular, if $\kappa=2+\varepsilon$, $0<\varepsilon<1$ the right hand side of (\ref{eqreg}) is a genuine function if $t>0$.

From (\ref{F-f}) and (\ref{IInt}) it follows that
$$
(\mathcal{I}\otimes {\rm Id}_V){\bf F}^-(f)(x)=\sum_{\tiny\begin{array}{c}\beta\in \M' :\\|\beta|<\kappa\end{array}}\tau_\beta({\rm Id}_{V_\beta}\otimes f_\beta(x))\in \D^{\kappa}(T),
$$
and hence by (\ref{Ff}) for any $\beta\in \M'$ the $V_\beta^*$--component of equation (\ref{eqreg}) is reduced to the identity $f_\beta=f_\beta$. 

For $\n\in \N^4$, $|\n|<\kappa$ the $V^*_{\delta_\n}$--component of (\ref{eqreg}) yields by (\ref{defPP}) for $\Pi_x$, by (\ref{Ff}) and by (\ref{F-f})
\be\lab{equcomp}
{\rm Id}_{W_{\delta_\n}}^V f_{\delta_\n}(x)= \frac{1}{\n!}\partial^\n\K_m(\mathcal{R}{\bf F}^-(f))(x)-
\ee
$$
-\frac{1}{\n!}\sum_{\tiny\begin{array}{c}\beta\in \M': \\ |\beta|\leq |\n|\end{array}}(\partial^\n\K\Pi_{x\beta}^-)(f_\beta(x))(x)+\frac{1}{\n!}\partial^\n\widehat{u}_0(x).
$$

Since by the choice of $\alpha$ in Section \ref{scaleh} and by (\ref{hom}) for $\beta\not\in \M_{pp}$  one has $|\beta|\not\in \N$, the sum over $\beta$ in (\ref{equcomp}) is actually over $\beta\in \M'$ such that $|\beta|< |\n|$, and hence by (\ref{Pder}) for such $\beta$ one has $\partial^\n\K\Pi_{x\beta}^-=\partial^\n\Pi_{x\beta}$. Thus (\ref{equcomp}) takes the form
\be\lab{equcomp1}
{\rm Id}_{W_{\delta_\n}}^V f_{\delta_\n}(x)= \frac{1}{\n!}\partial^\n\K_m(\mathcal{R}{\bf F}^-(f))(x)-
\ee
$$
-\frac{1}{\n!}\sum_{\tiny\begin{array}{c}\beta\in \M': \\ |\beta|\leq |\n|\end{array}}(\partial^\n\Pi_{x\beta})(f_\beta(x))(x)+\frac{1}{\n!}\partial^\n\widehat{u}_0(x).
$$

Observe that by (\ref{defPP}), (\ref{Pvan}), (\ref{bcoh}), (\ref{bcoh1}) and (\ref{rec})
\be\lab{Ar}
A(x):=(\mathcal{R}{\bf F}(f))(x)=\sum_{\tiny\begin{array}{c}\beta\in \M: \\ |\beta|\leq 0 \end{array}}\Pi_{x\beta}(f_\beta(x))(x)=
\ee
$$
=v(x)+\Pi_0(x)+g\Pi_{\delta_\bg}(x),~v(x):={\rm Id}_{W_{\delta_{\bf 0}}}^V f_{\delta_{\bf 0}}(x),
$$
as $\delta_{\bf 0}$, $0$ and $\delta_\bg$ are the only multi--indices with non--positive homogeneities.

Similarly, by (\ref{defPP}), (\ref{P-van}) and (\ref{rec})
\be\lab{Rr}
(\mathcal{R}{\bf F}^-(f))(x)=\sum_{\tiny\begin{array}{c}\beta\in \M': \\ |\beta|\leq 2 \end{array}}\Pi_{x\beta}^-(f_\beta(x))(x).
\ee 

Now for $\n={\bf 0}$ (\ref{equcomp1}) takes the form
\be\lab{equcomp2}
v(x)= \K_m(\mathcal{R}{\bf F}^-(f))(x)
-\sum_{\tiny\begin{array}{c}\beta\in \M': \\ |\beta|\leq 0\end{array}}\Pi_{x\beta}(f_\beta(x))(x)+\widehat{u}_0(x).
\ee
Rearranging terms and using (\ref{Ar}) we obtain
\be\lab{Ax}
A(x)=\K_m(\mathcal{R}{\bf F}^-(f))(x)+\widehat{u}_0(x).
\ee

By (\ref{equcomp2}) formula (\ref{equcomp1}) for $\n\neq {\bf 0}$ takes the form
\be\lab{equcomp3}
{\rm Id}_{W_{\delta_\n}}^V f_{\delta_\n}(x)= \frac{1}{\n!}\partial^\n v(x)-\frac{1}{\n!}\sum_{\tiny\begin{array}{c}\beta\in \M': \\ 0<|\beta|\leq |\n|\end{array}}(\partial^\n\Pi_{x\beta})(f_\beta(x))(x).
\ee

Note that for $\beta\in \M'$ the components $f_\beta$ are expressed in terms of $f_{\delta_\m}$, $\m\in \N^4$, $|\m|<|\beta|$ with the help of (\ref{bcoh2}). Using these expressions together with (\ref{equcomp3}) one can find all $f_\beta$, $\beta\in \M$ by induction once $f_{\delta_{\bf 0}}$ is known, at least when $|\beta|<2+\varepsilon$, $0<\varepsilon<1$ and the right hand side of (\ref{equcomp3}) is a genuine function as shown above, so that all products in (\ref{bcoh2}) are defined. 

To obtain a closed equation for $f_{\delta_{\bf 0}}$ we observe that applying the operator $L$ to (\ref{Ax}) and using (\ref{inc}) and the definition of $\K_m$ one gets for $t>0$
\be\lab{LA}
LA=\mathcal{R}{\bf F}^-(f).
\ee
We show that this equation coincides with equation (\ref{LangdTrmAB}) which is equivalent to (\ref{LangHcrm}).

Indeed, by (\ref{equcomp3}) with $|\n|=1$
$$
\partial^\n v(x)={\rm Id}_{W_{\delta_\n}}^V f_{\delta_\n}(x)+\sum_{\tiny\begin{array}{c}\beta\in \M': \\ 0<|\beta|\leq 1 \end{array}}(\partial^\n\Pi_{x\beta})(f_\beta(x))(x)=
$$
$$
=\sum_{\tiny\begin{array}{c}\beta\in \M: \\ 0<|\beta|\leq 1 \end{array}}(\partial^\n\Pi_{x\beta})(f_\beta(x))(x),
$$
and hence by (\ref{Ar})
$$
\partial^\n A(x)=\partial^\n v(x)+\sum_{\tiny\begin{array}{c}\beta\in \M: \\ |\beta|< 0 \end{array}}(\partial^\n \Pi_{x\beta})(f_\beta(x))(x)=
\sum_{\tiny\begin{array}{c}\beta\in \M: \\ |\beta|\leq 1 \end{array}}(\partial^\n\Pi_{x\beta})(f_\beta(x))(x),
$$
or
$$
\partial_i A(x)=
\sum_{\tiny\begin{array}{c}\beta\in \M: \\ |\beta|\leq 1 \end{array}}(\partial_i\Pi_{x\beta})(f_\beta(x))(x),~i=1,2,3.
$$
Substituting this expression and (\ref{Ar}) into the right hand side of (\ref{LangdTrmAB}), recalling that $\Pi_{x0}^-=\xi^\rho$ and using (\ref{Pvan}), (\ref{P-van}), (\ref{bcoh}) and (\ref{bcoh1}) we obtain
$$
g\A(A,A)(x) + g^2 \mB(A,A,A)(x)+cA(x)+\xi^\rho(x)=
$$
$$
=g\hspace{-0.8cm}\sum_{\tiny\begin{array}{c}\beta_{1,2}\in \M: \\ |\beta_1|\leq 0, |\beta_2|\leq 1 \end{array}}\hspace{-0.8cm}\A(\Pi_{x\beta_1},\Pi_{x\beta_2})(f_{\beta_1}(x)f_{\beta_2}(x))(x)+
$$
$$
+g^2\hspace{-0.5cm}\sum_{\tiny\begin{array}{c}\beta_{1,2,3}\in \M: \\ |\beta_{1,2,3}|\leq 0 \end{array}}\hspace{-0.5cm}\mB(\Pi_{x\beta_1},\Pi_{x\beta_2},\Pi_{x\beta_3})(f_{\beta_1}(x)f_{\beta_2}(x)f_{\beta_3}(x))(x)+
$$
$$
+\hspace{-1cm}\sum_{\tiny \begin{array}{c} k\in  \{1,2,3,4\},\beta_1\in \M: \\ |\beta_1|\leq 0\end{array}}\hspace{-1cm}c_kg^k\Pi_{x\beta_1}(f_{\beta_1}(x))(x)+\Pi_{x0}^-(f_0(x))(x)=
$$
$$
=\sum_{\tiny \begin{array}{c}\beta\in \Ml: \\ \beta\neq 0\end{array}}\left(\sum_{\tiny\begin{array}{c}\beta_{1},\beta_{2}\in \Ml: \\ \beta_1+\beta_2+\delta_\bg=\beta \end{array}}\hspace{-0.8cm}\A(\Pi_{x\beta_1},\Pi_{x\beta_2})(f_{\beta}(x))(x)+\right.
$$
$$
+\hspace{-0.5cm}\sum_{\tiny\begin{array}{c}\beta_{1},\beta_2,\beta_3\in \Ml: \\ \beta_1+\beta_2+\beta_3+2\delta_\bg=\beta\\  \end{array}}\hspace{-0.5cm}\mB(\Pi_{x\beta_1},\Pi_{x\beta_2},\Pi_{x\beta_3})(f_{\beta}(x))(x)+
$$
$$
\left. +\hspace{-1cm}\sum_{\tiny \begin{array}{c} k\in  \{1,2,3,4\},\beta_1\in \Ml: \\ k\delta_\bg+\beta_1=\beta \end{array}}\hspace{-1cm}c_k\Pi_{x\beta_1}(f_{\beta}(x))(x)\right)+\Pi_{x0}^-(f_0(x))(x).
$$
Finally by (\ref{defP-}) this identity takes the form
$$
g\A(A,A)(x) + g^2 \mB(A,A,A)(x)+cA(x)+\xi^\rho(x)=\sum_{\beta\in \Ml}\Pi_{x\beta}^-(f_{\beta}(x))(x)=
$$
$$
=\sum_{\beta\in \M':|\beta|\leq 2}\Pi_{x\beta}^-(f_{\beta}(x))(x)=\mathcal{R}{\bf F}^-(f)=LA,
$$
where at the last two steps we also used (\ref{P-van}), (\ref{Rr}) and (\ref{LA}). The last identity confirms that $A$ defined by (\ref{Ar}) satisfies (\ref{LangHcrm}).

In conclusion we note that by (\ref{Ar}) $f_{\delta_{\bf 0}}$, and hence all $f_\beta$, $\beta\in \M$, can be found once $A$ is known.


\section{Analytic tools}\label{Anal}

In this section we recall some analytic results required for the proof of stochastic estimates for a model in Section \ref{Stohe}. In Section \ref{Diffrec} we also obtain some results related to Malliavin (directional) derivatives for recentered maps used in this proof.


\subsection{Random variables and function spaces}

\setcounter{equation}{0}
\setcounter{theorem}{0}

Let $E$ be a finite-dimensional real vector space. Fix a norm $|\cdot|$ on $E$.
For any function $\varphi:\R^4\to E$ we denote by 
$$
\left\|\varphi\right\|_{L_\infty}:={\rm ess}-\sup_{x\in \R^4}|\varphi(x)|
$$
its (essential) supremum norm, and for any $r\in \N$ let
$$
\left\|\varphi\right\|_{C^r}:=\max_{\n\in \N^4:|\n|\leq r}\left\|\partial^\n\varphi\right\|_{L_\infty}.
$$
Note that we use definition (\ref{nhom}) of $|\n|$ in the last formula. 

For any $r\in \N$, $s\in \R$, let
$$
\mB^r:=\{\varphi\in \D(B_1):\left\|\varphi\right\|_{C^r}\leq 1\},
$$
and
$$
\mB^r_s:=\{\varphi\in \mB^r:\int_{\R^4}x^\n\varphi(x)dx=0~{\rm for~ all}~\n\in \N^4:|\n|\leq s \},
$$
where we assume that $\mB^r_s=\mB^r$ if $s<0$.

We shall need a parabolic version of the Morrey inequality which relates the H\"{o}lder and the Sobolev seminorms as follows (see \cite{BIN}, Theorem 16.10 for the most general scaled version of it, \cite{WYW}, Theorem 7.1.1 for the parabolic version, and \cite{DiNPV}, Theorem 8.2, \cite{GM}, Theorem 5.5, \cite{WYW}, Theorem 6.1.1 for more elementary proofs of the usual Euclidean version).

Let $\Omega\subset \R^4$ be a bounded open subset. Denote by ${\rm diam}(\Omega)$ the parabolic diameter of $\Omega$, ${\rm diam}(\Omega):=\sup_{x,y\in \Omega}|x-y|$. Assume that there exists $A>0$ such that for all $x\in \Omega$, $0\leq r<{\rm diam}(\Omega)$ one has 
$$
|\int_{B_r(x)\cap \Omega}\hspace{-1cm} 1 dx|\geq Ar^d.
$$
For instance, this condition is satisfied if the boundary of $\Omega$ is of class $C^1$. 
 
Let $1\leq p<\infty$, $\vartheta>d$. Then for $f:\Omega\to E$ and for all $x,y\in \Omega$
\be\lab{Morr}
|f(x)-f(y)|^p\lesssim \iint_{\Omega\times \Omega}\frac{|f(z)-f(z')|^p}{|z-z'|^{d+\vartheta}}dzdz'|x-y|^{\vartheta-d},
\ee 
where the constant in the inequality only depends on $\Omega$, $\vartheta$ and $p$. 

In particular, we shall use this inequality for $\Omega=B_R(y)$, $R>0$, $y \in \R^4$. Due to obvious translational invariance (\ref{Morr}) implies in this case by the triangle inequality
\be\lab{Morr1}
|f(x)|^p\lesssim |f(y)|^p+\iint_{B_R(y)\times B_R(y)}\frac{|f(z)-f(z')|^p}{|z-z'|^{d+\vartheta}}dzdz',~x\in B_R(y),
\ee 
where the constant in the inequality only depends on $R$, $p$ and $\vartheta$.

We shall also need the following elementary inequality which is valid for continuously differentiable $E$--valued functions $f$ on the interval $(0,\bar{\lambda}]$, $\bar{\lambda}>0$, $1<p<\infty$, $a\in \R$,
\be\lab{Sobt}
\sup_{0<t\leq \bar{\lambda}}|t^{-a+\frac 1p}f(t)|^p\lesssim \int_0^{\bar{\lambda}}(|t^{-a}f(t)|^p+|t^{-a+1}f'(t)|^p)dt,
\ee
where the constant in the inequality only depends on $a,p,\bar{\lambda}$.

This inequality can be found in \cite{BOS}, Appendix B. For completeness we give its proof. The proof is very elementary and relies on the Newton formula.

Indeed, we may assume that left hand side of (\ref{Sobt}) is finite, and hence $|t^{-a+1}f(t)|=|t^{-a+\frac 1p+1-\frac 1p}f(t)|=  |t^{-a+\frac 1p}f(t)|t^{1-\frac 1p}\to 0$ when $t\to 0$ as $1-\frac 1p>0$. Thus if we denote $h(t):=t^{-a+1}f(t)$ then we can extend $h$ to a continuous function on the interval $[0,\bar{\lambda}]$ in such a way that $h(0)=0$. Then, as $h$ is continuously differentiable on $(0,\bar{\lambda}]$, by the Newton formula we have for $t\in (0,\bar{\lambda}]$
$$
h(t)=\int_0^th'(s)ds.
$$
Since $h'(s)=(1-a)s^{-a}f(s)+s^{-a+1}f'(s)$ the previous identity and the H\"{o}lder inequality imply
$$
|h(t)|\lesssim \left(\int_0^t (|s^{-a}f(s)|^p+|s^{-a+1}f'(s)|^p)ds\right)^{\frac 1p}t^{1-\frac 1p}\lesssim 
$$
$$
\lesssim \left(\int_0^{\bar{\lambda}} (|s^{-a}f(s)|^p+|s^{-a+1}f'(s)|^p)ds\right)^{\frac 1p}t^{1-\frac 1p},
$$
or recalling that $h(t)=t^{-a+1}f(t)$, 
$$
|t^{-a+\frac 1p}f(t)|^p \lesssim \int_0^{\bar{\lambda}} (|s^{-a}f(s)|^p+|s^{-a+1}f'(s)|^p)ds,
$$
where the constant in the inequality only depends on $a,p,\bar{\lambda}$. The last inequality yields (\ref{Sobt}).

For an $E$--valued function $f$ of variable $x\in \R^4$, $1\leq p\leq \infty$, and any measurable subset $\Omega\subset \R^4$ we denote
$$
\left\|f(x)\right\|_{L_p(\Omega,x)}:=\left(\int_{\Omega}|f(x)|^pdx\right)^{\frac 1p},
$$
along with the usual convention that if $p=\infty$ the integral in the right hand side is replaced with the essential supremum norm.

If $\Omega=\R^4$ we simply write $L_p(\R^4)=L_p$ and $\left\|f(x)\right\|_{L_p(\R^4,x)}=\left\|f\right\|_{L_p}$.

For an $E$--valued random variable $f$ and $1\leq p\leq \infty$ we define
$$
\left\|f\right\|_{\L_p}:=\E^{\frac 1p}(|f|^p),
$$
along with the usual convention that if $p=\infty$ the corresponding norm is the almost sure supremum norm.

We shall call $E$--valued random distributions $\Sw'\otimes E$--valued random variables.
If an $E$--valued random distribution takes values in $C(\R^4)\otimes E$ we call it continuous and if it takes values in $C^\infty(\R^4)\otimes E$ we call it smooth.

For an $E$--valued random distribution $\zeta$ and a function $\varphi\in \Sw$ we mean by $\zeta(\varphi)\in E$ the evaluation of $\zeta$ on $\varphi$ with respect to the first factor in the tensor product $\Sw'\otimes E$. We keep the same notation in the case when $\zeta$ is actually continuous or smooth.

For any $1\leq p,q\leq \infty$ and for a family of $E$--valued random variables $f(x)$, $x\in \R^4$, which can be regarded as a random variable with values in the space of $E$--valued functions of $x$ on $\R^4$, we introduce
$$
\left\|f(x)\right\|_{pq\Omega x}:=\left\|\left\|f(x)\right\|_{\L_p}\right\|_{L_q(\Omega,x)}.
$$

The following stochastic version of Sobolev's inequality can be found in \cite{BOT}, Section 3.5 (see \cite{AF}, Section 4.16 or \cite{BIN}, Theorem 10.4  for the usual deterministic version of it; the latter reference contains the most general scaled version that we use in the case of a parabolic scaling):
\bp\lab{sob}
For any smooth $E$--valued random distribution $f$, $1\leq q<\infty$, $R>0$, and any $k>\frac d2$ one has
$$
\left\|f(y)\right\|_{\L_q}\lesssim \sum_{\tiny\begin{array}{c} \n\in \N^4:\\ |\n|\leq k\end{array}}R^{|\n|-\frac d2}\left\|(\partial^\n f)(x+y)\right\|_{q2B_Rx},
$$
where the constant in the inequality only depends on $k$.
\ep

For $1\leq p,q\leq \infty$, $s\in \R$, $r\in \N$, $r>-s$, $\bar{\lambda}\in \R$, $\bar{\lambda}>0$, any  $E$--valued random distribution $\zeta$  we define 
\be\lab{defbspq<}
\left\|\zeta\right\|_{\B_{pq}^s}:=\sup_{0<\lambda\leq \bar{\lambda}}\frac{\left\|\sup_{\varphi\in \mB^r}\left\|\zeta(\varphi_x^\lambda)\right\|_{\L_p}\right\|_{L_q(\R^4,x)}}{\lambda^s} ~{\rm if }~ s<0,
\ee
and
\be\lab{defbspq>}
\left\|\zeta\right\|_{\B_{pq}^s}:=\left\|\sup_{\varphi\in \mB^r}\left\|\zeta(\varphi_x)\right\|_{\L_p}\right\|_{L_q(\R^4,x)}+
\ee
$$
+\sup_{0<\lambda\leq \bar{\lambda}}\frac{\left\|\sup_{\varphi\in \mB^r_{\left\lfloor s \right\rfloor}}\left\|\zeta(\varphi_x^\lambda)\right\|_{\L_p}\right\|_{L_q(\R^4,x)}}{\lambda^s} ~{\rm if }~ s\geq 0.
$$

Let $\B_{pq}^s$ be the space of all $E$--valued random distributions $\zeta$ such that $\left\|\zeta\right\|_{\B_{pq}^s}<\infty$. Following the arguments in the proof of Proposition A5 in \cite{BL} one can show that this definition does not depend on the choice of $r$ and $\bar{\lambda}$. 

$\B_{pq}^s$ can be regarded as natural stochastic counterparts of Besov spaces (see, e.g., \cite{Tr,Tr1,Tr2} and Appendix A in \cite{BL} for similar deterministic versions of these spaces and their properties).
Note, however, that in the case when $E=\R$ the subspace of $\B_{pq}^s$ which consists of elements of $\Sw'\otimes E=\Sw'$ is $B_{q\infty}^{s}$ in the usual Besov space notation, so that in our notation the indices $p$ and $q$ do not correspond to the usual indices $p$ and $q$ in the definition of Besov spaces. 

If $s<0$ then $\B_{pq}^s$ can also be defined as the space of all $E$--valued random distributions $\zeta$ such that there exists $\varphi\in \D$ with 
$$
\int_{\R^4}\varphi(x)dx\neq 0,
$$
and 
\be\lab{defbspq<'}
\left\|\zeta\right\|_{\B_{pq}^s}':=\sup_{0<\lambda\leq \bar{\lambda}}\frac{\left\|\left\|\zeta(\varphi_x^\lambda)\right\|_{\L_p}\right\|_{L_q(\R^4,x)}}{\lambda^s} <\infty.
\ee
Moreover, $\left\|\cdot\right\|_{\B_{pq}^s}'$ is an equivalent quasi--norm on $\B_{pq}^s$. The last statement can be obtained by replacing the moduli of real numbers with the norm $\left\|~\cdot~\right\|_{\L_p}$ in the arguments in the proof of Proposition A5 (1) in \cite{BL} and by using these arguments verbatim. A similar statement for the usual Besov spaces can also be found in \cite{Tr}, Corollary 1.12.

We shall also need the so-called spectral gap inequality in the case when the underlying probability space is the white noise space modeled on $\Sw'\otimes V$, where $V$ is a Euclidean space, with the covariance induced by the scalar product on $V$, so that its Cameron--Martin space ${\rm CM}$ is $L_2(\R^4)\otimes V$ with the natural induced scalar product. In fact, in our case $V=\k\otimes\R^3$  where the scalar product is induced by $-(~\cdot~,~\cdot~ )$ on $\k$, and by the standard Euclidean scalar product on $\R^3$ fixed in the Introduction. Denote also by ${\rm WN}$ the corresponding space of $\Sw'\otimes V$--valued random variables. Note that the white noise space variable $\xi$ can be regarded as an element of ${\rm WN}$. 

For an $E$--valued random variable $f$ of random argument $\xi$  we denote by $\delta_\upsilon f$ or simply by $\delta f$, if it does not cause a confusion, its directional derivative with respect to an element $\upsilon\in {\rm CM}$, 
$$
\delta_\upsilon f(\xi)=\frac{d}{dt}|_{t=0}f(\xi+t\upsilon).
$$
If $\delta$ appears in a formula several times it is meant that this is the directional derivative with respect to the same element $\upsilon$ in all occurrences.

The spectral gap inequality will be used in the form of equation (187) in \cite{BOT},
\be\lab{SpG}
\left\|f\right\|_{\L_p}\lesssim |\E(f)|+\sup_{\tiny \begin{array}{c}\upsilon\in \L_{p^*}(L_2): \\ \left\|\left\|\upsilon\right\|_{L_2}\right\|_{\L_{p^*}}\leq 1\end{array}}\left\|\delta_\upsilon f\right\|_{\L_q},
\ee
where $\L_{p^*}(L_2):=\{F\in {\rm WN}: \left\|\left\|F\right\|_{L_2}\right\|_{\L_{p^*}}<\infty\}$, $1\leq q<p^*\leq 2\leq p$, $\frac 1p+\frac{1}{p^*}=1$, and the constant in the inequality depends on $p$. $\delta_\upsilon f$ in the right hand side of the spectral gap inequality makes sense if $f$ belongs to the Malliavin--Sobolev class ${\mathbb W}_p^1$ (see \cite{Bog}, Chapter 5 for their definition). Indeed, $\delta_\upsilon f$ should be viewed as the Malliavin derivative of $f$ evaluated at $\upsilon$. Recall that the Malliavin derivative is well defined as a closed operator from $\L_p$ to $\L_p(L_2)=\{F\in {\rm WN}: \left\|\left\|F\right\|_{L_2}\right\|_{\L_p}<\infty\}$ for $p\geq 2$, with the domain being the Malliavin--Sobolev space ${\mathbb W}_p^1$ (see \cite{Bog}, Chapter 5). Hence by the natural duality for $\L_p(L_2)$--spaces, for $p\geq 2$, $\E (\delta_\upsilon f)$ is finite for $\upsilon\in \L_{p^*}(L_2)$, $f\in {\mathbb W}_p^1$. Note that by the H\"{o}lder inequality $|\E(\delta_\upsilon f)|\leq \left\|\delta_\upsilon f\right\|_{\L_q}$ as $1\leq q$, so, generally speaking, if $q>1$ the right hand side of the spectral gap inequality in form (\ref{SpG}) is finite under a restriction stronger than just $f\in {\mathbb W}_p^1$.       

In fact, the spectral gap inequality with $q=1$ holds for bounded cylindrical functions provided that the underlying probability measure is centered Gaussian which is the case in the case of a white noise. Then it can be extended to the Malliavin--Sobolev classes ${\mathbb W}_p^1$  (see \cite {Bog}, Section 5.5 for the case $p=2$, and \cite{JO}, Lemma 3.1 for a generalization if $p>2$). 



\subsection{Schauder estimates and reconstruction theorem for germs}\lab{Shaud}

\setcounter{equation}{0}
\setcounter{theorem}{0}

Following the strategy of \cite{BOT,LOTT} for establishing stochastic estimates for models we shall use stochastic versions of Schauder estimates for the operator $\K$ acting on germs which are similar to those presented in \cite{BCZ} in the deterministic case (see also \cite{FH}, Chapter 14). These estimates are closely related to multilevel Schauder estimates which initially appeared in \cite{H}. 

Another important ingredient used in \cite{BOT,LOTT} is a stochastic version of the reconstruction theorem for germs. We state and use it in a form similar to Theorem 3.2 in \cite{BL}. This theorem is, in turn, a generalization of the reconstruction theorem suggested in \cite{H} for modelled distributions (see also \cite{CZ,HL} for other versions of the reconstruction theorem).

Note that a major difference between the Schauder estimates and the reconstruction theorem used in \cite{BOT} and analogous statements in this section is that in the latter case the rescaling parameter, which naturally appears in both cases, takes values in a bounded interval (compare with the discussion in Section \ref{stat} of the scaling properties of the white noise and of less rough noises utilized in \cite{BOT}).

Recall that a $p$--germ, $1\leq p\leq \infty$, is a family $F=(F_x)_{x\in \R^4}$ of $E$--valued random distributions $F_x$ such that for any $\varphi\in \D$ the map
$$
x\mapsto \left\|F_x(\varphi)\right\|_{\L_p} 
$$
is measurable.

If the random distributions $F_x$ are smooth, we call the germ $F$ smooth.

For a $p$--germ $F$, $1\leq p,q\leq \infty$, $a,b, c\in \R$, $a\leq c$, $\bar{\lambda}>0$, $r\in \N$, any measurable subset $\Omega\subset \R^4$, and $\varphi\in \D$ we define
$$
\left\|F\right\|_{pq\Omega \bar{\lambda} \varphi}^b:=\sup_{0<\lambda\leq \bar{\lambda}}\frac{\left\|F_x(\varphi_x^\lambda)\right\|_{pq\Omega x}}{\lambda^b},
$$
$$
\left\|F\right\|_{pq\Omega \bar{\lambda} \varphi}^{ac}:=\sup_{h\in \R^4:|h|\leq 2\bar{\lambda}}\sup_{0<\lambda\leq \bar{\lambda}}\frac{\left\|(F_{x+h}-F_x)(\varphi_x^\lambda)\right\|_{pq\Omega x}}{\lambda^a(\lambda+|h|)^{c-a}},
$$
\be\lab{bbnr}
\left\|F_x\right\|_{p\bar{\lambda}r}^b:=\sup_{0<\lambda\leq \bar{\lambda}}\sup_{\varphi\in \mB^r}\frac{\left\|F_x(\varphi^\lambda_x)\right\|_{\L_p}}{\lambda^b},
\ee
\be\lab{bcnr}
\left\|F_x\right\|_{p\bar{\lambda}r}^{\circ b}:=\sup_{0<\lambda\leq \bar{\lambda}}\sup_{\varphi\in \mB^r_b}\frac{\left\|F_x(\varphi^\lambda_x)\right\|_{\L_p}}{\lambda^b},
\ee
\be\lab{bonr}
\left\|F\right\|_{pq\bar{\lambda}\Omega r}^b:=\sup_{0<\lambda\leq \bar{\lambda}}\frac{\left\|\sup_{\varphi\in \mB^r}\left\|F_x(\varphi^\lambda_x)\right\|_{\L_p}\right\|_{L_q(\Omega,x)}}{\lambda^b},
\ee
\be\lab{bocnr}
\left\|F\right\|_{pq\bar{\lambda}\Omega r}^{\circ b}:=\sup_{0<\lambda\leq \bar{\lambda}}\frac{\left\|\sup_{\varphi\in \mB^r_b}\left\|F_x(\varphi^\lambda_x)\right\|_{\L_p}\right\|_{L_q(\Omega,x)}}{\lambda^b}.
\ee


We shall need the following pointwise and local Schauder type estimates for the action of the operator $\K$ on germs in Besov type seminorms.

\bp\lab{emb}
For any $1\leq p,q\leq \infty$, $b\in \R$, $r>|b+2|$, $\bar{\lambda}>0$, $x\in \R^4$, any compact subset $\mK\subset \R^4$, and any $p$--germ $F$ one has
\be\lab{emb1}
\left\|\K F_x\right\|_{p\bar{\lambda}r}^{\circ b+2}\lesssim \left\|F_x\right\|_{p\bar{\lambda}'r}^{\circ b}\leq \left\|F_x\right\|_{p\bar{\lambda}'r}^b,~\bar{\lambda}'=\max\{2\bar{\lambda},2\},
\ee
and
\be\lab{emb2}
\left\|\K F\right\|_{pq\bar{\lambda}\mK r}^{\circ b+2}\lesssim \left\|F\right\|_{pq\bar{\lambda}'\mK r}^{\circ b}\leq \left\|F\right\|_{pq\bar{\lambda}'\mK r}^b,
\ee
where $(\K F_x)(\varphi)=F_x(\K^-\varphi)$, $\varphi\in \D$, $\K^-$ is the integral operator with the kernel $K^-(x,y):=K(-x+y)$, the constants in the first inequalities in (\ref{emb1}) and (\ref{emb2}) only depend on $b$ and $r$. 
\ep

Note that the second inequalities in (\ref{emb1}) and (\ref{emb2}) are trivial, and the first ones can be obtained following the arguments in the proof of Theorem 14.17 in \cite{FH} verbatim with $\zeta=F_x$, where, in the discussion after the identity (14.7), instead of the modulus of real numbers one has to use the seminorm $\sup_{\varphi\in \mB^r_{b+2}}\left\|F_x(\varphi^\lambda_x)\right\|_{\L_p}$ in the case of (\ref{emb1}) and the seminorm $\left\|\sup_{\varphi\in \mB^r_{b+2}}\left\|F_x(\varphi^\lambda_x)\right\|_{\L_p}\right\|_{L_q(\mK,x)}$ in the case of (\ref{emb2}). A similar result is also stated and proved as Theorem 5.4 in \cite{BCZ}.

In view of definition (\ref{defP}) we shall also need the following weakened stochastic version of Theorem 5.10 in \cite{BCZ} which allows to estimate seminorms (\ref{bcnr}) and (\ref{bocnr}) in terms of (\ref{bbnr}) and (\ref{bonr}), respectively.
\bp\lab{bcest}
Let $F$ be a smooth $p$--germ, $1\leq p\leq \infty$. Then for any $1\leq q\leq \infty$, $b\in \R$, $b\geq 0$, $b\not\in \N$, $r>b$, $\bar{\lambda}>0$, $x\in \R^4$, any compact subset $\mK\subset \R^4$ one has
\be\lab{emb1'}
\left\|F^b_x\right\|_{p\bar{\lambda}r}^b\lesssim \left\|F_x\right\|_{p\bar{\lambda}r}^{\circ b},
\ee
and
\be\lab{emb2'}
\left\| F^b\right\|_{pq\bar{\lambda}\mK r}^b\lesssim \left\|F\right\|_{pq\bar{\lambda}\mK r}^{\circ b},
\ee
where the constants in the inequalities in (\ref{emb1'}) and (\ref{emb2'}) only depend on $b$ and $r$, and $F^b$ is the smooth germ defined by
$$
F^b_x(y):=F_x(y)-\sum_{\n\in \N^4:|\n|<b}\partial^\n F_x(x)\frac{(y-x)^\n}{\n!}.
$$

\ep

The proof of this proposition follows verbatim the arguments in the proof of Theorem 5.10 in \cite{BCZ} with the moduli of real numbers replaced by appropriate seminorms as described in the discussion after Proposition \ref{emb}.

From Propositions \ref{emb} and \ref{bcest} we immediately obtain the following corollary.
\bc
For any $1\leq p,q\leq \infty$, $b\not\in \mathbb{Z}$, $r>b+2\geq 0$, $\bar{\lambda}>0$, $x\in \R^4$, any compact subset $\mK\subset \R^4$, and any smooth $p$--germ $F$ one has
\be\lab{emb1c}
\left\|(\K F)_x^{b+2}\right\|_{p\bar{\lambda}r}^{ b+2}\lesssim \left\|F_x\right\|_{p\bar{\lambda}'r}^b,~\bar{\lambda}'=\max\{2\bar{\lambda},2\},
\ee
and
\be\lab{emb2c}
\left\|(\K F)^{b+2}\right\|_{pq\bar{\lambda}\mK r}^{ b+2}\lesssim  \left\|F\right\|_{pq\bar{\lambda}'\mK r}^b,
\ee
where the constants in the inequalities in (\ref{emb1c}) and (\ref{emb2c}) only depend on $b$ and $r$. 
\ec

Now we turn to the reconstruction theorem for germs.
The following proposition is an obvious stochastic version of Theorem 3.2 in \cite{BL}.
\bp
Let $F$ be a $p$--germ, $1\leq p,q\leq \infty$, $a,b, c\in \R$, $a\leq c$, $c\neq 0$, $\bar{\lambda}>0$. Assume that there exists a function $\varphi\in \D$ such that $\int_{\R^4}\varphi(x)dx\neq 0$, and for any compact subset $\mK\subset \R^4$
$$
\left\|F\right\|_{pq\mK \bar{\lambda} \varphi}^b<\infty.
$$
Suppose also that
\be\lab{cohcdn1}
\left\|F\right\|_{pq\mK \bar{\lambda} \varphi}^{ac}<\infty.
\ee

Then there exists an $E$--valued random distribution $\mathcal{R}(F)$ such that for any $r\in \N$, $r>\max\{-a,-b\}$ and any compact subset $\mK\subset \R^4$ one has
$$
\left\|\mathcal{R}(F)-F\right\|_{pq\bar{\lambda}\mK r}^c\lesssim \left\|F\right\|_{pq\mK_{2\bar{\lambda}} \bar{\lambda} \varphi}^{ac},
$$
where $\mathcal{R}(F)$ is regarded as a germ independent of the base point, i.e. $\mathcal{R}(F)_x=\mathcal{R}(F)$, for all $x\in \R^4$, the constant in the inequality is independent of $\mK$, and $\mK_{2\bar{\lambda}}:=\mK+B_{2\bar{\lambda}}=\{x+y:x\in \mK,y\in B_{2\bar{\lambda}} \}$. 

$\mathcal{R}(F)$ is unique if $c>0$, and if $F_x(x)$, $x\in \R^4$ is a continuous $E$--valued random distribution then  $\mathcal{R}(F)(x)=F_x(x)$.
\ep

The proof of this proposition can be obtained using the arguments in the proof of Theorem 3.2 in \cite{BL} verbatim with the moduli of real numbers replaced by appropriate $\L_p$--norms.

For a $p$--germ $F$, $1\leq p,q\leq \infty$, $a,b,b', c,c'\in \R$, $a\leq c$, $c'\geq 0$ $\bar{\lambda}>0$, $r\in \N$, any compact subset $\mK\subset \R^4$, and $\varphi\in \D$ we define
$$
\left\|F\right\|_{pq R \bar{\lambda} \varphi}^{acc'}:=\sup_{h\in \R^4:|h|\leq 2\bar{\lambda}}\sup_{0<\lambda\leq \bar{\lambda}}\frac{\left\|(F_{x+h}-F_x)(\varphi_{x+h}^\lambda)\right\|_{pq B_R x}}{\lambda^a(\lambda+|h|)^{c-a}(\lambda+R+|h|)^{c'}},
$$
$$
\left\|F\right\|_{pq\bar{\lambda} R r}^{bb'}:=\sup_{0<\lambda\leq \bar{\lambda}}\frac{\left\|\sup_{\varphi\in \mB^r}\left\|F_x(\varphi^\lambda_x)\right\|_{\L_p}\right\|_{L_q(B_R,x)}}{\lambda^b(\lambda+R)^{b'}}.
$$

The next proposition is another version of Theorem 3.2 in \cite{BL} the proof of which is obtained by a minor modification of the proof of the reconstruction bound in \cite{BL}, Section 4.4, replacing the moduli of real numbers by appropriate $\L_p$--norms.
\bp\lab{recprop}
Let $F$ be a $p$--germ, $1\leq p,q\leq \infty$, $a, c,c'\in \R$, $a\leq c$, $c\neq 0$, $c'\geq 0$, $\bar{\lambda}>0$. Let $F_x(x)$, $x\in \R^4$ be a continuous $E$--valued random distribution and $\varphi\in \D$ a function such that $\int_{\R^4}\varphi(x)dx\neq 0$. 

Assume that for some $R>0$ 
\be\lab{cohcdn2}
\left\|F\right\|_{pq R' \bar{\lambda} \varphi}^{acc'}\leq C, 0<R'\leq R+\bar{\lambda}
\ee
for some $C>0$ independent of $R'$ and $R$.

Then for any $r\in \N$, $r>-a$ one has
$$
\left\|\mathcal{R}(F)-F\right\|_{pq\bar{\lambda}R r}^{cc'}\lesssim C,
$$
where $\mathcal{R}(F)(x)=F_x(x)$ is regarded as a germ independent of the base point, i.e. $\mathcal{R}(F)_x=\mathcal{R}(F)$, for all $x\in \R^4$, and the constant in the inequality is independent of $R$. 

\ep

Conditions of type (\ref{cohcdn1}) and (\ref{cohcdn2}) are called coherence conditions.

From the previous proposition we immediately obtain the following obvious corollary.
\bc\lab{reccor}
Assume that the conditions of the previous proposition are satisfied, and $\mathcal{R}(F)(x)=F_x(x)=0$. Then for $r\in \N$, $r>-a$ one has
\be\lab{germ1c}
\left\|F\right\|_{pq\bar{\lambda}R r}^{cc'}\lesssim C,
\ee
where $C$ is defined in the previous proposition and the constant in the inequality is independent of $R$.

\ec

\subsection{Differentiation of the recentered maps}\lab{Diffrec}

\setcounter{equation}{0}
\setcounter{theorem}{0}

In this section we obtain some results related to Malliavin derivatives of the recentered maps required in the proof of stochastic estimates for them. We refer the reader to \cite{BOT}, Section 2.5 for motivations of the constructions in this section.

First note that by the definition the random variables $\Pi_{x\beta}$ and $\Pi_{x\beta}^-$ depend polynomially on the white noise space variable $\xi$ via the mollification $\xi^\rho(x)=\xi(\eta_x^\rho)$. Definition (\ref{Ph0}) and induction over $|\beta|$ using formulas (\ref{defP-}), (\ref{defP}) show that $\Pi_{x\beta}$ and $\Pi_{x\beta}^-$ have directional derivatives of polynomial growth along all directions in the Cameron--Martin space (see the discussion in \cite{HS}, Section 2.2.2. for more details). Therefore $\Pi_{x\beta}$ and $\Pi_{x\beta}^-$ are Malliavin differentiable, and the spectral gap inequality is satisfied for them (see Lemma 2.23 in \cite{HS}). We shall frequently use this fact without comments. 

Moreover, $\delta_\upsilon\Pi_{x\beta}$ and $\delta_\upsilon\Pi_{x\beta}^-$ are well defined for random $\upsilon$ with values in ${\rm CM}$. We shall use the symbols $\delta_\upsilon\Pi_{x\beta}$ and $\delta_\upsilon\Pi_{x\beta}^-$ in this sense without comments. 
In fact, since the random variables $\Pi_{x\beta}$ and $\Pi_{x\beta}^-$ depend polynomially on the white noise space variable $\xi$ via the mollification $\xi^\rho(x)$, $\delta_\upsilon\Pi_{x\beta}$ and $\delta_\upsilon\Pi_{x\beta}^-$ take values in $\O_\beta$ and $\O_\beta^-$, respectively, for all random $\upsilon$ with values in ${\rm CM}$.


Note also that by (\ref{defP}) 
$$
(\delta\Pi_{x\beta})(y)=(\K\delta \Pi_{x\beta}^-)(y)-\sum_{{\n}\in \N^4: |{\n}|<|\beta|}\partial^{\n}(\K\delta \Pi_{x\beta}^-)(x)\frac{(y-x)^{\n}}{{\n}!}.
$$

Below for linear endomorphisms of $\R[[\z]]$ and induced endomorphisms of $\Sm[[\z]]$  we use the notation and conventions of Section \ref{Frec}. 

We proceed with the following algebraic statement which will be used for the study of the properties of the Malliavin derivative of $\Pi_{0\beta}$. 
\bp\lab{F*ext}
For any $x\in \R^4$, $f_\n\in V\otimes\R[[\z]]$, $\n\in \N^4$, $|\n|<2-\e$, where $\e$ is as in the definition of $\alpha=-\frac12 -\e$, there exists a unique continuous linear endomorphism $F^*$ of $\R[[\z]]$ such that
\be\lab{cdn1}
F^*1=0,~F^*\z_\bg=0,
\ee
\be\lab{cdn2}
({\rm Id}_{W_{\delta_\n}}^V\otimes F^*)\z^{\delta_\n}=\left\{\begin{array}{ll} 0 &{\rm if}~ \n\in \N^4,~|\n|\geq 2-\e \\ f_\n &{\rm if}~ \n\in \N^4,~|\n|< 2-\e \end{array}\right. ,
\ee
\be\lab{F*hom}
F^*(ZZ')=F^*(Z)\G_{x0}^*(Z')+\G_{x0}^*(Z)F^*(Z')
\ee
for any $Z,Z'\in \R[[\z]]$.

$F^*$ satisfies
\be\lab{prp1}
(F)_0^\gamma=0~{\rm if}~\gamma\neq \delta_\n~{\rm for~ some}~\n\in \N^4,~|\n|< 2-\e,
\ee
and
\be\lab{prp2}
(F)_\beta^\gamma=0~{\rm if}~|\gamma|\geq |\beta|+\frac d2. 
\ee

\ep

We start the proof of this proposition with the following lemma which is similar to Lemma \ref{Gbgmult}.
\bl\lab{Fbg'}
For any linear map $F^*:\R[\z]\to \R[[\z]]$ which satisfies (\ref{cdn1}) and (\ref{F*hom}) with arbitrary $Z,Z'\in \R[\z]$ and any multi--index $\gamma=l\delta_\bg+\sum_{i=1}^p\delta_{{\m}_i}\neq 0$, $\m_i\in \N^4$, $l\in \N$  one has
\be\lab{GbF}
(F)_{\beta}^\gamma= \hspace{10cm}
\ee
$$
=  \sum_{k=1}^p\sum_{\tiny\begin{array}{c}\beta_1,\ldots,\beta_p\in \Ml: \\ \beta_1+\ldots+\beta_p+l\delta_\bg=\beta\end{array}} \hspace{-0.5cm}  (\G_{x0})_{\beta_1}^{\delta_{\m_1}}\ldots (\G_{x0})_{\beta_{k-1}}^{\delta_{\m_{k-1}}}(F)_{\beta_k}^{\delta_{\m_k}}(\G_{x0})_{\beta_{k+1}}^{\delta_{\m_{k+1}}}\ldots (\G_{x0})_{\beta_p}^{\delta_{\m_p}},
$$
where we assume that an empty sum is zero and an empty product is 1. 

Moreover, if $l>0$ then in the non--zero terms in the sum in (\ref{GbF}) $|\beta_i|<|\beta|$, $i=1,\ldots,p$.
\el

\bpr
Firstly, similarly to the proof of (\ref{Gb}), relation (\ref{GbF}) follows from (\ref{cdn1}) and (\ref{F*hom}).

From definition (\ref{hom}) it follows that $|\cdot|-\alpha$ is a linear function on $\Ml$ achieving its minimum equal to zero at $\beta=0$, and hence for each non--zero term in (\ref{GbF}) we have as in (\ref{bs})
\be\lab{bs'}
|\beta|-\alpha=\sum_{i=1}^p(|\beta_i|-\alpha)+l(\alpha+1),~|\beta_i|-\alpha\geq 0,~i=1,\ldots, p.
\ee

Since $(\G_{x0})_0^{\delta_{\m_i}}=0$, $i=1,\ldots,p$ by (\ref{triangG}) as $|\delta_{\m_i}|=|\m_i|\geq 0>|0|=\alpha$, we may assume that for each fixed $k$ and for any $i\neq k$ in the right hand side of (\ref{GbF}) one has $\beta_i\neq 0$, and hence by Lemma \ref{locf} $|\beta_i|-\alpha>0$. Note also that $l(\alpha+1)> 0$ if $l>0$. 
Therefore the term $l(\alpha+1)$ is strictly positive in the right hand side of (\ref{bs'}) which implies that for $p>1$ one has $|\beta_i|-\alpha<|\beta|-\alpha$, $i=1,\ldots,p$, or $|\beta_i|<|\beta|$, $i=1,\ldots,p$, and if $p=1$ then $|\beta|=|\beta_1|+l(\alpha+1)>|\beta_1|$.

\epr

\noindent
{\em Proof of Proposition \ref{F*ext}.} Firstly, conditions (\ref{cdn1}), (\ref{cdn2}) and formula (\ref{GbF}) uniquely define a linear map $F^*:\R[\z]\to \R[[\z]]$ which satisfies (\ref{F*hom}) with arbitrary $Z,Z'\in \R[\z]$. Using Lemma \ref{Flemc} (i) we show that this map can be extended to a continuous linear map $F^*:\R[[\z]]\to \R[[\z]]$, so that (\ref{F*hom}) will hold with arbitrary $Z,Z'\in \R[[\z]]$ by continuity. In order to verify condition (\ref{triangF'}) it suffices to show that $F^*:\R[\z]\to \R[[\z]]$ satisfies (\ref{prp1}) and (\ref{prp2}) which will also immediately imply that the same properties hold for the extension $F^*:\R[[\z]]\to \R[[\z]]$.

We start with (\ref{prp1}). If $\gamma=0$ then $(F)_0^0=(F^*1)_0=0$ by (\ref{cdn1}).

If $\gamma\neq 0$ then the right hand side of (\ref{GbF}) with $\beta=0$ does not vanish only if $l=0$ and $\beta_i=0$, $i=1,\ldots,p$. But in this case $(\G_{x0})_0^{\delta_{\m_i}}=0$, $i=1,\ldots,p$ by (\ref{triangG}) as $|\delta_{\m_i}|=|\m_i|\geq 0>|0|=\alpha$. Thus the right hand side of (\ref{GbF}) with $\beta=0$ does not vanish only if $l=0$ and $p=1$, i.e. if $\gamma=\delta_\n$ for some $\n\in \N^4$, and by (\ref{cdn2}) in this case $|\n|< 2-\e$.

Now we prove (\ref{prp2}).
If $\gamma=0$ property (\ref{prp2}) is obvious as by (\ref{cdn1})
$(F)_\beta^0=(F^*1)_\beta=0$ for all $\beta$.

If $\gamma=\delta_\n$ for some $\n\in \N^4$, then by (\ref{cdn2}) $(F)_\beta^\gamma=(F)_\beta^{\delta_\n}$ may be non--zero only if $|\n|< 2-\e$, and in this case for any $\beta$
$$
|\gamma|=|\delta_\n|=|\n|<2-\e=\alpha+\frac d2\leq |\beta|+\frac d2,
$$
where at the last step we used Lemma \ref{locf}. This proves (\ref{prp2}) for $\gamma\in \M_{pp}$.

Now for each fixed $\gamma\neq 0$, $\gamma\not\in \M_{pp}$, $\gamma(\bg)>0$ we prove (\ref{prp2}) by induction over $|\beta|$.

If $|\beta|=\alpha$ is minimal possible then $\beta=0$ by Lemma \ref{locf} and (\ref{prp2}) follows from (\ref{prp1}) as $\gamma\not\in \M_{pp}$.


Now assume that (\ref{prp2}) holds for all $\beta'\in \Ml$ with $|\beta'|<|\beta|$, $\beta\in \Ml$. Note that in the considered case $\gamma\not\in \M_{pp}$, $\gamma(\bg)>0$ and $\gamma\neq 0$. Therefore $\gamma=l\delta_\bg+\sum_{i=1}^p\delta_{{\m}_i}\neq 0$, $\m_i\in \N^4$, $l\in \N$, where $l>0$. Thus by Lemma \ref{Fbg'} in the right hand side of (\ref{GbF}) $|\beta_i|<|\beta|$, $i=1,\ldots,p$, and hence by (\ref{triangG}), by the induction assumption, by the already established (\ref{prp2}) in the case of polynomial $\gamma$, and by (\ref{bs'})
$$
|\beta|-\alpha=\sum_{i=1}^p(|\beta_i|-\alpha)+l(\alpha+1)>\sum_{i=1}^p(|\delta_{\m_i}|-\alpha)+|l\delta_\bg|-\alpha-\frac d2=|\gamma|-\frac d2-\alpha,
$$
where at the last step we used linearity of $|\cdot|-\alpha$. Thus if $(F)^\gamma_\beta\neq 0$ then $|\beta|>|\gamma|-\frac d2$.
This establishes the induction step and proves (\ref{prp2}) for $\gamma\neq 0$, $\gamma\not\in \M_{pp}$, $\gamma(\bg)>0$.

Finally, if $\gamma=\sum_{i=1}^p\delta_{{\m}_i}\neq 0$, $\m_i\in \N^4$ then each in non--zero term in the right hand side of (\ref{GbF}) labeled by $k$, for $i\neq k$ one has 
\be\lab{leq1F}
|\delta_{{\m}_i}|\leq |\beta_i|
\ee 
by (\ref{triangG}). By the already established (\ref{prp2}) in the case of polynomial $\gamma$ we obtain 
\be\lab{leq2F}
|\delta_{{\m}_k}|< |\beta_k|+\frac d2.
\ee

Adding (\ref{leq1F}) over all $i\neq k$ and (\ref{leq2F}) and using linearity of $|\cdot|-\alpha$ on $\Ml$ we obtain
$$
|\gamma|-\alpha=\sum_{i=1}^p(|\delta_{{\m}_i}|-\alpha)<\sum_{i=1}^p(|\beta_i|-\alpha)+\frac d2=|\beta|-\alpha+\frac d2,
$$
i.e.
$$
|\gamma|<|\beta|+\frac d2.
$$
This justifies (\ref{prp2}) for $\gamma=\sum_{i=1}^p\delta_{{\m}_i}\neq 0$, $\m_i\in \N^4$ and completes the proof.

\qed

Now we state the main result of this section. 
\bp\lab{G*ext}
For any random variable $\upsilon$ with values in ${\rm CM}$ and $x\in \R^4$ there exists a unique continuous linear endomorphism ${\rm d}\G_x^*$ of $\R[[\z]]$ which satisfies (\ref{cdn1}), (\ref{F*hom}), (\ref{prp1}), (\ref{prp2}), and
\be\lab{cdn30}
 {\rm d}\G_x^*\z^{\delta_{\n}}=0,~\n\in \N^4,~|\n|\geq 2-\e,
\ee
\be\lab{cdn3}
(\partial^\n\delta\Pi_0)(x)(\z)-({\rm d}\G_x^*\partial^\n\Pi_x(x))(\z)=0~{\rm for}~|\n|<2-\e,
\ee
where $\delta=\delta_\upsilon$, $\e$ is as in the definition of $\alpha=-\frac12 -\e$.
\ep

\bpr
Fix $x\in \R^4$. First observe that for any $\n\in \N^4$ by (\ref{Pvan}) one has
\be\lab{Pvan1}
\partial^\n(\Pi_{x\beta})(x)=0 ~{\rm for~ any }~ x\in \R^4,~\beta\in \Ml,~\n\in \N^4 ~{\rm such~ that }~|\n|<|\beta|,
\ee
and by (\ref{poph}) for $\beta\in \M_{\geq 0}\cup \M_{pp}$ and $|\beta|=|\n|$ one has $\beta=\delta_\m$ for some $\m\in \N^4$. Therefore using (\ref{Pi0}) and definition (\ref{Pipol}) we obtain in this case
\be\lab{Pipol1}
\frac{1}{\n!}\partial^\n\Pi_{x\beta}(y)=\delta_\beta^{\delta_\n}{\rm Id}_{W_{\delta_\n}}^V,
\ee
where for $\beta,\gamma\in \Ml$ we define
$$
\delta_\beta^\gamma=\left\{\begin{array}{ll} 1 & {\rm if}~\beta=\gamma \\ 0 & {\rm otherwise}\end{array}\right. .
$$

We can rewrite (\ref{Pvan1}) and (\ref{Pipol1}) as a single identity,
\be\lab{Yid}
(\frac{1}{\n!}\partial^\n\Pi_x(x)(\z)-({\rm Id}_{W_{\delta_\n}}^V\otimes {\rm Id}_{\R[[\z]]})\z^{\delta_\n})_\beta=0 
\ee
for any  $\beta\in \Ml$, $\n\in \N^4$ such that $|\n|\leq|\beta|$.

Note that if $\beta(\m)>0$ for some $\m\in \N^4$ then by (\ref{hom}) $|\beta|\geq |\m|$.  Thus  $|\n|>|\beta|$ implies $\beta(\m)=0$ for $\m\in \N^4$, $|\m|\geq |\n|$, and hence from (\ref{Yid}) we immediately obtain
\be\lab{Ydef}
Y_\n(x)(\z):=\frac{1}{\n!}\partial^\n\Pi_x(x)(\z)-({\rm Id}_{W_{\delta_\n}}^V\otimes {\rm Id}_{\R[[\z]]})\z^{\delta_\n}\in 
\ee
$$
\in V[[\z_\bg,\{\z_\m^i:i=1,\ldots,{\rm dim}(V),\m\in \N^4,|\m|<|\n|\}]].
$$

Now we define ${\rm d}\G_x^*$ using this property, Proposition \ref{F*ext} and induction over $r:=|\n|\in \N$ as follows.

First observe that if $\n\in \N^4$ and $|\n|=0$ then $\n=0$. Thus (\ref{Ydef}) reduces to
$$ 
Y_{\bf 0}(x)(\z)=\Pi_x(x)(\z)-({\rm Id}_{W_{\delta_{\bf 0}}}^V\otimes {\rm Id}_{\R[[\z]]})\z^{\delta_{\bf 0}}\in V[[\z_\bg]],
$$
and we can define a unique continuous endomorphism ${\rm d}\G_x^*$ of $V[[\z_\bg]]$ using (\ref{cdn1}). In fact, this way we obtain the zero endomorphism. Therefore by Proposition \ref{F*ext} we can also define a unique continuous endomorphism of $\R[[\z]]$, which we denote by the same symbol ${\rm d}\G_x^*$, by requiring that 
$$
({\rm Id}_{W_{\delta_{\bf 0}}}^V\otimes {\rm d}\G_x^*)\z^{\delta_{\bf 0}}=\delta\Pi_0(x)(\z)-({\rm d}\G_x^*Y_{\bf 0}(x))(\z)=\delta\Pi_0(x)(\z),
$$
and 
$$
({\rm Id}_{W_{\delta_{\n}}}^V \otimes {\rm d}\G_x^*)\z^{\delta_{\n}}=0,~\n\neq {\bf 0}.
$$
${\rm d}\G_x^*$ also satisfies (\ref{cdn1}), (\ref{F*hom}), (\ref{prp1}) and (\ref{prp2}).

Now assume that for $0<r<2-\e$ 
a continuous endomorphism ${\rm d}\G_x^*$ of $\R[[\z]]$ is defined and satisfies (\ref{cdn1}), (\ref{F*hom}), (\ref{prp1}), (\ref{prp2}),
\be\lab{Gind}
({\rm Id}_{W_{\delta_\m}}^V\otimes {\rm d}\G_x^*)\z^{\delta_\m}=
\ee
$$
=\frac{1}{\m!}(\partial^\m\delta\Pi_0)(x)(\z)-({\rm d}\G_x^*Y_\m(x))(\z),~\m\in \N^4,~|\m|<r,
$$
and 
$$
({\rm Id}_{W_{\delta_\m}}^V \otimes {\rm d}\G_x^*)\z^{\delta_{\m}}=0,~\m\in \N^4,~|\m|\geq r. 
$$

Then by Proposition \ref{F*ext} we can define another unique continuous endomorphism of $\R[[\z]]$, which we denote by the same symbol ${\rm d}\G_x^*$, by requiring that
$$
({\rm Id}_{W_{\delta_\m}}^V\otimes {\rm d}\G_x^*)\z^{\delta_\m}~\m\in \N^4,~|\m|<r 
$$
coincide with the values of the endomorphism ${\rm d}\G_x^*$ defined at the previous step of the induction, for $\n\in \N^4$, $|\n|=r$
$$
({\rm Id}_{W_{\delta_\n}}^V\otimes {\rm d}\G_x^*)\z^{\delta_\n}=\frac{1}{\n!}(\partial^\n\delta\Pi_0)(x)(\z)-({\rm d}\G_x^*Y_\n(x))(\z),
$$
where in the right hand side in the expression $({\rm d}\G_x^*Y_\n(x))(\z)$ we use the action of ${\rm d}\G_x^*$ on $\R[[z]]$ defined at the previous step of the induction, and 
$$
({\rm Id}_{W_{\delta_\m}}^V \otimes {\rm d}\G_x^*)\z^{\delta_{\m}}=0,~\m\in \N^4,~|\m|> r. 
$$

By construction the action of the newly defined ${\rm d}\G_x^*$ introduced this way on $V[[\z_\bg,\{\z_\m^i:i=1,\ldots,{\rm dim}(V),\m\in \N^4,|\m|< r\}]]$ 
coincides with the action of ${\rm d}\G_x^*$ introduced at the previous step of the induction, and hence the newly defined ${\rm d}\G_x^*$ still satisfies (\ref{Gind}) as by (\ref{Ydef}) for $\n\in \N^4$, $|\n|=r$
$Y_\n(x)(\z)\in  V[[\z_\bg,\{\z_\m^i:i=1,\ldots,{\rm dim}(V),\m\in \N^4,|\m|<|\n|=r\}]]$.
Thus by construction for the newly defined ${\rm d}\G_x^*$ we have
$$
({\rm Id}_{W_{\delta_\m}}^V\otimes {\rm d}\G_x^*)\z^{\delta_\m}=
$$
$$
=\frac{1}{\m!}(\partial^\m\delta\Pi_0)(x)(\z)-({\rm d}\G_x^*Y_\m(x))(\z),~\m\in \N^4,~|\m|\leq r, 
$$
$$
({\rm Id}_{W_{\delta_\m}}^V \otimes {\rm d}\G_x^*)\z^{\delta_{\m}}=0,~\m\in \N^4,~|\m|> r, 
$$

By Proposition \ref{F*ext} the newly defined ${\rm d}\G_x^*$ also satisfies (\ref{cdn1}), (\ref{F*hom}), (\ref{prp1}) and (\ref{prp2}).

Then we can proceed by induction until we reach $r$ corresponding to $\n\in \N^4$ with the largest possible $r=|\n|<2-\e$. By construction the continuous endomorphism ${\rm d}\G_x^*$ of $\R[[\z]]$ defined at the last step satisfies (\ref{cdn1}), (\ref{F*hom}), (\ref{prp1}), (\ref{prp2}),
\be\lab{YY}
({\rm Id}_{W_{\delta_\n}}^V\otimes {\rm d}\G_x^*)\z^{\delta_\n}=
\ee
$$
=\frac{1}{\n!}(\partial^\n\delta\Pi_0)(x)(\z)-({\rm d}\G_x^*Y_\n(x))(\z),~\n\in \N^4,~|\n|<2-\e,
$$
and 
$$
({\rm Id}_{W_{\delta_\n}}^V \otimes {\rm d}\G_x^*)\z^{\delta_{\n}}=0,~\n\in \N^4,~|\n|\geq 2-\e.
$$
Therefore for $\n\in\N^4$, $|\n|<2-\e$ one has from (\ref{YY})
$$
\frac{1}{\n!}(\partial^\n\delta\Pi_0)(x)(\z)=({\rm Id}_{W_{\delta_\n}}^V\otimes {\rm d}\G_x^*)\z^{\delta_\n}+({\rm d}\G_x^*Y_\n(x))(\z)=\frac{1}{\n!}({\rm d}\G_x^*\partial^\n\Pi_x(x))(\z),
$$
where at the last step we used (\ref{Ydef}). This is exactly identity (\ref{cdn3}). Note that by (\ref{Ydef}) properties (\ref{cdn3}) and (\ref{YY}) are equivalent.

Uniqueness of ${\rm d}\G_x^*$ follows from Proposition \ref{F*ext} and the inductive definition of ${\rm d}\G_x^*$. This completes the proof.

\epr

\br
Remark \ref{indepbas} and the definition of ${\rm d}\G_x^*$ imply that the linear maps $({\rm d}\G_x)_\beta^\gamma$, $\beta,\gamma\in \Ml$ do not depend on the choice of the basis $f_i$, $i=1,\ldots, {\rm dim}(V)$ in Section \ref{Frec}.
\er

We shall also need the following analogue of property (\ref{cdn3}) for $\Pi_0^-(x)(\z)$ and $\Pi_x^-(x)(\z)$.
\bp\lab{dGG-dif}
For any random variable $\upsilon$ with values in ${\rm CM}$, $x\in \R^4$ the endomorphism ${\rm d}\G_x^*$ satisfies 
\be\lab{cdn3-}
\delta\Pi_0^-(x)(\z)-({\rm d}\G_x^*\Pi_x^-(x))(\z)-\upsilon^\rho(x)=0,
\ee
where $\delta=\delta_\upsilon$, and $\upsilon^\rho(x)=\upsilon(\eta_x^\rho)$.
\ep

\bpr
First observe that (\ref{Ph0}) and (\ref{defP-}) for $\Pi_{x\beta}$, $\Pi_{x\beta}^-$ can be rewritten as a single identity,
\be\lab{defP-z}
\Pi_{x}^-(y)(\z)=\z_\bg \A(\Pi_{x}(y)(\z),\Pi_{x}(y)(\z))+
\ee
$$
+\z_\bg^2\mB(\Pi_{x}(y)(\z),\Pi_{x}(y)(\z),\Pi_{x}(y)(\z))+c(\z)\Pi_{x}(y)(\z)+\xi^\rho(y),
$$
where
$$
c(\z):=\sum_{k=1}^4c_k\z_\bg^k.
$$

Applying $\delta$ to (\ref{defP-z}) with $x=0$ and $y=x$ and using (\ref{cdn3}) with $\n={\bf 0}$, and $\n\in \N^4$, $|\n|=1$ we obtain
\be\lab{defP-z1}
\delta\Pi_0^-(x)(\z)=\z_\bg \A(({\rm d}\G_x^*\Pi_x(x))(\z),\Pi_0(x)(\z))+
\ee
$$
+\z_\bg \A(\Pi_0(x)(\z),({\rm d}\G_x^*\Pi_x(x))(\z))+\z_\bg^2\mB(({\rm d}\G_x^*\Pi_x(x))(\z),\Pi_0(x)(\z),\Pi_0(x)(\z))+ 
$$
$$
+\z_\bg^2\mB(\Pi_0(x)(\z),({\rm d}\G_x^*\Pi_x(x))(\z),\Pi_0(x)(\z))+\z_\bg^2\mB(\Pi_0(x)(\z),\Pi_0(x)(\z),({\rm d}\G_x^*\Pi_x(x))(\z))+
$$
$$
+c(\z)({\rm d}\G_x^*\Pi_x(x))(\z)+\upsilon^\rho(x).
$$

Using the fact that by (\ref{PG}) $\Pi_0(x)(\z)=(\G_{x0}^*\Pi_x(x))(\z)$ and properties (\ref{G1g}), (\ref{cdn1}), (\ref{F*hom}) we can reduce the right hand side of (\ref{defP-z1}) as follows
$$
\delta\Pi_0^-(x)(\z)={\rm d}\G_x^*\large\left( \z_\bg \A(\Pi_x(x)(\z),\Pi_x(x)(\z))+\right.
$$
$$
\left. +\z_\bg^2\mB(\Pi_x(x)(\z),\Pi_x(x)(\z),\Pi_x(x)(\z))
+c(\z)\Pi_x(x)(\z)\large\right)+\upsilon^\rho(x)=
$$
$$
=({\rm d}\G_x^*\Pi_x^-(x))(\z)+\upsilon^\rho(x),
$$
where at the last step we used (\ref{defP-z}) with $y=x$ and (\ref{cdn1}). This completes the proof.

\epr

Now we study some algebraic properties of the endomorphism ${\rm d}\G_x^*$.
\bp
For any random variable $\upsilon$ with values in ${\rm CM}$, and $x\in \R^4$ the endomorphism ${\rm d}\G_x^*$ satisfies the following properties:
\be\lab{tr1G*}
({\rm d}\G_x)_\beta^\gamma=0~{\rm for ~any}~\gamma\in \M,~\beta\in \Ml,~[\beta]<0,
\ee
\be\lab{tr2G*}
({\rm d}\G_x)_\beta^\gamma=0~{\rm for~ any}~\beta,\gamma\in \Ml,~|\gamma|_\prec\geq |\beta|_\prec.
\ee
\ep

\bpr
We start with (\ref{tr1G*}) which is obvious for $\gamma=0$ by the first identity in (\ref{cdn1}) and by the definition of $({\rm d}\G_x)_\beta^\gamma=({\rm d}\G_x^*\z^\gamma)_\beta$. 

For any $\gamma\in \M$ and $\beta\in \Ml$ with minimal possible $|\beta|=\alpha$ we have $\beta=0$, and the claim is void as $[0]=0$.

Now assume that for some $\beta\in \Ml$ the statement is true for all $\gamma\in \M$ and $\beta'\in \Ml$ with $|\beta'|<|\beta|$.

For $\gamma=l\delta_\bg+\sum_{i=1}^p\delta_{{\m}_i}\neq 0$ and $\gamma\not\in \M_{pp}$, i.e. $l>0$, we proceed by induction over $|\beta|$ the base of which is already established. 

We use (\ref{GbF}) with $F={\rm d}\G_x$, where by Lemma \ref{Fbg'} $|\beta_i|<|\beta|$, $i=1,\ldots,p$ as $l>0$. Thus by the induction assumption the k-th term in the sum over k in (\ref{GbF}) vanishes unless $[\beta_k]\geq 0$ in all summands which this term contains. By (\ref{tr3}) this term also vanishes unless $\beta_i\in \M$, $i=1,\ldots,k-1,k+1,\ldots,p$ in all summands which this term contains, and for all such $\beta_i$ one, in particular, has $[\beta_i]\geq -1$. Thus, by the additivity of $[~\cdot~]$, in this case $({\rm d}\G_x)_\beta^\gamma=0$ unless $[\beta]\geq l-p+1$. 

But by the definition of $\gamma$ one also has in this case $[\gamma]=l-p\geq -1$ as $\gamma\in \M$. Thus $({\rm d}\G_x)_\beta^\gamma=0$ unless $[\beta]\geq [\gamma]+1\geq 0$. This establishes the induction step if $\gamma\neq 0$ and $\gamma\not\in \M_{pp}$.

The case $\gamma=\delta_\n\in \M_{pp}$, $\n\in \N^4$ is trivial if $|\n|\geq 2-\e$ by (\ref{cdn30}), and for $|\n|<2-\e$ we use (\ref{cdn3}) and rewrite it in terms of components with the help of (\ref{Pi0}) as follows 
$$
\partial^\n\delta\Pi_{0\beta}(x)=\sum_{\mu\in \M_{\geq 0}\cup\M_{pp}}\partial^\n\Pi_{x\mu}(x)\circ ({\rm d}\G_x)_\beta^\mu,
$$
or, since 
$$
\partial^\n\Pi_{x\delta_\m}(x)=(\partial^\n(y-x)^\m)(x){\rm Id}_{W_{\delta_\m}}^V=\n!\delta_{\delta_\m}^{\delta_\n}{\rm Id}_{W_{\delta_\n}}^V,
$$
as
\be\lab{dGind}
\n!{\rm Id}_{W_{\delta_\n}}^V\circ({\rm d}\G_x)_\beta^{\delta_\n}=\partial^\n\delta\Pi_{0\beta}(x)-\sum_{\mu\in \M_{\geq 0}}\partial^\n\Pi_{x\mu}(x)\circ ({\rm d}\G_x)_\beta^\mu.
\ee

If $[\beta]<0$ then the first term in the right hand side of this formula vanishes as by (\ref{Pi0}) $\Pi_{0\beta}$ is non--zero in this case only if $\beta=\delta_\m$ for some $\m\in \N^4$, and by (\ref{Pipol}) in this case $\Pi_{0\delta_\m}=x^\m$ is independent of the white noise space variables, so $\delta\Pi_{0\beta}(x)=0$.

Also, for $[\beta]<0$ $({\rm d}\G_x)_\beta^\mu=0$ in the last sum in the right hand side of (\ref{dGind}) for all $\mu\in \M_{\geq 0}$ by the part of the induction step already proved as $\mu\in \M$ and $\mu\not\in \M_{pp}$. 

We deduce that if $[\beta]<0$ then the right hand side of (\ref{dGind}) vanishes, i.e. $({\rm d}\G_x)_\beta^{\delta_\n}=0$. This completes the proof of the induction step and of property (\ref{tr1G*}).

Now we turn to (\ref{tr2G*}) which trivially holds for $\gamma=l\delta_\bg$, $l\in \N$ by (\ref{cdn1}) and (\ref{F*hom}).

If $\gamma=\delta_\n$, $\n\in \N^4$ then ${\rm d}\G_x^*\z^{\delta_n}=0$ for $|\delta_n|_\prec=|\delta_n|=|\n|\geq 2-\e$ by (\ref{cdn30}). Thus $({\rm d}\G_x)_\beta^{\delta_\n}=0$ for any $\beta$ in this case.

If $|\delta_n|_\prec=|\delta_n|=|\n|< 2-\e$, by (\ref{tr1G*}) $({\rm d}\G_x)_\beta^{\delta_\n}=0$ for any $\beta$ with $[\beta]<0$, and if $[\beta]\geq 0$ we have by (\ref{homm})
$$
|\beta|_\prec=|\beta|+\frac{d}{2}([\beta]+1)\geq \alpha+\frac d2=2-\e,
$$
so that $|\beta|_\prec\geq 2-\e>|\delta_n|_\prec$ for possible non--zero $({\rm d}\G_x)_\beta^{\delta_\n}$. This proves (\ref{tr2G*}) for $\gamma\in \M_{pp}$.

If $\gamma=l\delta_\bg+\sum_{i=1}^p\delta_{{\m}_i}\neq 0$, $\m_i\in \N^4$, $l\in \N$ with $p>0$ 
we proceed using (\ref{GbF}) for $F={\rm d}\G_x$. Assume that $|\gamma|_\prec\geq |\beta|_\prec$ in (\ref{GbF}) with $F={\rm d}\G_x$. By the linearity of the function $|\cdot|_+:=|\cdot|_\prec-\frac d2 -\alpha$ this implies that for the k-th term in the sum in the right hand side of (\ref{GbF}) and for all summands which this term contains one has
$$
|l\delta_\bg|_+ +\sum_{i=1}^p|\delta_{{\m}_i}|_+=|\gamma|_+\geq |\beta|_+=|l\delta_\bg|_+ + \sum_{i=1}^p |\beta_i|_+,
$$
or
\be\lab{dbg}
\sum_{i=1}^p|\delta_{{\m}_i}|_\prec\geq \sum_{i=1}^p|\beta_i|_\prec,
\ee
where, if $p>1$, by (\ref{triangGm}) non--zero contributions to the right hand side of (\ref{GbF}) may only occur when $|\beta_i|_\prec\geq |\delta_{{\m}_i}|_\prec$, $i=1,\ldots,k-1,k+1,\ldots,p$. This implies
$$
\sum_{\tiny\begin{array}{c} i=1 \\ i\neq k\end{array}}^p|\delta_{{\m}_i}|_\prec\leq \sum_{\tiny\begin{array}{c} i=1 \\ i\neq k\end{array}}^p|\beta_i|_\prec
$$
which by (\ref{dbg}) yields $|\delta_{{\m}_k}|_\prec\geq |\beta_k|_\prec$ if $p>1$.

If $p=1$ then by (\ref{dbg}) the same inequality holds with $k=1$.

Since we already proved (\ref{tr2G*}) for $\gamma\in \M_{pp}$ this forces $({\rm d}\G_x)_{\beta_k}^{\delta_{{\m}_k}}=0$. Thus the right hand side of (\ref{GbF}) for $F={\rm d}\G_x$ vanishes. This completes the proof of (\ref{tr2G*}).

\epr

From (\ref{triangGm}), (\ref{tr3}), (\ref{tr1G*})  and (\ref{tr2G*}) we immediately deduce the following corollary.
\bc
Let $S_{xy}:={\rm d}\G_x-\G_{xy}{\rm d}\G_y$, $x,y\in \R^4$. Then
\be\lab{tr1S*}
(S_{xy})_\beta^\gamma=0~{\rm for~ any}~\gamma\in \M,~\beta\in \Ml,~[\beta]<0,
\ee
and
\be\lab{tr2S*}
(S_{xy})_\beta^\gamma=0~{\rm for~ any}~\beta,\gamma\in \Ml,~|\gamma|_\prec\geq |\beta|_\prec.
\ee
\ec

In conclusion we state some other algebraic properties of ${\rm d}\G_x$ and $S_{xy}$ which will be used in the proof of the stochastic estimates for the model.
\bl\lab{ltriang}
Let $\gamma=l\delta_\bg+\sum_{i=1}^p\delta_{{\m}_i}$, $\m_i\in \N^4$, $l\in \N$ be a multi--index with $[\gamma]\geq -1$, and $p>1$ or $l>0$. Then the following statements are true. 

(i) Let $F$ be an algebra endomorphism of $\R[[\z]]$ which satisfies (\ref{triang1}). Assume also that
\be\lab{tr1F}
\mu\in \M_{\geq 0}\cup \M_{pp},~\beta\not\in \M_{\geq 0}\cup \M_{pp}\Longrightarrow (F)_{\beta}^{\mu} =0,
\ee  
Then in the non--zero terms in the right hand side of (\ref{Gb}) one has $[\beta_i]\leq [\beta]$ and $|\beta_i|_\prec<|\beta|_\prec$, $i=1,\ldots,p$.

(ii) Let $F^*:\R[\z]\to \R[[\z]]$ be a linear map which satisfies (\ref{cdn1}) and (\ref{F*hom}) with arbitrary $Z,Z'\in \R[\z]$ and assume that for any $m\in \N^4$ one has $(F)_{\beta}^{\delta_{\m}}=0$ if $[\beta]<-1$. Then in the non--zero terms in the right hand side of (\ref{GbF}) one has $[\beta_i]\leq [\beta]$ and $|\beta_i|_\prec<|\beta|_\prec$, $i=1,\ldots,p$.

(iii) $F={\rm d}\G_x$ and $F=S_{xy}={\rm d}\G_x-\G_{xy}{\rm d}\G_y$, $x,y\in \R^4$ satisfy the conditions of part (ii), and hence the statement in part (ii) holds for them.

\el

\bpr
We prove (i). Part (ii) is established using the same arguments.

First observe that by (\ref{tr1F}) we may assume that $[\beta_i]\geq -1$, $i=1,\ldots, p$ in all terms in the sum in the right hand side of (\ref{Gb}). Therefore for every such term and for any $k=1,\ldots, p$ by the definition of $[~\cdot~]$ we have
$$
[\beta]=\sum_{\tiny\begin{array}{c} i=1 \\ i\neq k\end{array}}^p([\beta_i]+1)+[\beta_k]+l-p+1=\sum_{\tiny\begin{array}{c} i=1 \\ i\neq k\end{array}}^p([\beta_i]+1)+[\beta_k]+[\gamma]+1\geq [\beta_k]
$$
as $[\beta_i]+1\geq 0$, $i=1,\ldots, p$ and $[\gamma]+1\geq 0$.

Now by the definition of the modified homogeneity and by Lemma \ref{Gbgmult} we deduce
$$
|\beta|_\prec=|\beta|+\frac d2 ([\beta]+1)>|\beta_k|+\frac d2 ([\beta_k]+1)=|\beta_k|_\prec.
$$

(iii) Note that both $F^*$ in part (iii) satisfy (\ref{cdn1}) and (\ref{F*hom}) with arbitrary $Z,Z'\in \R[\z]$. For $F^*={\rm d}\G_x^*$ this is clear from the definition, and for $F^*={\rm d}\G_x^*-{\rm d}\G_y^*\G_{xy}^*$ this follows from (\ref{cdn1}) and (\ref{F*hom}) for ${\rm d}\G_x^*$ and from the fact that $\G_{xy}^*:\R[[\z]]\to\R[[\z]]$ is an algebra automorphism satisfying (\ref{G1g}).

Now part (iii) immediately follows from (ii) as by (\ref{tr1G*}) $({\rm d}\G_x)_\beta^\gamma=0$ for any $\gamma\in \M$, $\beta\in \Ml$, $[\beta]<0$, so, in particular, for any $\m\in \N^4$ one has $({\rm d}\G_x)_{\beta}^{\delta_{\m}}=0$ if $[\beta]<-1$, and also by (\ref{tr1})
$$
(\G_{xy}{\rm d}\G_y)^{\delta_{\m}}_\beta=\sum_{\mu\in\M_\geq 0\cup\M_{pp}}(\G_{xy})^{\delta_{\m}}_\mu({\rm d}\G_y)^\mu_\beta=0
$$
if $[\beta]<-1$ as $\mu\in \M_{\geq 0}\cup\M_{pp}\subset \M$ and $({\rm d}\G_y)^\mu_\beta=0$ in this case by part (ii).

\epr


\section{Stochastic estimates for the model}\label{Stohe}

In this section we state and prove the main result of this paper, stochastic estimates for the smooth models $(\Pi_x,\G_{xy}^T)$ constructed in Section \ref{regm}. The proof is quite complicated but follows the proof of a similar statement in \cite{BOT} with some modifications. As we mentioned in the Introduction, the main difference compared with \cite{BOT} is that we obtain the stochastic estimates in a bounded range of the rescaling parameter which is one of the ingredients of them. This leads to some modification of the BPHZ condition imposed on the model in the course of the proof. The arguments in all statements related to estimates involving the BPHZ condition are also adjusted accordingly. In fact, these arguments are closer to the ones which appear in the proof of Proposition 5.2 in \cite{HS}. 

After recalling the BPHZ condition in Section \ref{BPHZsubs} and the translational invariance in law property for models in Section \ref{Trinvsubs} we state the main result in Section \ref{Mainstat} and then proceed with its proof.


\subsection{A generalized BPHZ condition}\lab{BPHZsubs}

\setcounter{equation}{0}
\setcounter{theorem}{0}

To obtain stochastic estimates stable under the limit $\rho\to 0$ for the smooth model introduced in the previous section one has to fix the coefficients $c_k$, $k=1,2,3,4$ in (\ref{defP-}) in an appropriate way. This is achieved by imposing the so-called BPHZ condition. We say that a model $(\Pi_x,\G_{xy}^T)$ satisfies a generalized BPHZ condition with respect to a family of functions $\psi_\beta \in \Sw$, $\beta\in \M'$, $|\beta|<2$ if for all $\beta\in \M'$, $|\beta|<2$ one has 
\be\lab{BPHZ}
\E\Pi_{0\beta}^-(\psi_\beta)=0,
\ee
where
$$
\Pi_{0\beta}^-(\psi_\beta)=\int_{\R^4}\Pi_{0\beta}^-(x)\psi_\beta(x)dx.
$$

For any given $\beta\in \M'$ with $|\beta|<2$ we also say that $\Pi_{0\beta}^-$ satisfies the generalized BPHZ condition with respect to $\psi_\beta$ if (\ref{BPHZ}) holds.

\bp\lab{BPHZc}
For any family of functions $\psi_\beta \in \Sw$, $\beta\in \M'$, $|\beta|<2$, which are even in the variables $x_i$, $i=1,2,3$, invariant under all permutations of the coordinates $x_i$, $i=1,2,3$, and satisfy
$$
\int_{\R^4} \psi_\beta(x) dx\neq 0, ~\beta\in \M', ~|\beta|<2,
$$ 
there is a unique model of the form $(\Pi_x,\G_{xy}^T)$ which satisfies (\ref{BPHZ}).

In fact, the unique model is completely determined by $\psi_\beta$ for $\beta=k\delta_\bg+\delta_{\bf 0}$, $k=1,2,3,4$, and for the other $\beta\in \M'$, $|\beta|<2$ the corresponding generalized BPHZ conditions are satisfied automatically.
\ep

In the proof of this proposition we follow the arguments of Section 1.12. in \cite{BOT}.
We firstly describe the relevant multi--indices.

\bl\lab{ind<2}
The multi--indices $\beta\in \M'$, $|\beta|<2$ are of the form
\be\lab{I1}
\beta=k\delta_\bg,~|\beta|=k(\alpha+1)+\alpha=\frac{k-1}{2}-(k+1)\e,~k=0,1,2,3,4,5;
\ee
\be\lab{I2}
\beta=k\delta_\bg+\delta_{\bf 0},~|\beta|=k(\alpha+1)=\frac{k}{2}-k\e,~k=1,2,3,4;
\ee
\be\lab{I3}
\beta=k\delta_\bg+2\delta_{\bf 0},~|\beta|=k(\alpha+1)-\alpha=\frac{(k+1)}{2}-(k-1)\e,~k=1,2,3;
\ee
\be\lab{I4}
\beta=k\delta_\bg+\delta_{{\bf e}_i},~i=1,2,3,~|\beta|=k(\alpha+1)+1=\frac{(k+2)}{2}-k\e,~k=1,2,
\ee
where
$$
{\bf e}_1=(0,1,0,0),~{\bf e}_2=(0,0,1,0),~{\bf e}_3=(0,0,0,1).
$$
\el

\bpr
For $\e$ in the range (\ref{erest}) formula (\ref{hom}) with $\alpha=-\frac 12 -\e$ for any $\beta\in \Ml$ yields that either $\beta({\bf n})=0$ for all ${\bf n}\in \N^4$ or 
\be\lab{bnr}
|\beta|\geq \sum_{{\n}\in \N^4}\beta({\n})|{\n}|.
\ee

Thus, if $|\beta|<2$, in both cases $\beta({\bf n})=0$ for $|{\bf n}|\geq 2$. By (\ref{nhom}) the last condition implies that $\beta({\bf n})\neq 0$ only for ${\bf n}\in \{{\bf 0},{\bf e}_1,{\bf e}_2,{\bf e}_3\}$. 
Note that by (\ref{bnr}) $\beta$ can be nonzero on at most one of the ${\bf e}_i$, $i=1,2,3$ and then its value there is $1$. 

To estimate the value of $\beta\in \M'$ with $|\beta|<2$ at ${\bf 0}$ we observe that for $\beta\in \M'$ one has $[\beta]\geq -1$ which is equivalent by (\ref{pop}) to
\be\lab{pop1}
\beta(\bg)\geq \sum_{{\n}\in \N^4}\beta({\n})-1.
\ee

Now we obtain from (\ref{pop}) using (\ref{pop1}) that
$$
2> |\beta|=(\alpha+1)\beta({\bg})+\sum_{{\n}\in \N^4}\beta({\n})(|{\n}|-\alpha)+\alpha \geq 
$$
$$
\geq (\alpha+1)(\sum_{{\n}\in \N^4}\beta({\n})-1)+\sum_{{\n}\in \N^4}\beta({\n})(|{\n}|-\alpha)+\alpha=\sum_{{\n}\in \N^4}\beta({\n})(|{\n}|+1)-1\geq \sum_{{\n}\in \N^4}\beta({\n})-1
$$
which implies 
$$
\sum_{{\n}\in \N^4}\beta({\n})< 3.
$$
In particular, $\beta({\bf 0})< 3$.

Using the last condition, recalling that $\beta$ may not vanish on a unique ${\bf e}_i$, $i=1,2,3$, and it must take value $1$ on it, we deduce using (\ref{hom}) by examining all possible cases that any $\beta\in \M'$ with $|\beta|<2$ has one of the forms listed in (\ref{I1})--(\ref{I4}).

\epr

\noindent
{\em Proof of Proposition \ref{BPHZc}.}

Firstly we show that condition (\ref{BPHZ}) is satisfied for $\beta$ of types ({\ref{I1}), ({\ref{I3}) and ({\ref{I4}) for any coefficients $c_k$, $k=1,2,3,4$.

Recall that $\xi^\rho(x)=\xi(\eta_x^\rho)$, $\rho>0$, where $\eta\in \Sw$ satisfies $\int_{\R^4}\eta(y)dy=1$, $\eta(t,x_1,x_2,x_3)$ is even in the spatial coordinates $x_1,x_2,x_3$ on the $\R^3$--factor in $\R\times \R^3$ and is invariant under arbitrary permutations of these coordinates. Here, as before, by $\xi(\eta_x^\rho)$ we mean that scalar components of $\xi$ are evaluated at $\eta_x^\rho$.

Note that by the definition the covariance of the white noise $\xi=\sum_{i=1}^3\xi_i dx_i$ is invariant under translations of $\R\times \R^3$, orthogonal transformations of the $\R^3$--factor in $\R\times \R^3$, and the adjoint action of $\mathtt K$. These symmetries imply the following symmetries in law
\be\lab{sym1}
h\circ(\xi^\rho)=_{\rm law}\xi^\rho,~h\in {\mathtt K},
\ee
\be\lab{sym2}
R_i^*\xi^\rho=_{\rm law}\xi^\rho,~i=1,2,3,
\ee
\be\lab{sym3}
\pi^*\xi^\rho=_{\rm law}\xi^\rho,
\ee
where the group $\mathtt K$ acts on $\xi^\rho$ componentwise via the adjoint representation,
$$
h\circ(\xi^\rho)(x)=\sum_{i=1}^3{\rm Ad}h(\xi_i(x)) dx_i,~x\in \R^4, 
$$
$R_i:\R\times \R^3\to \R\times \R^3$ is the reflection in the plane defined by the equation $x_i=0$, $i=1,2,3$, and $R_i^*$ is the corresponding transformation of differential 1-forms, and for any permutation $\pi$ of the spatial coordinates $x_1,x_2,x_3$ $\pi^*$ stands for the corresponding transformation of differential 1-forms.

Note that the adjoint action of $\mathtt K$ on $\k$ naturally induces an action on each $V_\beta$ and $V_\beta^-$, $\beta\in \Ml$, and hence on $\O_\beta=\O\otimes V_\beta$ and on $\O_\beta^-=\O\otimes V_\beta^-$ which we denote by the same symbol $\circ$.

Applying this action to (\ref{Ph0}), (\ref{defP-}), (\ref{Pipol}), and (\ref{defP}) in the case of $\Pi_{0\beta}$ and $\Pi_{0\beta}^-$ we infer by induction over $\beta(\bg)$ that symmetry (\ref{sym1}) induces the following symmetries in law 
\be\lab{sym1p}
h\circ \Pi_{0\beta}=_{\rm law}\Pi_{0\beta},~h\circ \Pi_{0\beta}^-=_{\rm law}\Pi_{0\beta}^-,~h\in {\mathtt K}.
\ee  

Similarly, each $R_i$, $i=1,2,3$ and each permutation $\pi$ of the coordinates $x_1,x_2,x_3$ induces a linear map of each $V_\beta$ and $V_\beta^-$, $\beta\in \Ml$, and hence of $\O_\beta=\O\otimes V_\beta$ and of $\O_\beta^-=\O\otimes V_\beta^-$ which we denote by $R_i^*$ and by $\pi^*$, respectively.

Recall that by the assumption the function $\varsigma \in \D$ is invariant under the spatial reflections $R_i$, $i=1,2,3$, and hence, recalling formula (\ref{GrL}) for $K_m$ we deduce that the kernel $K(x)=K_m(x)\varsigma(x)$ is invariant under such reflections as well.

Applying $R_i^*$, $i=1,2,3$ to (\ref{Ph0}), (\ref{defP-}), (\ref{Pipol}), and (\ref{defP}) in the case of $\Pi_{0\beta}$ and $\Pi_{0\beta}^-$ we infer by induction over $\beta(\bg)$ that symmetry (\ref{sym2}) induces the following symmetries in law 
\be\lab{sym2p}
(-1)^{\sum_{{\bf n}\in \N^4} n_i\beta({\bf n})}R_i^*\Pi_{0\beta}=_{\rm law}\Pi_{0\beta},~(-1)^{\sum_{{\bf n}\in \N^4} n_i\beta({\bf n})}R_i^*\Pi_{0\beta}^-=_{\rm law}\Pi_{0\beta}^-.
\ee  

Finally, recall that by the assumption the function $\varsigma \in \D$ is invariant under each permutation $\pi$ of the coordinates $x_1,x_2,x_3$, and hence, recalling formula (\ref{GrL}) for $K_m$ we deduce that the kernel $K(x)=K_m(x)\varsigma(x)$ is invariant under $\pi$ well. Applying $\pi^*$ to (\ref{Ph0}), (\ref{defP-}), (\ref{Pipol}), and (\ref{defP}) in the case of $\Pi_{0\beta}$ and $\Pi_{0\beta}^-$ we infer by induction over $\beta(\bg)$ that symmetry (\ref{sym3}) induces the following symmetries in law 
\be\lab{sym3p}
\pi^*\Pi_{0\beta}=_{\rm law}\Pi_{0\beta},~\pi^*\Pi_{0\beta}^-=_{\rm law}\Pi_{0\beta}^-.
\ee 

For $\beta$ of type (\ref{I1}), (\ref{I2}) or (\ref{I3}) one has from the second symmetry in (\ref{sym2p}) by the definition of $R_p^*$, $p=1,2,3$
\be\lab{symIpi}
R_p^*\Pi_{0\beta}^-=_{\rm law}\Pi_{0\beta}^-.
\ee

If $\beta=k\delta_\bg$, $k=0,1,2,3,4,5$ is of type (\ref{I1})
then $V_\beta^-\simeq V\simeq \k\otimes \R^3$.  Thus $\Pi_{0\beta}^-\in \O_\beta^-\simeq \O\otimes \k\otimes \R^3$. Denote by $(\Pi_{0\beta}^-)_{p}$ the $p$-th component of the vector $\Pi_{0\beta}^-$ with respect to the basis $e^p$, $p=1,2,3$ of $\R^3$.

In terms of the components $(\Pi_{0\beta}^-)_{p}$ the definition of $R_p$, $p=1,2,3$ and the symmetry (\ref{symIpi}) yield
\be\lab{iiI1}
-(\Pi_{0\beta}^-)_{p}(R_p x)=_{\rm law}(\Pi_{0\beta}^-)_{p}(x),~p=1,2,3.
\ee

Since $\psi_\beta (t,x_1,x_2,x_3)$ are even in the spatial coordinates $x_1,x_2,x_3$, symmetry (\ref{iiI1}) implies (\ref{BPHZ}) for $\beta$ of type (\ref{I1}).

If $\beta=k\delta_\bg+2\delta_{\bf 0}$, $k=1,2,3$ is of type (\ref{I3}) then $V_\beta^-\simeq V\otimes V^*\otimes V^*\simeq \k\otimes \k^*\otimes \k^*\otimes \R^3\otimes \R^3\otimes \R^3$, where we identify $\R^3$ and ${\R^3}^*$ using the fixed scalar product on $\R^3$.  Then $\Pi_{0\beta}^-\in \O_\beta^-\simeq \O\otimes \k\otimes \k^*\otimes \k^*\otimes \R^3\otimes \R^3\otimes \R^3$. Denote by $(\Pi_{0\beta}^-)_{pql}$ the $(p,q,l)$--component of the tensor $\Pi_{0\beta}^-$ with respect to the basis $e^p\otimes e^q\otimes e^l$, $p,q,l=1,2,3$ of $\R^3\otimes \R^3\otimes \R^3$.

In terms of the components $(\Pi_{0\beta}^-)_{pql}$ the definition of $R_p$, $p=1,2,3$ and symmetry (\ref{symIpi}) yield
$$
(-1)^{\delta_{ip}+\delta_{iq}+\delta_{il}}(\Pi_{0\beta}^-)_{pql}(R_i x)=_{\rm law}(\Pi_{0\beta}^-)_{pql}(x),~p,q,l,i=1,2,3,
$$
and hence
\be\lab{iiI3}
-(\Pi_{0\beta}^-)_{pql}(R_i x)=_{\rm law}(\Pi_{0\beta}^-)_{pql}(x),
\ee
if all $p,q,l$ are different and $i=1,2,3$,
\be\lab{iiI3a}
-(\Pi_{0\beta}^-)_{pql}(R_i x)=_{\rm law}(\Pi_{0\beta}^-)_{pql}(x),
\ee
if two of the indices $p,q,l$ are equal and they are different from the third index which is equal to $i$,
\be\lab{iiI3b}
-(\Pi_{0\beta}^-)_{ppp}(R_p x)=_{\rm law}(\Pi_{0\beta}^-)_{ppp}(x),~p=1,2,3.
\ee

Since $\psi_\beta (t,x_1,x_2,x_3)$ are even in the spatial coordinates $x_1,x_2,x_3$, symmetries (\ref{iiI3}) (\ref{iiI3a}) and (\ref{iiI3b}) imply (\ref{BPHZ}) for $\beta$ of type (\ref{I3}).

For $\beta=k\delta_\bg+\delta_{{\bf e}_i}$, $i=1,2,3$, $k=1,2$ of type (\ref{I4}) one has from the second symmetry in (\ref{sym2p}) by the definition of $R_p^*$, $p=1,2,3$
\be\lab{sym2pi}
-R_i^*\Pi_{0\beta}^-=_{\rm law}\Pi_{0\beta}^-,~R_p^*\Pi_{0\beta}^-=_{\rm law}\Pi_{0\beta}^-,~p\neq i.
\ee

If $\beta=k\delta_\bg+\delta_{{\bf e}_i}$, $i=1,2,3$, $k=1,2$ then $V_\beta^-\simeq V\otimes V^*\simeq \k\otimes \k^*\otimes \R^3\otimes \R^3$, where we identify $\R^3$ and ${\R^3}^*$ using the fixed scalar product on $\R^3$.  Then $\Pi_{0\beta}^-\in \O_\beta^-\simeq \O\otimes \k\otimes \k^*\otimes \R^3\otimes \R^3$. Denote by $(\Pi_{0\beta}^-)_{pq}$ the $(p,q)$--component of the tensor $\Pi_{0\beta}^-$ with respect to the basis $e^p\otimes e^q$, $p,q=1,2,3$ of $\R^3\otimes \R^3$.

In terms of the components $(\Pi_{0\beta}^-)_{pq}$ the definition of $R_p$, $p=1,2,3$ and the first symmetry in (\ref{sym2pi}) yield
\be\lab{ii}
-(\Pi_{0\beta}^-)_{ii}(R_i x)=_{\rm law}(\Pi_{0\beta}^-)_{ii}(x),-(\Pi_{0\beta}^-)_{pq}(R_i x)=_{\rm law}(\Pi_{0\beta}^-)_{pq}(x),~p,q\neq i,
\ee
and from the second symmetry in (\ref{sym2pi}) we obtain
$$
-(\Pi_{0\beta}^-)_{pq}(R_p x)=_{\rm law}(\Pi_{0\beta}^-)_{pq}(x),~p\neq i,~p\neq q,
$$
$$
-(\Pi_{0\beta}^-)_{qp}(R_p x)=_{\rm law}(\Pi_{0\beta}^-)_{qp}(x),~p\neq i,~p\neq q,
$$
so, in particular,
\be\lab{pi}
-(\Pi_{0\beta}^-)_{pi}(R_p x)=_{\rm law}(\Pi_{0\beta}^-)_{pi}(x),~-(\Pi_{0\beta}^-)_{ip}(R_p x)=_{\rm law}(\Pi_{0\beta}^-)_{ip}(x),~p\neq i.
\ee
Since $\psi_\beta (t,x_1,x_2,x_3)$ are even in the spatial coordinates $x_1,x_2,x_3$, symmetries (\ref{ii}) and (\ref{pi}) imply (\ref{BPHZ}) for $\beta$ of type (\ref{I4}).

If $\beta=k\delta_\bg+\delta_{\bf 0}$, $k=1,2,3,4$ is of type (\ref{I2}) then $V_\beta^-\simeq V\otimes V^*\simeq \k\otimes \k^*\otimes \R^3\otimes \R^3$, where we identify $\R^3$ and ${\R^3}^*$ using the fixed scalar product on $\R^3$.  Then $\Pi_{0\beta}^-\in \O_\beta^-\simeq \O\otimes \k\otimes \k^*\otimes \R^3\otimes \R^3$. Denote by $(\Pi_{0\beta}^-)_{pq}$ the $(p,q)$--component of the tensor $\Pi_{0\beta}^-$ with respect to the basis $e^p\otimes e^q$, $p,q=1,2,3$ of $\R^3\otimes \R^3$.

In terms of the components $(\Pi_{0\beta}^-)_{pq}$ the definition of $R_p$, $p=1,2,3$ and symmetry (\ref{symIpi}) yield
\be\lab{iiI2}
-(\Pi_{0\beta}^-)_{ip}(R_i x)=_{\rm law}(\Pi_{0\beta}^-)_{ip}(x),
\ee
$$
-(\Pi_{0\beta}^-)_{pi}(R_i x)=_{\rm law}(\Pi_{0\beta}^-)_{pi}(x),~p,i=1,2,3,~p\neq i.
$$

Since $\psi_\beta (t,x_1,x_2,x_3)$ are even in the spatial coordinates $x_1,x_2,x_3$, symmetries (\ref{iiI2}) imply for $\beta$ of type (\ref{I2})
\be\lab{I2a}
\E\int_{\R^4}(\Pi_{0 k\delta_\bg+\delta_{\bf 0}}^-)_{ij}(x)\psi_{k\delta_\bg+\delta_{\bf 0}}(x)dx=0,~i,j=1,2,3,~i\neq j.
\ee

From the second symmetry in (\ref{sym1p}) we also have
$$
h\circ \Pi_{0 k\delta_\bg+\delta_{\bf 0}}^-=_{\rm law}\Pi_{0 k\delta_\bg+\delta_{\bf 0}}^-,~h\in {\mathtt K}.
$$
This implies that for $i=1,2,3$ the element
$$
\E\int_{\R^4}(\Pi_{0 k\delta_\bg+\delta_{\bf 0}}^-)_{ii}(x)\psi_{k\delta_\bg+\delta_{\bf 0}}(x)dx\in \k\otimes \k^*
$$
is invariant under the action of ${\mathtt K}$ induced by the adjoint action. The subspace of such invariants is one--dimensional. It is generated by the Casimir element which corresponds to the identity endomorphism ${\rm Id}_\k$ of $\k$ under the linear isomorphism $\k\otimes \k^*\simeq {\rm End}(\k)$. This fact together with (\ref{I2a}) implies that for $i,j=1,2,3$ and $k=1,2,3,4$
\be\lab{I2a'}
\E\int_{\R^4}(\Pi_{0 k\delta_\bg+\delta_{\bf 0}}^-)_{ij}(x)\psi_{k\delta_\bg+\delta_{\bf 0}}(x)dx=c_{ik}\delta_{ij}{\rm Id}_\k,~c_{ik}\in \R. 
\ee

Symmetry (\ref{sym3p}) and the fact that $\psi_{k\delta_\bg+\delta_{\bf 0}} (t,x_1,x_2,x_3)$, $k=1,2,3,4$ are invariant under permutations of the coordinates $x_1,x_2,x_3$ imply that in fact for each $k=1,2,3,4$ the constants $c_{ik}$ are equal for different $i$. Hence (\ref{I2a'}) takes the form
\be\lab{I2b}
\E\int_{\R^4}\Pi_{0 k\delta_\bg+\delta_{\bf 0}}^-(x)\psi_{k\delta_\bg+\delta_{\bf 0}}(x)dx=c_k'{\rm Id}_{W_{\delta_{\bf 0}}}^V,~c_k'\in \R. 
\ee

It remains to show that the inductive algorithm for defining $\Pi_{x\beta}$ described in Section \ref{recm} along with condition (\ref{BPHZ}) uniquely determine the model. Since for given coefficients $c_p$, $p=1,2,3,4$ this algorithm combined with the results of Section \ref{StrG} uniquely determines the model by the definition, it suffices to show that this algorithm is compatible with condition (\ref{BPHZ}) and that they uniquely determine the coefficients $c_p$, $p=1,2,3,4$. 

Indeed, formula (\ref{defP-}) in the case of $\Pi_{0\beta}^-$ implies that for any given $\beta$ the recentered map $\Pi_{0\beta}^-$ depends on $\Pi_{0\gamma}$ with $\gamma(\bg)<\beta(\bg)$ and $c_p$ with $p\leq \beta(\bg)$. Therefore for non--purely polynomial $\beta$ with $|\beta|<2$, $\beta\neq 0$ one can determine $\Pi_{0\beta}^-$ and $\Pi_{0\beta}$ by induction over $k:=\beta(\bg)$ with the help of (\ref{defP-}) and (\ref{defP}).

In the course of application of this inductive algorithm over $k=1,2,3,4,5$, assuming that $\Pi_{0\gamma}^-$, $\Pi_{0\gamma}$ and $c_p$ are already known for $p<k$ and $\gamma$ with $\gamma(\bg)<k$, we can firstly define $\Pi_{0 k\delta_\bg+\delta_{\bf 0}}^-$ and $\Pi_{0 k\delta_\bg+\delta_{\bf 0}}$ (if $k<5$) and then determine $\Pi_{0\beta}^-$ and $\Pi_{0\beta}$ for all the other $\beta$ with $\beta(\bg)=k$.

Note that in the case of $\Pi_{0 k\delta_\bg+\delta_{\bf 0}}^-$, $1\leq k<5$ the term containing $c_p$, $p=1,2,3,4$ in the right hand side of formula (\ref{defP-}) takes the form
$$
\sum_{l=1}^kc_l\Pi_{0 (k-l)\delta_\bg+\delta_{\bf 0}}=c_k{\rm Id}_{W_{\delta_{\bf 0}}}^V+\sum_{l=1}^{k-1}c_l\Pi_{0 (k-l)\delta_\bg+\delta_{\bf 0}} 
$$
as $\Pi_{0\delta_{\bf 0}}={\rm Id}_{W_{\delta_{\bf 0}}}^V$ by the definition. Therefore $\Pi_{0 k\delta_\bg+\delta_{\bf 0}}^--c_k{\rm Id}_{W_{\delta_{\bf 0}}}^V$ depends on $c_l$, $l<k$ and $\Pi_{0\gamma}$ with $\gamma(\bg)<k$, which are already known by the induction assumption. Hence, recalling (\ref{I2b}), condition (\ref{BPHZ}) holds for $\beta= k\delta_\bg+\delta_{\bf 0}$ provided that $c_k$ is determined from the relation
$$
c_k{\rm Id}_{W_{\delta_{\bf 0}}}^V=- \frac{\E\int_{\R^4}(\Pi_{0 k\delta_\bg+\delta_{\bf 0}}^-(x)-c_k{\rm Id}_{W_{\delta_{\bf 0}}}^V)\psi_{k\delta_\bg+\delta_{\bf 0}}(x)dx}{\int_{\R^4}\psi_{k\delta_\bg+\delta_{\bf 0}}(x)dx}=
$$
$$
=-\frac{c_k'-c_k\int_{\R^4}\psi_{k\delta_\bg+\delta_{\bf 0}}(x)dx}{\int_{\R^4}\psi_{k\delta_\bg+\delta_{\bf 0}}(x)dx}{\rm Id}_{W_{\delta_{\bf 0}}}^V,
$$ 
where the expression $c_k'-c_k\int_{\R^4}\psi_{k\delta_\bg+\delta_{\bf 0}}(x)dx$ is known by the induction assumption.

Thus the inductive algorithm for defining $\Pi_{x\beta}$ described in Section \ref{recm}  is compatible with condition (\ref{BPHZ}) and they uniquely determine the coefficients $c_p$, $p=1,2,3,4$ and the model as well. This completes the proof.

\qed


\subsection{Translational invariance in law}\lab{Trinvsubs}

\setcounter{equation}{0}
\setcounter{theorem}{0}

Another important property of the recentered maps $\Pi_{x\beta}(y)$, $\Pi_{x\beta}^-(y)$ and of the structure group maps $\G_{xy}$, which is crucial in the proof of the stochastic estimates for them, is translational invariance in law stated in the next proposition.
\bp
The recentered maps $\Pi_{x\beta}(y)$, $\Pi_{x\beta}^-(y)$ and the structure group maps $\G_{xy}$ are invariant in law with respect to translations in the sense that for all $x,y,h\in \R^4$, $\beta\in  \Ml$ one has
\be\lab{transl1}
\Pi_{x\beta}(y)=_{\rm law}\Pi_{x+h\beta}(y+h),~\Pi_{x\beta}^-(y)=_{\rm law}\Pi_{x+h\beta}^-(y+h),
\ee
\be\lab{transl2}
\G_{x+hy+h}=_{\rm law}\G_{xy}.
\ee
\ep

\bpr
By the definition the covariance of the white noise, $\xi=\sum_{i=1}^3\xi_i dx_i$ is invariant under translations in $\R\times \R^3$. This symmetry implies the following symmetry in law
\be\lab{sym1h}
\xi^\rho(x+h)=_{\rm law}\xi^\rho(x),~x,h\in \R^4.
\ee

Applying translation by $h$ to both arguments $x$ and $y$ in (\ref{Ph0}), (\ref{defP-}), (\ref{Pipol}), and (\ref{defP}) in the case of $\Pi_{x\beta}(y)$ and $\Pi_{x\beta}^-(y)$ and using (\ref{sym1h}) we infer by induction over $\beta(\bg)$, as in Section 1.12 of \cite{BOT}, that symmetry (\ref{sym1h}) induces symmetries (\ref{transl1}).

To justify (\ref{transl2}) we restrict $\G_{xy}^*$ to $\R_n:=\R[[\z_\bg,\{\z_\m^i:i=1,\ldots,{\rm dim}(V),\m\in \N^4,|\m|<n\}]]$, $n\in \N$ and use induction over $n$.

If $n=0$, $\R_0=\R[[\z_\bg]]$. By (\ref{G1g}) the restriction of $\G_{xy}^*$ to $\R_0=\R[[\z_\bg]]$ satisfies  
$$
\G_{x+hy+h}^*=\G_{xy}^*,~x,y,h\in \R^4 .
$$
This establishes the base of the induction.

Now assume that the restriction of $\G_{xy}^*$ to $\R_n$ satisfies (\ref{transl2}) for some $n\in \N$. Let $\n\in \N^4$ be such that $|\n|=n$.
Observe that by (\ref{transl1}) and by (\ref{Ydef}) $Y_\n(x+h)(\z)=_{\rm law}Y_\n(x)(\z)$, and hence by (\ref{Ydef}), by (\ref{PG}) and by the induction assumption
$$
\G_{xy}^*({\rm Id}_{W_{\delta_\n}}^V\otimes {\rm Id}_{\R[[\z]]})\z^{\delta_\n}=
\frac{1}{\n!}\G_{xy}^*\partial^\n\Pi_x(x)(\z)-\G_{xy}^*Y_\n(x)(\z)=
$$
$$
=\frac{1}{\n!}\partial^\n\Pi_y(x)(\z)-\G_{xy}^*Y_\n(x)(\z)=_{\rm law}\frac{1}{\n!}\partial^\n\Pi_{y+h}(x+h)(\z)-\G_{x+hy+h}^*Y_\n(x+h)(\z)=
$$
$$
=\frac{1}{\n!}\G_{x+hy+h}^*\partial^\n\Pi_{x+h}(x+h)(\z)-\G_{x+hy+h}^*Y_\n(x+h)(\z)=
\G_{x+hy+h}^*({\rm Id}_{W_{\delta_\n}}^V\otimes {\rm Id}_{\R[[\z]]})\z^{\delta_\n},
$$
where at the last step we used (\ref{Ydef}) again.

Repeating the same arguments for all $\n\in \N^4$ with $|\n|=n$ we deduce that 
the restriction of $\G_{xy}^*$ to $\R_{n+1}$ satisfies (\ref{transl2}). This establishes the induction step and completes the proof.

\epr


\subsection{The main statement}\lab{Mainstat}

\setcounter{equation}{0}
\setcounter{theorem}{0}

Now we are in a position to state the main result of the paper.
For this purpose observe that by Lemma \ref{locf} the set $\{|\beta|_\prec:\beta \in \M_{\geq 0}\}\subset I+\frac 52 \N$ is locally finite and bounded from below, so that the set $\M_N:=\{|\beta|_\prec:\beta \in \M_{\geq 0},~|\beta|_\prec<N\}$ is finite for any $N\in \R$. 

Let $n_0<n_1<\ldots<n_{M(N)}$ be the elements of the set $\M_N$, and for $\bar{\lambda}_0>0$, $k=0,\ldots, M(N)$ denote $\M_k:=\{\beta:\beta \in \M_{\geq 0},~|\beta|_\prec=n_k\}$, $\bar{\lambda}_k:=\frac{\bar{\lambda}_0}{32^k}$. 

\bt\lab{mainTes}
For any given mollifying parameter $\rho>0$, any $N\geq 17$ and $\bar{\lambda}_0>0$  satisfying $\bar{\lambda}_{M(N)}=\frac{\bar{\lambda}_0}{32^{M(N)}}\geq 1$ there exists a smooth model $(\Pi_x,\G_{xy}^T)$ which satisfies a generalized BPHZ condition and such that for all $1\leq p<\infty$, $\gamma,\mu \in \M$ the following inequality holds 
\be\lab{GestM}
\left\|(\G_{xy})^\gamma_\mu\right\|_{\L_p}\lesssim |x-y|^{|\mu|-|\gamma|},
\ee
for all $x,y\in \R^4$ if $\mu\not\in \M_{\geq 0}$, and if $\mu\in \M_k$, $k=0,\ldots, M(N)$ it holds for $x,y\in \R^4$, $|x-y|\leq \bar{\lambda}_k$.

Also, for all $1\leq p<\infty$, $\mu \in \M$, $\varphi\in \mB^r$, $r>{\rm max}\{2-\alpha, |\mu|\}$ the following inequalities hold
\be\lab{EphiM} 
\left\|\Pi_{0\mu}^-(\varphi^{\lambda}_x)\right\|_{\L_p}\lesssim \lambda^{\alpha-2}(\lambda+|x|)^{|\mu|-\alpha},
\ee
\be\lab{EPPsin'>M} 
\left\|\Pi_{0\mu}(\varphi_x^\lambda)\right\|_{\L_p}\lesssim \lambda^{\alpha}(\lambda+|x|)^{|\mu|-\alpha}
\ee
for all $x\in \R^4$, $0<\lambda<\infty$ if $\mu \not\in \M_{\geq 0}$, and if $\mu \in \M_k$, $k=0,\ldots, M(N)$ they hold for all $x\in \R^4$, $|x|\leq \bar{\lambda}_k$, $0<\lambda\leq \bar{\lambda}_k$.

The constants in inequalities (\ref{GestM}), (\ref{EphiM}) and (\ref{EPPsin'>M}) do not depend on $\lambda$, $\varphi$, $x$, $y$ and $\rho$.

If $N\geq 17$ the set $\{\beta\in \M: |\beta|_\prec<N \}$ contains all $\beta\in \M$ such that $|\beta|\leq 2$.
\et

This theorem is proved by induction over the modified homogeneity. The proof will be presented in Section \ref{Indstep}. In the next section we obtain the desired estimates in the case of the base of this induction when $\mu=0$ which will be used in the proof of Theorem \ref{mainTes}.


\subsection{The base case}\lab{Baseind}

\setcounter{equation}{0}
\setcounter{theorem}{0}


Firstly, for completeness and for referencing purposes we recall some algebraic properties of $\G_{xy}$, ${\rm d}\G_{x}$ and $S_{x+y x}$ which are nearly trivial in the base case.

\bp\lab{Gindpr0}
For all $x,y\in \R^4$ one has $(\G_{xy})_0^\gamma=\delta_0^\gamma{\rm Id}_{W_0}$.
\ep

\bpr
As by (\ref{triangG}) and (\ref{triangG1}) $\G_{xy}:W_0\to W_0$ is the identity transformation, and $W_0=\R$, we deduce $(\G_{xy})_0^\gamma=\delta_0^\gamma{\rm Id}_{W_0}$. 

\epr

\bp\lab{dGindprp0}
For $\gamma\neq \delta_\n$, $\n\in \N^4$, $|\n|< 2-\e$ one has $({\rm d}\G_{x})_0^\gamma=0$.
\ep

\bpr
By (\ref{prp1}) $({\rm d}\G_{x})_0^\gamma=0$ if $\gamma\neq \delta_\n$ for some $\n\in \N^4$, $|\n|< 2-\e$. 

\epr

\bp\lab{GestP0}
For $\gamma\neq \delta_\n$, $\n\in \N^4$, $|\n|< 2-\e$ one has $(S_{x+y x})^\gamma_0=0$.
\ep

\bpr
Note that by (\ref{prp1}) $({\rm d}\G_{x})_0^\gamma=0$ if $\gamma\neq \delta_\n$ for some $\n\in \N^4$, $|\n|< 2-\e$. Using this fact and (\ref{Gn}) we deduce that $(S_{x+y x})^\gamma_0=0$ if $\gamma\neq \delta_\n$ for some $\n\in \N^4$, $|\n|< 2-\e$. 

\epr

Next, we proceed with a series of estimates for Malliavin derivatives. We start with the following elementary basic result.
\bp\lab{dPi-00}
For any $1\leq q< 2$, $\bar{\lambda}_0>0$ and $\upsilon\in \L_q(L_2)$ satisfying $\left\|\left\|\upsilon(x)\right\|_{L_2(\R^4,x)}\right\|_{\L_q}\leq 1$ one has
\be\lab{Pi00est}
\left\|\delta\Pi_{00}^-\right\|_{q\bar{\lambda}_0 0}^{\alpha-2}=\sup_{0<\lambda\leq \bar{\lambda}_0}\sup_{\varphi\in \mB^0}\frac{\left\|\delta\Pi_{00}^-(\varphi^\lambda_0)\right\|_{\L_q}}{\lambda^{\alpha-2}}\lesssim 1,
\ee
where $\delta=\delta_\upsilon$ and the constant in the inequality only depends on $\bar{\lambda}_0$.
\ep

\bpr
By the definition of $\Pi_{00}^-(x)=\xi^\rho(x)=\xi(\eta^\rho_x)$, so by the H\"{o}lder and by the Young convolution inequalities we have for any $\varphi\in \mB^0$
$$
\left\|\delta\Pi_{00}^-(\varphi^\lambda_0)\right\|_{\L_q}= \left\|\upsilon^\rho(\varphi^\lambda_0)\right\|_{\L_q}\leq \left\|\left\|\upsilon(x)\right\|_{L_2(\R^4,x)}\left\|\eta^\rho_0(x)\right\|_{L_1(\R^4,x)}\left\|\varphi^\lambda_0(x)\right\|_{L_2(\R^4,x)}\right\|_{\L_q}\leq
$$
\be\lab{dPi00est}
\leq \lambda^{-\frac d2}\left\|\left\|\upsilon(x)\right\|_{L_2(\R^4,x)}\right\|_{\L_q}\left\|\varphi(x)\right\|_{L_2(\R^4,x)}\lesssim \lambda^{-\frac d2}\lesssim \lambda^{\alpha-2}, 
\ee
for $0<\lambda \leq \bar{\lambda}_0$ as
$$
\left\|\left\|\upsilon(x)\right\|_{L_2(\R^4,x)}\right\|_{\L_q}\leq 1
$$
by the assumption,
$$
\left\|\eta^\rho_0(x)\right\|_{L_1(\R^4,x)}=1,
$$
$\left\|\varphi(x)\right\|_{L_2(\R^4,x)}\lesssim 1$, and $-\frac d2>\alpha-2=-\frac d2-\e$.
This completes the proof.

\epr

Now using Schauder estimates for $\K$ we can derive an estimate for $\delta\Pi_{00}$.
\bp\lab{dPi0best0}
For any $1\leq q< 2$, $\bar{\lambda}_0\geq 1$, $r\in \N$, $r\geq 1$ and $\upsilon\in \L_{q}(L_2)$ satisfying $\left\|\left\|\upsilon(x)\right\|_{L_2(\R^4,x)}\right\|_{\L_q}\leq 1$ one has
$$
\left\|\delta\Pi_{00}\right\|_{q\bar{\lambda}_0 r}^{\alpha}=\sup_{0<\lambda\leq \bar{\lambda}_0}\sup_{\varphi\in \mB^{r}}\frac{\left\|\delta\Pi_{00}(\varphi^\lambda_0)\right\|_{\L_q}}{\lambda^{\alpha}}\lesssim 1,
$$
where $\delta=\delta_\upsilon$ and the constant in the inequality only depends on $\bar{\lambda}_0$ and $r$.
\ep

\bpr
Since by (\ref{Ph0}) $\Pi_{00}$ depends on $\xi$ polynomially, one obviously has from (\ref{Ph0})
$$
\delta\Pi_{00}=\K\delta\Pi_{00}^-.
$$

Observe also that $r\geq 1>|\alpha|=\e+\frac 12$.
Hence, the estimate in the statement follows from (\ref{Pi00est}), and from (\ref{emb1}) with $x=0$, and $F=\delta\Pi_{00}^-$.

\epr

Next, we proceed with a Besov type estimate for a germ.
\bp\lab{Fgerm0est0}
For any $\bar{\lambda}_0\geq 1$, $R>0$, $r\in \N$, $r\geq 2$, $1\leq q<2$, $\upsilon\in \L_{q}(L_2)$ with $\left\|\left\|\upsilon(x)\right\|_{L_2(\R^4,x)}\right\|_{\L_q}\leq 1$ the germ
$$
F_x(y):=\delta\Pi_{00}(y)-\sum_{\gamma\in \Ml}\Pi_{x\gamma}(y)\circ ({\rm d}\G_x)_0^\gamma
$$ 
satisfies
\be\lab{PGest=0}
\left\|F\right\|_{q2\bar{\lambda}_0 R r}^{\alpha+\frac d2}=\sup_{0<\lambda\leq \bar{\lambda}_0}\frac{\left\|\sup_{\varphi\in \mB^r}\left\|F_x(\varphi^\lambda_x)\right\|_{\L_q}\right\|_{L_2(B_R,x)}}{\lambda^{\alpha+\frac d2}}\lesssim 1,
\ee
where $\delta=\delta_\upsilon$, and the constant in the inequality only depends on $\bar{\lambda}_0$ and $r$.

\ep

\bpr
By (\ref{cdn30})
$$
F_x(y)=\delta\Pi_{00}(y)-\sum_{\tiny\begin{array}{c} \n\in \N^4: \\ |\n|<2-\e\end{array}}\Pi_{x\delta_\n}(y)\circ ({\rm d}\G_x)_0^{\delta_\n}.
$$ 

From this observation and (\ref{cdn3}) we deduce that for each $x\in \R^4$ 
$F_x(y)$ is obtained from $\delta\Pi_{00}(y)$ by subtracting its Taylor polynomial at $x$ of parabolic degree $1$. Thus $F_x=(\delta\Pi_{00})_x^{2-\e}=(\K \upsilon^\rho)_x^{2-\e}$, and we can apply (\ref{emb2c}) with $b=-\e$, $r>2-\e$ to get
\be\lab{emb2c'0}
\left\|F \right\|_{q2\bar{\lambda}_0 B_R r}^{2-\e}\lesssim  \left\|\upsilon^\rho\right\|_{q2~2\bar{\lambda}_0 B_R r}^{-\e}\leq \left\|\upsilon^\rho\right\|_{q2~2\bar{\lambda}_0 \R^4 r}^{-\e}\lesssim \left\|\upsilon^\rho\right\|_{\B_{q2}^{-\e}},
\ee
where at the last step we used definition (\ref{defbspq<}).

Now we estimate the right hand side using definition (\ref{defbspq<'}) of the space  $\B_{q2}^{-\e}$ according to which
\be\lab{emb2c'01}
\left\|\upsilon^\rho\right\|_{\B_{q2}^{-\e}}\lesssim \sup_{0<\lambda\leq 2\bar{\lambda}_0}\frac{\left\|\left\|\upsilon^\rho(\varphi_x^\lambda)\right\|_{\L_q}\right\|_{L_2(\R^4,x)}}{\lambda^{-\e}},
\ee
where $\varphi\in \D$ is such that $\int_{\R^4}\varphi(x)dx\neq 0$.
Since $\e>0$, this definition implies 
\be\lab{emb2c'02}
\left\|\upsilon^\rho\right\|_{q2~2\bar{\lambda}_0 B_R r}^{-\e}\lesssim \left\|\upsilon^\rho\right\|_{\B_{q2}^{-\e}}\lesssim \sup_{0<\lambda\leq 2\bar{\lambda}_0}\left\|\left\|\upsilon^\rho(\varphi_x^\lambda)\right\|_{\L_q}\right\|_{L_2(\R^4,x)}.
\ee

As $1\leq q<2$ we can proceed further by the Minkowski and by the Young convolution inequalities,
$$
\left\|\upsilon^\rho\right\|_{q2~2\bar{\lambda}_0 B_R r}^{-\e}\lesssim \sup_{0<\lambda\leq 2\bar{\lambda}_0}\left\|\left\|\upsilon^\rho(\varphi_x^\lambda)\right\|_{L_2(\R^4,x)}\right\|_{\L_q}\leq 
$$
\be\lab{emb2c'00}
\leq \sup_{0<\lambda\leq 2\bar{\lambda}_0}\left\|\left\|\upsilon\right\|_{L_2(\R^4)}\left\|\eta^\rho\right\|_{L_1(\R^4)}\left\|\varphi^\lambda\right\|_{L_1(\R^4)}\right\|_{\L_q}\lesssim\left\|\left\|\upsilon\right\|_{L_2(\R^4)}\right\|_{\L_q}\leq 1,
\ee
where we also used the fact that by the definition $\left\|\eta^\rho\right\|_{L_1(\R^4)}=\left\|\eta\right\|_{L_1(\R^4)}=1$, $\left\|\varphi^\lambda\right\|_{L_1(\R^4)}=\left\|\varphi\right\|_{L_1(\R^4)}$ does not depend on $\lambda$, and that by the assumption $\left\|\left\|\upsilon\right\|_{L_2(\R^4)}\right\|_{\L_q}\leq 1$.

Note that the last chain of inequalities resembles the Besov space inclusion $L_2\subset B_{2\infty}^{-\e}$.

Finally (\ref{PGest=0}) follows from (\ref{emb2c'0}) and from (\ref{emb2c'00}) given that $2-\e=\alpha+\frac d2$. This completes the proof.

\epr

Now we can obtain an estimate for $\Pi_{x0}^-$ using the previously obtained estimate for its Malliavin derivative and the spectral gap inequality.
\bp\lab{Ephi0}
For any $\bar{\lambda}_0>0$, and any $1\leq p<\infty$
\be\lab{Epsin0p} 
\left\|\Pi_{x0}^-\right\|_{p\bar{\lambda}_0 0}^{\alpha-2}=\sup_{0<\lambda\leq \bar{\lambda}_0}\sup_{\varphi\in \mB^0}\frac{\left\|\Pi_{x0}^-(\varphi^\lambda_x)\right\|_{\L_p}}{\lambda^{\alpha-2}}=
\ee
$$
=\sup_{0<\lambda\leq \bar{\lambda}_0}\sup_{\varphi\in \mB^0}\frac{\left\|\Pi_{00}^-(\varphi^\lambda_x)\right\|_{\L_p}}{\lambda^{\alpha-2}}\lesssim 1
$$
where the constant in the inequality only depends on $\bar{\lambda}_0$ and $p$.
\ep

\bpr
By translational invariance in law it suffices to consider again the case $x=0$, and by the H\"{o}lder inequality in probability one can restrict to the case when $p> 2$. Also, the equality in the right hand side of (\ref{Epsin0p}) is obvious as $\Pi_{x0}^-=\xi^\rho$ does not depend on $x$.

We shall use spectral gap inequality (\ref{SpG}) for $f=\Pi_{00}^-(\varphi^\lambda_0)$, $\varphi\in \mB^0$. 
Since the underlying Gaussian measure is centered, $\E(\Pi_{00}^-(\varphi^\lambda_0))=\E(\xi^\rho(\varphi^\lambda_0))=0$, and the spectral gap inequality takes the form
$$
\left\|\Pi_{00}^-(\varphi^\lambda_0)\right\|_{\L_p}\lesssim \sup_{\tiny \begin{array}{c}\upsilon\in \L_{p^*}(L_2): \\ \left\|\left\|\upsilon\right\|_{L_2(\R^4)}\right\|_{\L_{p^*}}\leq 1\end{array}}\left\|\delta_\upsilon \Pi_{00}^-(\varphi^\lambda_0)\right\|_{\L_q}\leq
$$
$$
\leq \sup_{\tiny \begin{array}{c}\upsilon\in \L_{p^*}(L_2): \\ \left\|\left\|\upsilon\right\|_{L_2(\R^4)}\right\|_{\L_{p^*}}\leq 1\end{array}}\left\|\delta_\upsilon \Pi_{00}^-(\varphi^\lambda_0)\right\|_{\L_{p^*}},
$$
where $1\leq q< p^*< 2< p$, $\frac 1p+\frac{1}{p^*}=1$, and at the last step we used the H\"{o}lder inequality in probability.

Thus by  Proposition \ref{dPi-00} with $q=p^*$ for $0<\lambda\leq \bar{\lambda}_0$ one has
$$
\left\|\Pi_{00}^-(\varphi^\lambda_0)\right\|_{\L_p}\lesssim \lambda^{\alpha-2},
$$
where the constant in the inequality only depends on $\bar{\lambda}_0$ and $p$. This implies (\ref{Epsin0p}).

\epr

Finally, using Schauder estimates for $\K$ we can derive an estimate for $\Pi_{x0}$.
\bp\lab{Epsin0p+0}
For any $\bar{\lambda}_0\geq 1$, $r\geq 1$ and any $1\leq p<\infty$
\be\lab{Epsin0p+} 
\left\|\Pi_{x0}\right\|_{p\bar{\lambda}_0 r}^{\alpha}=\sup_{0<\lambda\leq \bar{\lambda}_0}\sup_{\varphi\in \mB^r}\frac{\left\|\Pi_{x0}(\varphi^\lambda_x)\right\|_{\L_p}}{\lambda^{\alpha}}=
\ee
$$
=\sup_{0<\lambda\leq \bar{\lambda}_0}\sup_{\varphi\in \mB^r}\frac{\left\|\Pi_{00}(\varphi^\lambda_x)\right\|_{\L_p}}{\lambda^{\alpha}}\lesssim 1,
$$
where the constant in the inequality only depends on $\bar{\lambda}_0$, $p$ and $r$.

\ep

\bpr
By translational invariance in law, the estimate in the statement follows from (\ref{Epsin0p}) and from (\ref{emb1}) with $b=|0|-2=\alpha-2$ and $F_x=\Pi_{x0}^-$ as in this case by (\ref{defP})
$\Pi_{x0}=\K\Pi_{x0}^-$.
 
\epr

The last two propositions in this section stated below will be used at the next step of the inductive proof of Theorem \ref{mainTes}.
\bp\lab{dGpp0}
For any $1\leq q<2$, $\bar{\lambda}_0>0$, $0<R\leq\bar{\lambda}_0$ and $\upsilon\in \L_{q}(L_2)$ with $\left\|\left\|\upsilon(x)\right\|_{L_2(\R^4,x)}\right\|_{\L_q}\leq 1$  one has
\be\lab{G0estn0}
\left\|({\rm d}\G_{x})^{\delta_\n}_0\right\|_{q2B_Rx}\lesssim R^{\frac d2+\alpha-|\delta_\n|},
\ee
where ${\rm d}\G_{x}$ is defined with the help of $\upsilon$, and the constant in the inequality only depends on $\bar{\lambda}_0$. The left hand side of the inequality vanishes if $|\n|\geq 2-\e$.

\ep

The vanishing property in Proposition \ref{dGpp0} follows from Proposition \ref{dGindprp0}.
The proof of estimate (\ref{G0estn0}) will be given after the proof of Proposition \ref{dGpp}.

\bp\lab{Sxypp0}
For any $1\leq q<2$, $\bar{\lambda}_0>0$, $R\in \R$, $y\in \R^4$, $0<R\leq\bar{\lambda}_0$, $0<|y|\leq\bar{\lambda}_0$, $\n\in \N^4$ and $\upsilon\in \L_{q}(L_2)$ with $\left\|\left\|\upsilon(x)\right\|_{L_2(\R^4,x)}\right\|_{\L_q}\leq 1$  one has
\be\lab{S0estn0}
\left\|(S_{x+y x})^{\delta_\n}_0\right\|_{q2B_Rx}\lesssim |y|^{\frac d2+\alpha-|\delta_\n|},
\ee
where $S_{x+y x}$ is defined with the help of $\upsilon$, and the constant in the inequality only depends on $\bar{\lambda}_0$. The left hand side of the inequality vanishes if $|\n|\geq 2-\e$.

\ep

The vanishing property in Proposition \ref{Sxypp0} follows from Proposition \ref{GestP0}.
The proof of estimate (\ref{S0estn0}) will be given after the proof of Proposition \ref{Sxypp}.


\subsection{The induction step statement}\lab{Indstep}

In this section we formulate the induction step statement required for the proof of Theorem \ref{mainTes}.
We start by introducing a function which will be used to impose a BPHZ condition.

For any function $\varphi:\R^4\to \R$ and any $\lambda>0$ we denote by $\varphi^\lambda(y):=\varphi^\lambda_0(-y)$ the function $\varphi$ rescaled by $\lambda$. We also denote by $\ast$ the convolution of functions.

Let $\omega\in \D$ be an even function invariant under permutations of the coordinates $x_1,x_2,x_3$, and such that ${\rm supp}(\omega)\subset B_{\frac 16}$, and
\be\lab{omegapr*}
\int_{\R^4}\omega(x)dx=1,~\int_{\R^4}x^\m\omega(x)dx=0,~0<|\m|\leq r-1,~\m \in \N^4,
\ee 
where $r\in \N$, $r>2-\alpha$. Such function exists, e.g., by Lemma 8.1 in \cite{CZ}.

Define the function 
\be\lab{defpsi}
\psi:=\omega^2\ast \omega,
\ee 
where $\omega^2$ is $\omega$ rescaled by $2$. Note that $\psi$ is also even, compactly supported and invariant under permutations of the coordinates $x_1,x_2,x_3$.

\bp\lab{indstepprp}
Let $\rho>0$ be the mollification parameter in the definition of smooth models $(\Pi_x,\G_{xy}^T)$ in Section \ref{regm} with recentered maps $\Pi_{x\mu}$, $\Pi_{x\mu}^-$, and structure group maps $\G_{xy}$. Let $\bar{\lambda}_\prec\in \R$ be such that $\frac{\bar{\lambda}_\prec}8 \geq 1$, and $\beta\in \M_{\geq 0}$, $\beta\neq 0$. 

Suppose that for all $1\leq p<\infty$, $\gamma,\mu \in \M$ with $|\mu|_\prec<|\beta|_\prec$  one has 
\be\lab{Gest*}
\left\|(\G_{xy}\right)^\gamma_\mu\|_{\L_p}\lesssim |x-y|^{|\mu|-|\gamma|},~x,y\in \R^4,~|x-y|\leq \bar{\lambda}_\prec,
\ee
where the constant in the inequality is independent of $x$ and $y$ and $\rho$.

Assume furthermore that for any $1\leq p<\infty$, $\mu \in \M$ with $|\mu|_\prec<|\beta|_\prec$, and any $\varphi\in \mB^r$ for some (resp. for any) $r>{\rm max}\{2-\alpha, |\mu|\}$ one has
\be\lab{Ephi*} 
\left\|\Pi_{0\mu}^-(\varphi^{\lambda}_x)\right\|_{\L_p}\lesssim \lambda^{\alpha-2}(\lambda+|x|)^{|\mu|-\alpha},~x\in \R^4,~|x|\leq \bar{\lambda}_\prec,~0<\lambda\leq \bar{\lambda}_\prec,
\ee
and
\be\lab{EPPsin'>*} 
\left\|\Pi_{0\mu}(\varphi_x^\lambda)\right\|_{\L_p}\lesssim \lambda^{\alpha}(\lambda+|x|)^{|\mu|-\alpha},~x\in \R^4,~|x|\leq \bar{\lambda}_\prec,~0<\lambda\leq \bar{\lambda}_\prec,
\ee
where the constants in the inequalities do not depend on $\lambda$, $\varphi$, $x$ and $\rho$.

If $|\beta|<2$, assume in addition that the recentered map $\Pi_{0\beta}^-$ satisfies the generalized BPHZ condition with the function $\psi^{\bar{\lambda}_\prec}$, and that for all $\gamma,\mu\in \M$ with $|\mu|_\prec<|\beta|_\prec$, $1\leq q<p^*<2$, ${\rm d}\G_{x}$ and $S_{x+y x}$ defined with the help of $\upsilon\in \L_{p^*}(L_2)$ such that $\left\|\left\|\upsilon\right\|_{L_2}\right\|_{\L_{p^*}}\leq 1$ one has 
\be\lab{G0est*}
\left\|({\rm d}\G_{x})^\gamma_\mu\right\|_{q2B_Rx}\lesssim R^{\frac d2+|\mu|-|\gamma|},~0<R\leq \bar{\lambda}_\prec,
\ee
and  
\be\lab{G0estS*}
\left\|(S_{x+y x})^\gamma_\mu\right\|_{q2B_Rx}\lesssim
\ee
$$
\lesssim \left\{\begin{array}{ll} |y|^{\frac d2-|\gamma|_p+\alpha}(|y|+R)^{|\mu|-|\gamma|+|\gamma|_p-\alpha} & {\rm if}~\frac d2-|\gamma|_p+\alpha>0 \\ (|y|+R)^{|\mu|-|\gamma|+\frac d2} &{\rm else}\end{array}\right.,~0<R,|y|\leq \frac{\bar{\lambda}_\prec}2,
$$
where the constants in the inequalities are independent of $\rho$, $R$, $\upsilon$ and $y$.

Then if $|\beta|<2$ (resp. $|\beta|\geq 2$) estimates (\ref{Gest*}), (\ref{Ephi*}), (\ref{EPPsin'>*}), (\ref{G0est*}) and (\ref{G0estS*}) (resp. (\ref{Gest*}), (\ref{Ephi*}), (\ref{EPPsin'>*})) hold with $\mu$ replaced by $\beta$, $\bar{\lambda}_\prec$ replaced by $\bar{\lambda}=\frac{\bar{\lambda}_\prec}{32}$, and with the same other conditions as in the assumptions, where $\mu$ is replaced by $\beta$, and $\bar{\lambda}_\prec$ is replaced by $\bar{\lambda}$.


\ep

\br
By Lemma 8.1 in \cite{CZ}, given any smooth compactly supported function $\varphi$ on $\R^4$ with non--vanishing integral over $\R^4$ one can find a linear combination of a finite number of rescaled $\varphi$ which satisfies (\ref{omegapr*}). Therefore the BPHZ condition in Proposition \ref{indstepprp} can be imposed with respect to a function from a larger class. To simplify exposition we prefer to fix $\omega$ satisfying (\ref{omegapr*}) from the very beginning.
\er

\br
In the statement of Proposition \ref{indstepprp} when $|\beta|\geq 2$ one can actually put $\bar{\lambda}=\frac{\bar{\lambda}_\prec}{8}$. For uniformity of the statement we define $\bar{\lambda}=\frac{\bar{\lambda}_\prec}{32}$ for all $\beta\in \M_{\geq 0}$.
\er

\noindent
{\em Strategy of the proof of Proposition \ref{indstepprp}.}
We outline the strategy of the proof which will be presented as a series of statements in Sections \ref{indstGGS}, \ref{indstdP<2}, \ref{indstepPb} and \ref{indstGGSpp}. The conditions in each statement in these sections depend on the results obtained in the preceding statements only and on the induction assumptions imposed in this proposition. However, claims in Sections \ref{indstGGS}, \ref{indstdP<2}, \ref{indstepPb} and \ref{indstGGSpp} are often stated under weaker assumptions and sometimes results stronger than those required for the induction step are proved.

Estimates (\ref{Gest*}), (\ref{G0est*}) and (\ref{G0estS*}) for $\bar{\lambda}_\prec$ replaced by $\bar{\lambda}$, $\mu=\beta$ and $\gamma\in\M\setminus\M_{pp}=\M'$ are established in Propositions \ref{Gindpr}, \ref{dGindprp} and \ref{GestP}, respectively, (see also Corollary \ref{Gnonpp} for the exact statement), and for $\mu=\beta$ and $\gamma\in\M_{pp}$ they are obtained in Propositions \ref{Gestn}, \ref{dGpp}, and \ref{Sxypp}, respectively (see Corollaries \ref{Gestcor} and \ref{dGScor} for the exact statement). 

Estimate (\ref{Ephi*}) for $\bar{\lambda}_\prec$ replaced by $\bar{\lambda}$, $\mu=\beta$ follows from Lemma \ref{psilem0}, Corollary \ref{psicor} and from Proposition \ref{Pi-b<} for $|\beta|<2$, and it is derived in Proposition \ref{Pi-b>} for $|\beta|>2$ (see Corollary \ref{P-onpp} for the exact statement).

Estimate (\ref{EPPsin'>*}) for $\bar{\lambda}_\prec$ replaced by $\bar{\lambda}$, $\mu=\beta$ is obtained in Proposition \ref{Piint} (see Corollary \ref{Pestcor} for the exact statement).

\vskip 0.3cm

\qed

\vskip 0.5cm

\noindent
{\em Proof of Theorem \ref{mainTes}.}
First note that by Lemma \ref{Gprop} (iv) and (v) for any smooth model estimate (\ref{GestM}) holds for all $\mu\in \M_{pp}$ and for all $\mu\in \M'\setminus \M_{\geq 0}$, $x,y\in \R^4$ as by this lemma in the former case we can assume that $\gamma\in \M_{pp}$, and in the latter case we can assume that $\gamma\in \M'\setminus \M_{\geq 0}$, and estimate (\ref{GestM}) holds with any $x,y\in \R^4$ by (\ref{Gn}) and (\ref{m'}), respectively.  

Using this observation, (\ref{Pipol}), (\ref{Pi0}), (\ref{Pi-0}), (\ref{defP-}), Lemma \ref{lF} (i), (ii), (iii), (iv), (\ref{Pibnd2}) we deduce that for any smooth model estimates (\ref{EphiM}), (\ref{EPPsin'>M}) also hold for all $\mu\not\in \M_{\geq 0}$ and for any $0<\lambda<\infty$, $x\in \R^4$. For instance, if $\mu\not\in \M_{\geq 0}$ then $\Pi_{0\mu}\neq 0$ only if $\mu\in \M_{pp}$, and in this case by Lemma \ref{lF} (i), (iii), (iv), by already established (\ref{GestM}) for $\mu\in \M_{pp}$, $x,y\in \R^4$ and by (\ref{Pipol}) we have for all $0<\lambda<\infty$, $x\in \R^4$ using the H\"{o}lder inequality
\be\lab{Piunif}
\left\|\Pi_{0\mu}(\varphi^{\lambda}_x)\right\|_{\L_p}=\left\|\sum_{\tiny\begin{array}{c}\gamma\in \M_{pp} \\ |\gamma|<|\mu|\end{array}}\Pi_{x\gamma}(\varphi^{\lambda}_x)(\G_{x0})_\mu^\gamma+\Pi_{x\mu}(\varphi^{\lambda}_x)\right\|_{\L_p}\leq 
\ee
$$
\leq\sum_{\tiny\begin{array}{c}\gamma\in \M_{pp} \\ |\gamma|<|\mu|\end{array}} \left\|\Pi_{x\gamma}(\varphi^{\lambda}_x)\right\|_{\L_{2p}}\left\|(\G_{x0})_\mu^\gamma\right\|_{\L_{2p}}+\left\|\Pi_{x\mu}(\varphi^{\lambda}_x)\right\|_{\L_{2p}}\lesssim 
$$
$$
\lesssim\lambda^{\alpha}\sum_{\tiny\begin{array}{c}\gamma\in \M_{pp} \\ |\gamma|\leq |\mu|\end{array}}\lambda^{|\gamma|-\alpha}|x|^{|\mu|-|\gamma|}\lesssim \lambda^{\alpha}(\lambda+|x|)^{|\mu|-\alpha},
$$
where at the last step we used (\ref{min}).

The case of $\Pi_{0\mu}^-$ with $\mu\not\in \M_{\geq 0}$ is treated similarly noting that in this case $\Pi_{0\mu}^-\neq 0$ only if $\mu\in \M'\setminus \M_{\geq 0}$ by Lemma \ref{lF} (iii).

By (\ref{tr1G*}) and  (\ref{tr1S*}) estimates (\ref{G0est*}) and (\ref{G0estS*}) trivially hold for any smooth model, for all $\mu\not\in \M_{\geq 0}$ and for any $\bar{\lambda}_\prec$. 

Therefore it suffices to establish estimates (\ref{GestM}), (\ref{EphiM}) and (\ref{EPPsin'>M}) for $\mu\in \M_{\geq 0}$, and we can assume in the proof that estimates (\ref{Gest*}), (\ref{G0est*}), (\ref{G0estS*}), (\ref{Ephi*}) and (\ref{EPPsin'>*}) hold for all $\mu\in \M\setminus \M_{\geq 0}$ and for any $\bar{\lambda}_\prec$.

Next, using Proposition \ref{BPHZc} we fix a smooth model for which estimates (\ref{GestM}), (\ref{EphiM}) and (\ref{EPPsin'>M}) will be established.

Recall that by Lemma \ref{locf} the set $\{|\beta|_\prec:\beta \in \M_{\geq 0}\}\subset I+\frac 52 \N$ is locally finite and bounded from below, so that the set $\M_N=\{|\beta|_\prec:\beta \in \M_{\geq 0},~|\beta|_\prec<N\}$ is finite. 

By checking the modified homogeneities of the elements $\beta\in \M_{\geq 0}$, $|\beta|<2$ from the list in Lemma \ref{ind<2} we deduce that the condition $N\geq 17$ implies 
\be\lab{MNb}
\M_N\supset \{|\beta|_\prec:\beta\in \M_{\geq 0},|\beta|<2\}.
\ee 

As before, let $n_0<n_1<\ldots<n_{M(N)}$ be the elements of the set $\M_N$, and for $\bar{\lambda}_0>0$, $k=0,\ldots, M(N)$ denote $\bar{\lambda}_k=\frac{\bar{\lambda}_0}{32^k}$. Fix $\bar{\lambda}_0>0$ such that $\bar{\lambda}_{M(N)}\geq 1$.

Note that by the definition the function $\psi$ introduced in (\ref{defpsi}) is even and 
$$
\int_{\R^4}\psi(x)dx=1.
$$
Thus if for $k\delta_\bg+\delta_{\bf 0}\in \M_{\geq 0}$, $k=1,2,3,4$ we write $|k\delta_\bg+\delta_{\bf 0}|_\prec=n_{i_k}$ then the functions $\psi^{\bar{\lambda}_{n_{i_k-1}}}$, $k=1,2,3,4$ satisfy the conditions of Proposition \ref{BPHZc}.
 
Now for any given mollification parameter $\rho>0$ consider the smooth model $(\Pi_x,\G_{xy}^T)$ which is uniquely defined as in Proposition \ref{BPHZc} by the generalized BPHZ conditions
\be\lab{BPHZfix}
\E\Pi_{0~k\delta_\bg+\delta_{\bf 0}}^-(\psi^{\bar{\lambda}_{n_{i_k-1}}})=0, 
\ee
We claim that for this model estimates (\ref{GestM}), (\ref{EphiM}) and (\ref{EPPsin'>M}) hold. The proof is by induction over the modified homogeneity.

Firstly, one can verify straightforwardly that $\mu=0$ is the only element of $\M_{\geq 0}$ with minimal possible modified homogeneity, $|0|_\prec=2-\e$. For $\mu=0$ and arbitrary $\bar{\lambda}_\prec=\bar{\lambda}_0>0$ estimates (\ref{Gest*}), (\ref{G0est*}), (\ref{G0estS*}), (\ref{Ephi*}) and (\ref{EPPsin'>*}) follow from the results of Section \ref{Baseind}, see Proposition \ref{Gindpr0} for (\ref{Gest*}), Propositions \ref{dGindprp0} and \ref{dGpp0} for (\ref{G0est*}), Propositions \ref{GestP0} and \ref{Sxypp0} for (\ref{G0estS*}), Proposition \ref{Ephi0} for (\ref{Ephi*}), and Proposition \ref{Epsin0p+0} for (\ref{EPPsin'>*}). Therefore for $\beta\in \M_{\geq 0}$, such that $\beta\neq 0$ and $|\beta|_\prec$ is minimal possible, the induction assumptions are satisfied, and we can begin the induction procedure over the modified homogeneity the induction step of which is established in Proposition \ref{indstepprp}. Conditions (\ref{BPHZfix}}) ensure that the generalized BPHZ conditions used in the induction procedure are satisfied. 

Note that by the condition $\bar{\lambda}_{M(N)}\geq 1$ the induction step can be iterated until we exhaust the set of all multi-indices $\beta\in \M_{\geq 0}$ with $|\beta|_\prec<N$, and after completing the last iteration we obtain estimates (\ref{GestM}), (\ref{EphiM}) and (\ref{EPPsin'>M}) for all $\mu \in \M_{\geq 0}$ with $|\mu|_\prec<N$.

Finally we prove the last claim of the theorem. Observing that by the definition $[\beta]=-1$ for all $\beta\in \M_{pp}\cup (\M'\setminus \M_{\geq 0})$ we deduce that for such $\beta$ one has $|\beta|=|\beta|_\prec$. Now from (\ref{poph}), (\ref{MNb}), from formulas (\ref{n1}) and (\ref{n2}) for the homogeneities of the elements of the set $\M'\setminus \M_{\geq 0}$ and from the definition of the homogeneity which implies $|\delta_\n|=|\n|$, $\n\in \N^4$ we infer that $\{\beta\in \M:|\beta|_\prec<N\}\supset \{\beta\in \M:|\beta|\leq 2\}$ if $N\geq 17$. 
This completes the proof.

\qed


\subsection{The induction step for $(\G_{xy})^\gamma_\beta$, $({\rm d}\G_x)^\gamma_\beta$, $(S_{xy})^\gamma_\beta$, $\gamma\not\in \M_{pp}$}\lab{indstGGS}

\setcounter{equation}{0}
\setcounter{theorem}{0}

In this section we establish the induction step for $(\G_{xy})^\gamma_\beta$, $({\rm d}\G_x)^\gamma_\beta$, $(S_{xy})^\gamma_\beta$, $\gamma\not\in \M_{pp}$ using the so--called algebraic argument following the terminology of \cite{BOT,LOTT}.

\bp\lab{Gindpr}
Suppose that for some $\beta\in \M$, $\beta\neq 0$ and for all $1\leq p<\infty$, $\gamma\in\M_{pp}$, $\mu \in \M$ with $|\mu|_\prec<|\beta|_\prec$ one has
\be\lab{Gest}
\left\|(\G_{xy}\right)^\gamma_\mu\|_{\L_p}\lesssim |x-y|^{|\mu|-|\gamma|},~x,y\in \R^4,~|x-y|\leq \bar{\lambda}_\prec, \bar{\lambda}_\prec>0,
\ee
where the constant in the inequality is independent of $x$ and $y$ and $\rho$.
Then (\ref{Gest}) holds for $\mu=\beta$, for all $1\leq p<\infty$ and $\gamma\in\M\setminus\M_{pp}$. 
\ep

\bpr
For $\gamma=0$ the statement is void by the first identity in (\ref{G1g}).

If $\gamma\neq 0$, $\gamma\in\M\setminus\M_{pp}$ then $\gamma=l\delta_\bg+\sum_{i=1}^r\delta_{{\m}_i}$, $\m_i\in \N^4$, $l\in \N$, $l>0$.

If $r=0$, i.e. $\gamma=l\delta_\bg$, $l>0$ the statement obviously holds by the second identity in (\ref{G1g}).

If $r>0$ we apply $\left\|\cdot\right\|_{\L_p}$ to (\ref{Gb}) with $F=\G_{xy}$ and use the H\"{o}lder inequality together with (\ref{tr1}) to get
$$
\left\|(\G_{xy}\right)^\gamma_\beta\|_{\L_p}\leq \hspace{-1cm}  \sum_{\tiny\begin{array}{c}\beta_1,\ldots,\beta_r\in \M_{\geq 0}\cup \M_{pp}: \\ \beta_1+\ldots+\beta_r+l\delta_\bg=\beta\end{array}} \hspace{-1cm}  \left\|(\G_{xy})_{\beta_1}^{\delta_{\m_1}}\right\|_{\L_{p_1}}\ldots \left\|(\G_{xy})_{\beta_r}^{\delta_{\m_r}}\right\|_{\L_{p_r}},
$$
where $\frac{1}{p_1}+\ldots+\frac{1}{p_r}=\frac{1}{p}$.  

Since $l>0$, by Lemma \ref{ltriang} (i) $|\beta_i|_\prec<|\beta|_\prec$, $i=1,\ldots,r$ in the right hand side of the previous formula, so we can use the induction assumption. This yields
\be\lab{Gxyind}
\left\|(\G_{xy}\right)^\gamma_\beta\|_{\L_p}\lesssim \hspace{-1cm}  \sum_{\tiny\begin{array}{c}\beta_1,\ldots,\beta_r\in \M_{\geq 0}\cup \M_{pp}: \\ \beta_1+\ldots+\beta_r+l\delta_\bg=\beta\end{array}} \hspace{-1cm}|x-y|^{\sum_{i=1}^r(|\beta_i|-|\delta_{\m_i}|)}.
\ee

By the linearity of the function $|\cdot|-\alpha$ we obtain
\be\lab{bd-}
\sum_{i=1}^r(|\beta_i|-|\delta_{\m_i}|)=|\beta|-|\gamma|.
\ee
Substituting this into (\ref{Gxyind}) yields the result.
 
\epr

\bp\lab{dGindprp}
Suppose the conditions of Proposition \ref{Gindpr} are satisfied, that $|\beta|<2$, and that for all $\gamma\in \M_{pp}$, $\mu \in \M$ with $|\mu|_\prec<|\beta|_\prec$, $1\leq q<p^*<2$, ${\rm d}\G_{x}$ defined with the help of $\upsilon\in \L_{p^*}(L_2)$ such that $\left\|\left\|\upsilon\right\|_{L_2}\right\|_{\L_{p^*}}\leq 1$ one has
\be\lab{G0est}
\left\|({\rm d}\G_{x})^\gamma_\mu\right\|_{q2B_Rx}\lesssim R^{\frac d2+|\mu|-|\gamma|},~0<R\leq \bar{\lambda}_\prec,
\ee
where  the constant in the inequality is independent of $\rho$, $R$ and $\upsilon$.

Then (\ref{G0est}) holds for $\mu=\beta$, for all $1\leq q<p^*<2$, $\gamma\in\M\setminus\M_{pp}=\M'$, ${\rm d}\G_{x}$ defined with the help of $\upsilon\in \L_{p^*}(L_2)$ such that $\left\|\left\|\upsilon\right\|_{L_2}\right\|_{\L_{p^*}}\leq 1$. The constant in the inequality is independent of $\rho$, $R$ and $\upsilon$.
\ep

\bpr
Fix $1\leq q<p^*<2$ and assume that ${\rm d}\G_{x}$ is defined with the help of $\upsilon\in \L_{p^*}(L_2)$ such that $\left\|\left\|\upsilon\right\|_{L_2}\right\|_{\L_{p^*}}\leq 1$.

For $\gamma=l\delta_\bg$, $l\in \N$ the statement is obvious by the second identity in (\ref{cdn1}) and by (\ref{F*hom}).

If $\gamma\neq l\delta_\bg$, $l\in \N$, $\gamma\in\M\setminus\M_{pp}$ then $\gamma=l\delta_\bg+\sum_{i=1}^r\delta_{{\m}_i}$, $\m_i\in \N^4$, $l\in \N$, $l>0$, $r\geq 1$, and we apply $\left\|\cdot\right\|_{q2B_Rx}$ to (\ref{GbF}) with $F={\rm d}\G_{x}$ and use the H\"{o}lder inequality together with (\ref{tr1}) and (\ref{tr1G*}) to get for $\upsilon\in \L_{p^*}(L_2)$ satisfying $\left\|\left\|\upsilon\right\|_{L_2}\right\|_{\L_{p^*}}\leq 1$
\be\lab{dGest}
\left\|({\rm d}\G_{x})^\gamma_\beta\right\|_{q2B_Rx}\leq \hspace{7cm}
\ee
$$
\leq \sum_{k=1}^r\sum_{\tiny\begin{array}{c}\beta_1,\ldots, \beta_{k-1}, \\ \beta_{k+1},\ldots ,\beta_r\in \M_{\geq 0}\cup \M_{pp}, \\ \beta_k\in \M_{\geq 0}: \\ \beta_1+\ldots+\beta_r+l\delta_\bg=\beta\end{array}} \hspace{-0.5cm}  \left\|(\G_{x0})_{\beta_1}^{\delta_{\m_1}}\right\|_{q_1\infty B_Rx}\ldots \left\|(\G_{x0})_{\beta_{k-1}}^{\delta_{\m_{k-1}}}\right\|_{q_{k-1}\infty B_Rx}\times
$$
$$
\times\left\|({\rm d}\G_{x})_{\beta_k}^{\delta_{\m_k}}\right\|_{q_k2B_Rx}\left\|(\G_{x0})_{\beta_{k+1}}^{\delta_{\m_{k+1}}}\right\|_{q_{k+1}\infty B_Rx}\ldots \left\|(\G_{x0})_{\beta_r}^{\delta_{\m_r}}\right\|_{q_r\infty B_Rx},
$$
where $\frac{1}{q_1}+\ldots+\frac{1}{q_r}=\frac{1}{q}$, in the $k$-th term $1\leq q\leq q_k<p^*<2$. 

Since $l>0$, by Lemma \ref{ltriang} (iii) $|\beta_i|_\prec<|\beta|_\prec$, $i=1,\ldots,r$ in the right hand side of the previous formula, so we can use (\ref{Gest}) and (\ref{G0est}). Since $|\beta_i|\geq |\delta_{\m_i}|$, $i\neq k$ in the $k$-th term (see (\ref{triangG})) we obtain by (\ref{Gest}) that $\left\|(\G_{x0})_{\beta_{i}}^{\delta_{\m_{i}}}\right\|_{q_{i}\infty B_Rx}\lesssim R^{|\beta_i|- |\delta_{\m_i}|}$, as well as $\left\|({\rm d}\G_{x})_{\beta_k}^{\delta_{\m_k}}\right\|_{q_k2B_Rx}\lesssim R^{\frac d2 +|\beta_k|- |\delta_{\m_k}|}$ by (\ref{G0est}). Using these estimates in (\ref{dGest}) together with (\ref{bd-}) yields the result.
 
\epr

For any $\gamma\in \Ml$ we denote $|\gamma|_p:=\sum_{\n\in\N^4}|\n|\gamma(\n)$.

\bp\lab{GestP}
Suppose the conditions of Proposition \ref{Gindpr} are satisfied, that $|\beta|<2$,  and that for all $\gamma\in \M_{pp}$, $\mu \in \M$ with $|\mu|_\prec<|\beta|_\prec$, $0<R\leq \frac{\bar{\lambda}_\prec}2$, $y\in \R^4$, $|y|\leq \frac{\bar{\lambda}_\prec}2$, $1\leq q<p^*<2$, $S_{x+y x}$ defined with the help of $\upsilon\in \L_{p^*}(L_2)$ such that $\left\|\left\|\upsilon\right\|_{L_2}\right\|_{\L_{p^*}}\leq 1$ one has 
\be\lab{G0estS}
\left\|(S_{x+y x})^\gamma_\mu\right\|_{q2B_Rx}\lesssim
\ee
$$
\lesssim \left\{\begin{array}{ll} |y|^{\frac d2-|\gamma|_p+\alpha}(|y|+R)^{|\mu|-|\gamma|+|\gamma|_p-\alpha} & {\rm if}~\frac d2-|\gamma|_p+\alpha>0 \\ (|y|+R)^{|\mu|-|\gamma|+\frac d2} &{\rm else}\end{array}\right. ,
$$
where the constant in the inequality is independent of $\rho$, $R$, $\upsilon$ and $y$.

Then (\ref{G0estS}) holds for $\mu=\beta$, for all $\gamma\in\M\setminus\M_{pp}$, $0<R\leq \frac{\bar{\lambda}_\prec}2$, $y\in \R^4$, $|y|\leq \frac{\bar{\lambda}_\prec}2$, $1\leq q<p^*<2$, $S_{x+y x}$ defined with the help of $\upsilon\in \L_{p^*}(L_2)$ such that $\left\|\left\|\upsilon\right\|_{L_2}\right\|_{\L_{p^*}}\leq 1$, and the constant in the inequality is independent of $y$, $\rho$, $R$, $\upsilon$. 

Moreover, the left hand side of inequality (\ref{G0estS}) for $\mu=\beta$ vanishes if $|\beta|-|\gamma|+|\gamma|_p-\alpha<0$ or $\sum_{\n\in\N^4}\gamma(\n)=0$.
\ep

\bpr
Fix $1\leq q<p^*<2$ and assume that $S_{x+y x}$ is defined with the help of $\upsilon\in \L_{p^*}(L_2)$ such that $\left\|\left\|\upsilon\right\|_{L_2}\right\|_{\L_{p^*}}\leq 1$.

If $\sum_{\n\in\N^4}\gamma(\n)=0$ then $\gamma=l\delta_\bg$, and in this case $(S_{x+y x})_{\beta}^{l\delta_\bg}=0$ by the definition of $S_{x+y x}={\rm d}\G_{x+y}-\G_{x+y x}{\rm d}\G_x$, by (\ref{G1g}) and by (\ref{cdn1}) as these properties imply
$$
(S_{x+y x})_{\beta}^{l\delta_\bg}\z^\beta=(S_{x+y x}^*\z_\bg^l)_\beta=
(({\rm d}\G_{x+y}^*-{\rm d}\G_x^*\G_{x+y x}^*)\z_\bg^l)_\beta=-({\rm d}\G_x^*\z_\bg^l)_\beta=0.
$$
In particular, (\ref{G0estS}) holds for $\gamma=0$.

If $\gamma\neq 0$, $\gamma\in\M\setminus\M_{pp}$ then $\gamma=l\delta_\bg+\sum_{i=1}^r\delta_{{\m}_i}$, $\m_i\in \N^4$, $l\in \N$, $l>0$, and we apply $\left\|\cdot\right\|_{q2B_Rx}$ to (\ref{GbF}) with $F=S_{x+y x}$ and use the H\"{o}lder inequality together with (\ref{tr1}) and (\ref{tr1S*}) to get
\be\lab{SGest}
\left\|(S_{x+y x})^\gamma_\beta\right\|_{q2B_Rx}\leq \hspace{7cm}
\ee
$$
\leq \sum_{k=1}^r\sum_{\tiny\begin{array}{c}\beta_1,\ldots, \beta_{k-1}, \\ \beta_{k+1},\ldots ,\beta_r\in \M_{\geq 0}\cup \M_{pp}, \\ \beta_k\in \M_{\geq 0}: \\ \beta_1+\ldots+\beta_r+l\delta_\bg=\beta\end{array}} \hspace{-0.5cm}  \left\|(\G_{x+y0})_{\beta_1}^{\delta_{\m_1}}\right\|_{q_1\infty B_Rx}\ldots \left\|(\G_{x+y0})_{\beta_{k-1}}^{\delta_{\m_{k-1}}}\right\|_{q_{k-1}\infty B_Rx}\times
$$
$$
\times\left\|(S_{x+y x})_{\beta_k}^{\delta_{\m_k}}\right\|_{q_k2B_Rx}\left\|(\G_{x+y0})_{\beta_{k+1}}^{\delta_{\m_{k+1}}}\right\|_{q_{k+1}\infty B_Rx}\ldots \left\|(\G_{x+y0})_{\beta_r}^{\delta_{\m_r}}\right\|_{q_r\infty B_Rx}.
$$
Here $\frac{1}{q_1}+\ldots+\frac{1}{q_r}=\frac{1}{q}$, in the $k$-th term $1\leq q\leq q_k<p^*<2$, and $\upsilon\in \L_{p^*}(L_2)$ satisfies $\left\|\left\|\upsilon\right\|_{L_2}\right\|_{\L_{p^*}}\leq 1$.

Since $l>0$, by Lemma \ref{ltriang} (iii) $|\beta_i|_\prec<|\beta|_\prec$, $i=1,\ldots,r$ in the right hand side of the previous formula, so we can use (\ref{Gest}) and (\ref{G0estS}). Since $|\beta_i|\geq |\delta_{\m_i}|$, $i\neq k$ in the $k$-th term (see (\ref{triangG})) we obtain by (\ref{Gest}) that $\left\|(\G_{x+y0})_{\beta_{i}}^{\delta_{\m_{i}}}\right\|_{q_{i}\infty B_Rx}\lesssim (|y|+R)^{|\beta_i|- |\delta_{\m_i}|}$, $|y|,R\leq \frac{\bar{\lambda}_\prec}2$. Using these estimates in (\ref{SGest}) together with estimate (\ref{G0estS}) for $\left\|(S_{x+y x})_{\beta_k}^{\delta_{\m_k}}\right\|_{q_k2B_Rx}$ and (\ref{bd-}) we obtain  (\ref{G0estS}) for $\mu=\beta$.

Note that for the $k$-th term in (\ref{SGest}) we have by (\ref{bd-})
$$
\sum_{i=1}^r|\beta_i|=|\beta|-|\gamma|+\sum_{i=1}^r|\delta_{\m_i}|=|\beta|-|\gamma|+|\gamma|_p.
$$
Since by (\ref{triangG}) in the $k$-th term in sum (\ref{SGest}) one has $|\beta_i|\geq |\delta_{\m_i}|=|\m_i|\geq 0$, $i\neq k$, and $|\beta_k|\geq \alpha$ we deduce
$$
\alpha\leq\sum_{i=1}^r|\beta_i|=|\beta|-|\gamma|+|\gamma|_p,
$$
i.e. the left hand side of inequality (\ref{G0estS}) for $\mu=\beta$ vanishes if
$|\beta|-|\gamma|+|\gamma|_p-\alpha<0$.
This completes the proof.

\epr

Since if $\beta\in \M_{\geq 0}$ the assumptions of Proposition \ref{indstepprp} are stronger than those in Propositions \ref{Gindpr}, \ref{dGindprp} and \ref{GestP} we immediately obtain the following corollary from them.
\bc\lab{Gnonpp}
Under the induction step assumptions in Proposition \ref{indstepprp} estimates (\ref{Gest}), (\ref{G0est}) and (\ref{G0estS}) hold with $\mu$ replaced by $\beta$ and with the same conditions as in Propositions \ref{Gindpr}, \ref{dGindprp}, and \ref{GestP}, respectively, and estimates (\ref{Gest*}), (\ref{G0est*}) and (\ref{G0estS*}) hold with $\mu$ replaced by $\beta$, $\bar{\lambda}_\prec$ replaced by $\bar{\lambda}=\frac{\bar{\lambda}_\prec}{32}$, for all $\gamma\in \M'$, and with the same other conditions as in the assumptions of Proposition \ref{indstepprp}, where $\mu$ is replaced by $\beta$, and $\bar{\lambda}_\prec$ is replaced by $\bar{\lambda}$.
\ec


\subsection{Induction step estimates for $\delta\Pi_\beta$ and $\delta\Pi_\beta^-$, $|\beta|<2$}\lab{indstdP<2}

\setcounter{equation}{0}
\setcounter{theorem}{0}

This section contains the most difficult part in the proof of the induction step. As in the base case to obtain estimates for $\Pi_\beta$ and $\Pi_\beta^-$, $|\beta|<2$ using the spectral gap inequality we have to derive some estimates for their Malliavin derivatives which is done in this section. For this purpose we use the reconstruction theorem and Schauder estimates for germs several times.

In this section it is always assumed that $|\beta|<2$. 
We start by establishing a coherence condition for a germ.
\bl\lab{reclem}
Assume that (\ref{G0estS}) holds for $\mu=\beta\in \M_{\geq 0}$, where $\beta\neq 0$, $|\beta|<2$, for all $\gamma\in\M\setminus\M_{pp}=\M'$, $0<R\leq \frac{\bar{\lambda}_\prec}2$, $y\in \R^4$, $|y|\leq \frac{\bar{\lambda}_\prec}2$, $1\leq q<p^*<2$, $S_{x+y x}$ defined with the help of $\upsilon\in \L_{p^*}(L_2)$ such that $\left\|\left\|\upsilon\right\|_{L_2}\right\|_{\L_{p^*}}\leq 1$. Suppose also that for some $\varphi\in \D$ with $\int_{\R^4}\varphi(x)dx\neq 0$ and all $1\leq p<\infty$, $0<R\leq \bar{\lambda}_\prec$, $\gamma\in\M'$with $|\gamma|_\prec<|\beta|_\prec$ one has
$$
\left\|\Pi_{x\gamma}^-\right\|_{p\infty B_R \bar{\lambda}_\prec \varphi}^{|\gamma|-2}=\sup_{0<\lambda\leq \bar{\lambda}_\prec}\frac{\left\|\Pi_{x\gamma}^-(\varphi_x^\lambda)\right\|_{p\infty B_R x}}{\lambda^{|\gamma|-2}}\lesssim 1,
$$
where the constant in the inequality does not depend on $\rho$ and $R$.

Then for all $1\leq q<p^*<2$ the germ $F_x(y):=\sum_{\gamma\in \Ml}\Pi_{x\gamma}^-(y)\circ ({\rm d}\G_x)_\beta^\gamma$, where ${\rm d}\G_x$ is defined with the help of $\upsilon \in \L_{p^*}(L_2)$ such that $\left\|\left\|\upsilon\right\|_{L_2}\right\|_{\L_{p^*}}\leq 1$, satisfies for any $0<R\leq 2\bar{\lambda}':=\frac{\bar{\lambda}_\prec}{2}$
$$
\left\|F\right\|_{q2 R \bar{\lambda}' \varphi}^{\alpha-1~2\alpha+\frac d2-1~ |\beta|-1-2\alpha}=
$$
\be\lab{gestM}
=\sup_{h\in \R^4:|h|\leq 2\bar{\lambda}'}\sup_{0<\lambda\leq \bar{\lambda}'}\frac{\left\|(F_{x+h}-F_x)(\varphi_{x+h}^\lambda)\right\|_{q2 B_R x}}{\lambda^{\alpha-1}(\lambda+|h|)^{\alpha+\frac d2}(\lambda+R+|h|)^{|\beta|-1-2\alpha}}\lesssim 1,
\ee
and the constant in the inequality is independent of $\rho$, $R$ and $\upsilon$.
\el

\bpr
Fix $1\leq q<p^*<2$ and assume that $S_{x+y x}$ is defined with the help of $\upsilon\in \L_{p^*}(L_2)$ such that $\left\|\left\|\upsilon\right\|_{L_2}\right\|_{\L_{p^*}}\leq 1$.

First observe that by (\ref{Pi-0}), (\ref{PG}), (\ref{prp2}), and (\ref{tr2S*}) 
$$
(F_{x+h}-F_x)(\varphi_{x+h}^\lambda)=\sum_{\tiny\begin{array}{c}\gamma\in \M': \\ |\gamma|_\prec<|\beta|_\prec \\ |\gamma|< |\beta|+\frac d2 \end{array}}\Pi_{x+h\gamma}^-(\varphi_{x+h}^\lambda)\circ(S_{x+h x})_\beta^\gamma,
$$
and hence by H\"{o}lder's inequality
\be\lab{germe1}
\left\|(F_{x+h}-F_x)(\varphi_{x+h}^\lambda)\right\|_{q2 R x}\leq 
\ee
$$
\leq \sum_{\tiny\begin{array}{c}\gamma\in \M': \\ |\gamma|_\prec<|\beta|_\prec \\ |\gamma|< |\beta|+\frac d2 \end{array}}\left\|(S_{x+h x})^\gamma_\beta\right\|_{q'2B_Rx}\left\|\Pi_{x+h\gamma}^-(\varphi_{x+h}^\lambda)\right\|_{p\infty B_R x},
$$
where $\frac 1q=\frac{1}{q'}+\frac 1p$, $1\leq q< q'<p^*<2$, $1\leq p<\infty$.

By the assumption 
\be\lab{p-est}
\left\|\Pi_{x+h\gamma}^-(\varphi_{x+h}^\lambda)\right\|_{p\infty B_R x}\lesssim \lambda^{|\gamma|-2},~0<\lambda\leq \bar{\lambda}_\prec,
\ee
where the constant in the inequality is independent on $\rho$ and $R$, and by translational invariance in law for $\Pi_{x+h\gamma}^-$ it is also independent of $x+h$. 

Now suppose that $0<R\leq 2\bar{\lambda}'$. Then by the conditions of the lemma we can use (\ref{G0estS}) in (\ref{germe1}) for $h\in \R^4$, $|h|\leq 2\bar{\lambda}'$, $\mu=\beta\in \M$, $1\leq q<q'<p^*<2$ and $\gamma\in\M\setminus\M_{pp}$.

We now estimate each term in the sum in the right hand side of (\ref{germe1}).

Observe that in (\ref{G0estS}) we can assume that $\sum_{\n\in\N^4}\gamma(\n)>0$, for otherwise its left hand side is zero, and $\gamma(\bg)>0$ for $\gamma\in \M'$. Thus in this case
\be\lab{case1'}
|\gamma|-|\gamma|_p=(\alpha+1)\gamma(\bg)-\alpha\sum_{\n\in\N^4}\gamma(\n)+\alpha\geq\alpha+1.
\ee 
Therefore, since we always have $|\gamma|_p\geq 0$,
$$
|\gamma|\geq |\gamma|_p+\alpha+1\geq \alpha+1,
$$
and hence
\be\lab{case1''}
|\gamma|-2\geq \alpha-1.
\ee

If 
\be\lab{case1}
\frac d2-|\gamma|_p+\alpha>0,
\ee
using inequality (\ref{case1''}), (\ref{p-est}) and (\ref{G0estS}) in the case when $\frac d2-|\gamma|_p+\alpha>0$ and for $h\in \R^4$, $|h|\leq 2\bar{\lambda}'$, $\mu=\beta\in \M$, $1\leq q< q'<p^*<2$ we obtain for $0<R\leq 2\bar{\lambda}'$
$$
\left\|(S_{x+h x})^\gamma_\beta\right\|_{q'2B_Rx}\left\|\Pi_{x+h\gamma}^-(\varphi_{x+h}^\lambda)\right\|_{p\infty B_R x}\lesssim \lambda^{|\gamma|-2}|h|^{\frac d2-|\gamma|_p+\alpha}(|h|+R)^{|\beta|-|\gamma|+|\gamma|_p-\alpha}=
$$
\be\lab{trest}
=\lambda^{\alpha-1}\lambda^{|\gamma|-2-\alpha+1}|h|^{\frac d2-|\gamma|_p+\alpha}(|h|+R)^{|\beta|-|\gamma|+|\gamma|_p-\alpha}.
\ee

Now recall that for any $v,w,a,b\geq 0$
\be\lab{min}
v^aw^b\leq (v+w)^{a+b}.
\ee 

By (\ref{case1}) and (\ref{case1''}) we can apply (\ref{min}) to the second and to the third factor in (\ref{trest}) and obtain
$$
\left\|(S_{x+h x})^\gamma_\beta\right\|_{q'2B_Rx}\left\|\Pi_{x+h\gamma}^-(\varphi_{x+h}^\lambda)\right\|_{p\infty B_R x}\lesssim
$$
$$
\lesssim \lambda^{\alpha-1}(\lambda+|h|)^{\frac d2-|\gamma|_p+\alpha+|\gamma|-2-\alpha+1}(|h|+R)^{|\beta|-|\gamma|+|\gamma|_p-\alpha}=
$$
\be\lab{trest1}
=\lambda^{\alpha-1}(\lambda+|h|)^{\frac d2+\alpha}(\lambda+|h|)^{|\gamma|-|\gamma|_p-(\alpha+1)}(|h|+R)^{|\beta|-|\gamma|+|\gamma|_p-\alpha}.
\ee

Note that when using (\ref{G0estS}) we can only consider the case $|\beta|-|\gamma|+|\gamma|_p-\alpha\geq 0$ as otherwise by Proposition \ref{GestP} the left hand side of (\ref{G0estS}) with $\mu=\beta$ is zero. Recalling also (\ref{case1'}) we can apply (\ref{min}) to the last two factors in the right hand side of (\ref{trest1}). This yields
$$
\left\|(S_{x+h x})^\gamma_\beta\right\|_{q'2B_Rx}\left\|\Pi_{x+h\gamma}^-(\varphi_{x+h}^\lambda)\right\|_{p\infty B_R x}\lesssim
$$
\be\lab{trest2}
\lesssim \lambda^{\alpha-1}(\lambda+|h|)^{\frac d2+\alpha}(\lambda+|h|+R)^{|\beta|-2\alpha-1}.
\ee

If 
\be\lab{case2}
\frac d2-|\gamma|_p+\alpha\leq 0
\ee
then by (\ref{case1''}), (\ref{case2}), and by the obvious inequality $|\gamma|_p\geq 0$ we have
\be\lab{case2'}
|\gamma|-2=|\gamma|-|\gamma|_p+|\gamma|_p-2\geq (\alpha+1)-2+|\gamma|_p\geq 2\alpha-1+\frac d2.
\ee

Using this inequality, (\ref{p-est}) and (\ref{G0estS}) in the case $\frac d2-|\gamma|_p+\alpha\leq 0$ for $h\in \R^4$, $|h|\leq 2\bar{\lambda}'$, $\mu=\beta\in \M$, $1\leq q< q'<p^*<2$ we obtain for $0<R\leq 2\bar{\lambda}'$
$$
\left\|(S_{x+h x})^\gamma_\beta\right\|_{q'2B_Rx}\left\|\Pi_{x+h\gamma}^-(\varphi_{x+h}^\lambda)\right\|_{p\infty B_R x}\lesssim \lambda^{|\gamma|-2}(|h|+R)^{|\beta|-|\gamma|+\frac d2}=
$$
\be\lab{trest3}
=\lambda^{2\alpha-1+\frac d2}\lambda^{|\gamma|-1-2\alpha-\frac d2}(|h|+R)^{|\beta|-|\gamma|+\frac d2}.
\ee

Recalling that $|\beta|-|\gamma|+\frac d2>0$ in (\ref{germe1}) and using (\ref{case2'}) we can apply (\ref{min}) to the second and to the third factor in (\ref{trest3}) and obtain
$$
\left\|(S_{x+h x})^\gamma_\beta\right\|_{q'2B_Rx}\left\|\Pi_{x+h\gamma}^-(\varphi_{x+h}^\lambda)\right\|_{p\infty B_R x}\lesssim
$$
\be\lab{trest4}
\lesssim \lambda^{2\alpha-1+\frac d2}(\lambda+|h|+R)^{|\beta|-2\alpha-1}.
\ee
Finally note that $\alpha+\frac d2 >0$, and hence we can rewrite (\ref{trest4}) as follows
$$
\left\|(S_{x+h x})^\gamma_\beta\right\|_{q'2B_Rx}\left\|\Pi_{x+h\gamma}^-(\varphi_{x+h}^\lambda)\right\|_{p\infty B_R x}\lesssim \lambda^{\alpha-1}(\lambda+|h|)^{\alpha+\frac d2}(\lambda+|h|+R)^{|\beta|-2\alpha-1}.
$$
Substituting this inequality and (\ref{trest2}) into (\ref{germe1}) we obtain (\ref{gestM}).
This completes the proof

\epr

Now using the coherence condition established in the previous lemma we can obtain the following result.
\bl\lab{dPidGPi}
Assume that the conditions of Lemma \ref{reclem} are satisfied. Then for all $r\in \N$, $r>1-\alpha$, $0<R\leq \bar{\lambda}'=\frac{\bar{\lambda}_\prec}{4}$, $1\leq q<p^*<2$ the germ $F_x(y):=\delta\Pi_{0\beta}^-(y)-\sum_{\gamma\in \Ml}\Pi_{x\gamma}^-(y)\circ ({\rm d}\G_x)_\beta^\gamma$, where ${\rm d}\G_x$ is defined with the help of $\upsilon \in \L_{p^*}(L_2)$ such that $\left\|\left\|\upsilon\right\|_{L_2}\right\|_{\L_{p^*}}\leq 1$, satisfies
\be\lab{dPcoh}
\left\|F\right\|_{q2\bar{\lambda}'R r}^{2\alpha+\frac d2-1~ |\beta|-1-2\alpha}
=\sup_{0<\lambda\leq \bar{\lambda}'}\frac{\left\|\sup_{\varphi\in \mB^r}\left\|F_x(\varphi^\lambda_x)\right\|_{\L_q}\right\|_{L_2(B_R,x)}}{\lambda^{2\alpha+\frac d2-1}(\lambda+R)^{|\beta|-1-2\alpha}}\lesssim 1,
\ee
and the constant in the inequality is independent of $\rho$, $\upsilon$ and $R$.
\el

\bpr
Note that by examining the list of $\beta\in \M'$ with $|\beta|<2$ from Lemma \ref{ind<2} one can deduce that the condition $\beta\neq 0$ ensures that $|\beta|-1-2\alpha\geq 0$. Thus
by Lemma \ref{reclem} and by Proposition \ref{dGG-dif} the germ $F_x$ satisfies the conditions of Corollary \ref{reccor}, and the statement of this lemma immediately follows from Corollary \ref{reccor}.

\epr

Finally, using Lemma \ref{dPidGPi} we can obtain the desired estimate for the Malliavin derivative of $\Pi_{0\beta}^-$.
\bp\lab{dPiestprp}
Assume that (\ref{dPcoh}) holds for some $\beta\in \M_{\geq 0}$, $\beta\neq 0$, $|\beta|<2$, for some $r\in \N$, $r>1-\alpha$, for all $1\leq p< \infty$, $0<R\leq \bar{\lambda}'=\frac{\bar{\lambda}_\prec}{4}$, $1\leq q<p^*<2$, and with ${\rm d}\G_x$  defined with the help of $\upsilon \in \L_{p^*}(L_2)$ such that $\left\|\left\|\upsilon\right\|_{L_2}\right\|_{\L_{p^*}}\leq 1$, $\delta=\delta_\upsilon$, and that (\ref{G0est}) holds for $\mu=\beta$, for all $1\leq q<p^*<2$, and with ${\rm d}\G_x$  defined with the help of $\upsilon \in \L_{p^*}(L_2)$ such that $\left\|\left\|\upsilon\right\|_{L_2}\right\|_{\L_{p^*}}\leq 1$, and for arbitrary $\gamma\in\M\setminus\M_{pp}=\M'$. Suppose also that for all $1\leq p<\infty$, $\gamma\in\M'$ with $|\gamma|_\prec<|\beta|_\prec$ one has
\be\lab{Pigest}
\left\|\Pi_{x\gamma}^-\right\|_{p\bar{\lambda}_\prec r}^{|\gamma|-2}=\sup_{0<\lambda\leq \bar{\lambda}_\prec}\sup_{\varphi\in \mB^r}\frac{\left\|\Pi_{x\gamma}^-(\varphi^\lambda_x)\right\|_{\L_p}}{\lambda^{|\gamma|-2}}\lesssim 1,
\ee
where the constant in the inequality does not depend on $\rho$.

Then for $r'=r+3$, any $0<\lambda\leq \bar{\lambda}'=\frac{\bar{\lambda}_\prec}{4}$, $y\in \R^4$, $|y|+\lambda\leq \bar{\lambda}'$, $0<R\leq \bar{\lambda}'$, $1\leq q<p^*<2$, $\upsilon \in \L_{p^*}(L_2)$ such that $\left\|\left\|\upsilon\right\|_{L_2}\right\|_{\L_{p^*}}\leq 1$ one has 
\be\lab{dPiestR}
\left\|\delta\Pi_{0\beta}^-\right\|_{q2\bar{\lambda}'Rr}=\sup_{0<\lambda\leq \bar{\lambda}'}\frac{\left\|\sup_{\varphi\in \mB^r}\left\|\delta\Pi_{0\beta}^-(\varphi^\lambda_{x})\right\|_{\L_q}\right\|_{L_2(B_R,x)}}{\lambda^{\alpha-2}(\lambda+R)^{|\beta|-\alpha+\frac d2}}\lesssim  1,
\ee
and
\be\lab{dPiest}
\sup_{\varphi\in \mB^{r'}}\left\|\delta\Pi_{0\beta}^-(\varphi^\lambda_y)\right\|_{\L_q}\lesssim \lambda^{\alpha-2-\frac d2}(\lambda+|y|)^{|\beta|-\alpha+\frac d2},
\ee
where $\delta=\delta_\upsilon$, and the constants in the inequalities do not depend on $\upsilon$, $\rho$, $R$, $y$ and $\lambda$.
In particular,
\be\lab{Pi0best}
\left\|\delta\Pi_{0\beta}^-\right\|_{q\bar{\lambda}'r'}^{|\beta|-2}=\sup_{0<\lambda\leq \bar{\lambda}'}\sup_{\varphi\in \mB^{r'}}\frac{\left\|\delta\Pi_{0\beta}^-(\varphi^\lambda_0)\right\|_{\L_q}}{\lambda^{|\beta|-2}}\lesssim 1.
\ee

\ep

\bpr
Fix $1\leq q<p^*<2$ and assume that ${\rm d}\G_{x}$ is defined with the help of $\upsilon\in \L_{p^*}(L_2)$ such that $\left\|\left\|\upsilon\right\|_{L_2}\right\|_{\L_{p^*}}\leq 1$.

First note that by (\ref{Pi-0}), (\ref{prp2}), and (\ref{tr2G*}) the germ $F_x$ defined in Lemma \ref{dPidGPi} can be written in the form 
$$
F_x(y)=\delta\Pi_{0\beta}^-(y)-\sum_{\tiny\begin{array}{c}\gamma\in \M': \\ |\gamma|_\prec<|\beta|_\prec \\ \alpha\leq |\gamma|< |\beta|+\frac d2 \end{array}}\Pi_{x\gamma}^-(y)\circ ({\rm d}\G_x)_\beta^\gamma,
$$
and hence, taking into account that by translational invariance in law $\left\|\Pi_{x\gamma}^-(\varphi^\lambda_x)\right\|_{\L_p}$ does not depend on $x$, by the triangle and by the H\"{o}lder inequalities we obtain
$$
\left\|\sup_{\varphi\in \mB^r}\left\|\delta\Pi_{0\beta}^-(\varphi^\lambda_{x+y})\right\|_{\L_q}\right\|_{L_2(B_R,x)}
\leq \left\|\sup_{\varphi\in \mB^r}\left\|\delta\Pi_{0\beta}^-(\varphi^\lambda_{x})\right\|_{\L_q}\right\|_{L_2(B_{R+|y|},x)}\leq
$$
$$
\leq \left\|\sup_{\varphi\in \mB^r}\left\|F_x(\varphi^\lambda_{x})\right\|_{\L_q}\right\|_{L_2(B_{R+|y|},x)} + \sum_{\tiny\begin{array}{c}\gamma\in \M': \\ |\gamma|_\prec<|\beta|_\prec \\ \alpha\leq |\gamma|< |\beta|+\frac d2 \end{array}}\left\|({\rm d}\G_x)_\beta^\gamma\right\|_{q_1 2 B_{R+|y|}x}\sup_{\varphi\in \mB^r}\left\|\Pi_{x\gamma}^-(\varphi^\lambda_x)\right\|_{\L_{q_2}},
$$
where $\frac 1q=\frac 1{q_1}+\frac 1{q_2}$, $1\leq q< q_1< p^*<2$, $1\leq q_2<\infty$.

Now, for $0<\lambda\leq \bar{\lambda}'$, $|y|+R\leq \bar{\lambda}'$, by (\ref{dPcoh}), (\ref{G0est}) and (\ref{Pigest}) we can estimate the right hand side of the last inequality as follows
$$
\left\|\sup_{\varphi\in \mB^r}\left\|\delta\Pi_{0\beta}^-(\varphi^\lambda_{x+y})\right\|_{\L_q}\right\|_{L_2(B_R,x)}\lesssim \lambda^{2\alpha+\frac d2 -1}(\lambda+R+|y|)^{|\beta|-2\alpha-1}+
$$
$$
+\sum_{\tiny\begin{array}{c}\gamma\in \M': \\ \alpha\leq |\gamma|_\prec<|\beta|_\prec \\ |\gamma|< |\beta|+\frac d2 \end{array}}\lambda^{|\gamma|-2}(R+|y|)^{|\beta|-|\gamma|+\frac d2}=
$$
$$
= \lambda^{\alpha-2}\lambda^{\alpha+\frac d2 +1}(\lambda+R+|y|)^{|\beta|-2\alpha-1}+\sum_{\tiny\begin{array}{c}\gamma\in \M': \\ \alpha\leq |\gamma|_\prec<|\beta|_\prec \\ |\gamma|< |\beta|+\frac d2 \end{array}}\lambda^{\alpha-2}\lambda^{|\gamma|-\alpha}(R+|y|)^{|\beta|-|\gamma|+\frac d2}\lesssim
$$
$$
\lesssim \lambda^{\alpha-2}(\lambda+R+|y|)^{|\beta|-\alpha+\frac d2},
$$
where the last inequality is obtained by applying (\ref{min}) to the second and the third factor in each term of the sum. This is possible since by the definition of $\alpha$ one has $\alpha+\frac d2 +1>0$, by the definition of the homogeneity $|\beta|-2\alpha-1\geq 0$ for $\beta\neq 0$, $\beta \in \M'$ as $\beta(\bg)\neq 0$, and $|\gamma|-\alpha\geq 0$, $|\beta|-|\gamma|+\frac d2> 0$ by the choice of the limits in the sum over $\gamma$. 

For $y=0$ we obtain (\ref{dPiestR}) from the previous inequality, and for $R=\lambda$ we get from it for any $\varphi\in \mB^{r}$
$$
\left\|\delta\Pi_{0\beta}^-(\varphi^\lambda_{x+y})\right\|_{q2B_\lambda x}\lesssim \lambda^{\alpha-2}(\lambda+|y|)^{|\beta|-\alpha+\frac d2},~0<\lambda\leq \bar{\lambda}', ~|y|+\lambda\leq \bar{\lambda}'.
$$

Note that for $\varphi\in \mB^{r'}$ and $|\n|\leq 3$ we have $\partial^\n\varphi \in \mB^r$, and we obtain from the previous estimate
$$
\left\|\partial^\n\delta\Pi_{0\beta}^-(\varphi^\lambda_{x+y})\right\|_{q2B_\lambda x}=\lambda^{-|\n|}\left\|\delta\Pi_{0\beta}^-((\partial^\n\varphi)^\lambda_{x+y})\right\|_{q2B_\lambda x}\lesssim
$$
$$
\lesssim \lambda^{\alpha-2-|\n|}(\lambda+|y|)^{|\beta|-\alpha+\frac d2},~0<\lambda\leq \bar{\lambda}', ~|y|+\lambda\leq \bar{\lambda}'.
$$

Now (\ref{dPiest}) follows from the last estimate and Proposition \ref{sob} with $k=3$, $R=\lambda$ applied to $f(y)=\delta\Pi_{0\beta}^-(\varphi^\lambda_{y})$, and (\ref{Pi0best}) is (\ref{dPiest}) with $y=0$. This completes the proof.

\epr

From Proposition \ref{dPiestprp} we also immediately obtain estimates for the Malliavin derivative of $\Pi_{0\beta}$ using Schauder estimates.
\bp\lab{dPi0best}
Assume that the conditions of Proposition \ref{dPiestprp} hold and that $\bar{\lambda}'=\frac{\bar{\lambda}_\prec}{4}\geq 2$. 
Then for any $r'\in \N$, $r'\geq r+3$, $1\leq q<p^*<2$ and $\upsilon \in \L_{p^*}(L_2)$ such that $\left\|\left\|\upsilon\right\|_{L_2}\right\|_{\L_{p^*}}\leq 1$ one has
\be\lab{dPindest}
\left\|\delta\Pi_{0\beta}\right\|_{q\bar{\lambda}'' r'}^{|\beta|}=\sup_{0<\lambda\leq \bar{\lambda}''}\sup_{\varphi\in \mB^{r'}}\frac{\left\|\delta\Pi_{0\beta}(\varphi^\lambda_0)\right\|_{\L_q}}{\lambda^{|\beta|}}\lesssim 1, ~\bar{\lambda}'':=\frac{\bar{\lambda}'}2\geq 1,
\ee
where $\delta=\delta_\upsilon$, and the constant in the inequality is independent of $\upsilon$ and $\rho$.
\ep

\bpr
Since $\Pi_{0\beta}$ depends on $\xi$ polynomially, one obviously has by (\ref{defP}) using the notation as in (\ref{emb1c})
$$
\delta\Pi_{0\beta}(y)=(\K\delta\Pi_{0\beta}^-)(y)-\sum_{{\n}\in \N^4: |{\n}|<|\beta|}\partial^{\n}(\K\delta\Pi_{0\beta}^-)(0)\frac{y^{\n}}{{\n}!}=(\K\delta\Pi_{0\beta}^-)_0^{|\beta|}(y).
$$
Observe also that $r'\geq r+3>1-\alpha+3=\frac 92 +\e>||\beta||$ as $||\beta||<2$, where $||\beta||$ is the modulus of $|\beta|$.
Hence, the statement follows from (\ref{Pi0best}), and from (\ref{emb1}) (resp. (\ref{emb1c})) with $x=0$, $F=\delta\Pi_{0\beta}^-$ if $|\beta|<0$ (resp. if $|\beta|>0$).

\epr

Since if $\beta\in \M_{\geq 0}$ the assumptions of Proposition \ref{indstepprp} are stronger than those in Proposition \ref{GestP}, which implies (\ref{G0estS}), and, as for $|\beta|<2$ one has $2-\alpha=2\frac 12 +\e> |\beta|$, the other conditions in Lemmas \ref{reclem}, \ref{dPidGPi} and in Propositions \ref{dPiestprp}, \ref{dPi0best} are also weaker than the assumptions of Proposition \ref{indstepprp}, we immediately obtain the following corollary from them.
\bc\lab{dP-onpp}
Under the induction step assumptions in Proposition \ref{indstepprp} estimates (\ref{dPiestR}), (\ref{dPiest}), (\ref{Pi0best}) hold with the same conditions as in Proposition \ref{dPiestprp}, and estimate (\ref{dPindest}) holds with the same conditions as in Proposition \ref{dPi0best}.
\ec

The remaining results in this section will be used in Section \ref{indstGGSpp} to establish the induction step estimates for $(\G_{xy})^\gamma_\beta$, $({\rm d}\G_x)^\gamma_\beta$, $(S_{xy})^\gamma_\beta$, $\gamma\in \M_{pp}$.

\bl\lab{PidGest}
Assume that for some $\beta\in \M_{\geq 0}$, $\beta\neq 0$, $|\beta|<2$ (\ref{G0est}) holds for $\mu=\beta$, for all $\gamma\in\M\setminus\M_{pp}=\M'$, $1\leq q<p^*<2$ and ${\rm d}\G_x$ defined with the help of $\upsilon \in \L_{p^*}(L_2)$ such that $\left\|\left\|\upsilon\right\|_{L_2}\right\|_{\L_{p^*}}\leq 1$. 

Then the following statements are true.

(i) If for all $1\leq p<\infty$, $\gamma\in\M'$ with $|\gamma|_\prec<|\beta|_\prec$, $x\in \R^4$ one has for some $r\in \N$
\be\lab{Pigest-}
\left\|\Pi_{x\gamma}^-\right\|_{p\bar{\lambda}_\prec r}^{|\gamma|-2}=\sup_{0<\lambda\leq \bar{\lambda}_\prec}\sup_{\varphi\in \mB^r}\frac{\left\|\Pi_{x\gamma}^-(\varphi^\lambda_x)\right\|_{\L_p}}{\lambda^{|\gamma|-2}}\lesssim 1,
\ee
where the constant in the inequality does not depend on $\rho$, then for any $0<R\leq \bar{\lambda}_\prec$,  $1\leq q<p^*<2$ and ${\rm d}\G_x$ defined with the help of $\upsilon \in \L_{p^*}(L_2)$ such that $\left\|\left\|\upsilon\right\|_{L_2}\right\|_{\L_{p^*}}\leq 1$ the germ 
$$
F_x(y):=\sum_{\tiny\begin{array}{c}\gamma\in \Ml: \\ 2-\e\leq |\gamma|\end{array}}\Pi_{x\gamma}^-(y)\circ ({\rm d}\G_x)_\beta^\gamma
$$ 
satisfies
\be\lab{PG-est}
\left\|F\right\|_{q2\bar{\lambda}_\prec R r}^{\alpha-2+\frac d2~|\beta|-\alpha}=\sup_{0<\lambda\leq \bar{\lambda}_\prec}\frac{\left\|\sup_{\varphi\in \mB^r}\left\|F_x(\varphi^\lambda_x)\right\|_{\L_q}\right\|_{L_2(B_R,x)}}{\lambda^{\alpha-2+\frac d2}(\lambda+R)^{|\beta|-\alpha}}\lesssim 1,
\ee
where the constant in the inequality does not depend on $R$, $\rho$ and $\upsilon$.

(ii) If for all $1\leq p<\infty$, $\gamma\in\M'$ with $|\gamma|_\prec<|\beta|_\prec$, $x\in \R^4$ one has for some $r\in \N$
\be\lab{Pigest+}
\left\|\Pi_{x\gamma}\right\|_{p\bar{\lambda}_\prec r}^{|\gamma|}=\sup_{0<\lambda\leq \bar{\lambda}_\prec}\sup_{\varphi\in \mB^r}\frac{\left\|\Pi_{x\gamma}(\varphi^\lambda_x)\right\|_{\L_p}}{\lambda^{|\gamma|}}\lesssim 1,
\ee
where the constant in the inequality does not depend on $\rho$, then for any $0<R\leq \bar{\lambda}_\prec$,  $1\leq q<p^*<2$ and ${\rm d}\G_x$ defined with the help of $\upsilon \in \L_{p^*}(L_2)$ such that $\left\|\left\|\upsilon\right\|_{L_2}\right\|_{\L_{p^*}}\leq 1$ the germ 
$$
F_x(y):=\sum_{\tiny\begin{array}{c}\gamma\in \Ml: \\ 2-\e\leq |\gamma|\end{array}}\Pi_{x\gamma}(y)\circ ({\rm d}\G_x)_\beta^\gamma
$$ 
satisfies
\be\lab{PGest}
\left\|F\right\|_{q2\bar{\lambda}_\prec R r}^{\alpha+\frac d2~|\beta|-\alpha}=\sup_{0<\lambda\leq \bar{\lambda}_\prec}\frac{\left\|\sup_{\varphi\in \mB^r}\left\|F_x(\varphi^\lambda_x)\right\|_{\L_q}\right\|_{L_2(B_R,x)}}{\lambda^{\alpha+\frac d2}(\lambda+R)^{|\beta|-\alpha}}\lesssim 1,
\ee
where the constant in the inequality does not depend on $R$, $\rho$ and $\upsilon$.

\el

\bpr
Fix $1\leq q<p^*<2$ and assume that ${\rm d}\G_{x}$ is defined with the help of $\upsilon\in \L_{p^*}(L_2)$ such that $\left\|\left\|\upsilon\right\|_{L_2}\right\|_{\L_{p^*}}\leq 1$.

(i) By (\ref{Pi-0}), (\ref{tr2G*}) and (\ref{prp2}) we can rewrite the definition of $F_x$ as follows
$$
F_x(y)=\sum_{\tiny\begin{array}{c}\gamma\in \M': \\ |\gamma|_\prec <|\beta|_\prec \\ 2-\e\leq |\gamma|<|\beta|+\frac d2\end{array}}\Pi_{x\gamma}^-(y)\circ ({\rm d}\G_x)_\beta^\gamma,
$$
Note that by the translational invariance in law $\left\|\Pi_{x\gamma}^-(\varphi^\lambda_x)\right\|_{\L_p}$ does not depend on $x$. Hence by (\ref{G0est}), (\ref{Pigest-}), by the triangle and by the H\"{o}lder inequalities we have for $0<R\leq \bar{\lambda}_\prec$, $0<\lambda\leq \bar{\lambda}_\prec$, $\frac 1p +\frac 1{q'}=\frac 1q$, $1\leq q<q'<p^*<2$, $1\leq p<\infty$
$$
\left\|\sup_{\varphi\in \mB^r}\left\|F_x(\varphi^\lambda_x)\right\|_{\L_q}\right\|_{L_2(B_R,x)}\leq \sum_{\tiny\begin{array}{c}\gamma\in \M': \\ |\gamma|_\prec <|\beta|_\prec \\ 2-\e\leq |\gamma|<|\beta|+\frac d2\end{array}}\sup_{\varphi\in \mB^r}\left\|\Pi_{x\gamma}^-(\varphi^\lambda_x)\right\|_{\L_p}\left\|({\rm d}\G_{x})^\gamma_\beta\right\|_{q'2B_Rx}\lesssim
$$
$$
\lesssim \sum_{\tiny\begin{array}{c}\gamma\in \M': \\ |\gamma|_\prec <|\beta|_\prec \\ 2-\e\leq |\gamma|<|\beta|+\frac d2\end{array}}\lambda^{|\gamma|-2}R^{\frac d2+|\beta|-|\gamma|}=\sum_{\tiny\begin{array}{c}\gamma\in \M': \\ |\gamma|_\prec <|\beta|_\prec \\ \frac d2+\alpha\leq |\gamma|<|\beta|+\frac d2\end{array}}\lambda^{\alpha+\frac d2 -2}\lambda^{|\gamma|-\frac d2-\alpha}R^{\frac d2+|\beta|-|\gamma|}\leq
$$
\be\lab{eql}
\leq\sum_{\tiny\begin{array}{c}\gamma\in \M': \\ |\gamma|_\prec <|\beta|_\prec \\ \frac d2+\alpha\leq |\gamma|<|\beta|+\frac d2\end{array}}\lambda^{\alpha+\frac d2 -2}(\lambda+R)^{|\beta|-\alpha}\lesssim \lambda^{\alpha+\frac d2 -2}(\lambda+R)^{|\beta|-\alpha},
\ee
where to obtain (\ref{eql}) from the preceding expression we applied (\ref{min}) to the last two factors in each term of the sum over $\gamma$. Estimate (\ref{eql}) justifies (\ref{PG-est}).

(ii)
By (\ref{Pi0}) $\gamma\in \M_{\geq 0}\cup \M_{pp}$ in the sum in the definition of $F_x$.
By (\ref{cdn30}) ${\rm d}\G_x^*\z^{\delta_{\n}}=0$ for $\n\in \N^4$, $|\n|\geq 2-\e$. Thus for such $\n$
\be\lab{dGdvan}
({\rm d}\G_x)_\beta^{\delta_\n}=({\rm d}\G_x^*\z^{\delta_{\n}})_\beta=0,~|\n|\geq 2-\e,
\ee
and using also (\ref{tr2G*}) and (\ref{prp2}) we can rewrite the definition of $F_x$ as follows
$$
F_x(y)=\sum_{\tiny\begin{array}{c}\gamma\in \M_{\geq 0}: \\ |\gamma|_\prec <|\beta|_\prec \\ 2-\e\leq |\gamma|<|\beta|+\frac d2\end{array}}\Pi_{x\gamma}(y)\circ ({\rm d}\G_x)_\beta^\gamma,
$$
Note that by the translational invariance in law $\left\|\Pi_{x\gamma}(\varphi^\lambda_x)\right\|_{\L_p}$ does not depend on $x$. Hence by (\ref{G0est}), (\ref{Pigest+}), by the triangle and by the H\"{o}lder inequalities we have for $0<R\leq \bar{\lambda}_\prec$, $0<\lambda\leq \bar{\lambda}_\prec$, $\frac 1p +\frac 1{q'}=\frac 1q$, $1\leq q<q'<p^*<2$, $1\leq p<\infty$
$$
\left\|\sup_{\varphi\in \mB^r}\left\|F_x(\varphi^\lambda_x)\right\|_{\L_q}\right\|_{L_2(B_R,x)}\leq \sum_{\tiny\begin{array}{c}\gamma\in \M_{\geq 0}: \\ |\gamma|_\prec <|\beta|_\prec \\ 2-\e\leq |\gamma|<|\beta|+\frac d2\end{array}}\sup_{\varphi\in \mB^r}\left\|\Pi_{x\gamma}(\varphi^\lambda_x)\right\|_{\L_p}\left\|({\rm d}\G_{x})^\gamma_\beta\right\|_{q'2B_Rx}\lesssim
$$
$$
\lesssim \sum_{\tiny\begin{array}{c}\gamma\in \M_{\geq 0}: \\ |\gamma|_\prec <|\beta|_\prec \\ 2-\e\leq |\gamma|<|\beta|+\frac d2\end{array}}\lambda^{|\gamma|}R^{\frac d2+|\beta|-|\gamma|}=\sum_{\tiny\begin{array}{c}\gamma\in \M_{\geq 0}: \\ |\gamma|_\prec <|\beta|_\prec \\ \frac d2+\alpha\leq |\gamma|<|\beta|+\frac d2\end{array}}\lambda^{\alpha+\frac d2}\lambda^{|\gamma|-\frac d2-\alpha}R^{\frac d2+|\beta|-|\gamma|}\leq
$$
\be\lab{eql1}
\leq\sum_{\tiny\begin{array}{c}\gamma\in \M_{\geq 0}: \\ |\gamma|_\prec <|\beta|_\prec \\ \frac d2+\alpha\leq |\gamma|<|\beta|+\frac d2\end{array}}\lambda^{\alpha+\frac d2}(\lambda+R)^{|\beta|-\alpha}\lesssim \lambda^{\alpha+\frac d2}(\lambda+R)^{|\beta|-\alpha},
\ee
where to obtain (\ref{eql1}) from the preceding expression we applied (\ref{min}) to the last two factors in each term of the sum over $\gamma$. Estimate (\ref{eql1}) justifies (\ref{PGest}). This completes the proof.

\epr

The next lemma is a preparation for the proof of Proposition \ref{Fgerm0est}.
\bl
Assume that for some $\beta\in \M_{\geq 0}$, $\beta\neq 0$, $|\beta|<2$ estimate (\ref{G0est}) (resp. (\ref{G0estS})) holds for $\mu=\beta$, and for all $\gamma\in\M\setminus\M_{pp}=\M'$, $0<R\leq \bar{\lambda}_\prec$ (resp. for all $0<R\leq \frac{\bar{\lambda}_\prec}2$, $y\in \R^4$, $|y|\leq \frac{\bar{\lambda}_\prec}2$), $1\leq q<p^*<2$ and ${\rm d}\G_x$ (resp. $S_{x+y x}$) defined with the help of $\upsilon \in \L_{p^*}(L_2)$ such that $\left\|\left\|\upsilon\right\|_{L_2}\right\|_{\L_{p^*}}\leq 1$. Suppose also that for all $1\leq p<\infty$, $\gamma\in\M'$ with $|\gamma|_\prec<|\beta|_\prec$, $x\in \R^4$, and for some $r\in \N$, $r>4-\alpha$ (\ref{Pigest-}) holds, and that $\bar{\lambda}''=\frac{\bar{\lambda}'}2 =\frac{\bar{\lambda}_\prec}8\geq 1$.

Then for any $0<R\leq \bar{\lambda}''=\frac{\bar{\lambda}_\prec}8$,  $1\leq q<p^*<2$ and ${\rm d}\G_x$ defined with the help of $\upsilon \in \L_{p^*}(L_2)$ such that $\left\|\left\|\upsilon\right\|_{L_2}\right\|_{\L_{p^*}}\leq 1$ the germ 
$$
F_x(y):=\delta\Pi_{0\beta}(y)-\sum_{\tiny\begin{array}{c}\gamma\in \Ml: \\  |\gamma|<2-\e \end{array}}\Pi_{x\gamma}(y)\circ ({\rm d}\G_x)_\beta^\gamma
$$ 
satisfies
\be\lab{PGest<}
\left\|F\right\|_{q2\frac{\bar{\lambda}''}2 R r}^{\alpha+\frac d2~|\beta|-\alpha}=\sup_{0<\lambda\leq \frac{\bar{\lambda}''}2}\frac{\left\|\sup_{\varphi\in \mB^r}\left\|F_x(\varphi^\lambda_x)\right\|_{\L_q}\right\|_{L_2(B_R,x)}}{\lambda^{\alpha+\frac d2}(\lambda+R)^{|\beta|-\alpha}}\lesssim 1,
\ee
where $\delta=\delta_\upsilon$, and the constant in the inequality does not depend on $\upsilon$, $\rho$ and $R$.
\el

\bpr
Fix $1\leq q<p^*<2$ and assume that ${\rm d}\G_{x}$ and $S_{x+y x}$ are defined with the help of $\upsilon\in \L_{p^*}(L_2)$ such that $\left\|\left\|\upsilon\right\|_{L_2}\right\|_{\L_{p^*}}\leq 1$.

Consider the germ
$$
G_x(y):=\delta\Pi_{0\beta}^-(y)-\sum_{\tiny\begin{array}{c}\gamma\in \Ml: \\  |\gamma|<2-\e \end{array}}\Pi_{x\gamma}^-(y)\circ ({\rm d}\G_x)_\beta^\gamma=
$$ 
$$
=\delta\Pi_{0\beta}^-(y)-\sum_{\tiny\begin{array}{c}\gamma\in \M': \\  |\gamma|<2-\e \end{array}}\Pi_{x\gamma}^-(y)\circ ({\rm d}\G_x)_\beta^\gamma.
$$

Checking the list of multi--indices $\gamma\in \M'$ with $|\gamma|<2$ from Lemma \ref{ind<2} shows that the condition $|\gamma|<2-\e$ for multi--indices $\gamma \in \M'$ is equivalent to $|\gamma|<2$. From this observation and (\ref{defP}) we deduce that for each $x\in \R^4$ the function
$F_x(y)$ is obtained from $(\K G)_x(y)$ by subtracting a polynomial in $y$ of parabolic degree less than $2$ (or equivalently less than $2-\e$). By (\ref{cdn3}) this polynomial must be the Taylor polynomial of $(\K G)_x(y)$ at $x$ of parabolic degree $1$. Thus $F=(\K G)^{2-\e}$, and we can apply (\ref{emb2c}) with $b=-\e$, $r>4-\alpha$ to get
\be\lab{emb2c'}
\left\|F \right\|_{q2\frac{\bar{\lambda}''}2 B_R r}^{2-\e}\lesssim  \left\|G\right\|_{q2\bar{\lambda}'' B_R r}^{-\e}.
\ee

Now we estimate $\left\|G\right\|_{q2 \bar{\lambda}'' B_R r}^{-\e}$ for $0<R\leq \bar{\lambda}''$. We split the sum in the definition of $G$ as follows
$$
G_x(y):=\large\left(\delta\Pi_{0\beta}^-(y)-\sum_{\gamma\in \Ml}\Pi_{x\gamma}^-(y)\circ ({\rm d}\G_x)_\beta^\gamma\large\right)+\sum_{\tiny\begin{array}{c}\gamma\in \Ml: \\  |\gamma|\geq 2-\e \end{array}}\Pi_{x\gamma}^-(y)\circ ({\rm d}\G_x)_\beta^\gamma
$$ 
and apply Lemma \ref{dPidGPi} to the expression in the large brackets and Lemma \ref{PidGest} to the remaining sum over $\gamma$. This yields for $0<\lambda\leq \bar{\lambda}''$ with the help of the triangle inequality
$$
\left\|\sup_{\varphi\in \mB^r}\left\|G_x(\varphi^\lambda_x)\right\|_{\L_q}\right\|_{L_2(B_R,x)}\lesssim \lambda^{2\alpha+\frac d2-1}(\lambda+R)^{|\beta|-1-2\alpha}+\lambda^{\alpha-2+\frac d2}(\lambda+R)^{|\beta|-\alpha}\lesssim
$$
$$
\lesssim \lambda^{\alpha-2+\frac d2}(\lambda+R)^{|\beta|-\alpha}\large\left(\frac{\lambda}{\lambda+R}\large\right)^{\alpha+1}+ \lambda^{\alpha-2+\frac d2}(\lambda+R)^{|\beta|-\alpha}\lesssim \lambda^{\alpha-2+\frac d2}(\lambda+R)^{|\beta|-\alpha}.
$$
For $\lambda\leq R$ we obtain from the last estimate
\be\lab{G1est}
\left\|\sup_{\varphi\in \mB^r}\left\|G_x(\varphi^\lambda_x)\right\|_{\L_q}\right\|_{L_2(B_R,x)}\lesssim \lambda^{\alpha-2+\frac d2}R^{|\beta|-\alpha}=\lambda^{-\e}R^{|\beta|-\alpha}.
\ee

If $\lambda\geq R$ we write by (\ref{tr2G*}) 
$$
G_x(y)=\delta\Pi_{0\beta}^-(y)-\sum_{\tiny\begin{array}{c}\gamma\in \M': \\ |\gamma|_\prec <|\beta|_\prec \\ \alpha \leq |\gamma|<2-\e\end{array}}\Pi_{x\gamma}^-(y)\circ ({\rm d}\G_x)_\beta^\gamma,
$$
and hence obtain by the triangle and the H\"{o}lder inequalities using translational invariance in law of $\Pi_{x\gamma}^-(\varphi^\lambda_x)$
$$
\left\|\sup_{\varphi\in \mB^r}\left\|G_x(\varphi^\lambda_x)\right\|_{\L_q}\right\|_{L_2(B_R,x)}\leq \left\|\sup_{\varphi\in \mB^r}\left\|\delta\Pi_{0\beta}^-(\varphi^\lambda_x)\right\|_{\L_q}\right\|_{L_2(B_R,x)}+
$$
$$
+\sum_{\tiny\begin{array}{c}\gamma\in \M': \\ |\gamma|_\prec <|\beta|_\prec \\ \alpha \leq |\gamma|<2-\e\end{array}}\sup_{\varphi\in \mB^r}\left\|\Pi_{x\gamma}^-(\varphi^\lambda_x)\right\|_{\L_p} \left\|({\rm d}\G_x)_\beta^\gamma\right\|_{q'2B_Rx},
$$
where $\frac 1p +\frac 1{q'}=\frac 1q$, $1\leq q< q'<p^*<2$, $1\leq p<\infty$. 

Using (\ref{dPiest}) in the first term in the right hand side and (\ref{G0est}), (\ref{Pigest-}) in the sum over $\gamma$ we obtain for $R\leq\lambda\leq \bar{\lambda}''$ with the help $\frac{R}{\lambda}\leq 1$, $\frac d2-|\beta|+\alpha=-|\beta|+2-\e>0$ and $-|\gamma|+2-\e>0$  
$$
\left\|\sup_{\varphi\in \mB^r}\left\|G_x(\varphi^\lambda_x)\right\|_{\L_q}\right\|_{L_2(B_R,x)} \lesssim \lambda^{\alpha-2-\frac d2}(\lambda+R)^{|\beta|-\alpha+\frac d2}R^{\frac d2} +
$$
$$
+\sum_{\tiny\begin{array}{c}\gamma\in \M': \\ |\gamma|_\prec <|\beta|_\prec \\ \alpha \leq |\gamma|<2-\e\end{array}}\lambda^{|\gamma|-2}R^{\frac d2+|\beta|-|\gamma|}\leq
$$
$$
\leq \lambda^{|\beta|-2}R^{\frac d2}+\lambda^{-\e}R^{|\beta|-\alpha}\sum_{\tiny\begin{array}{c}\gamma\in \M': \\ |\gamma|_\prec <|\beta|_\prec \\ \alpha \leq |\gamma|<2-\e\end{array}}\large\left(\frac{R}{\lambda}\large\right)^{-|\gamma|+2-\e}\lesssim
$$
$$
\lesssim \lambda^{-\e}R^{|\beta|-\alpha}\large\left(\frac{R}{\lambda}\large\right)^{\frac d2-|\beta|+\alpha}+   \lambda^{-\e}R^{|\beta|-\alpha}\lesssim \lambda^{-\e}R^{|\beta|-\alpha}.
$$

The last inequality and (\ref{G1est}) imply that for $0<R\leq \bar{\lambda}''$ one has
$$
\left\|G\right\|_{q2 \bar{\lambda}'' B_R r}^{-\e}\lesssim R^{|\beta|-\alpha}.
$$

Now by (\ref{emb2c'}) for $0<R\leq \bar{\lambda}''$ and  $0<\lambda\leq \frac{\bar{\lambda}''}2$
$$
\left\|\sup_{\varphi\in \mB^r}\left\|F_x(\varphi^\lambda_x)\right\|_{\L_q}\right\|_{L_2(B_R,x)}\lesssim \lambda^{2-\e}R^{|\beta|-\alpha}\lesssim \lambda^{\alpha+\frac d2}(\lambda+R)^{|\beta|-\alpha}
$$
as $|\beta|-\alpha\geq 0$. This proves (\ref{PGest<}).

\epr

From the previous lemma and from Lemma \ref{PidGest} (ii) we immediately obtain the following statement.
\bp\lab{Fgerm0est}
Assume that for some $\beta\in \M_{\geq 0}$, $\beta\neq 0$, $|\beta|<2$ estimates (\ref{G0est}) (resp. (\ref{G0estS})) holds for $\mu=\beta$, for all $\gamma\in\M\setminus\M_{pp}=\M'$, $0<R\leq \bar{\lambda}_\prec$ (resp. $0<R\leq \frac{\bar{\lambda}_\prec}2$, $y\in \R^4$, $|y|\leq \frac{\bar{\lambda}_\prec}2$), $1\leq q<p^*<2$ and ${\rm d}\G_x$ (resp. $S_{x+y x}$) defined with the help of $\upsilon \in \L_{p^*}(L_2)$ such that $\left\|\left\|\upsilon\right\|_{L_2}\right\|_{\L_{p^*}}\leq 1$. Suppose also that for all $1\leq p<\infty$, $\gamma\in\M'$ with $|\gamma|_\prec<|\beta|_\prec$, $x\in \R^4$, and for some $r\in \N$, $r>4-\alpha$ (\ref{Pigest-}) and (\ref{Pigest+}) hold, and that $\bar{\lambda}''=\frac{\bar{\lambda}'}2 =\frac{\bar{\lambda}_\prec}8\geq 1$.

Then for any $0<R\leq \bar{\lambda}''=\frac{\bar{\lambda}_\prec}8$,  $1\leq q<p^*<2$ and ${\rm d}\G_x$ defined with the help of $\upsilon \in \L_{p^*}(L_2)$ such that $\left\|\left\|\upsilon\right\|_{L_2}\right\|_{\L_{p^*}}\leq 1$ the germ 
$$
F_x(y):=\delta\Pi_{0\beta}(y)-\sum_{\gamma\in \Ml}\Pi_{x\gamma}(y)\circ ({\rm d}\G_x)_\beta^\gamma
$$ 
satisfies
\be\lab{PGest=}
\left\|F\right\|_{q2\frac{\bar{\lambda}''}2 R r}^{\alpha+\frac d2~|\beta|-\alpha}=\sup_{0<\lambda\leq \frac{\bar{\lambda}''}2}\frac{\left\|\sup_{\varphi\in \mB^r}\left\|F_x(\varphi^\lambda_x)\right\|_{\L_q}\right\|_{L_2(B_R,x)}}{\lambda^{\alpha+\frac d2}(\lambda+R)^{|\beta|-\alpha}}\lesssim 1,
\ee
where $\delta=\delta_\upsilon$, and the constant in the inequality does not depend on $\upsilon$, $\rho$ and $R$.

\ep

\bpr
We can write
$$
F_x(y)=\large\left(\delta\Pi_{0\beta}(y)-\sum_{\tiny\begin{array}{c}\gamma\in \Ml: \\  |\gamma|<2-\e \end{array}}\Pi_{x\gamma}(y)\circ ({\rm d}\G_x)_\beta^\gamma\large\right)-\sum_{\tiny\begin{array}{c}\gamma\in \Ml: \\  |\gamma|\geq 2-\e \end{array}}\Pi_{x\gamma}(y)\circ ({\rm d}\G_x)_\beta^\gamma
$$
and then estimate $\left\|F\right\|_{q2\frac{\bar{\lambda}''}2 R r}^{\alpha+\frac d2~|\beta|-\alpha}$ using the triangle inequality in the right hand side, the previous lemma for the first term in the large brackets and Lemma \ref{PidGest} (ii) for the last sum over $\gamma\in \Ml$, $|\gamma|\geq 2-\e$. This yields the desired estimate.

\epr


\subsection{The induction step for $\Pi_{x\beta}$, $\Pi_{x\beta}^-$}\lab{indstepPb}

\setcounter{equation}{0}
\setcounter{theorem}{0}

In this section we establish the induction step estimates for $\Pi_{x\beta}$, $\Pi_{x\beta}^-$.

We firstly consider the case when $|\beta|<2$ when an estimate for $\Pi_{x\beta}^-$ is obtained using the spectral gap inequality. We start with an estimate for the expectation of $\Pi_{x\beta}^-$ which appears in the right hand side of the spectral gap inequality. To estimate the expectation we shall have to impose a BPHZ condition in a special way (compare the following proposition and its proof with Proposition 5.2 in \cite{HS} and with its proof). 

\bp\lab{Pi-psi}
Let $\omega\in \D$ be an even function supported in $B_{\frac 16}$, invariant under permutations of the coordinates $x_1,x_2,x_3$, and satisfying the first condition in (\ref{omegapr*}). Assume that for some $\beta\in \M'$, $|\beta|<2$, $\beta\neq 0$ the recentered map $\Pi_{x\beta}^-$ satisfies the generalized BPHZ condition with respect to the function $\psi^{\bar{\lambda}_\prec}$,  $\bar{\lambda}_\prec>0$, where $\psi:=\omega^2\ast \omega$, and $\omega^2$ is $\omega$ rescaled by $2$. Suppose also that for all $1\leq p<\infty$, $\gamma \in\M'$ one has
\be\lab{Gest11}
\left\|(\G_{xy}\right)^\gamma_\beta\|_{\L_p}\lesssim |x-y|^{|\beta|-|\gamma|},~x,y\in \R^4,~|x-y|\leq \bar{\lambda}_\prec,
\ee
where the constant in the inequality does not depend on $x,y$ and $\rho$, and that for all $1\leq p<\infty$, $\gamma\in\M'$ with $|\gamma|_\prec<|\beta|_\prec$ and $|\gamma|<|\beta|$ one has
\be\lab{Pigest11}
\sup_{0<\lambda\leq \bar{\lambda}_\prec}\frac{\left\|\Pi_{x\gamma}^-(\omega^\lambda_x)\right\|_{\L_p}}{\lambda^{|\gamma|-2}}\lesssim 1,
\ee
where the constant in the inequality does not depend on $\rho$.
Then for any $x\in \R^4$ and $n\in \N$
\be\lab{Epsin} 
|\E(\Pi_{x\beta}^-(\psi^{\e_n}_x))|\lesssim {\e_n}^{|\beta|-2},
\ee
where $\e_n=\bar{\lambda}_\prec 2^{-n}$, $n\in \N$, and the constant in the inequality does not depend on $n$, $x$ and $\rho$.
\ep

\bpr
First note that from the definition of $\psi$ it follows that $\int_{\R^4}\psi(x)dx=1$, and 
\be\lab{psisupp}
{\rm supp}(\psi)\subset B_{\frac 12}\subset B_1. 
\ee

We also have the following identity
$$
\psi^{\frac 12}-\psi=\omega\ast (\omega^{\frac 12}-\omega^2),
$$
and ${\rm supp}(\omega^{\frac 12}-\omega^2)\subset B_{\frac 23}$.

Using the previous identity we can write for $n\in \N$
$$
\psi^{\e_{n+1}}-\psi^{\e_n}=\omega^{\e_n}\ast (\omega^{\frac 12}-\omega^2)^{\e_n},
$$
and hence
\be\lab{tele}
\psi^{\e_n}=\psi^{\e_0}+\sum_{k=0}^{n-1}(\psi^{\e_{k+1}}-\psi^{\e_k})=\psi^{\bar{\lambda}_\prec}+\sum_{k=0}^{n-1}\omega^{\e_k}\ast(\omega^{\frac 12}-\omega^2)^{\e_k}.
\ee

By translational invariance in law it suffices to establish (\ref{Epsin}) for $x=0$. Using (\ref{tele}) and the BPHZ condition $\Pi_{0\beta}^-(\psi^{\bar{\lambda}_\prec})=0$ we obtain
$$
\E(\Pi_{0\beta}^-(\psi^{\e_n}))=\sum_{k=0}^{n-1}\int_{\R^4}\E(\Pi_{0\beta}^-(\omega^{\e_k}_x))(\omega^{\frac 12}-\omega^2)^{\e_k}(x)dx=
$$
$$
=\sum_{k=0}^{n-1}\sum_{\gamma\in \M'}\int_{\R^4}\E(\Pi_{x\gamma}^-(\omega^{\e_k}_x)\circ (\G_{x0})_\beta^\gamma) (\omega^{\frac 12}-\omega^2)^{\e_k}(x)dx,
$$
where at the last step we used (\ref{PG}). 

Since by the translational invariance in law $\E(\Pi_{x\gamma}^-(\omega^{\e_k}_x))$ does not depend on $x$ and 
$$
\int_{\R^4}(\omega^{\frac 12}-\omega^2)^{\e_k}(x)dx=0
$$
we can rewrite the last identity as follows
$$
\E(\Pi_{0\beta}^-(\psi^{\e_n}))=\sum_{k=0}^{n-1}\sum_{\gamma\in \M'}\int_{\R^4}\E(\Pi_{x\gamma}^-(\omega^{\e_k}_x)\circ (\G_{x0}-{\rm Id})_\beta^\gamma) (\omega^{\frac 12}-\omega^2)^{\e_k}(x)dx.
$$

Now by (\ref{triangG}) and (\ref{triangGm}) we can further reduce this identity to
$$
\E(\Pi_{0\beta}^-(\psi^{\e_n}))=\sum_{k=0}^{n-1}\sum_{\tiny\begin{array}{c}\gamma\in \M': \\ |\gamma|_\prec<|\beta|_\prec \\ \alpha\leq |\gamma|<|\beta|\end{array}}\int_{\R^4}\E(\Pi_{x\gamma}^-(\omega^{\e_k}_x)\circ (\G_{x0}-{\rm Id})_\beta^\gamma) (\omega^{\frac 12}-\omega^2)^{\e_k}(x)dx.
$$
Since $(\omega^{\frac 12}-\omega^2)^{\e_k}$, $k\in \N$ is supported in $B_{\frac{2\bar{\lambda}_\prec}3}$ the domain of integration is reduced to $B_{\frac{2\bar{\lambda}_\prec}3}$, and using the H\"{o}lder inequality in probability with H\"{o}lder conjugate exponents $p$, $p'$, $1<p,p'<\infty$ we can estimate the terms in the sum with the help of the assumptions of the proposition as follows
$$
|\E(\Pi_{0\beta}^-(\psi^{\e_n}))|\leq
$$
$$
\leq\sum_{k=0}^{n-1}\sum_{\tiny\begin{array}{c}\gamma\in \M': \\ |\gamma|_\prec<|\beta|_\prec \\ \alpha\leq |\gamma|<|\beta|\end{array}}\int_{B_{\frac{2\bar{\lambda}_\prec}3}}\left\|\Pi_{x\gamma}^-(\omega^{\e_k}_x)\right\|_{\L_p} \left\|(\G_{x0})_\beta^\gamma)\right\|_{\L_{p'}} |(\omega^{\frac 12}-\omega^2)^{\e_k}(x)|dx\lesssim
$$
$$
\lesssim \sum_{k=0}^{n-1}\sum_{\tiny\begin{array}{c}\gamma\in \M': \\ |\gamma|_\prec<|\beta|_\prec \\ \alpha\leq |\gamma|<|\beta|\end{array}}\int_{B_{\frac{2\bar{\lambda}_\prec}3}}\e_k^{|\gamma|-2}|x|^{|\beta|-|\gamma|}|(\omega^{\frac 12}-\omega^2)^{\e_k}(x)|dx=
$$
$$
=\sum_{k=0}^{n-1}\sum_{\tiny\begin{array}{c}\gamma\in \M': \\ |\gamma|_\prec<|\beta|_\prec \\ \alpha\leq |\gamma|<|\beta|\end{array}}\int_{\R^4}\e_k^{|\beta|-|\gamma|+|\gamma|-2}|y|^{|\beta|-|\gamma|}(\omega^{\frac 12}-\omega^2)(y)dy\lesssim
$$
$$
\lesssim \sum_{k=0}^{n-1}\e_k^{|\beta|-2}=\frac{\e_n^{|\beta|-2}-1}{\e_1^{|\beta|-2}-1}\lesssim \e_n^{|\beta|-2},
$$
where we also used the substitution of variables $x=R_{\e_k}y$ in the integral in the $k$-th term in the sum over $k$. This completes the proof.

\epr

Now we can apply the spectral gap inequality to estimate $\L_p$--norms of $\Pi_{x\beta}^-$.
\bl\lab{psilem0}
Assume that the conditions of Propositions \ref{dPiestprp} and \ref{Pi-psi} are satisfied. Then for any $1\leq p<\infty$, $x\in \R^4$ and $n\in \N$
\be\lab{Epsinp} 
\left\|\Pi_{x\beta}^-(\psi^{\e_n}_x)\right\|_{\L_p} \lesssim {\e_n}^{|\beta|-2},
\ee
where $\e_n=\bar{\lambda}' 2^{-n}=\bar{\lambda}_\prec 2^{-(n+2)}$, $n\in \N$, and the constant in the inequality does not depend on $n$, $x$ and $\rho$.
\el

\bpr
By translational invariance in law it suffices to consider the case $x=0$, and by the H\"{o}lder inequality in probability one can restrict to the case when $p> 2$. 

We shall use spectral gap inequality (\ref{SpG}) for $f=\Pi_{0\beta}^-(\psi^{\e_n}_0)$, $n\in \N$,
$$
\left\|\Pi_{0\beta}^-(\psi^{\e_n}_0)\right\|_{\L_p}\leq |\E(\Pi_{0\beta}^-(\psi^{\e_n}_0))|+\sup_{\tiny \begin{array}{c}\upsilon\in \L_{p^*}(L_2): \\ \left\|\left\|\upsilon\right\|_{L_2(\R^4)}\right\|_{\L_{p^*}}\leq 1\end{array}}\left\|\delta_\upsilon \Pi_{0\beta}^-(\psi^{\e_n}_0)\right\|_{\L_q},
$$
where $1\leq q<p^*< 2< p$, and $\frac 1p+\frac{1}{p^*}=1$.

Now we can apply (\ref{Epsin}) in the first term and (\ref{Pi0best}) in the second one to obtain from the previous estimate for $n\in \N$
$$
\left\|\Pi_{0\beta}^-(\psi^{\e_n}_0)\right\|_{\L_p}\lesssim {\e_n}^{|\beta|-2}
$$
which completes the proof.

\epr

Next, we shall derive more general estimates for $\Pi_{x\beta}^-$ depending on the base point $x$ and for any value of the rescaling parameter $\lambda$ in the domain $0<\lambda\leq 2\bar{\lambda}'$. We start with the following useful lemma. 
\bl\lab{Pi-int}
Assume that for some $\beta\in \M'$, $\lambda>0$, $\varphi\in \D$ and for all $1\leq p<\infty$, $\gamma\in \M'$ with $|\gamma|_\prec<|\beta|_\prec$ and for $\gamma=\beta$ one has
\be\lab{Pi-x}
\left\|\Pi_{x\gamma}^-(\varphi_x^\lambda))\right\|_{\L_p}\leq \lambda^{|\gamma|-2},~ x\in \R^4.
\ee

Suppose also that for all $\gamma\in\M'$ and $1\leq p<\infty$ one has
\be\lab{Gest'}
\left\|(\G_{xy}\right)^\gamma_\beta\|_{\L_p}\lesssim |x-y|^{|\beta|-|\gamma|},~x,y\in \R^4,~|x-y|\leq \bar{\lambda}_\prec, \bar{\lambda}_\prec>0,
\ee
where the constant in the inequality does not depend on $x,y$ and $\rho$.
Then for any $x\in \R^4$, $|x|\leq \bar{\lambda}_\prec$ 
\be\lab{Epsin'} 
\left\|\Pi_{0\beta}^-(\varphi_x^\lambda)\right\|_{\L_p}\lesssim \lambda^{\alpha-2}(\lambda+|x|)^{|\beta|-\alpha},
\ee
where the constant in the inequality does not depend on $\lambda$, $x$ and $\rho$.
\el

\bpr
By (\ref{PG}), (\ref{triangG}), (\ref{triangG1}) and (\ref{triangGm}) we can write
$$
\Pi_{0\beta}^-(\varphi_x^\lambda)=\sum_{\tiny\begin{array}{c} \gamma\in \M': \\ \alpha\leq |\gamma|<|\beta| \\ |\gamma|_\prec<|\beta|_\prec \end{array}}\Pi_{x\gamma}^-(\varphi_x^\lambda)\circ (\G_{x0})^\gamma_\beta +\Pi_{x\beta}^-(\varphi_x^\lambda).
$$

Thus by the H\"{o}lder inequality in probability we get with the help of (\ref{Pi-x}) and (\ref{Gest'})
$$
\left\|\Pi_{0\beta}^-(\varphi_x^\lambda)\right\|_{\L_p}\leq\sum_{\tiny\begin{array}{c} \gamma\in \M': \\ \alpha\leq |\gamma|<|\beta| \\ |\gamma|_\prec<|\beta|_\prec \end{array}}\left\|\Pi_{x\gamma}^-(\varphi_x^\lambda)\right\|_{\L_{2p}} \left\|(\G_{x0})^\gamma_\beta\right\|_{\L_{2p}} +\left\|\Pi_{x\beta}^-(\varphi_x^\lambda)\right\|_{\L_p}\lesssim
$$
$$
\lesssim \sum_{\tiny\begin{array}{c} \gamma\in \M': \\ \alpha\leq |\gamma|\leq|\beta| \\ |\gamma|_\prec\leq |\beta|_\prec \end{array}}\lambda^{|\gamma|-2}|x|^{|\beta|-|\gamma|}\lesssim
\lambda^{\alpha-2}(\lambda+|x|)^{|\beta|-\alpha},
$$
where at the last step we used (\ref{min}). This completes the proof.

\epr

\bc\lab{psicor}
Assume that for some $\beta\in \M'$ and for all $x\in \R^4$, $1\leq p<\infty$, $\gamma\in \M'$ with $|\gamma|_\prec<|\beta|_\prec$
$$ 
\left\|\Pi_{x\gamma}^-(\psi^{\e_n}_x)\right\|_{\L_p} \lesssim {\e_n}^{|\gamma|-2},
$$
and
$$
\left\|\Pi_{x\beta}^-(\psi^{\e_n}_x)\right\|_{\L_p} \lesssim {\e_n}^{|\beta|-2},
$$
where $\e_n=\bar{\lambda}' 2^{-n}=\bar{\lambda}_\prec 2^{-(n+2)}$, $\bar{\lambda}_\prec>0$, $n\in \N$, and the constants in the inequalities do not depend on $n$, $x$ and $\rho$.

Suppose also that for all $\gamma\in\M'$ and $1\leq p<\infty$ one has
$$
\left\|(\G_{xy}\right)^\gamma_\beta\|_{\L_p}\lesssim |x-y|^{|\beta|-|\gamma|},~x,y\in \R^4,~|x-y|\leq \bar{\lambda}_\prec,
$$
where the constant in the inequality does not depend on $x,y$ and $\rho$.

Then for any $x\in \R^4$, $|x|\leq \bar{\lambda}_\prec$, $n\in \N$
\be\lab{Epsin''} 
\left\|\Pi_{0\beta}^-(\psi^{\e_n}_x)\right\|_{\L_p}\lesssim \e_n^{\alpha-2}(\e_n+|x|)^{|\beta|-\alpha},
\ee
where the constant in the inequality does not depend on $x$, $n$ and $\rho$.
\ec

\bpr
The proof is obtained by applying the previous lemma to $\varphi=\psi$, $\lambda={\e_n}=\bar{\lambda}' 2^{-n}=\bar{\lambda}_\prec 2^{-(n+2)}$, $n\in \N$, $x\in \R^4$, $|x|\leq \bar{\lambda}_\prec$.

\epr

Now the desired induction step estimate (\ref{Ephi*}), where $\mu=\beta$, $\bar{\lambda}_\prec$ is replaced by $\bar{\lambda}$, in full generality in the case when $|\beta|<2$ follows from Lemma \ref{psilem0}, the previous corollary and from the next proposition.
\bp\lab{Pi-b<}
Let $\e_n=\bar{\lambda}' 2^{-n}=\bar{\lambda}_\prec 2^{-(n+2)}$, $n\in \N$, for some $\bar{\lambda}'=\frac{\bar{\lambda}_\prec}4>0$, $\psi=\omega^2\ast \omega$, and $\omega\in \D$  an even function supported in $B_{\frac 16}$ and satisfying (\ref{omegapr*}) with some $r\in \N$, $r>2-\alpha$.

Assume that some $\beta\in \M'$ and for all $1\leq p<\infty$, $x\in \R^4$, $|x|\leq \bar{\lambda}_\prec$, $n\in \N$
\be\lab{Epsin''*} 
\left\|\Pi_{0\beta}^-(\psi^{\e_n}_x)\right\|_{\L_p}\lesssim \e_n^{\alpha-2}(\e_n+|x|)^{|\beta|-\alpha},
\ee
where the constant in the inequality does not depend on $x$, $n$ and $\rho$.

Suppose also that for all $\varphi\in \mB^r$ and for all $1\leq p<\infty$, $\gamma\in \M'$ with $|\gamma|_\prec<|\beta|_\prec$ and $0<\lambda\leq \bar{\lambda}_\prec$ one has
\be\lab{Pi-x1}
\left\|\Pi_{x\gamma}^-(\varphi_x^\lambda))\right\|_{\L_p}\lesssim \lambda^{|\gamma|-2},
\ee
and that for all $\gamma\in\M'$ and $1\leq p<\infty$ one has
\be\lab{Gest'1}
\left\|(\G_{xy}\right)^\gamma_\beta\|_{\L_p}\lesssim |x-y|^{|\beta|-|\gamma|},~x,y\in \R^4,~|x-y|\leq \bar{\lambda}_\prec,
\ee
where the constants in the inequalities do not depend on $x,y$, $\lambda$, $\rho$ and $\varphi$.

Then for any $1\leq p<\infty$, $\varphi\in \mB^r$, $x\in \R^4$, $|x|\leq \bar{\lambda}_\prec$, $0<\lambda\leq 2\bar{\lambda}'=\frac{\bar{\lambda}_\prec}{2}$ one has
\be\lab{Ephi} 
\left\|\Pi_{0\beta}^-(\varphi^{\lambda}_x)\right\|_{\L_p}\lesssim \lambda^{\alpha-2}(\lambda+|x|)^{|\beta|-\alpha},
\ee
where the constant in the inequality does not depend on $\lambda$, $\varphi$, $x$ and $\rho$.
\ep

\bpr
As in the proof of Proposition \ref{Pi-psi} we introduce
\be\lab{cps}
\check{\psi}:=\psi^{\frac 12}-\psi^2.
\ee
By (\ref{psisupp}), by the definition of $\psi=\omega^2\ast \omega$ and by the second property in (\ref{omegapr*}) we infer
\be\lab{cpsupp}
{\rm supp}(\check{\psi})\subset B_1,~\int_{\R^4}x^\m\check{\psi}(x)dx=0,~0\leq|\m|\leq r-1,~\m \in \N^4.
\ee

We shall also need the function
$$
\sigma:=\psi^2\ast \psi.
$$

The first property in (\ref{omegapr*}) and the definition of $\psi=\omega^2\ast \omega$ imply
$$
\int_{\R^4}\sigma(x)dx=1,
$$
and hence $\sigma$ can be used as a mollifier, so that for $\e_k:=\bar{\lambda}'2^{-k}$, $k\in \mathbb{Z}$ and any $\varphi\in \D$ we have
\be\lab{moll}
\varphi=\lim_{k\to \infty}\varphi\ast \sigma^{\e_k}
\ee

Using (\ref{cps}) we can write for $n\in \N$
$$
\sigma^{\e_{n+1}}-\sigma^{\e_n}=\psi^{\e_n}\ast (\psi^{\frac 12}-\psi^2)^{\e_n}=\psi^{\e_n}\ast \check{\psi}^{\e_n},
$$
and hence
\be\lab{teleps}
\sigma^{\e_k}=\sigma^{\e_n}+\sum_{m=n}^{k-1}(\sigma^{\e_{m+1}}-\sigma^{\e_m})=\psi^{\e_{n-1}}\ast \psi^{\e_n}+\sum_{m=n}^{k-1}\psi^{\e_m}\ast \check{\psi}^{\e_m}.
\ee 

Thus by (\ref{moll}) and by the previous presentation for $\sigma^{\e_k}$ we can write for any $\varphi\in \mB^r$, $n\in \N$
\be\lab{P-appr}
\Pi_{0\beta}^-(\varphi^{\e_{n-1}}_0)=(\Pi_{0\beta}^-\ast \varphi^{\e_{n-1}})(0)=\lim_{k\to \infty}(\Pi_{0\beta}^-\ast \varphi^{\e_{n-1}}\ast \sigma^{\e_k})(0)=
\ee
$$
=(\Pi_{0\beta}^-\ast \varphi^{\e_{n-1}}\ast \psi^{\e_{n-1}}\ast \psi^{\e_n})(0)+\sum_{m=n}^{\infty}(\Pi_{0\beta}^-\ast \varphi^{\e_{n-1}}\ast \psi^{\e_m}\ast \check{\psi}^{\e_m})(0)=
$$
$$
=\int_{\R^4}\Pi_{0\beta}^-(\psi^{\e_n}_z) (\varphi^{\e_{n-1}}\ast \psi^{\e_{n-1}})(-z)dz+\sum_{m=n}^{\infty }\int_{\R^4}\Pi_{0\beta}^-(\psi^{\e_m}_z) (\varphi^{\e_{n-1}}\ast \check{\psi}^{\e_m})(-z)dz.
$$

By (\ref{psisupp}) for $n\in \N$ one has ${\rm supp}(\varphi^{\e_{n-1}}\ast \psi^{\e_{n-1}})\subset B_{\e_{n-2}}\subset B_{\bar{\lambda}_\prec}$ and by the first property in (\ref{cpsupp}) ${\rm supp}(\varphi^{\e_{n-1}}\ast \check{\psi}^{\e_m})\subset B_{\e_{n-2}}\subset B_{\bar{\lambda}_\prec}$ for $m\geq n$.
Using these observations, the Minkowski and the H\"{o}lder inequalities we have for $1\leq p<\infty$
\be\lab{Pi0bn}
\left\|\Pi_{0\beta}^-(\varphi^{\e_{n-1}}_0)\right\|_{\L_p}\leq \int_{B_{\e_{n-2}}}\left\|\Pi_{0\beta}^-(\psi^{\e_n}_z)\right\|_{\L_p} dz \left\|\varphi^{\e_{n-1}}\ast \psi^{\e_{n-1}}\right\|_{L_\infty}+
\ee
$$
+\sum_{m=n}^{\infty }\int_{B_{\e_{n-2}}}\left\|\Pi_{0\beta}^-(\psi^{\e_m}_z)\right\|_{\L_p} dz \left\|\varphi^{\e_{n-1}}\ast \check{\psi}^{\e_m}\right\|_{L_\infty}.
$$

Now, since $\left\|\varphi\right\|_{L_\infty}\leq \left\|\varphi\right\|_{C^r}\leq 1$, we have  
\be\lab{est1aux}
\left\|\varphi^{\e_{n-1}}\ast \psi^{\e_{n-1}}\right\|_{L_\infty}\lesssim 2^{nd}\left\|\varphi\right\|_{L_\infty} \left\|\psi\right\|_{L_1}\lesssim 2^{nd},
\ee
where the constant in the inequality only depends on $\bar{\lambda}'$ and $\omega$.

By Lemma 9.2 in \cite{CZ} we also obtain
\be\lab{est2aux}
\left\|\varphi^{\e_{n-1}}\ast \check{\psi}^{\e_m}\right\|_{L_\infty}\lesssim 2^{n(r+d)-mr}\left\|\varphi\right\|_{C^r}\left\|\check{\psi}\right\|_{L_1}\lesssim 2^{n(r+d)-mr},
\ee
where the constant in the inequality only depends on $\bar{\lambda}'$, $r$ and $\omega$.

Using the last two estimates and (\ref{Epsin''}) in (\ref{Pi0bn}) we obtain
\be\lab{Pib-psf}
\left\|\Pi_{0\beta}^-(\varphi^{\e_{n-1}}_0)\right\|_{\L_p}\leq 2^{nd}\int_{B_{\e_{n-2}}}\e_n^{\alpha-2}(\e_n+|z|)^{|\beta|-\alpha} dz +
\ee
$$
+\sum_{m=n}^{\infty }2^{n(r+d)-mr}\int_{B_{\e_{n-2}}}\e_m^{\alpha-2}(\e_m+|z|)^{|\beta|-\alpha} dz \lesssim
$$
$$
\lesssim 2^{nd-nd}\e_n^{\alpha-2+|\beta|-\alpha}+\e_n^{\alpha-2+|\beta|-\alpha}\sum_{m=n}^{\infty}2^{n(r+d)-mr-m(\alpha-2)+n(\alpha-2)-nd}\lesssim
$$
$$
\lesssim \e_{n-1}^{|\beta|-2}\sum_{m=n}^{\infty}2^{-(m-n)(r+\alpha-2)}\lesssim \e_{n-1}^{|\beta|-2},
$$ 
where the series in the last expression converges due to the condition $r+\alpha-2>0$ and 
the constant in the inequality does not depend on $n$, $\varphi$ and $\rho$.

Now for any $0<\lambda \leq 2\bar{\lambda}'$, let $n\in \N$ be the largest possible such that $\e_{n-1}\geq \lambda$, so that $\e_n\leq \lambda$, and for $\lambda':=\e_{n-1}^{-1}\lambda$ we have $\frac 12\leq \lambda'\leq 1$. 

Also if $\varphi\in \mB^r$ then $\varphi^{\lambda'}\in C\mB^r$, where $C>0$ is a constant independent of $\lambda'$ and $\varphi$.

Since by the definition $\varphi^\lambda=(\varphi^{\lambda'})^{\e_{n-1}}$ we obtain from the last observation and from (\ref{Pib-psf}) that
$$
\left\|\Pi_{0\beta}^-(\varphi_0^{\lambda})\right\|_{\L_p}\lesssim \e_{n-1}^{|\beta|-2}\leq \lambda^{|\beta|-2},
$$
where the constant in the inequality does not depend on $\lambda$, $\varphi$ and $\rho$.


Finally, the last estimate, (\ref{Pi-x1}), (\ref{Gest'1}) and Lemma \ref{Pi-int} imply (\ref{Ephi}). This completes the proof.

\epr

Next, in the case when $|\beta|>2$ we obtain the induction step estimate (\ref{Ephi*}), where $\mu=\beta$, and $\bar{\lambda}_\prec$ is replaced by $\bar{\lambda}$. This is done by applying the reconstruction theorem for germs.
\bp\lab{Pi-b>}
Assume that for some $\beta\in \M'$, $|\beta|>2$, $\bar{\lambda}_\prec>0$, $r\in \N$, $r>2-\alpha$ one has
\be\lab{Pi-x>}
\left\|\Pi_{x\gamma}^-\right\|_{p\bar{\lambda}_\prec r}^{|\gamma|-2}=\sup_{0<\lambda\leq \bar{\lambda}_\prec}\sup_{\varphi\in \mB^r}\frac{\left\|\Pi_{x\gamma}^-(\varphi_x^\lambda))\right\|_{\L_p}}{\lambda^{|\gamma|-2}}\lesssim 1
\ee
for all $x\in \R^4$, $1\leq p<\infty$, $\gamma\in \M'$ with $|\gamma|_\prec<|\beta|_\prec$, where the constant in the inequality does not depend on $\rho$ and $x$.


Suppose also that for all $\gamma\in\M'$ and $1\leq p<\infty$ one has
\be\lab{Gest'>}
\left\|(\G_{xy}\right)^\gamma_\beta\|_{\L_p}\lesssim |x-y|^{|\beta|-|\gamma|},~x,y\in \R^4,~|x-y|\leq \bar{\lambda}_\prec,
\ee
where the constant in the inequality does not depend on $x,y$ and $\rho$.
Then for any $1\leq p<\infty$, $x\in \R^4$, $|x|\leq \bar{\lambda}_\prec$, $0<\lambda\leq 2\bar{\lambda}'=\frac{\bar{\lambda}_\prec}{2}$, $\varphi\in \mB^r$ 
\be\lab{Epsin'>} 
\left\|\Pi_{0\beta}^-(\varphi_x^\lambda)\right\|_{\L_p}\lesssim \lambda^{\alpha-2}(\lambda+|x|)^{|\beta|-\alpha},
\ee
where the constant in the inequality does not depend on $\lambda$, $\varphi$, $x$ and $\rho$.
\ep

\bpr
Let $F_x:=\Pi_{x\beta}^-$. Then for any $\varphi$ as in (\ref{Pi-x>}) and any $R>0$ we have by (\ref{Pi-0}), (\ref{triangG}), (\ref{triangG1}), (\ref{triangGm}), and (\ref{PG}) 
$$
\left\|(F_{x+h}-F_x)(\varphi_{x+h}^\lambda)\right\|_{p\infty B_R x}=\|\hspace{-0.5cm} \sum_{\tiny\begin{array}{c} \gamma\in \M': \\ |\gamma|_\prec<|\beta|_\prec \\ \alpha\leq |\gamma|<|\beta|\end{array}}\hspace{-0.5cm} \Pi_{x+h\gamma}^-(\varphi_{x+h}^\lambda)\circ ({\rm Id}-\G_{x+h x})^\gamma_\beta \|_{p\infty B_R x}\leq
$$
$$
\leq \sum_{\tiny\begin{array}{c} \gamma\in \M': \\ |\gamma|_\prec<|\beta|_\prec \\ \alpha\leq |\gamma|<|\beta|\end{array}}\left\|\Pi_{x+h\gamma}^-(\varphi_{x+h}^\lambda)\right\|_{p'\infty B_R x} \left\|(\G_{x+h x})^\gamma_\beta\right\|_{p''\infty B_R x},
$$
where $\frac{1}{p'}+\frac{1}{p''}=\frac{1}{p}$, $1\leq p',p''<\infty$, and at the last step we used the H\"{o}lder inequality.

Now using (\ref{Pi-x>}) and (\ref{Gest'>})  we obtain for $|h|\leq \bar{\lambda}_\prec$, $0<\lambda\leq \bar{\lambda}_\prec$
$$
\left\|(F_{x+h}-F_x)(\varphi_{x+h}^\lambda)\right\|_{p\infty B_R x}\lesssim\sum_{\tiny\begin{array}{c} \gamma\in \M': \\ |\gamma|_\prec<|\beta|_\prec \\ \alpha\leq |\gamma|<|\beta|\end{array}}|h|^{|\beta|-|\gamma|}\lambda^{|\gamma|-2}\lesssim \lambda^{\alpha-2}(\lambda+|h|)^{|\beta|-\alpha},
$$
where at the last step we used (\ref{min}).

Thus
$$
\left\|F\right\|_{p\infty R \frac{\bar{\lambda}_\prec}2 \varphi}^{\alpha-2~|\beta|-2~0}=\sup_{h\in \R^4:|h|\leq \bar{\lambda}_\prec}\sup_{0<\lambda\leq \frac{\bar{\lambda}_\prec}2}\frac{\left\|(F_{x+h}-F_x)(\varphi_{x+h}^\lambda)\right\|_{p\infty B_R x}}{\lambda^{\alpha-2}(\lambda+|h|)^{|\beta|-\alpha}}\lesssim 1,
$$
where the constant in the inequality is independent of $R$. Hence by (\ref{P-van}) with $\n={\bf 0}$ and by Corollary \ref{reccor} we obtain for $r>2-\alpha$
$$
\left\|F\right\|_{p\infty 2\bar{\lambda}' R r}^{|\beta|-2~0}=\sup_{0<\lambda\leq 2\bar{\lambda}'}\frac{\left\|\sup_{\varphi\in \mB^r}\left\|F_x(\varphi^\lambda_x)\right\|_{\L_p}\right\|_{L_\infty(B_R,x)}}{\lambda^{|\beta|-2}}\lesssim 1,
$$
where the constant in the inequality is independent of $R$. Finally, the last estimate and Lemma \ref{Pi-int} imply (\ref{Epsin'>}). This completes the proof.

\epr

Suppose that $\beta\in \M_{\geq 0}$. Then the assumptions of Proposition \ref{indstepprp} imply those in Lemma \ref{psilem0} for $|\beta|<2$ (see Corollary \ref{dP-onpp}), and together with the results of Lemma \ref{psilem0} they imply the assumptions in Corollary \ref{psicor}. In turn, the results of Corollary \ref{psicor} together with the assumptions of Proposition \ref{indstepprp} imply those in Proposition \ref{Pi-b<} if $|\beta|<2$. The assumptions of Proposition \ref{indstepprp} are also stronger than those in Proposition \ref{Pi-b>} if $|\beta|>2$. Therefore the induction step estimate (\ref{Ephi*}), where $\mu=\beta$, $\bar{\lambda}_\prec$ is replaced by $\bar{\lambda}=\frac{\bar{\lambda}_\prec}{32}$, follows from Proposition \ref{Pi-b<} if $|\beta|<2$ and from Proposition \ref{Pi-b>} if $|\beta|>2$. More precisely, we have the following corollary which implies the induction step estimate for $\Pi_{0\beta}^-$ in Proposition \ref{indstepprp}.
\bc\lab{P-onpp}
Under the induction step assumptions in Proposition \ref{indstepprp} estimate (\ref{Ephi}) holds with the same conditions as in Proposition \ref{Pi-b<} if $|\beta|<2$, and estimate (\ref{Epsin'>})  holds with the same conditions as in Proposition \ref{Pi-b>} if $|\beta|>2$. In particular, the induction step statement for the estimate (\ref{Ephi*}), where $\mu=\beta$, $\bar{\lambda}_\prec$ is replaced by $\bar{\lambda}=\frac{\bar{\lambda}_\prec}{32}$, in Proposition \ref{indstepprp} holds.
\ec

Finally we can obtain a preliminary estimate for $\Pi_{x\beta}$ using Lemma \ref{psilem0}, Corollary \ref{psicor}, Propositions \ref{Pi-b<}, \ref{Pi-b>} and Schauder estimates for germs. The proper induction step estimate for $\Pi_{x\beta}$ will be obtained in the next section as it is based on the induction step estimates for $(\G_{xy})^{\delta_\n}_\beta$, $\n\in \N^4$ derived there.
\bp\lab{Pibestprel}
Assume that $\beta\in \M_{\geq 0}$, $\beta\neq 0$ and that the conditions of Lemma \ref{psilem0}, Corollary \ref{psicor} and of Proposition \ref{Pi-b<} are satisfied if $|\beta|<2$, and that the conditions of Proposition \ref{Pi-b>} are satisfied if $|\beta|>2$. Suppose also that $\bar{\lambda}'=\frac{\bar{\lambda}_\prec}{4}\geq 1$. Then for any $r>\max\{2-\alpha, |\beta|\}$, $1\leq p<\infty$ one has
\be\lab{Pibestp}
\left\|\Pi_{x\beta}\right\|_{p\bar{\lambda}'r}^{|\beta|}=\sup_{0<\lambda\leq \bar{\lambda}'}\frac{\sup_{\varphi\in \mB^r}\left\|\Pi_{x\beta}(\varphi^\lambda_x)\right\|_{\L_p}}{\lambda^{|\beta|}}=
\ee
$$
=\sup_{0<\lambda\leq \bar{\lambda}'}\frac{\sup_{\varphi\in \mB^r}\left\|\Pi_{0\beta}(\varphi^\lambda_0)\right\|_{\L_p}}{\lambda^{|\beta|}}\lesssim 1,
$$
where the constant in the inequality does not depend on $\rho$.


\ep

\bpr
If $|\beta|<0$, by translational invariance in law, the estimate in the statement follows from Lemma \ref{psilem0}, Corollary \ref{psicor}, Proposition \ref{Pi-b<} and from (\ref{emb1}) with $b=|\beta|-2$ and $F_x=\Pi_{x\beta}^-$ as in this case by (\ref{defP})
$\Pi_{x\beta}=\K\Pi_{x\beta}^-$.

If $|\beta|>0$, by translational invariance in law, the estimate in the statement follows from Lemma \ref{psilem0}, Corollary \ref{psicor}, Proposition \ref{Pi-b<} for $|\beta|<2$ or from Proposition \ref{Pi-b>} for $|\beta|>2$, and from (\ref{emb1c}) with $b=|\beta|-2$, $F_x=\Pi_{x\beta}^-$ along with the definition of $\Pi_{x\beta}$ in (\ref{defP}).

\epr

The arguments in the discussion before Corollary \ref{P-onpp}, along with the fact that the inequality $\bar{\lambda}'=\frac{\bar{\lambda}_\prec}{4}\geq 1$ is also ensured by the assumptions in Proposition \ref{indstepprp}, and Proposition \ref{Pibestprel} immediately lead to the following statement.
\bc\lab{Ponpprl}
Under the induction step assumptions in Proposition \ref{indstepprp} estimate (\ref{Pibestp}) holds with the same conditions as in Proposition \ref{Pibestprel}.
\ec


\subsection{The induction step for $(\G_{xy})^\gamma_\beta$, $({\rm d}\G_x)^\gamma_\beta$, $(S_{xy})^\gamma_\beta$, $\gamma\in \M_{pp}$}\lab{indstGGSpp}

\setcounter{equation}{0}
\setcounter{theorem}{0}


In this section we complete the proof of the induction step statement by establishing estimates for $(\G_{xy})^\gamma_\beta$, $({\rm d}\G_x)^\gamma_\beta$, $(S_{xy})^\gamma_\beta$, $\gamma\in \M_{pp}$. We start with 
the following auxiliary result which can be found in \cite{BCZ}, Section 5.5, proof of Lemma 3.15.
\bl\lab{phi0lem}
For any $N>0$, and any $\m\in \N^4$ with $|\m|\leq N$ there exists a function $\phi\in \D$ supported in $B_{\frac 14}$ such that for $\n\in \N^4$, $0\leq |\n|\leq N$
$$
\int_{\R^4}(-x)^\n\phi(x)dx=\left\{\begin{array}{ll} 1 & {\rm if}~\n=\m \\ 0 & {\rm otherwise}\end{array}\right. .
$$
\el

We shall also need the following lemma which will be used several times in the proofs in this section.
\bl\lab{philem}
Let $\phi\in \D(B_{\frac{1}{2k}})\cap \mB^r$ for some $r\in \N$, $k\geq 1$. Then for any $z\in \R^4$ with $|z|\leq \frac{1}{2k}$ the function $\tilde{\phi}(y):=\phi(y+z)$ is an element of $\D(B_{\frac 1k})\cap\mB^r$, and $\left\|\tilde{\phi}\right\|_{C^{r}}=\left\|\phi\right\|_{C^{r}}$ is independent of $z$.
\el

\bpr
It follows from the triangle inequality in the form $|y|\leq |z|+|y+z|$ that for $z\in \R^4$ with $|z|\leq \frac{1}{2k}$ and $|y+z|\leq \frac{1}{2k}$, $k\geq 1$ one has $|y|\leq \frac 1k$. Thus the function $\tilde{\phi}(y)$ is supported in $B_{\frac 1k}$ as $\phi\in \D$ is supported in $B_{\frac{1}{2k}}$.  

The remaining part of the statement is obvious by the translational invariance of the $C^r$--norm.

\epr

In Propositions \ref{Gestn}, \ref{dGpp}, \ref{Sxypp} below we shall obtain the induction step estimates for $(\G_{xy})^\gamma_\beta$, $({\rm d}\G_x)^\gamma_\beta$, $(S_{xy})^\gamma_\beta$, $\gamma\in \M_{pp}$, respectively. These proofs are similar and use a trick called the three--point argument in \cite{BOT}. 
\bp\lab{Gestn}
Assume that for some $\beta\in \M_{\geq 0}$, $\beta\neq 0$, estimates (\ref{Gest}) hold for $\mu=\beta$, for all $1\leq p<\infty$ and all $\gamma\in\M\setminus\M_{pp}=\M'$, $x,y\in \R^4$, $|x-y|\leq \bar{\lambda}_\prec$. Suppose also that $\bar{\lambda}' =\frac{\bar{\lambda}_\prec}4\geq 1$ and that for all $1\leq p<\infty$, $\gamma\in\M_{\geq 0}$ with $|\gamma|_\prec<|\beta|_\prec$, and for some $r\in \N$, $r>{\rm max}\{2-\alpha, |\beta|\}$  (\ref{Pigest+}) holds and (\ref{Pibestp}) also holds.

Then for any $1\leq p<\infty$, $x,y\in \R^4$, $0<|x-y|\leq \bar{\lambda}''=\frac{\bar{\lambda}'}2$, $\n\in \N^4$ one has
\be\lab{Gxyestn}
\left\|(\G_{xy})^{\delta_\n}_\beta\right\|_{\L_p}\lesssim |x-y|^{|\beta|-|\delta_\n|},
\ee
where the constant in the inequality does not depend on $x,y$ and $\rho$.


\ep

\bpr
Note that by translational invariance in law it suffices to consider the case $y=0$.

By (\ref{PG}), (\ref{Pi0}), (\ref{triangG}), (\ref{triangG1}) and (\ref{triangGm})  we have
$$
\Pi_{0\beta}(y)=\sum_{\gamma\in \Ml}\Pi_{x\gamma}(y)\circ (\G_{x0})_\beta^\gamma=\sum_{\tiny\begin{array}{c}\gamma\in \M_{\geq 0}\cup\M_{pp}: \\ |\gamma|_\prec<|\beta|_\prec \\ |\gamma|<|\beta| \end{array}}\Pi_{x\gamma}(y)\circ (\G_{x0})_\beta^\gamma+\Pi_{x\beta}(y) =
$$
$$
=\Pi_{x\beta}(y)+\sum_{\tiny\begin{array}{c}\gamma\in \M_{\geq 0}: \\ |\gamma|_\prec<|\beta|_\prec \\ |\gamma|<|\beta| \end{array}}\Pi_{x\gamma}(y)\circ (\G_{x0})_\beta^\gamma+\sum_{\tiny\begin{array}{c}\n\in \N^4: \\ |\n|<|\beta|_\prec \end{array}}\Pi_{x\delta_\n}(y)\circ (\G_{x0})_\beta^{\delta_\n}.
$$
Using the definition $\Pi_{x\delta_\n}(y)=(y-x)^\n{\rm Id}^V_{W_{\delta_\n}}$ we can rewrite the previous identity as follows
$$
\sum_{\tiny\begin{array}{c}\n\in \N^4: \\ |\n|<|\beta|_\prec \end{array}}(y-x)^\n{\rm Id}^V_{W_{\delta_\n}} (\G_{x0})_\beta^{\delta_\n}=\sum_{\tiny\begin{array}{c}\n\in \N^4: \\ |\n|<|\beta|_\prec \end{array}}\Pi_{x\delta_\n}(y)\circ (\G_{x0})_\beta^{\delta_\n}=
$$
$$
=\Pi_{0\beta}(y)-\sum_{\tiny\begin{array}{c}\gamma\in \M_{\geq 0}: \\ |\gamma|_\prec<|\beta|_\prec \\ |\gamma|<|\beta| \end{array}}\Pi_{x\gamma}(y)\circ (\G_{x0})_\beta^\gamma-\Pi_{x\beta}(y).
$$

Now we multiply this identity by the function $\phi^\lambda_x(y)$ and integrate with respect to  $y$ over $\R^4$, where $\phi$ is given in Lemma \ref{phi0lem} with $N=|\beta|_\prec$ and arbitrary fixed $\m\in \N^4$, $|\m|<N$. This yields
$$
\lambda^{|\m|}{\rm Id}^V_{W_{\delta_\m}} (\G_{x0})_\beta^{\delta_\m}=\Pi_{0\beta}(\phi^\lambda_x)-\Pi_{x\beta}(\phi^\lambda_x)-\sum_{\tiny\begin{array}{c}\gamma\in \M_{\geq 0}: \\ |\gamma|_\prec<|\beta|_\prec \\|\gamma|<|\beta| \end{array}}\Pi_{x\gamma}(\phi^\lambda_x)\circ (\G_{x0})_\beta^\gamma.
$$

Thus using the H\"{o}lder inequality we obtain
\be\lab{Gm}
\lambda^{|\m|}\left\|(\G_{x0})^{\delta_\m}_\beta\right\|_{\L_p}\leq \left\|\Pi_{0\beta}(\phi^\lambda_x)\right\|_{\L_p}+\left\|\Pi_{x\beta}(\phi^\lambda_x)\right\|_{\L_p}+
\ee
$$
+\sum_{\tiny\begin{array}{c}\gamma\in \M_{\geq 0}: \\ |\gamma|_\prec<|\beta|_\prec \\|\gamma|<|\beta| \end{array}}\left\|\Pi_{x\gamma}(\phi^\lambda_x)\right\|_{\L_{p'}} \left\|(\G_{x0})_\beta^\gamma\right\|_{\L_{p''}},
$$
where $\frac{1}{p'}+\frac{1}{p''}=\frac 1p$, $1<p',p''<\infty$.

Now assume that $0<\lambda \leq \bar{\lambda}'$.
 
By Lemma \ref{philem}, if $\frac{|x|}\lambda\leq \frac 12$ then for any $r\in \N$ $\phi^\lambda_x=\tilde{\phi}^\lambda_0$, where $\tilde{\phi}(y):=\phi(y+R_{\frac 1\lambda}x)\in C_r\mB^r$, $C_r\geq 0$, and $\left\|\tilde{\phi}\right\|_{C^r}$ is independent of $x$ and $\lambda$. 
Therefore for $|x| = \frac{\lambda}2$ we can apply (\ref{Pibestp}) in the first and in the second terms. We can also use  (\ref{Gest}), (\ref{Pigest+}) to estimate the factors in the terms of the sum over $\gamma$. This way we obtain with $\lambda =2|x|$
$$
\lambda^{|\m|}\left\|(\G_{x0})^{\delta_\m}_\beta\right\|_{\L_p}\lesssim \lambda^{|\beta|}+\lambda^{|\beta|}+\sum_{\tiny\begin{array}{c}\gamma\in \M_{\geq 0}: \\ |\gamma|_\prec<|\beta|_\prec \\ |\gamma|<|\beta| \end{array}}|x|^{|\beta|-|\gamma|}\lambda^{|\gamma|}\lesssim |x|^{|\beta|},
$$
or
$$
\left\|(\G_{x0})^{\delta_\m}_\beta\right\|_{\L_p}\lesssim |x|^{|\beta|-|\m|},
$$
which is equivalent to (\ref{Gxyestn}).

\epr

Under the induction step assumptions in Proposition \ref{indstepprp}, which ensure, in particular, that (\ref{Pigest+}) holds, from Corollaries \ref{Gnonpp}, \ref{Ponpprl} and from Proposition \ref{Gestn} we obtain the following corollary.
\bc\lab{Gestcor}
Under the induction step assumptions in Proposition \ref{indstepprp} estimate (\ref{Gxyestn}) holds with the same conditions as in Proposition \ref{Gestn}, and estimate (\ref{Gest*}) holds for all $\gamma\in \M_{pp}$ with $\mu$ replaced by $\beta$, $\bar{\lambda}_\prec$ replaced by $\bar{\lambda}=\frac{\bar{\lambda}_\prec}{32}$, and with the same other conditions as in the assumptions in Proposition \ref{indstepprp}, where $\mu$ is replaced by $\beta$, and $\bar{\lambda}_\prec$ is replaced by $\bar{\lambda}$.
\ec

At this point, based on the result of the previous proposition, we can finally justify the induction step estimate for $\Pi_{0\beta}$.
\bp\lab{Piint}
Assume that for some $\beta\in \M_{\geq 0}$, $\beta\neq 0$ estimates (\ref{Gest}) hold for $\mu=\beta$, for all $1\leq p<\infty$ and all $\gamma\in\M$, $x,y\in \R^4$, $|x-y|\leq \bar{\lambda}''=\frac{\bar{\lambda}'}2$. Suppose also that $\bar{\lambda}' =\frac{\bar{\lambda}_\prec}4\geq 1$ and that for all $1\leq p<\infty$, $\gamma\in\M$ with $|\gamma|_\prec<|\beta|_\prec$, and for some $r\in \N$, $r>{\rm max}\{2-\alpha, |\beta|\}$  (\ref{Pigest+}) holds and (\ref{Pibestp}) also holds.

Then for any $1\leq p<\infty$, $x\in \R^4$, $|x|\leq \bar{\lambda}''$, $0<\lambda\leq \bar{\lambda}'=\frac{\bar{\lambda}_\prec}{4}$, $\varphi\in \mB^r$ 
\be\lab{EPPsin'>} 
\left\|\Pi_{0\beta}(\varphi_x^\lambda)\right\|_{\L_p}\lesssim \lambda^{\alpha}(\lambda+|x|)^{|\beta|-\alpha},
\ee
where the constant in the inequality does not depend on $\lambda$, $\varphi$, $x$ and $\rho$.

\ep

\bpr
By (\ref{PG}), (\ref{triangG}), (\ref{triangG1}) and (\ref{triangGm}) we can write for $\varphi\in \mB^r$
$$
\Pi_{0\beta}(\varphi_x^\lambda)=\sum_{\tiny\begin{array}{c} \gamma\in \M: \\ \alpha\leq |\gamma|<|\beta| \\ |\gamma|_\prec<|\beta|_\prec \end{array}}\Pi_{x\gamma}(\varphi_x^\lambda)\circ (\G_{x0})^\gamma_\beta +\Pi_{x\beta}(\varphi_x^\lambda).
$$

Thus by the H\"{o}lder inequality in probability we get for $x\in \R^4$, $|x|\leq \bar{\lambda}''$, $0<\lambda\leq \bar{\lambda}'$ with the help of (\ref{Pigest+}), (\ref{Pibestp}) and (\ref{Gest})
$$
\left\|\Pi_{0\beta}(\varphi_x^\lambda)\right\|_{\L_p}\leq\sum_{\tiny\begin{array}{c} \gamma\in \M: \\ \alpha\leq |\gamma|<|\beta| \\ |\gamma|_\prec<|\beta|_\prec \end{array}}\left\|\Pi_{x\gamma}(\varphi_x^\lambda)\right\|_{\L_{2p}} \left\|(\G_{x0})^\gamma_\beta\right\|_{\L_{2p}} +\left\|\Pi_{x\beta}(\varphi_x^\lambda)\right\|_{\L_p}\lesssim
$$
$$
\lesssim \sum_{\tiny\begin{array}{c} \gamma\in \M: \\ \alpha\leq |\gamma|\leq|\beta| \\ |\gamma|_\prec\leq |\beta|_\prec \end{array}}\lambda^{|\gamma|}|x|^{|\beta|-|\gamma|}\lesssim
\lambda^{\alpha}(\lambda+|x|)^{|\beta|-\alpha},
$$
where at the last step we used (\ref{min}). This completes the proof.

\epr

From Corollary \ref{Gestcor} and Proposition \ref{Piint} we immediately derive the following statement.
\bc\lab{Pestcor}
Under the induction step assumptions in Proposition \ref{indstepprp} estimate (\ref{EPPsin'>}) holds with the same conditions as in Proposition \ref{Piint}, and estimate (\ref{EPPsin'>*}) holds with $\mu$ replaced by $\beta$, $\bar{\lambda}_\prec$ replaced by $\bar{\lambda}=\frac{\bar{\lambda}_\prec}{32}$, and with the same other conditions as in the assumptions in Proposition \ref{indstepprp}, where $\mu$ is replaced by $\beta$, and $\bar{\lambda}_\prec$ is replaced by $\bar{\lambda}$.
\ec

As promised, to complete the proof of the induction step we proceed with the induction step estimates for $({\rm d}\G_x)^\gamma_\beta$, $(S_{xy})^\gamma_\beta$, $\gamma\in \M_{pp}$. 
\bp\lab{dGpp}
Assume that for some $\beta\in \M_{\geq 0}$, $\beta\neq 0$, $|\beta|<2$ estimate (\ref{G0est}) (resp. (\ref{G0estS})) holds for $\mu=\beta$, for all $\gamma\in\M\setminus\M_{pp}=\M'$, $0<R\leq \bar{\lambda}_\prec$ (resp. $0<R\leq \frac{\bar{\lambda}_\prec}2$, $y\in \R^4$, $|y|\leq \frac{\bar{\lambda}_\prec}2$), $1\leq q<p^*<2$ and ${\rm d}\G_x$ (resp. $S_{x+y x}$) defined with the help of $\upsilon \in \L_{p^*}(L_2)$ such that $\left\|\left\|\upsilon\right\|_{L_2}\right\|_{\L_{p^*}}\leq 1$. Suppose also that for all $1\leq p<\infty$, $\gamma\in\M'$ with $|\gamma|_\prec<|\beta|_\prec$, $x\in \R^4$, and for some $r\in \N$, $r>4-\alpha$ (\ref{Pigest-}) and (\ref{Pigest+}) hold, and that $\bar{\lambda}''=\frac{\bar{\lambda}'}2 =\frac{\bar{\lambda}_\prec}8\geq 1$.

Then for any $0<R\leq \frac{\bar{\lambda}''}4$, $\n\in \N^4$, $1\leq q<p^*<2$, ${\rm d}\G_{x}$ defined with the help of $\upsilon\in \L_{p^*}(L_2)$ such that $\left\|\left\|\upsilon\right\|_{L_2}\right\|_{\L_{p^*}}\leq 1$ one has
\be\lab{G0estn}
\left\|({\rm d}\G_{x})^{\delta_\n}_\beta\right\|_{q2B_Rx}\lesssim R^{\frac d2+|\beta|-|\delta_\n|},
\ee
where the constant in the inequality is independent of $R$, $\upsilon$, $\rho$.
\ep

\bpr
By (\ref{Pi0}), (\ref{tr2G*})  we have
$$
\sum_{\gamma\in \Ml}\Pi_{x\gamma}(y)\circ ({\rm d}\G_x)_\beta^\gamma=\sum_{\tiny\begin{array}{c}\gamma\in \M_{\geq 0}\cup\M_{pp}: \\ |\gamma|_\prec<|\beta|_\prec \end{array}}\Pi_{x\gamma}(y)\circ ({\rm d}\G_x)_\beta^\gamma=
$$
$$
=\sum_{\tiny\begin{array}{c}\gamma\in \M_{\geq 0}: \\ |\gamma|_\prec<|\beta|_\prec  \end{array}}\Pi_{x\gamma}(y)\circ ({\rm d}\G_x)_\beta^\gamma+\sum_{\tiny\begin{array}{c}\n\in \N^4: \\ |\n|<|\beta|_\prec \end{array}}\Pi_{x\delta_\n}(y)\circ ({\rm d}\G_x)_\beta^{\delta_\n}.
$$
Using the germ
$$
F_x(y):=\delta\Pi_{0\beta}(y)-\sum_{\gamma\in \Ml}\Pi_{x\gamma}(y)\circ ({\rm d}\G_x)_\beta^\gamma
$$
and the definition $\Pi_{x\delta_\n}(y)=(y-x)^\n{\rm Id}^V_{W_{\delta_\n}}$ we can rewrite the previous identity as follows
$$
\sum_{\tiny\begin{array}{c}\n\in \N^4: \\ |\n|<|\beta|_\prec \end{array}}(y-x)^\n{\rm Id}^V_{W_{\delta_\n}} ({\rm d}\G_x)_\beta^{\delta_\n}=\sum_{\tiny\begin{array}{c}\n\in \N^4: \\ |\n|<|\beta|_\prec \end{array}}\Pi_{x\delta_\n}(y)\circ ({\rm d}\G_x)_\beta^{\delta_\n}=
$$
$$
=\delta\Pi_{0\beta}(y)-F_x(y)-\sum_{\tiny\begin{array}{c}\gamma\in \M_{\geq 0}: \\ |\gamma|_\prec<|\beta|_\prec  \end{array}}\Pi_{x\gamma}(y)\circ ({\rm d}\G_x)_\beta^\gamma.
$$

Now we multiply this identity by the function $\phi^\lambda_x(y)$ and integrate with respect to  $y$ over $\R^4$, where $\phi$ is given in Lemma \ref{phi0lem} with $N=|\beta|_\prec$ and arbitrary fixed $\m\in \N^4$, $|\m|<N$. This yields
$$
\lambda^{|\m|}{\rm Id}^V_{W_{\delta_\m}} ({\rm d}\G_x)_\beta^{\delta_\m}=\delta\Pi_{0\beta}(\phi^\lambda_x)-F_x(\phi^\lambda_x)-\sum_{\tiny\begin{array}{c}\gamma\in \M_{\geq 0}: \\ |\gamma|_\prec<|\beta|_\prec  \end{array}}\Pi_{x\gamma}(\phi^\lambda_x)\circ ({\rm d}\G_x)_\beta^\gamma.
$$

Thus
\be\lab{dGm}
\lambda^{|\m|}\left\|({\rm d}\G_{x})^{\delta_\m}_\beta\right\|_{q2B_Rx}\leq \left\|\delta\Pi_{0\beta}(\phi^\lambda_x)\right\|_{q2B_Rx}+\left\|F_x(\phi^\lambda_x)\right\|_{q2B_Rx}+
\ee
$$
+\sum_{\tiny\begin{array}{c}\gamma\in \M_{\geq 0}: \\ |\gamma|_\prec<|\beta|_\prec  \end{array}}\left\|\Pi_{x\gamma}(\phi^\lambda_x)\right\|_{q'\infty B_Rx}\left\|({\rm d}\G_x)_\beta^\gamma\right\|_{q''2B_Rx},
$$
where $\frac{1}{q'}+\frac{1}{q''}=\frac 1q$, $1\leq q'<\infty$, $1\leq q< q''<p^*<2$, and $\upsilon$ satisfies $\left\|\left\|\upsilon\right\|_{L_2}\right\|_{\L_{p^*}}\leq 1$.

Now assume that $0<\lambda \leq \frac{\bar{\lambda}''}2$, $R=\frac{\lambda}2$.
 
By Lemma \ref{philem}, if $\frac{|x|}\lambda\leq \frac 12$ then for any $r\in \N$ $\phi^\lambda_x=\tilde{\phi}^\lambda_0$, where $\tilde{\phi}(y):=\phi(y+R_{\frac 1\lambda}x)\in C_r\mB^r$, $C_r\geq 0$, and $\left\|\tilde{\phi}\right\|_{C^r}$ is independent of $x$ and $\lambda$. 
Therefore for $|x|\leq R$ we can apply Proposition \ref{dPi0best}, conditions of which are guaranteed by the conditions of this proposition, in the first term. We can also use (\ref{PGest=}), applicability of which is ensured by the conditions of this proposition as well, to estimate the second term and (\ref{G0est}), (\ref{Pigest+}) to estimate the factors in the terms of the sum over $\gamma$. This way we obtain
$$
\lambda^{|\m|}\left\|({\rm d}\G_{x})^{\delta_\m}_\beta\right\|_{q2B_Rx}\lesssim |\lambda|^{|\beta|}R^{\frac d2}+\lambda^{\alpha+\frac d2}(\lambda+R)^{|\beta|-\alpha}+\sum_{\tiny\begin{array}{c}\gamma\in \M_{\geq 0}: \\ |\gamma|_\prec<|\beta|_\prec  \end{array}}R^{\frac d2 +|\beta|-|\gamma|}\lambda^{|\gamma|}.
$$
Now (\ref{G0estn}) is obtained from this estimate by substituting $\lambda =2R$.

\epr

\vskip 0.5cm

\noindent
{\em Proof of Proposition \ref{dGpp0}.}
The proof of Proposition \ref{dGpp0} is similar to that of Proposition \ref{dGpp} and even simpler as, if $\beta=0$, from the definition of the modified homogeneity it follows that the sum over $\gamma$ in the right hand side of (\ref{dGm}) is empty. To estimate the first and the second terms in the right hand side of (\ref{dGm}) when $\beta=0$ we use Proposition \ref{dPi0best0}  for the first term and (\ref{PGest=0}) for the second one. The remaining arguments are the same as in the proof of Proposition \ref{dGpp}.

\qed

\bp\lab{Sxypp}
Assume that for some $\beta\in \M_{\geq 0}$, $\beta\neq 0$, $|\beta|<2$ estimate (\ref{Gest}) holds for all $1\leq p<\infty$ and all $\gamma\in\M_{pp}$, $\mu\in \M$ with $|\mu|_\prec<|\beta|_\prec$,  and estimates (\ref{Gest}), (\ref{G0est}) (resp. (\ref{G0estS})) hold for $\mu=\beta$, for all $1\leq p<\infty$ and all $\gamma\in\M\setminus\M_{pp}=\M'$, $0<R\leq \bar{\lambda}_\prec$ (resp. $0<R\leq \frac{\bar{\lambda}_\prec}2$, $y\in \R^4$, $|y|\leq \frac{\bar{\lambda}_\prec}2$), $1\leq q<p^*<2$ and ${\rm d}\G_x$ (resp. $S_{x+y x}$) defined with the help of $\upsilon \in \L_{p^*}(L_2)$ such that $\left\|\left\|\upsilon\right\|_{L_2}\right\|_{\L_{p^*}}\leq 1$. Suppose also that for all $1\leq p<\infty$, $\gamma\in\M'$ with $|\gamma|_\prec<|\beta|_\prec$, $x\in \R^4$, and for some $r\in \N$, $r>4-\alpha$ (\ref{Pigest-}) and (\ref{Pigest+}) hold, and that $\bar{\lambda}''=\frac{\bar{\lambda}'}2 =\frac{\bar{\lambda}_\prec}8\geq 1$.

Then for any $R\in \R$, $y\in \R^4$, $0<R, |y|\leq \frac{\bar{\lambda}''}4$, $\n\in \N^4$,  $1\leq q<p^*<2$, $S_{x+y x}$ defined with the help of $\upsilon\in \L_{p^*}(L_2)$ such that $\left\|\left\|\upsilon\right\|_{L_2}\right\|_{\L_{p^*}}\leq 1$ one has
\be\lab{S0estn}
\left\|(S_{x+y x})^{\delta_\n}_\beta\right\|_{q2B_Rx}\lesssim
\ee
$$
\lesssim \left\{\begin{array}{ll} |y|^{\frac d2-|\n|+\alpha}(|y|+R)^{|\beta|-\alpha} & {\rm  if}~\frac d2-|\n|+\alpha>0 \\ (|y|+R)^{|\beta|-|\n|+\frac d2} &{\rm else}\end{array}\right. ,
$$
where the constant in the inequality is independent of $R$, $\upsilon$, $\rho$ and $y$.
\ep

\bpr
By (\ref{Pi0}), (\ref{tr2S*})  we have
$$
\sum_{\gamma\in \Ml}\Pi_{x\gamma}(z)\circ (S_{x+y x})_\beta^\gamma=\sum_{\tiny\begin{array}{c}\gamma\in \M_{\geq 0}\cup\M_{pp}: \\ |\gamma|_\prec<|\beta|_\prec \end{array}}\Pi_{x+y\gamma}(z)\circ (S_{x+y x})_\beta^\gamma=
$$
$$
=\sum_{\tiny\begin{array}{c}\gamma\in \M_{\geq 0}: \\ |\gamma|_\prec<|\beta|_\prec  \end{array}}\Pi_{x+y\gamma}(z)\circ (S_{x+y x})_\beta^\gamma+\sum_{\tiny\begin{array}{c}\n\in \N^4: \\ |\n|<|\beta|_\prec \end{array}}\Pi_{x+y\delta_\n}(z)\circ (S_{x+y x})_\beta^{\delta_\n}.
$$
Using the germ
$$
F_x(z):=\delta\Pi_{0\beta}(z)-\sum_{\gamma\in \Ml}\Pi_{x\gamma}(z)\circ ({\rm d}\G_x)_\beta^\gamma,
$$
property (\ref{PG}) and the definitions of $\Pi_{x\delta_\n}(z)=(z-x)^\n{\rm Id}^V_{W_{\delta_\n}}$ and $S_{x+y x}={\rm d}\G_{x+y}-\G_{x+y x}{\rm d}\G_x$ we can rewrite the previous identity as follows
$$
\sum_{\tiny\begin{array}{c}\n\in \N^4: \\ |\n|<|\beta|_\prec \end{array}}(z-x-y)^\n{\rm Id}^V_{W_{\delta_\n}} (S_{x+y x})_\beta^{\delta_\n}=\sum_{\tiny\begin{array}{c}\n\in \N^4: \\ |\n|<|\beta|_\prec \end{array}}\Pi_{x+y\delta_\n}(z)\circ (S_{x+y x})_\beta^{\delta_\n}=
$$
$$
=F_x(z)-F_{x+y}(z)-\sum_{\tiny\begin{array}{c}\gamma\in \M_{\geq 0}: \\ |\gamma|_\prec<|\beta|_\prec  \end{array}}\Pi_{x+y\gamma}(z)\circ (S_{x+y x})_\beta^\gamma.
$$

Now we multiply this identity by the function $\phi^\lambda_{x+y}(z)$ and integrate with respect to  $z$ over $\R^4$, where $\phi$ is given in Lemma \ref{phi0lem} with $N=|\beta|_\prec$ and arbitrary fixed $\m\in \N^4$, $|\m|<N$. This yields
$$
\lambda^{|\m|}{\rm Id}^V_{W_{\delta_\m}} (S_{x+y x})_\beta^{\delta_\m}=F_x(\phi^\lambda_{x+y})-F_{x+y}(\phi^\lambda_{x+y})-\sum_{\tiny\begin{array}{c}\gamma\in \M_{\geq 0}: \\ |\gamma|_\prec<|\beta|_\prec  \end{array}}\Pi_{x+y\gamma}(\phi^\lambda_{x+y})\circ (S_{x+y x})_\beta^\gamma.
$$

Thus
\be\lab{Sm}
\lambda^{|\m|}\left\|(S_{x+y x})^{\delta_\m}_\beta\right\|_{q2B_Rx}\leq \left\|F_x(\phi^\lambda_{x+y})\right\|_{q2B_Rx}+\left\|F_{x+y}(\phi^\lambda_{x+y})\right\|_{q2B_Rx}
\ee
$$
+\sum_{\tiny\begin{array}{c}\gamma\in \M_{\geq 0}: \\ |\gamma|_\prec<|\beta|_\prec  \end{array}}\left\|\Pi_{x+y\gamma}(\phi^\lambda_{x+y})\right\|_{q'\infty B_Rx}\left\|(S_{x+y x})_\beta^\gamma\right\|_{q''2B_Rx}\lesssim
$$
$$
\lesssim \left\|F_x(\tilde{\phi}^\lambda_x)\right\|_{q2B_Rx}+\left\|F_{x}(\phi^\lambda_{x})\right\|_{q2B_{R+|y|}x}+
$$
$$
+\sum_{\tiny\begin{array}{c}\gamma\in \M_{\geq 0}: \\ |\gamma|_\prec<|\beta|_\prec  \end{array}}\left\|\Pi_{x+y\gamma}(\phi^\lambda_{x+y})\right\|_{q'\infty B_Rx}\left\|(S_{x+y x})_\beta^\gamma\right\|_{q''2B_Rx},
$$
where $\tilde{\phi}(z):=\phi(z+R_{\frac 1\lambda}y)$, $\frac{1}{q'}+\frac{1}{q''}=\frac 1q$ and $1\leq q'<\infty$, $1\leq q< q''<p^*<2$, and $\upsilon$ satisfies $\left\|\left\|\upsilon\right\|_{L_2}\right\|_{\L_{p^*}}\leq 1$.

Now assume that $0<\lambda \leq \frac{\bar{\lambda}''}2$, $0<R\leq \frac{\bar{\lambda}''}4$.
By Lemma \ref{philem}, if $\frac{|y|}\lambda\leq \frac 12$ then for any $r\in \N$ $\phi^\lambda_{x+y}=\tilde{\phi}^\lambda_x$, and $\tilde{\phi}(z)=\phi(z+R_{\frac 1\lambda}y)\in C_r\mB^r$, $C_r\geq 0$, and $\left\|\tilde{\phi}\right\|_{C^r}$ is independent of $y$ and $\lambda$. 
Therefore for $|y|\leq \frac{\lambda}2\leq \frac{\bar{\lambda}''}4$ we can use (\ref{PGest=}), applicability of which is ensured by the conditions of this proposition, to estimate the first and the second terms and (\ref{G0estS}), (\ref{Pigest+}) to estimate the factors in the terms of the sum over $\gamma$. This way we obtain
\be\lab{S1est}
\lambda^{|\m|}\left\|(S_{x+y x})^{\delta_\m}_\beta\right\|_{q2B_Rx}\lesssim \lambda^{\alpha+\frac d2}(\lambda+R)^{|\beta|-\alpha}+\lambda^{\alpha+\frac d2}(\lambda+R+|y|)^{|\beta|-\alpha}+
\ee
$$
+\sum_{\tiny\begin{array}{c}\gamma\in \M_{\geq 0}: \\ |\gamma|_\prec<|\beta|_\prec \\ \frac d2-|\gamma|_p+\alpha>0 \\ \sum_{\n\in\N^4}\gamma(\n)\neq 0 \end{array}}|y|^{\frac d2-|\gamma|_p+\alpha}(|y|+R)^{|\beta|-|\gamma|+|\gamma|_p-\alpha}\lambda^{|\gamma|}+
$$
$$
+\sum_{\tiny\begin{array}{c}\gamma\in \M_{\geq 0}: \\ |\gamma|_\prec<|\beta|_\prec \\ \frac d2-|\gamma|_p+\alpha\leq0 \\ \sum_{\n\in\N^4}\gamma(\n)\neq 0\end{array}}(|y|+R)^{|\beta|-|\gamma|+\frac d2}\lambda^{|\gamma|}.
$$

Now observe that for $\gamma\in \M_{\geq 0}$ with $\sum_{\n\in\N^4}\gamma(\n)\neq 0$ 
$$
|\gamma|=(\alpha+1)\gamma(\bg)+\sum_{\n\in \N^4}\gamma(\n)(|\n|-\alpha)+\alpha> \sum_{\n\in \N^4}\gamma(\n)|\n|=|\gamma|_p
$$
as $(\alpha+1)\gamma(\bg)+\sum_{\n\in \N^4}\gamma(\n)(-\alpha)+\alpha>0$ in this case.
Using this observation we obtain from (\ref{S1est}) for $\lambda =2|y|$
\be\lab{S1est1}
\left\|(S_{x+y x})^{\delta_\m}_\beta\right\|_{q2B_Rx}\lesssim |y|^{\alpha+\frac d2-|\m|}(|y|+R)^{|\beta|-\alpha}+
\ee
$$
+|y|^{\alpha+\frac d2-|\m|}(|y|+R)^{|\beta|-\alpha}\sum_{\tiny\begin{array}{c}\gamma\in \M_{\geq 0}: \\ |\gamma|_\prec<|\beta|_\prec \\ \frac d2-|\gamma|_p+\alpha>0 \\ \sum_{\n\in\N^4}\gamma(\n)\neq 0 \end{array}}\large\left(\frac{|y|}{|y|+R}\large\right)^{|\gamma|-|\gamma|_p}+
$$
$$
+|y|^{\alpha+\frac d2-|\m|}(|y|+R)^{|\beta|-\alpha}\sum_{\tiny\begin{array}{c}\gamma\in \M_{\geq 0}: \\ |\gamma|_\prec<|\beta|_\prec \\ \frac d2-|\gamma|_p+\alpha\leq 0 \\ \sum_{\n\in\N^4}\gamma(\n)\neq 0\end{array}}\large\left(\frac{|y|}{|y|+R}\large\right)^{-\frac d2+|\gamma|-\alpha}\lesssim
$$
$$
\lesssim |y|^{\alpha+\frac d2-|\m|}(|y|+R)^{|\beta|-\alpha}
$$
as in the last sum $-\frac d2+|\gamma|-\alpha>-\frac d2+|\gamma|_p-\alpha\geq 0$.

Now the first case in (\ref{S0estn}) is obtained from this estimate if $\frac d2-|\m|+\alpha>0$.

If $\frac d2-|\n|+\alpha\leq 0$ one can obtain a better estimate as in the second case in (\ref{S0estn}). This condition implies $|\n|\geq \frac d2+\alpha=2-\e$, and hence by (\ref{dGdvan})
$$
({\rm d}\G_x)_\beta^{\delta_\n}=({\rm d}\G_x^*\z^{\delta_{\n}})_\beta=0,~|\n|\geq 2-\e.
$$
Thus 
$$
(S_{x+y x})_\beta^{\delta_\n}=-(\G_{x+y x}{\rm d}\G_x)_\beta^{\delta_\n}=-\sum_{\gamma\in \Ml}(\G_{x+y x})_\gamma^{\delta_\n}({\rm d}\G_x)_\beta^\gamma.
$$

For $\gamma=\delta_\m$, $\m\in \N^4$ in this sum one must have by (\ref{dGdvan}) in the non--zero terms $|\m|< 2-\e$. On the other hand in this case $(\G_{x+y x})_{\delta_\m}^{\delta_\n}=0$ by (\ref{triangG}) as $|\n|\geq 2-\e$, and conditions (\ref{triangG}), (\ref{tr1}) together with (\ref{prp2}) and (\ref{tr2G*}) imply that the expression for $(S_{x+y x})_\beta^{\delta_\n}$ is reduced to
$$
(S_{x+y x})_\beta^{\delta_\n}=-\sum_{\tiny \begin{array}{c}\gamma\in \M_{\geq 0}: \\ |\n|\leq |\gamma|< |\beta|+\frac d2 \\ |\gamma|_\prec< |\beta|_\prec\end{array}}(\G_{x+y x})_\gamma^{\delta_\n}({\rm d}\G_x)_\beta^\gamma.
$$
Thus we obtain by (\ref{Gest}) and by Proposition \ref{dGindprp}  for $0<R, |y|\leq \frac{\bar{\lambda}''}4$
$$
\left\|(S_{x+y x})_\beta^{\delta_\n}\right\|_{q2B_Rx}\lesssim\sum_{\tiny \begin{array}{c}\gamma\in \M_{\geq 0}: \\ |\n|\leq |\gamma|< |\beta|+\frac d2 \\ |\gamma|_\prec< |\beta|_\prec\end{array}}\left\|(\G_{x+y x})_\gamma^{\delta_\n}\right\|_{q'\infty B_Rx}\left\|({\rm d}\G_x)_\beta^\gamma\right\|_{q''2B_Rx}\lesssim
$$
$$
\lesssim \sum_{\tiny \begin{array}{c}\gamma\in \M_{\geq 0}: \\ |\n|\leq |\gamma|< |\beta|+\frac d2 \\ |\gamma|_\prec< |\beta|_\prec\end{array}}|y|^{-|\n|+|\gamma|}R^{\frac d2+|\beta|-|\gamma|}\lesssim
$$
$$
\lesssim (|y|+R)^{\frac d2+|\beta|-|\n|},
$$ 
where $\frac{1}{q'}+\frac{1}{q''}=\frac 1q$, $1< q'<\infty$, $1\leq q< q''<p^*<2$, $\upsilon$ satisfies $\left\|\left\|\upsilon\right\|_{L_2}\right\|_{\L_{p^*}}\leq 1$, and at the last step we used (\ref{min}). This completes the proof.

\epr

\vskip 0.5cm

\noindent
{\em Proof of Proposition \ref{Sxypp0}}
The proof of Proposition \ref{Sxypp0} is similar to that of Proposition \ref{Sxypp} and even simpler as, if $\beta=0$, the case $\frac d2-|\n|+\alpha\leq 0$ is trivial by (\ref{tr2S*}). Indeed, in this case 
$(S_{x+y x})_0^{\delta_\n}=0$ since $|\n|=|\delta_\n|_\prec\geq |0|_\prec=2-\e=\frac d2+\alpha$.

If $\frac d2-|\n|+\alpha> 0$ from the definition of the modified homogeneity it follows that the sum over $\gamma$ in the right hand side of (\ref{Sm}) is empty. To estimate the first and the second terms in the right hand side of (\ref{Sm}) when $\beta=0$ we use (\ref{PGest=0}). The remaining arguments are the same as in the proof of Proposition \ref{Sxypp}.

\qed

As we mentioned when describing the strategy of the proof of Proposition \ref{indstepprp}, in Sections \ref{indstGGS}, \ref{indstdP<2}, \ref{indstepPb} and \ref{indstGGSpp} the conditions in each statement depend on the results obtained in the preceding statements only and on the induction assumptions imposed in Proposition \ref{indstepprp}. Therefore Propositions \ref{dGpp} and \ref{Sxypp} imply the following corollary.
\bc\lab{dGScor}
Under the induction step assumptions in Proposition \ref{indstepprp} estimates (\ref{G0est*}) and (\ref{G0estS*}) hold for all $\gamma\in \M_{pp}$ and with $\mu$ replaced by $\beta$, $\bar{\lambda}_\prec$ replaced by $\bar{\lambda}=\frac{\bar{\lambda}_\prec}{32}$ and with the same conditions as in  Proposition \ref{indstepprp}.
\ec


\section{Convergence of models}\lab{conv}

In this section we show that the smooth models $(\Pi_x(\rho),\G_{xy}^T(\rho))$ converge in a certain sense when $\rho\to 0$ and the mollification is removed. The proof of this statement is based on a Cauchy property for models which we state in the next section. 


\subsection{The Cauchy property for models}\lab{Mainstatconv}

\setcounter{equation}{0}
\setcounter{theorem}{0}

The proof of the Cauchy property for models closely follows the same steps as the proof of Theorem \ref{mainTes}. We begin by stating this property in this section.

For any function $F(\rho)$ of a real variable $\rho>0$ we denote $\Delta_{\rho, \rho'}F:=F(\rho)-F(\rho')$. In the statement below we keep the notation of Theorem \ref{mainTes}.
\bt\lab{mainTconv}
For $\rho>0$, fixed $N\geq 17$ and $\bar{\lambda}_0>0$ satisfying $\bar{\lambda}_{M(N)}=\frac{\bar{\lambda}_0}{32^{M(N)}}\geq 1$ let $(\Pi_x(\rho),\G_{xy}^T(\rho))$ be the smooth models defined in Theorem \ref{mainTes}.

Then for all $\rho,\rho'>0$, $1\leq p<\infty$, $\gamma,\mu \in \M$ the following inequality holds  
\be\lab{GestMconv}
\left\|(\Delta_{\rho, \rho'}\G_{xy}\right)^\gamma_\mu\|_{\L_p}\lesssim ({\rm max}\{\rho,\rho'\})^\e |x-y|^{|\mu|-|\gamma|},
\ee
for $x,y\in \R^4$, $|x-y|\leq \bar{\lambda}_k$ if $\mu\in \M_k$, $k=0,\ldots, M(N)$.
The left hand side of this inequality is zero if $\mu\not\in \M_{\geq 0}$.

Also, for all $\rho,\rho'>0$, $1\leq p<\infty$, $\gamma,\mu \in \M$, $\varphi\in \mB^r$, $r>{\rm max}\{2-\alpha, |\mu|\}$ the following inequalities hold
\be\lab{EphiMconv} 
\left\|\Delta_{\rho, \rho'}\Pi_{0\mu}^-(\varphi^{\lambda}_x)\right\|_{\L_p}\lesssim ({\rm max}\{\rho,\rho'\})^\e\lambda^{\alpha-2}(\lambda+|x|)^{|\mu|-\alpha},
\ee
\be\lab{EPPsin'>Mconv} 
\left\|\Delta_{\rho, \rho'}\Pi_{0\mu}(\varphi_x^\lambda)\right\|_{\L_p}\lesssim ({\rm max}\{\rho,\rho'\})^\e\lambda^{\alpha}(\lambda+|x|)^{|\mu|-\alpha},
\ee
for all $x\in \R^4$, $|x|\leq \bar{\lambda}_k$, $0<\lambda\leq \bar{\lambda}_k$ if $\mu \in \M_k$, $k=0,\ldots, M(N)$. The left hand sides of inequalities (\ref{EphiMconv}) and (\ref{EPPsin'>Mconv}) are equal to zero if $\mu\not\in \M_{\geq 0}$.

The constants in inequalities (\ref{GestMconv}), (\ref{EphiMconv}) and (\ref{EPPsin'>Mconv}) do not depend on $\lambda$, $\varphi$, $x$, $y$ and $\rho$, $\rho'$.
\et

This theorem is proved by induction over the modified homogeneity similarly to Theorem \ref{mainTes}. The proof will be presented in Section \ref{Indstepconv}. In the next section we obtain the desired estimates in the case of the base of this induction when $\mu=0$.


\subsection{The base case}\lab{Baseindconv}

\setcounter{equation}{0}
\setcounter{theorem}{0}


In the base case when $\mu=0$ the proof of estimates in Theorem \ref{mainTes} differs considerably from the proofs of analogous statements in Section \ref{Baseind}. We start with the following lemma.
\bl\lab{estlem}
For any $0\leq \nu\leq 1$, $1\leq p<\infty$ and $\varphi\in \D$ the function
$$
\phi(x)=\int_{\R^4}\eta(y)\varphi(x-R_\nu y)dy-\varphi(x).
$$
satisfies
$$
\left\|\phi\right\|_{L_p}\lesssim \nu ~\max_{i=0,1,2,3}\left\|\frac{\partial \varphi}{\partial x_i}\right\|_{L_p},
$$
where the constant in the inequality only depends on $\eta$.
\el

\bpr
Since $\int_{\R^4}\eta(y)=1$ we can write
$$
\phi(x)=\int_{\R^4}\eta(y)(\varphi(x-R_\nu y)-\varphi(x))dy.
$$

By the Newton formula and by the chain rule we have
$$
\left|\varphi(x-R_\nu y)-\varphi(x)\right|=\left|\int_0^1\frac{d}{dt}\varphi(x-tR_\nu y)dt \right|\leq \int_0^1\left|\frac{d}{dt}\varphi(x-tR_\nu y)\right|dt\leq
$$
$$
\leq
\int_0^1\sum_{i=0}^3 \left|\frac{\partial \varphi(x-tR_\nu y)}{\partial x_i}\right|\left|(R_\nu y)_i\right|dt\leq \nu \int_0^1\sum_{i=0}^3 \left|\frac{\partial \varphi(x-tR_\nu y)}{\partial x_i}\right|\left|y_i\right|dt,
$$
where at the last step we used the fact that $0\leq \nu\leq 1$, so that $\nu^2\leq \nu$, and hence for $R_\nu y=(\nu^2 y_0,\nu y_1,\nu y_2, \nu y_3)$ each component satisfies $\left|(R_\nu y)_i\right|\leq \nu \left|y_i\right|$.

Now from the previous estimate we obtain using the Minkowski inequality and translational invariance of the $L_p$--norm
$$
\left\|\phi\right\|_{L_p}\leq \nu \int_0^1\left\||\eta(y)|\sum_{i=0}^3 \left\|\frac{\partial \varphi(x-tR_\nu y)}{\partial x_i}\right\|_{L_p(\R^4,x)}\left|y_i\right|\right\|_{L_1(\R^4,y)}dt \leq
$$
$$
\leq \nu ~\max_{i=0,1,2,3}\left\|\frac{\partial \varphi}{\partial x_i}\right\|_{L_p}\left\||\eta(y)|\sum_{i=0}^3\left|y_i\right|\right\|_{L_1(\R^4,y)}\lesssim \nu ~\max_{i=0,1,2,3}\left\|\frac{\partial \varphi}{\partial x_i}\right\|_{L_p},
$$
where at the last step we used the fact that $\eta\in \Sw$. This completes the proof.

\epr

The remaining statements in this section have counterparts in Section \ref{Baseind}. As in that section we start with estimates for Malliavin derivatives.
\bp\lab{dPi-00conv}
For any $1\leq q< 2$, $\bar{\lambda}_0>0$ and $\upsilon\in \L_{q}(L_2)$ satisfying $\left\|\left\|\upsilon(x)\right\|_{L_2(\R^4,x)}\right\|_{\L_q}\leq 1$ one has
\be\lab{Pi00estconv}
\left\|\delta\Delta_{\rho, \rho'}\Pi_{00}^-\right\|_{q\bar{\lambda}_0 1}^{\alpha-2}=\sup_{0<\lambda\leq \bar{\lambda}_0}\sup_{\varphi\in \mB^1}\frac{\left\|\delta\Delta_{\rho, \rho'}\Pi_{00}^-(\varphi^\lambda_0)\right\|_{\L_q}}{\lambda^{\alpha-2}}\lesssim ({\rm max}\{\rho,\rho'\})^\e,
\ee
where $\delta=\delta_\upsilon$ and the constant in the inequality only depends on $\bar{\lambda}_0$.
\ep

\bpr
Without loss of generality we can assume that $\rho\geq \rho'$. Let $0<\lambda\leq \bar{\lambda}_0$. 

Consider first the case when $\nu:=\frac{\rho}{\lambda}\geq 1$. Then under the assumption $0<\lambda \leq \bar{\lambda}_0$ by (\ref{dPi00est}) we have for any $\varphi\in \mB^1$
$$
\left\|\delta\Delta_{\rho, \rho'}\Pi_{00}^-(\varphi^\lambda_0)\right\|_{\L_q}\leq 
\left\|\delta\Pi_{00}^-(\rho)(\varphi^\lambda_0)\right\|_{\L_q}+\left\|\delta\Pi_{00}^-(\rho')(\varphi^\lambda_0)\right\|_{\L_q}\lesssim 
$$
\be\lab{dPi00estconv0}
\lesssim \lambda^{-\frac d2}\leq \lambda^{-\frac d2}\left(\frac{\rho}{\lambda}\right)^\e=({\rm max}\{\rho,\rho'\})^\e \lambda^{\alpha-2}.
\ee

If $\nu=\frac{\rho}{\lambda}\leq 1$ by the definition of $\Pi_{00}^-(x)=\xi^\rho(x)=\xi(\eta^\rho_x)$ we have for any $\varphi\in \mB^1$
\be\lab{dPi00estconv2}
\left\|\delta\Delta_{\rho, \rho'}\Pi_{00}^-(\varphi^\lambda_0)\right\|_{\L_q}= \left\|(\upsilon^\rho-\upsilon+\upsilon-\upsilon^{\rho'})(\varphi^\lambda_0)\right\|_{\L_q}\leq
\ee
$$
\leq \left\|(\upsilon^\rho-\upsilon)(\varphi^\lambda_0)\right\|_{\L_q}+\left\|(\upsilon-\upsilon^{\rho'})(\varphi^\lambda_0)\right\|_{\L_q}.
$$

Now we estimate the first term in in the right hand side. The second one is estimated in a similar way.

We proceed as in (\ref{dPi00est}) using the H\"{o}lder and the Young convolution inequalities, and the assumption
$$
\left\|\left\|\upsilon(x)\right\|_{L_2(\R^4,x)}\right\|_{\L_q}\leq 1
$$
which yields
$$
\left\|(\upsilon^\rho-\upsilon\right)(\varphi^\lambda_0)\|_{\L_q}\leq \left\|\left\|\upsilon(x)\right\|_{L_2(\R^4,x)}\left\|\eta^\rho(\varphi^\lambda_x)-\varphi^\lambda(x)\right\|_{L_2(\R^4,x)}\right\|_{\L_q}\leq 
$$
$$
\leq\left\|\eta^\rho(\varphi^\lambda_x)-\varphi^\lambda(x)\right\|_{L_2(\R^4,x)}=\left\|\phi^\lambda\right\|_{L_2(\R^4,x)}=\lambda^{-\frac d2}\left\|\phi\right\|_{L_2(\R^4,x)},
$$
where
$$
\phi(x)=\int_{\R^4}\eta(y)\varphi(x-R_\nu y)dy-\varphi(x).
$$

Now since $\varphi\in \mB^1$ we have by Lemma \ref{estlem} with $p=2$

\be\lab{dPi00estconv1}
\left\|(\upsilon^\rho-\upsilon)(\varphi^\lambda_0)\right\|_{\L_q} \lesssim \lambda^{-\frac d2}\nu \leq \nu^\e \lambda^{-\frac d2}= \rho^\e \lambda^{\alpha-2}, 
\ee
as for $i=0,1,2,3$ $\left\|\frac{\partial\varphi}{\partial x_i}(x)\right\|_{L_2(\R^4,x)}\lesssim 1$, and $\nu\leq \nu^e$.
Substituting the last estimate and a similar estimate with $\rho$ replaced by $\rho'$ into (\ref{dPi00estconv2}) we obtain
$$
\left\|\delta\Delta_{\rho, \rho'}\Pi_{00}^-(\varphi^\lambda_0)\right\|_{\L_q}\lesssim ({\rm max}\{\rho,\rho'\})^\e \lambda^{\alpha-2} 
$$
This estimate together with (\ref{dPi00estconv0}) imply the statement of this proposition.
This completes the proof.

\epr

Now, as in Proposition \ref{dPi0best0} we proceed with an estimate $\delta\Delta_{\rho, \rho'}\Pi_{00}$ using Schauder estimates for $\K$.

\bp\lab{dPi0best0conv}
For any $1\leq q< 2$, $\bar{\lambda}_0\geq 1$, $r\in \N$, $r\geq 1$ and $\upsilon\in \L_{q}(L_2)$ satisfying $\left\|\left\|\upsilon(x)\right\|_{L_2(\R^4,x)}\right\|_{\L_q}\leq 1$ one has
$$
\left\|\delta\Delta_{\rho, \rho'}\Pi_{00}\right\|_{q\bar{\lambda}_0 r}^{\alpha}=\sup_{0<\lambda\leq \bar{\lambda}_0}\sup_{\varphi\in \mB^{r}}\frac{\left\|\delta\Delta_{\rho, \rho'}\Pi_{00}(\varphi^\lambda_0)\right\|_{\L_q}}{\lambda^{\alpha}}\lesssim ({\rm \max}\{\rho,\rho'\})^\e,
$$
where $\delta=\delta_\upsilon$ and the constant in the inequality only depends on $\bar{\lambda}_0$ and $r$.
\ep

\bpr
Since by (\ref{Ph0}) $\Pi_{00}$ depends on $\xi$ polynomially, one obviously has from (\ref{Ph0})
$$
\delta\Delta_{\rho, \rho'}\Pi_{00}=\K\delta\Delta_{\rho, \rho'}\Pi_{00}^-.
$$

Observe also that $r\geq 1>|\alpha|=\e+\frac 12$.
Hence, the estimate in the statement follows from (\ref{Pi00estconv}), and from (\ref{emb1}) with $x=0$, and $F=\delta\Delta_{\rho, \rho'}\Pi_{00}^-$.

\epr

The next statement is an analogue of Proposition \ref{Fgerm0est0} for differences of germs.
\bp\lab{Fgerm0est0conv}
For any $\bar{\lambda}_0\geq 1$, $R>0$, $r\in \N$, $r\geq 2$, $1\leq q<2$, $\upsilon\in \L_{q}(L_2)$ with $\left\|\left\|\upsilon(x)\right\|_{L_2(\R^4,x)}\right\|_{\L_q}\leq 1$ the germ
$$
F_x(y):=\delta\Delta_{\rho, \rho'}\Pi_{00}(y)-\sum_{\gamma\in \Ml}\Delta_{\rho, \rho'}(\Pi_{x\gamma}(y)\circ ({\rm d}\G_x)_0^\gamma)
$$ 
satisfies
\be\lab{PGest=0conv}
\left\|F\right\|_{q2\bar{\lambda}_0 R r}^{\alpha+\frac d2}=\sup_{0<\lambda\leq \bar{\lambda}_0}\frac{\left\|\sup_{\varphi\in \mB^r}\left\|F_x(\varphi^\lambda_x)\right\|_{\L_q}\right\|_{L_2(B_R,x)}}{\lambda^{\alpha+\frac d2}}\lesssim ({\rm max}\{\rho,\rho'\})^\e,
\ee
where $\delta=\delta_\upsilon$, and the constant in the inequality only depends on $\bar{\lambda}_0$ and $r$.

\ep

\bpr
By (\ref{cdn30})
$$
F_x(y)=\delta\Delta_{\rho, \rho'}\Pi_{00}(y)-\sum_{\tiny\begin{array}{c} \n\in \N^4: \\ |\n|<2-\e\end{array}}\Delta_{\rho, \rho'}(\Pi_{x\delta_\n}(y)\circ ({\rm d}\G_x)_0^{\delta_\n}).
$$ 

From this observation and (\ref{cdn3}) we deduce that for each $x\in \R^4$ 
$F_x(y)$ is obtained from $\delta\Delta_{\rho, \rho'}\Pi_{00}(y)$ by subtracting its Taylor polynomial at $x$ of parabolic degree $1$. Thus $F_x=(\delta\Delta_{\rho, \rho'}\Pi_{00})_x^{2-\e}=(\K (\upsilon^\rho-\upsilon^{\rho'}))_x^{2-\e}$, and we can apply (\ref{emb2c}) with $b=-\e$, $r>2-\e$ to get
\be\lab{emb2c'0conv}
\left\|F \right\|_{q2\bar{\lambda}_0 B_R r}^{2-\e}\lesssim  \left\|\upsilon^\rho-\upsilon^{\rho'}\right\|_{q2~2\bar{\lambda}_0 B_R r}^{-\e}\leq
\ee
$$
\leq \left\|\upsilon^\rho-\upsilon^{\rho'}\right\|_{q2~2\bar{\lambda}_0 \R^4 r}^{-\e}\lesssim \left\|\upsilon^\rho-\upsilon^{\rho'}\right\|_{\B_{q2}^{-\e}},
$$
where at the last step we used definition (\ref{defbspq<}).

Now we estimate the right hand side using definition (\ref{defbspq<'}) of the space  $\B_{q2}^{-\e}$ according to which
\be\lab{emb2c'01conv}
\left\|\upsilon^\rho-\upsilon^{\rho'}\right\|_{\B_{q2}^{-\e}}\lesssim \sup_{0<\lambda\leq 2\bar{\lambda}_0}\frac{\left\|\left\|(\upsilon^\rho-\upsilon^{\rho'})(\varphi_x^\lambda)\right\|_{\L_q}\right\|_{L_2(\R^4,x)}}{\lambda^{-\e}},
\ee
where $\varphi\in \D$ is such that $\int_{\R^4}\varphi(x)dx\neq 0$. Thus we obtain
\be\lab{emb2c'02conv}
\left\|\upsilon^\rho-\upsilon^{\rho'}\right\|_{q2~2\bar{\lambda}_0 B_R r}^{-\e}\lesssim \left\|\upsilon^\rho-\upsilon^{\rho'}\right\|_{\B_{q2}^{-\e}}\lesssim 
\ee
$$
\lesssim
\sup_{0<\lambda\leq 2\bar{\lambda}_0}\frac{\left\|\left\|(\upsilon^\rho-\upsilon^{\rho'})(\varphi_x^\lambda)\right\|_{\L_q}\right\|_{L_2(\R^4,x)}}{\lambda^{-\e}}.
$$

As $1\leq q<2$ we can proceed further by the Minkowski inequality,
\be\lab{emb2c'00conv13}
\left\|\upsilon^\rho-\upsilon^{\rho'}\right\|_{q2~2\bar{\lambda}_0 B_R r}^{-\e}\lesssim \sup_{0<\lambda\leq 2\bar{\lambda}_0}\frac{\left\|\left\|(\upsilon^\rho-\upsilon^{\rho'})(\varphi_x^\lambda)\right\|_{L_2(\R^4,x)}\right\|_{\L_q}}{\lambda^{-\e}}. 
\ee

Now assume without loss of generality that $\rho\geq \rho'$, and as before $0<\lambda\leq 2\bar{\lambda}_0$.
If $\nu:=\frac{\rho}{\lambda}\geq 1$ we use the fact that $\e>0$, the triangle inequality and (\ref{emb2c'00})  to obtain
\be\lab{emb2c'00conv1}
\left\|\left\|(\upsilon^\rho-\upsilon^{\rho'})(\varphi_x^\lambda)\right\|_{L_2(\R^4,x)}\right\|_{\L_q}\leq 
\ee
$$
\leq \left\|\left\|\upsilon^\rho(\varphi_x^\lambda)\right\|_{L_2(\R^4,x)}\right\|_{\L_q}+\left\|\left\|\upsilon^{\rho'}(\varphi_x^\lambda)\right\|_{L_2(\R^4,x)}\right\|_{\L_q}\lesssim 1\leq \nu^\e=({\rm max}\{\rho,\rho'\})^\e\lambda^{-\e}.
$$

Now consider the case when $\frac{\rho}{\lambda}\leq 1$. Then we proceed from (\ref{emb2c'00conv13}) by the triangle inequality as follows
\be\lab{emb2c'00conv14}
\left\|\left\|(\upsilon^\rho-\upsilon^{\rho'})(\varphi_x^\lambda)\right\|_{L_2(\R^4,x)}\right\|_{\L_q}=\left\|\left\|(\upsilon^\rho-\upsilon+\upsilon-\upsilon^{\rho'})(\varphi_x^\lambda)\right\|_{L_2(\R^4,x)}\right\|_{\L_q}\leq
\ee
$$
\leq \left\|\left\|(\upsilon^\rho-\upsilon)(\varphi_x^\lambda)\right\|_{L_2(\R^4,x)}\right\|_{\L_q}+\left\|\left\|(\upsilon-\upsilon^{\rho'})(\varphi_x^\lambda)\right\|_{L_2(\R^4,x)}\right\|_{\L_q}.
$$

We shall estimate the first term in the right hand side of the previous inequality. The second one can be estimated in a similar way. We have by the the Young convolution inequality
\be\lab{emb2c'00conv2}
\left\|\left\|(\upsilon^\rho-\upsilon)(\varphi_x^\lambda)\right\|_{L_2(\R^4,x)}\right\|_{\L_q}\leq \left\|\left\|\upsilon\right\|_{L_2(\R^4)}\left\|\eta^\rho(\varphi^\lambda_x)-\varphi^\lambda(x)\right\|_{L_1(\R^4,x)}\right\|_{\L_q} =
\ee
$$
= \left\|\left\|\upsilon\right\|_{L_2(\R^4)}\right\|_{\L_q}\left\|\eta^\rho(\varphi^\lambda_x)-\varphi^\lambda(x)\right\|_{L_1(\R^4,x)}\leq \left\|\phi^\lambda\right\|_{L_1(\R^4)}=\left\|\phi\right\|_{L_1(\R^4)},
$$
where 
$$
\phi(x)=\int_{\R^4}\eta(y)\varphi(x-R_\nu y)dy-\varphi(x),
$$
and we also used the fact that
$$
\left\|\left\|\upsilon\right\|_{L_2(\R^4)}\right\|_{\L_q}\leq 1.
$$

Now since $\varphi\in \mB^1$ we have by Lemma \ref{estlem} with $p=1$ $\left\|\phi\right\|_{L_1(\R^4)}\lesssim \nu$. Thus from (\ref{emb2c'00conv2}) we deduce
$$
\left\|\left\|(\upsilon^\rho-\upsilon)(\varphi_x^\lambda)\right\|_{L_2(\R^4,x)}\right\|_{\L_q}\lesssim \nu \leq \nu^\e\leq ({\rm max}\{\rho, \rho'\})^\e\lambda^{-e}.
$$

This estimate and a similar estimate with $\rho$ replaced by $\rho'$ together with (\ref{emb2c'00conv14}) yield
$$
\left\|\left\|(\upsilon^\rho-\upsilon^{\rho'})(\varphi_x^\lambda)\right\|_{L_2(\R^4,x)}\right\|_{\L_q}\lesssim ({\rm max}\{\rho, \rho'\})^\e\lambda^{-e}.
$$

Finally from the last estimate and (\ref{emb2c'00conv1}) we deduce that
$$
\left\|\left\|(\upsilon^\rho-\upsilon^{\rho'})(\varphi_x^\lambda)\right\|_{L_2(\R^4,x)}\right\|_{\L_q}\lesssim ({\rm max}\{\rho, \rho'\})^\e\lambda^{-e}
$$
for all $\rho,\rho'>0$ and $0<\lambda\leq 2\bar{\lambda}_0$, where the constant in the inequality only depends on $\bar{\lambda}_0$ and $r$. Using this estimate, (\ref{emb2c'0conv}) and (\ref{emb2c'02conv}) we obtain the desired estimate. This completes the proof.

\epr

Next, as in Proposition \ref{Epsin0p} we can get an estimate for $\Delta_{\rho, \rho'}\Pi_{x0}^-$ with the help of the spectral gap inequality.
\bp\lab{Ephi0conv}
For any $\bar{\lambda}_0>0$, and any $1\leq p<\infty$
\be\lab{Epsin0pconv} 
\left\|\Delta_{\rho, \rho'}\Pi_{x0}^-\right\|_{p\bar{\lambda}_0 0}^{\alpha-2}=\sup_{0<\lambda\leq \bar{\lambda}_0}\sup_{\varphi\in \mB^0}\frac{\left\|\Delta_{\rho, \rho'}\Pi_{x0}^-(\varphi^\lambda_x)\right\|_{\L_p}}{\lambda^{\alpha-2}}=
\ee
$$
=\sup_{0<\lambda\leq \bar{\lambda}_0}\sup_{\varphi\in \mB^0}\frac{\left\|\Delta_{\rho, \rho'}\Pi_{00}^-(\varphi^\lambda_x)\right\|_{\L_p}}{\lambda^{\alpha-2}}\lesssim ({\rm max}\{\rho, \rho'\})^\e
$$
where the constant in the inequality only depends on $\bar{\lambda}_0$ and $p$.
\ep

\bpr
By translational invariance in law it suffices to consider again the case $x=0$, and by the H\"{o}lder inequality in probability one can restrict to the case when $p> 2$. Also, the equality in the right hand side of (\ref{Epsin0pconv}) is obvious as $\Pi_{x0}^-=\xi^\rho$ does not depend on $x$.

We shall use spectral gap inequality (\ref{SpG}) for $f=\Delta_{\rho, \rho'}\Pi_{00}^-(\varphi^\lambda_0)$, $\varphi\in \mB^0$. 
Since the underlying Gaussian measure is centered, $\E(\Delta_{\rho, \rho'}\Pi_{00}^-(\varphi^\lambda_0))=\E((\xi^\rho-\xi^{\rho'})(\varphi^\lambda_0))=0$, and the spectral gap inequality takes the form
$$
\left\|\Delta_{\rho, \rho'}\Pi_{00}^-(\varphi^\lambda_0)\right\|_{\L_p}\lesssim \sup_{\tiny \begin{array}{c}\upsilon\in \L_{p^*}(L_2): \\ \left\|\left\|\upsilon\right\|_{L_2(\R^4)}\right\|_{\L_{p^*}}\leq 1\end{array}}\left\|\delta_\upsilon \Delta_{\rho, \rho'}\Pi_{00}^-(\varphi^\lambda_0)\right\|_{\L_q}\leq
$$
$$
\leq \sup_{\tiny \begin{array}{c}\upsilon\in \L_{p^*}(L_2): \\ \left\|\left\|\upsilon\right\|_{L_2(\R^4)}\right\|_{\L_{p^*}}\leq 1\end{array}}\left\|\delta_\upsilon \Delta_{\rho, \rho'}\Pi_{00}^-(\varphi^\lambda_0)\right\|_{\L_{p^*}},
$$
where $1\leq q\leq p^*< 2< p$, $\frac 1p+\frac{1}{p^*}=1$, and at the last step we used the H\"{o}lder inequality in probability.

Thus by  Proposition \ref{dPi-00conv} with $q=p^*$ for $0<\lambda\leq \bar{\lambda}_0$ one has
$$
\left\|\Delta_{\rho, \rho'}\Pi_{00}^-(\varphi^\lambda_0)\right\|_{\L_p}\lesssim ({\rm max}\{\rho, \rho'\})^\e\lambda^{\alpha-2},
$$
where the constant in the inequality only depends on $\bar{\lambda}_0$ and $p$. This implies (\ref{Epsin0pconv}).

\epr

Finally, using Schauder estimates we can obtain an estimate for $\Delta_{\rho, \rho'}\Pi_{x0}$ similarly to Proposition \ref{Epsin0p+0}.
\bp\lab{Epsin0p+0conv}
For any $\bar{\lambda}_0\geq 1$, $r\in \N$, $r\geq 1$ and any $1\leq p<\infty$
\be\lab{Epsin0p+conv} 
\left\|\Delta_{\rho, \rho'}\Pi_{x0}\right\|_{p\bar{\lambda}_0 r}^{\alpha}=\sup_{0<\lambda\leq \bar{\lambda}_0}\sup_{\varphi\in \mB^r}\frac{\left\|\Delta_{\rho, \rho'}\Pi_{x0}(\varphi^\lambda_x)\right\|_{\L_p}}{\lambda^{\alpha}}=
\ee
$$
=\sup_{0<\lambda\leq \bar{\lambda}_0}\sup_{\varphi\in \mB^r}\frac{\left\|\Delta_{\rho, \rho'}\Pi_{00}(\varphi^\lambda_x)\right\|_{\L_p}}{\lambda^{\alpha}}\lesssim ({\rm max}\{\rho, \rho'\})^\e,
$$
where the constant in the inequality only depends on $\bar{\lambda}_0$, $p$ and $r$.

\ep

\bpr
By translational invariance in law, the estimate in the statement follows from (\ref{Epsin0pconv}) and from (\ref{emb1}) with $b=|0|-2=\alpha-2$ and $F_x=\Delta_{\rho, \rho'}\Pi_{x0}^-$, as in this case by (\ref{defP})
$\Delta_{\rho, \rho'}\Pi_{x0}=\K\Delta_{\rho, \rho'}\Pi_{x0}^-$.
 
\epr

The last two propositions in this section will be only used in the next induction step in the proof of Theorem \ref{mainTconv}.
\bp\lab{dGpp0conv}
For any $1\leq q<2$, $\bar{\lambda}_0>0$, $0<R\leq\bar{\lambda}_0$ and $\upsilon\in \L_{q}(L_2)$ with $\left\|\left\|\upsilon(x)\right\|_{L_2(\R^4,x)}\right\|_{\L_q}\leq 1$  one has
\be\lab{G0estn0conv}
\left\|(\Delta_{\rho, \rho'}{\rm d}\G_{x})^{\delta_\n}_0\right\|_{q2B_Rx}\lesssim ({\rm max}\{\rho, \rho'\})^\e R^{\frac d2+\alpha-|\delta_\n|},
\ee
where $\Delta_{\rho, \rho'}{\rm d}\G_{x}$ is defined with the help of $\upsilon$, and the constant in the inequality only depends on $\bar{\lambda}_0$. The left hand side of the inequality vanishes if $|\n|\geq 2-\e$.

\ep

The vanishing property follows from Proposition \ref{dGindprp0}.
The proof of estimate (\ref{G0estn0conv}) is similar to the proof of Proposition \ref{dGpp0}.

\bp\lab{Sxypp0conv}
For any $1\leq q<2$, $\bar{\lambda}_0>0$, $R\in \R$, $y\in \R^4$, $0<R, |y|\leq\bar{\lambda}_0$, $\n\in \N^4$ and $\upsilon\in \L_{q}(L_2)$ with $\left\|\left\|\upsilon(x)\right\|_{L_2(\R^4,x)}\right\|_{\L_q}\leq 1$  one has
\be\lab{S0estn0conv}
\left\|(\Delta_{\rho, \rho'}S_{x+y x})^{\delta_\n}_0\right\|_{q2B_Rx}\lesssim ({\rm max}\{\rho, \rho'\})^\e|y|^{\frac d2+\alpha-|\delta_\n|},
\ee
where $\Delta_{\rho, \rho'}S_{x+y x}$ is defined with the help of $\upsilon$, and the constant in the inequality only depends on $\bar{\lambda}_0$. The left hand side of the inequality vanishes if $|\n|\geq 2-\e$.

\ep

The vanishing property follows from Proposition \ref{GestP0}.
The proof of estimate (\ref{S0estn0conv}) is similar to the proof of Proposition \ref{Sxypp0}.


\subsection{The induction step statement}\lab{Indstepconv}

\setcounter{equation}{0}
\setcounter{theorem}{0}

The induction step statement in the proof of Theorem \ref{mainTconv} resembles Proposition \ref{indstepprp} used in the proof of Theorem \ref{mainTes}.

Let $\omega\in \D$ and $\psi$ be the functions fixed in Section \ref{Indstep}.
\bp\lab{indstepprpconv}
Let $\rho, \rho'>0$ and $\bar{\lambda}_\prec\in \R$ be such that $\frac{\bar{\lambda}_\prec}8 \geq 1$, and $\beta\in \M_{\geq 0}$, $\beta\neq 0$. Suppose that the conditions of Proposition \ref{indstepprp} are satisfied for the recentered maps and for the structure group maps of the smooth model as in Proposition \ref{indstepprp} with both mollification parameters $\rho, \rho'>0$, and that for all $1\leq p<\infty$, $\gamma,\mu \in \M$ with $|\mu|_\prec<|\beta|_\prec$  one has 
\be\lab{Gest*conv}
\left\|(\Delta_{\rho, \rho'}\G_{xy}\right)^\gamma_\mu\|_{\L_p}\lesssim ({\rm max}\{\rho, \rho'\})^\e|x-y|^{|\mu|-|\gamma|},~x,y\in \R^4,~|x-y|\leq \bar{\lambda}_\prec,
\ee
where the constant in the inequality is independent of $x$, $y$ $\rho$ and $\rho'$.

Assume also that for any $1\leq p<\infty$, $\mu \in \M$ with $|\mu|_\prec<|\beta|_\prec$, and any $\varphi\in \mB^r$ for some (or for any) $r>{\rm max}\{2-\alpha, |\mu|\}$ one has
\be\lab{Ephi*conv} 
\left\|\Delta_{\rho, \rho'}\Pi_{0\mu}^-(\varphi^{\lambda}_x)\right\|_{\L_p}\lesssim 
\ee
$$
\lesssim ({\rm max}\{\rho, \rho'\})^\e \lambda^{\alpha-2}(\lambda+|x|)^{|\mu|-\alpha},~x\in \R^4,~|x|\leq \bar{\lambda}_\prec,~0<\lambda\leq \bar{\lambda}_\prec,
$$
and
\be\lab{EPPsin'>*conv} 
\left\|\Delta_{\rho, \rho'}\Pi_{0\mu}(\varphi_x^\lambda)\right\|_{\L_p}\lesssim 
\ee
$$
\lesssim ({\rm max}\{\rho, \rho'\})^\e \lambda^{\alpha}(\lambda+|x|)^{|\mu|-\alpha},~x\in \R^4,~|x|\leq \bar{\lambda}_\prec,~0<\lambda\leq \bar{\lambda}_\prec,
$$
where the constants in the inequalities do not depend on $\lambda$, $\varphi$, $x$, $\rho$ and $\rho'$.

If $|\beta|<2$, assume in addition that the recentered map $\Pi_{0\beta}^-(\rho)$, $\Pi_{0\beta}^-(\rho')$ satisfy the generalized BPHZ condition with the function $\psi^{\bar{\lambda}_\prec}$, and that for all $1\leq q<2$, $\gamma,\mu\in \M$ with $|\mu|_\prec<|\beta|_\prec$, $q<p^*<2$, $\Delta_{\rho, \rho'}{\rm d}\G_{x}$ and $\Delta_{\rho, \rho'}S_{x+y x}$ defined with the help of $\upsilon\in \L_{p^*}(L_2)$ such that $\left\|\left\|\upsilon\right\|_{L_2}\right\|_{\L_{p^*}}\leq 1$ one has 
\be\lab{G0est*conv}
\left\|(\Delta_{\rho, \rho'}{\rm d}\G_{x})^\gamma_\mu\right\|_{q2B_Rx}\lesssim ({\rm max}\{\rho, \rho'\})^\e R^{\frac d2+|\mu|-|\gamma|},~0<R\leq \bar{\lambda}_\prec,
\ee
and  
\be\lab{G0estS*conv}
\left\|(\Delta_{\rho, \rho'}S_{x+y x})^\gamma_\mu\right\|_{q2B_Rx}\lesssim ({\rm max}\{\rho, \rho'\})^\e\times
\ee
$$
\times \left\{\begin{array}{ll}  |y|^{\frac d2-|\gamma|_p+\alpha}(|y|+R)^{|\mu|-|\gamma|+|\gamma|_p-\alpha} & {\rm if}~\frac d2-|\gamma|_p+\alpha>0 \\  (|y|+R)^{|\mu|-|\gamma|+\frac d2} &{\rm else}\end{array}\right.,~0<R,|y|\leq \frac{\bar{\lambda}_\prec}2,
$$
where the constants in the inequalities are independent of $\rho$, $\rho'$, $R$, $\upsilon$ and $y$.

Then if $|\beta|<2$ (resp. $|\beta|\geq 2$) estimates (\ref{Gest*conv}), (\ref{Ephi*conv}), (\ref{EPPsin'>*conv}), (\ref{G0est*conv}) and (\ref{G0estS*conv}) (resp. (\ref{Gest*conv}), (\ref{Ephi*conv}), (\ref{EPPsin'>*conv})) hold with $\mu$ replaced by $\beta$, $\bar{\lambda}_\prec$ replaced by $\bar{\lambda}=\frac{\bar{\lambda}_\prec}{32}$, and under the same conditions as in the assumptions, where $\mu$ is replaced by $\beta$, and $\bar{\lambda}_\prec$ is replaced by $\bar{\lambda}$.


\ep

\noindent
{\em Strategy of the proof.}
The strategy of the proof of this proposition is the same as that of Proposition \ref{indstepprp}: it is based on analogues of statements in Sections \ref{indstGGS}, \ref{indstdP<2}, \ref{indstepPb} and \ref{indstGGSpp} for the differences $\Delta_{\rho, \rho'}\G_{xy}$, $\Delta_{\rho, \rho'}\Pi_{0\mu}^-$, $\Delta_{\rho, \rho'}\Pi_{0\mu}$, $\Delta_{\rho, \rho'}{\rm d}\G_{x}$, $\Delta_{\rho, \rho'}S_{x+y x}$. In the proof we can assume that Proposition \ref{indstepprp} is already proved, and estimates (\ref{Gest*}), (\ref{Ephi*}), (\ref{EPPsin'>*}), (\ref{G0est*}) and (\ref{G0estS*}) for $\bar{\lambda}_\prec$ replaced by $\bar{\lambda}$, $\mu=\beta$ and $\gamma\in\M$ hold.

Then one can prove the same statements as in Sections \ref{indstGGS}, \ref{indstdP<2}, \ref{indstepPb} and \ref{indstGGSpp}, where in the statements all functions of $\rho$ are replaced with $\Delta_{\rho, \rho'}$ applied to these functions and the right hand sides in all estimates are multiplied by $({\rm max}\{\rho, \rho'\})^\e$. 

The order of the statements remains the same as in Sections \ref{indstGGS}, \ref{indstdP<2}, \ref{indstepPb} and \ref{indstGGSpp}, and the arguments in the proofs are also almost the same if we make the following simple observation. All estimates, arguments and identities which appear in the proofs in Sections \ref{indstGGS}, \ref{indstdP<2}, \ref{indstepPb} and \ref{indstGGSpp} remain valid if we apply $\Delta_{\rho, \rho'}$ to their right hand sides and left hand sides using the following formal rules. $\Delta_{\rho, \rho'}(1)=({\rm max}\{\rho, \rho'\})^\e$, $\Delta_{\rho, \rho'}$ is $\R$--linear, commutes with all norms and satisfies the product rule: for any functions $F_1,\ldots, F_n$, $n\in \N$ of $\rho>0$ 
\be\lab{telesc}
\Delta_{\rho, \rho'}(F_1\ldots F_n)=\sum_{i=1}^nF_1(\rho')\ldots F_{i-1}(\rho')\Delta_{\rho, \rho'}(F_i)F_{i+1}(\rho)\ldots F_n(\rho),
\ee
where the meaning of the product $F_1\ldots F_n$ depends on the context. E.g., for real valued functions, like norms of recentered maps, this is the usual product, for the recentered maps and the structure group maps it is an appropriate composition $\circ$, etc. We stress that these rules, except for (\ref{telesc}), have nothing to do with any proper algebraic identities, and their applicability  stems from the fact that all expressions in Sections \ref{indstGGS}, \ref{indstdP<2}, \ref{indstepPb} and \ref{indstGGSpp} depend on products of functions of $\rho$ estimates for which are already known, and in each statement the estimates for the differences of the functions at different points $\rho$ and $\rho'$, which are used in the proof, are known either by the induction assumption or by the assumption in the statement. So to estimate differences of products of such functions at different points $\rho$ and $\rho'$ one should use the same telescoping trick as in the proof of the product rule for differentiation. Then to estimate a norm of an expression like $\Delta_{\rho, \rho'}(F_1\ldots F_n)$ in (\ref{telesc}) one can use the known estimates for the factors in the products in the right hand side. Note that in each such product only one difference $\Delta_{\rho, \rho'}(F_i)$ appears, and hence the final estimate will only contain one factor $({\rm max}\{\rho, \rho'\})^\e$ in the right hand side, as by the induction assumption only estimates of differences $\Delta_{\rho, \rho'}(F_i)$ can give factors $({\rm max}\{\rho, \rho'\})^\e$.

Estimates (\ref{Gest*conv}), (\ref{G0est*conv}) and (\ref{G0estS*conv}) for $\bar{\lambda}_\prec$ replaced by $\bar{\lambda}$, $\mu=\beta$ and $\gamma\in\M\setminus\M_{pp}=\M'$ can be established similarly to Propositions \ref{Gindpr}, \ref{dGindprp} and \ref{GestP}, respectively, and for $\mu=\beta$ and $\gamma\in\M_{pp}$ they can be obtained similarly to Propositions \ref{Gestn}, \ref{dGpp}, and \ref{Sxypp}, respectively. 

Estimate (\ref{Ephi*conv}) for $\bar{\lambda}_\prec$ replaced by $\bar{\lambda}$, $\mu=\beta$ can be obtained using the telescoping trick and similarly to estimate (\ref{Ephi*}) in Lemma \ref{psilem0}, Corollary \ref{psicor}, Proposition \ref{Pi-b<} for $|\beta|<2$, and to Proposition \ref{Pi-b>} for $|\beta|>2$.

Estimate (\ref{EPPsin'>*conv}) for $\bar{\lambda}_\prec$ replaced by $\bar{\lambda}$, $\mu=\beta$ can be obtained similarly to Proposition \ref{Piint}.

\qed

\vskip 0.3cm

\noindent
{\em Proof of Theorem \ref{mainTconv}.}
Firstly, as in the proof of Theorem \ref{mainTes} note that by Lemma \ref{Gprop} (iv) and (v) estimate (\ref{Gest*conv}) holds for all $\mu\in \M_{pp}$ and for all $\mu\in \M'\setminus \M_{\geq 0}$ with any $\bar{\lambda}_\prec$, and by (\ref{Pipol}), (\ref{Pi0}), (\ref{Pi-0}), (\ref{defP-}), Lemma \ref{lF} (iii), (\ref{Pibnd2}) estimates (\ref{Ephi*conv}), (\ref{EPPsin'>*conv}) also hold for all $\mu\not\in \M_{\geq 0}$ and for any $\bar{\lambda}_\prec$. Actually in all cases the left hand sides in (\ref{Gest*conv}), (\ref{Ephi*conv}) and (\ref{EPPsin'>*conv}) are equal to zero as the expressions in the left hand sides of (\ref{Gest*}), (\ref{Ephi*}), (\ref{EPPsin'>*}) do not depend on $\rho$.

By (\ref{tr1G*}) and  (\ref{tr1S*}) estimates (\ref{G0est*conv}) and (\ref{G0estS*conv}) trivially hold for all $\mu\not\in \M_{\geq 0}$ and for any $\bar{\lambda}_\prec$. 

Therefore it suffices to establish estimates (\ref{GestMconv}), (\ref{EphiMconv}) and (\ref{EPPsin'>Mconv}) for the model defined in Theorem \ref{mainTes}, and for $\mu\in \M_{\geq 0}$, since for $\mu\not\in \M_{\geq 0}$ the left hand sides of these inequalities are equal to zero. We can also assume in the proof that estimates (\ref{Gest*conv}), (\ref{Ephi*conv}), (\ref{EPPsin'>*conv}), (\ref{G0est*conv}) and (\ref{G0estS*conv}) hold for all $\mu\in \M\setminus \M_{\geq 0}$ and for any $\bar{\lambda}_\prec$.

We claim that for the model defined in Theorem \ref{mainTes} estimates (\ref{GestMconv}), (\ref{EphiMconv}) and (\ref{EPPsin'>Mconv}) hold for $\mu\in \M_{\geq 0}$, $|\mu|_\prec <N$. The proof, which is similar to that of Theorem \ref{mainTes},  is by induction over the modified homogeneity.

Firstly, recall that $\mu=0$ is the only element of $\M_{\geq 0}$ with minimal possible modified homogeneity, $|0|_\prec=2-\e$. For $\mu=0$ and arbitrary $\bar{\lambda}_\prec=\bar{\lambda}_0>0$ estimates (\ref{Gest*conv}), (\ref{G0est*conv}), (\ref{G0estS*conv}), (\ref{Ephi*conv}) and (\ref{EPPsin'>*conv}) follow from the results of Sections \ref{Baseind} and \ref{Baseindconv}, see Proposition \ref{Gindpr0} for (\ref{Gest*conv}), Propositions \ref{dGindprp0} and \ref{dGpp0conv} for (\ref{G0est*conv}), Propositions \ref{GestP0} and \ref{Sxypp0conv} for (\ref{G0estS*conv}), Proposition \ref{Ephi0conv} for (\ref{Ephi*conv}), and Proposition \ref{Epsin0p+0conv} for (\ref{EPPsin'>*conv}). Therefore for $\beta\in \M_{\geq 0}$, such that $\beta\neq 0$ and $|\beta|_\prec$ is minimal possible, the induction assumptions are satisfied, and we can begin the induction procedure over the modified homogeneity the induction step of which is established in Proposition \ref{indstepprpconv}. Conditions (\ref{BPHZfix}}) ensure that the generalized BPHZ conditions used in the induction procedure are satisfied. 

Note that by the condition $\bar{\lambda}_{M(N)}\geq 1$, where $M$ is the same as in the proof of Theorem \ref{mainTes}, the induction step can be iterated until we exhaust the set of all multi-indices $\beta\in \M_{\geq 0}$ with $|\beta|_\prec<N$, and after completing the last iteration we obtain estimates (\ref{GestMconv}), (\ref{EphiMconv}) and (\ref{EPPsin'>Mconv}) for all $\mu \in \M_{\geq 0}$ with $|\mu|_\prec<N$ and $\bar{\lambda}=\bar{\lambda}_{M(N)}$. 
This completes the proof.

\qed


\subsection{Convergence of models}\lab{convspace}

\setcounter{equation}{0}
\setcounter{theorem}{0}

In this section we confirm that the smooth models $(\Pi_x(\rho),\G_{xy}^T(\rho))$ converge in a certain sense when $\rho\to 0$. To this end we have to define a space of models where this convergence takes place. In view of the definition of $(\Pi_x(\rho),\G_{xy}^T(\rho))$ and of Theorems \ref{mainTes} and \ref{mainTconv} this task naturally leads to a wider class of models which we are going to introduce now.
 
We keep the notation used in Theorem \ref{mainTes}.
By a model $(\Pi_x,\G_{xy}^T)$ of level $N\in \R$ at scale $\bar{\lambda}_0$ satisfying $\bar{\lambda}_{M(N)}\geq 1$ for the regularity structure with model space $T$ we mean a family of random variables $\Pi_{x\beta}$, $\Pi_{x\beta}^-$ and $(\G_{xy})_\beta^\gamma$ defined for $x,y\in \R^4$, $\beta, \gamma\in \M$, $|\beta|_\prec,|\gamma|_\prec<N$ with values in $\Sw_\beta':=\Sw'\otimes V_\beta$, ${\Sw_\beta'}^-:=\Sw'\otimes V_\beta^-\simeq \Sw_\beta'$, and ${\rm Hom}_{\R}(W_\beta, W_\gamma)$, respectively, and such that the following properties are satisfied for all $x,y,z\in \R^4$.

\begin{enumerate}
\item 
$$
\Pi_{x\delta_{\n}}(y)=(y-x)^{\n}{\rm Id}_{W_{\delta_\n}}^V.
$$

\item
If $\beta\in \M'\setminus \M_{\geq 0}$, $|\beta|_\prec<N$ then $\Pi_{x\beta}^-\in \P_{|\beta|-2}\otimes V_{\beta}^-$ is defined by formula (\ref{defP-}), where $\Phi_\beta^-$ is replaced by $\Pi_{x\beta}^-$, and  $\Phi_{\beta_i}$ are replaced by $\Pi_{x\beta_i}$, and $\beta_i\in \M_{pp}$, $i=1,2,3$ in this case.

\item
$$
\Pi_{x\beta}^-=0~{\rm for }~\beta \not\in \M',
$$
$$
\Pi_{x\beta}=0~{\rm for }~[\beta]<0,~\beta\neq \delta_{\n},~{\n}\in \N^4.
$$

\item
For all $\beta,\gamma\in \M$ with $|\beta|_\prec,|\gamma|_\prec<N$
$$
(\G_{xy})_\beta^\gamma=0~{\rm if}~|\gamma|>|\beta|~{\rm or}~|\gamma|_\prec>|\beta|_\prec,
$$
and $(\G_{xy})_\beta^\beta={\rm Id}_{W_\beta}$.

\item
$$
(\G_{xy})_\beta^0=\delta_\beta^0{\rm Id}_{W_0},~(\G_{xy})_\beta^{\delta_\bg}=\delta_\beta^{\delta_\bg}{\rm Id}_{W_{\delta_\bg}}.
$$

\item
For all $\beta,\gamma\in \M$ with $|\beta|_\prec,|\gamma|_\prec<N$
$$
(\G_{xz})_\beta^\gamma=\sum_{\tiny\begin{array}{c}\delta\in \M: \\ |\gamma|_\prec\leq |\delta|_\prec\leq |\beta|_\prec \\ |\gamma|\leq |\delta|\leq |\beta|\end{array}}(\G_{xy})_\delta^\gamma(\G_{yz})_\beta^\delta.
$$

\item
For any multi--indices $\beta,\gamma\in \M$, $|\beta|_\prec,|\gamma|_\prec<N$, $\gamma=l\delta_\bg+\sum_{i=1}^p\delta_{{\m}_i}\neq 0\in \M$, $\m_i\in \N^4$, $l\in \N$, $l>0$ one has
\be\lab{GbN}
(\G_{xy})_{\beta}^\gamma= \hspace{-1cm}  \sum_{\tiny\begin{array}{c}\beta_1,\ldots,\beta_p\in \M: \\ \beta_1+\ldots+\beta_p+l\delta_\bg=\beta \\ |\beta_i|<|\beta|,i=1,\ldots,p \\ |\beta_i|_\prec<|\beta|_\prec,i=1,\ldots,p \end{array}} \hspace{-1cm}  (\G_{xy})_{\beta_1}^{\delta_{\m_1}}\ldots (\G_{xy})_{\beta_p}^{\delta_{\m_p}},
\ee
where we assume that an empty sum is zero and an empty product is 1. 

Moreover, in the non--zero terms in the sum in the right hand side of (\ref{GbN}) $[\beta_i]<[\beta]$, $i=1,\ldots,p$.

\item
\be\lab{Gn'}
(\G_{xy})_{\delta_\m}^\gamma=\left\{\begin{array}{cc} {\m \choose \n}(x-y)^{\m-\n}{\rm Id}_{W_{\delta_{\m}}}^{W_{\delta_{\n}}} & {\rm if }~ \gamma =\delta_\n,~\n\in \N^4 \\ 0 & {\rm else} \end{array}\right. .
\ee

\item
$$
\gamma\in \M_{\geq 0}\cup \M_{pp},~\beta\not\in \M_{\geq 0}\cup \M_{pp}\Longrightarrow (\G_{xy})_{\beta}^\gamma =0,
$$
$$
\gamma\in \M'\setminus \M_{\geq 0},~\beta\not\in \M'\Longrightarrow(\G_{xy})_{\beta}^\gamma =0.
$$

\item
For $\beta, \gamma \in \M'\setminus \M_{\geq 0}$ the only nonzero $(\G_{xy})_\beta^\gamma\neq 0$ occur if for $\m:=\sum_{i=0}^\infty\gamma(\n_i)\n_i$, ${\bf k}:=\sum_{i=0}^\infty\beta(\n_i)\n_i\in \N^4$ one has $m_i\leq k_i$ for $i=0,1,2,3$, $W_\beta$ is naturally isomorphic to $W_\gamma$, and in this case $(\G_{xy})_\beta^\gamma$ is given by formula (\ref{m'}),
\be\lab{m''}
(\G_{xy})_\beta^\gamma=A_\beta^\gamma (x-y)^{{\bf k}-\m}{\rm Id}_{W_\beta}^{W_\gamma}.
\ee

\item
For all $\beta\in \M$, $|\beta|_\prec<N$
\be\lab{PiGlN}
\Pi_{y\beta}=\Pi_{x\beta}+\sum_{\tiny\begin{array}{c}\delta\in \M: \\ |\delta|_\prec < |\beta|_\prec \\ |\delta|< |\beta|\end{array}}\Pi_{x\delta}\circ (\G_{xy})^\delta_\beta,
\ee
\be\lab{Pi-GlN}
\Pi_{y\beta}^-=\Pi_{x\beta}^-+\sum_{\tiny\begin{array}{c}\delta\in \M: \\ |\delta|_\prec < |\beta|_\prec \\ |\delta|< |\beta|\end{array}}\Pi_{x\delta}^-\circ (\G_{xy})^\delta_\beta.
\ee

\item
For any $h\in \R^4$, $\beta, \gamma\in \M$ with $|\beta|_\prec,|\gamma|_\prec<N$
$$
\Pi_{x\beta}(y)=_{\rm law}\Pi_{x+h\beta}(y+h),~\Pi_{x\beta}^-(y)=_{\rm law}\Pi_{x+h\beta}^-(y+h),
$$
$$
(\G_{x+hy+h})_\beta^\gamma=_{\rm law}(\G_{xy})_\beta^\gamma.
$$

\item
If for $\mu\in \Ml$ we denote by $r_\mu$ the minimal possible natural number such that $r>{\rm max}\{2-\alpha, |\mu|\}$, $\M_k':=\M_k\cup \{\beta \in \M\setminus \M_{\geq 0}:|\beta|_\prec<N\}$ then for all $1\leq p<\infty$   
\be\lab{GestMlN}
\left\|\G\right\|_{Np}:=\max_{0\leq k\leq M}\max_{\tiny\begin{array}{c} \gamma\in \M, |\gamma|_\prec<N \\ \mu \in \M_k'  \end{array}}\sup_{\tiny\begin{array}{c} x,y\in \R^4 \\ x\neq y \\ |x-y|\leq \bar{\lambda}_k\end{array}}\frac{\left\|(\G_{xy})^\gamma_\mu\right\|_{\L_p}}{|x-y|^{|\mu|-|\gamma|}}<\infty,
\ee
\be\lab{EphiMlN} 
\left\|\Pi^-\right\|_{Np}:=\max_{0\leq k\leq M}\max_{\mu \in \M_k'}\sup_{0<\lambda\leq \bar{\lambda}_k}\sup_{\varphi\in \mB^{r_\mu}}\sup_{x\in \R^4}\frac{\left\|\Pi_{x\mu}^-(\varphi^{\lambda}_x)\right\|_{\L_p}}{\lambda^{|\mu|-2}}=
\ee
$$
=\max_{0\leq k\leq M}\max_{\mu \in \M_k'}\sup_{0<\lambda\leq \bar{\lambda}_k}\sup_{\varphi\in \mB^{r_\mu}}\frac{\left\|\Pi_{0\mu}^-(\varphi^{\lambda}_0)\right\|_{\L_p}}{\lambda^{|\mu|-2}}<\infty,
$$
and
\be\lab{EPPsin'>MlN} 
\left\|\Pi\right\|_{Np}:=\max_{0\leq k\leq M}\max_{\mu \in \M_k'}\sup_{0<\lambda\leq \bar{\lambda}_k}\sup_{\varphi\in \mB^{r_\mu}}\sup_{x\in \R^4}\frac{\left\|\Pi_{x\mu}(\varphi^{\lambda}_x)\right\|_{\L_p}}{\lambda^{|\mu|}}=
\ee
$$
=\max_{0\leq k\leq M}\max_{\mu \in \M_k'}\sup_{0<\lambda\leq \bar{\lambda}_k}\sup_{\varphi\in \mB^{r_\mu}}\frac{\left\|\Pi_{0\mu}(\varphi^{\lambda}_0)\right\|_{\L_p}}{\lambda^{|\mu|}}<\infty.
$$
Note that the identities in the left hand sides of the last two formulas follow from item (12) of this definition.

\end{enumerate}

\br\lab{detbound}
By (5), (8), (9), (10) for $\mu \in \M\setminus \M_{\geq 0}$, $|\mu|_\prec<N$ and $\mu=0$ one has the deterministic bound
\be\lab{Gxyleq0}
\max_{\tiny\begin{array}{c} \gamma\in \M \\ |\gamma|_\prec<N \end{array}}\sup_{\tiny\begin{array}{c} x,y\in \R^4 \\ x\neq y \end{array}}\frac{|(\G_{xy})^\gamma_\mu|}{|x-y|^{|\mu|-|\gamma|}}<\infty.
\ee
Indeed, by properties (8), (9) and (10) in the definition of models of level $N$ if $\mu\in \M_{pp}$ then we can assume that $\gamma\in \M_{pp}$ in the right hand side of (\ref{GestMlN}), and if $\mu\in \M'\setminus \M_{\geq 0}$ then we can assume that $\gamma\in \M'\setminus \M_{\geq 0}$. Hence the deterministic bound follows from (\ref{Gn'}) and (\ref{m''}) in these cases. 

If $\mu=0$ the bound trivially follows from (4) and the first property in (5). 

By (1), (2) and (3) for $\mu \in \M\setminus \M_{\geq 0}$, $|\mu|_\prec<N$ one has also the deterministic bounds
$$
\sup_{x\in \R^4}\sup_{0<\lambda<\infty }\sup_{\varphi\in \mB^{r_\mu}}\frac{|\Pi_{x\mu}^-(\varphi^{\lambda}_x)|}{\lambda^{|\mu|-2}}<\infty,
$$
$$
\sup_{x\in \R^4}\sup_{0<\lambda<\infty }\sup_{\varphi\in \mB^{r_\mu}}\frac{|\Pi_{x\mu}(\varphi^{\lambda}_x)|}{\lambda^{|\mu|}}<\infty.
$$
\er

\br
Properties (5) and (7) imply that for any multi--indices $\beta,\gamma\in \M$, $|\beta|_\prec,|\gamma|_\prec, |\gamma_1|_\prec,|\gamma_2|_\prec<N$, $\gamma=\gamma_1+\gamma_2\neq 0\in \M$, one has
\be\lab{GbNN}
(\G_{xy})_{\beta}^\gamma= \hspace{-1cm}  \sum_{\tiny\begin{array}{c}\beta_1,\beta_2\in \M: \\ \beta_1+\beta_2=\beta \\ |\beta_i|<|\beta|,i=1,2 \\ |\beta_i|_\prec<|\beta|_\prec,i=1,2 \end{array}} \hspace{-1cm}  (\G_{xy})_{\beta_1}^{\gamma_1} (\G_{xy})_{\beta_p}^{\gamma_2},
\ee
where we assume that an empty sum is zero and an empty product is 1. Clearly, iterating (\ref{GbNN}) and using (5) one can obtain (\ref{GbN}) as in the proof of identity (\ref{Gb}). We prefer to use (\ref{GbN}) in the definition of models of level $N$ to keep it transparent for future references.
\er

Observe that arguments similar to those in the proof of Lemma \ref{Pi-int} and Proposition \ref{Piint} imply that actually estimates in item (13) in the previous definition lead to the following estimates which hold for all $x\in \R^4$, $|x|\leq \bar{\lambda}_k$, $0<\lambda\leq \bar{\lambda}_k$, $\mu \in \M_k'$, $k=0,\ldots, M(N)$, and any $\varphi\in \mB^{r_\mu}$:
\be\lab{EphiMlN1} 
\left\|\Pi_{0\mu}^-(\varphi^{\lambda}_x)\right\|_{\L_p}\lesssim \lambda^{\alpha-2}(\lambda+|x|)^{|\mu|-\alpha},
\ee
\be\lab{EPPsin'>MlN1} 
\left\|\Pi_{0\mu}(\varphi_x^\lambda)\right\|_{\L_p}\lesssim \lambda^{\alpha}(\lambda+|x|)^{|\mu|-\alpha},
\ee
where the constants in the inequalities do not depend on $\lambda$, $\varphi$, and $x$.

Now we can introduce a space where the convergence of models $(\Pi_x(\rho),\G_{xy}^T(\rho))$ takes place.
Let ${\rm Md}_{N\bar{\lambda}_0}$ be the set of all models of level $N$ at scale $\bar{\lambda}_0$ with the model space $T$. If $(\Pi,\G^T)$, $(\Pi',{\G'}^T)$ are two such models we introduce a family of pseudo--distances ${\rm d}_{Np}$, $1\leq p<\infty$ between them by
$$
{\rm d}_{Np}((\Pi,\G^T),(\Pi',{\G'}^T))=\left\|\G-\G'\right\|_{Np}+\left\|\Pi^--{\Pi'}^-\right\|_{Np}+\left\|\Pi-\Pi'\right\|_{Np}.
$$

Equip ${\rm Md}_{N\bar{\lambda}_0}$ with a topology defined by the family of pseudo--distances ${\rm d}_{Np}$, $1\leq p<\infty$, where a sequence of models in ${\rm Md}_{N\bar{\lambda}_0}$ converges if it converges with respect to all pseudo--distances ${\rm d}_{Np}$, $1\leq p<\infty$. Denote by $\overline{{\rm Md}}_{N\bar{\lambda}_0}$ the corresponding topological space obtained by the completion of the metric identification of ${\rm Md}_{N\bar{\lambda}_0}$. Its topology is defined by the family of metrics corresponding to ${\rm d}_{Np}$, $1\leq p<\infty$. We denote these new metrics by the same symbols. When we speak about a model in the space $\overline{{\rm Md}}_{N\bar{\lambda}_0}$ we mean its equivalence class in $\overline{{\rm Md}}_{N\bar{\lambda}_0}$.

Note that the models in Theorems \ref{mainTes}, \ref{mainTconv} belong to ${\rm Md}_{N\bar{\lambda}_0}$, where $N\geq 17$. This follows from the estimates obtained in Theorem \ref{mainTes}, from the definition of the recentered maps in Section \ref{recm} (see especially properties (\ref{Pipol}), (\ref{Pi-0}), (\ref{Pi0}) and Lemma \ref{lF}(iii)), from the properties of the structure group maps stated in Proposition \ref{Gprop} (i)-(v), Lemma \ref{Gbgmult} and from (\ref{transl1}), (\ref{transl2}).

Now from Theorem \ref{mainTconv} using the Cauchy criterion and from formulas (\ref{Ph0})
we immediately obtain the following result
\bt\lab{mainmod}
For $N\geq 17$ the models $(\Pi(\rho),\G^T(\rho))\in \overline{{\rm Md}}_{N\bar{\lambda}_0}$, $\rho>0$ defined in Theorem \ref{mainTes} converge to a model $(\Pi,\G^T)\in \overline{{\rm Md}}_{N\bar{\lambda}_0}$ when $\rho\to 0$. This model, called a model of level $N$  at scale $\bar{\lambda}_0$ for the regularity structure with model space $T$ associated to equation (\ref{LangdT}), also satisfies $\Pi_{x0}^-=\xi$, $\Pi_{x0}=\K\Pi_{x0}^-=\K\xi$, where the last formula should be understood in the sense of generalized functions, and  
it defines a random variable with values in $\Sw_0'$.
\et


\section{Pointwise estimates for models}\lab{ptwise}

In this section we show that the recentered maps and the structure group maps for models of level $N$ at scale $\bar{\lambda}_0$ introduced in the previous section almost surely satisfy pointwise estimates similar to (\ref{GestMlN}), (\ref{EphiMlN1}) and (\ref{EPPsin'>MlN1}).


\subsection{Preliminary pointwise estimates for recentered maps with a fixed test function}\lab{prestrec}

\setcounter{equation}{0}
\setcounter{theorem}{0}
 
We start with pointwise estimates for the recentered maps.
As in the previous section we keep the notation used in Theorem \ref{mainTes}. For any $\beta\in \Ml$ we also denote by $|\beta|_-$ the corrected homogeneity which is defined by formula (\ref{hom}) where $\alpha$ is replaced with $\alpha_-:=\alpha-\e_-$ for some fixed arbitrary $\e_-$ satisfying $0<\e_-<\frac{\e}{100}$.

By a weight we mean a function $w\in C^\infty(\R^4)$ such that $w(x)>0$ for all $x\in \R^4$, 
\be\lab{sumw}
\sum_{\m\in \mathbb{Z}^4}\frac{1}{w^p(\m)}<\infty,
\ee 
for $p>\frac{12}{\e_-}$, and for all $x,y\in \R^4$, $\left\|x-y\right\|\leq 2$ or $|x-y|\leq 2$ 
\be\lab{bwcd}
\frac{1}{w(x)}\leq C\frac{1}{w(y)}.,
\ee 
where $C>0$ is independent of $x$ and $y$.

We shall also denote $B_1^e(x):=\{y\in \R^4:\left\|y-x\right\|\leq 1\}$.

An example of a weight is the function $w(x)=(1+\left\|x\right\|^2)^{\frac{5\e_-}{24}}$. Indeed, condition (\ref{sumw}) can be verified for it using the multi--dimensional generalization of the integral convergence test for series, and condition (\ref{bwcd}) follows from the obvious estimate for $x,y\in \R^4$ obtained with the help of the triangle inequality,
$$
1+\left\|y\right\|^2\leq 1+2(\left\|x\right\|^2+\left\|x-y\right\|^2)
$$
which for $\left\|x-y\right\|\leq 2$ yields
$$
1+\left\|y\right\|^2\leq 9+2\left\|x\right\|^2\leq 9(1+\left\|x\right\|^2),
$$
and for $|x-y|\leq 2$ it gives 
$$
1+\left\|y\right\|^2\leq 1+2A^2+2\left\|x\right\|^2\leq (1+2A^2+2)(1+\left\|x\right\|^2),
$$
where $A=\max_{z\in \R^4, |z|\leq 2}\left\|z\right\|$. The last two displayed inequalities imply (\ref{bwcd}).

Consider a model of level $N\in \R$ at scale $\bar{\lambda}_0$ satisfying $\bar{\lambda}_{M(N)}\geq 8$ for the regularity structure with model space $T$, with recentered maps $\Pi_{x\beta}$, $\Pi_{x\beta}^-$ and structure group maps $(\G_{xy})_\beta^\gamma$.

We shall obtain pointwise estimates for $\Pi_{x\beta}$, $\Pi_{x\beta}^-$ and $(\G_{xy})_\beta^\gamma$ which hold almost surely. 
\bp\lab{Pipmptwsprel}
Let $\Pi_{x\beta}$, $\Pi_{x\beta}^-$ be the recentered maps of a model of level $N\in \R$ at scale $\bar{\lambda}_0$ satisfying $\bar{\lambda}_{M(N)}\geq 8$. Then for any weight $w$, any $\varphi\in \D(B_{\frac 14})$, for all $\beta \in \M_{\geq 0}$, $|\beta|_\prec<N$, $x,h\in \R^4$, $0\leq |h|\leq \frac 14$, $0<\lambda\leq 8$ one has almost surely 
\be\lab{EphiMlNloc} 
|\Pi_{x\beta}^-({^h}\varphi^{\lambda}_x)|\lesssim \lambda^{|\beta|_--2}w(x),
\ee
\be\lab{EPPsin'>MlNloc} 
|\Pi_{x\beta}(({^h}\varphi)_x^\lambda)|\lesssim \lambda^{|\beta|_-}w(x),
\ee
where ${^h}\varphi(z)=\varphi(z+h)$, the constants in the inequalities do not depend on $h$, $\lambda$, and $x$, and they depend on the choice of $w$, and on the model via quantities (\ref{GestMlN}), (\ref{EphiMlN}) and (\ref{EPPsin'>MlN}), the last dependence being linear for each of them.
\ep

\bpr
We prove (\ref{EPPsin'>MlNloc}). Estimate (\ref{EphiMlNloc}) is established in a similar way. 

We shall apply the usual Kolmogorov continuity trick using a Sobolev type inequality (compare with \cite{DZ}, Proof of Theorem 3.3 or Appendix B in \cite{BOS}). Firstly we shall derive some preliminary estimates.

Fix $\beta \in \M_{\geq 0}$, $|\beta|_\prec<N$. Firstly we obtain a local estimate similar to (\ref{EPPsin'>MlNloc}).
Let $\phi\in \D(B_{\frac 12})$, $v\in \R^4$, $x,y\in B_1(v)$, so that $|x-y|\leq 2$, and $0<\lambda\leq 8$. 

For $1\leq p<\infty$ we start by estimating $\left\|\Pi_{y\beta}(\phi_x^\lambda) -\Pi_{y\beta}(\phi_y^\lambda)\right\|_{\L_p}$.

If $\frac{|x-y|}{\lambda}\leq \frac 12$ we use the Taylor formula and write for any $z\in \R^4$  
$$
\phi_x^\lambda(z)=\frac{1}{\lambda^d}\phi(R_{\frac{1}{\lambda}}(x-z))=\frac{1}{\lambda^d}\phi(R_{\frac{1}{\lambda}}(y-z))+
$$
$$
+\frac{1}{\lambda^d}\sum_{i=0}^3\int_0^1(\partial_i\phi)(R_{\frac{1}{\lambda}}(y-z+t(x-y)))(R_{\frac{1}{\lambda}}(x-y))_i dt=
$$
$$
=\phi_y^\lambda(z)+\sum_{i=0}^3\int_0^1{^i\phi}_y^\lambda(z)(R_{\frac{1}{\lambda}}(x-y))_i dt,
$$
where, due to the condition $t\frac{|x-y|}{\lambda}\leq \frac 12$ one has for $i=0,1,2,3$ that $^i\phi(z):=(\partial_i\phi)(z+tR_{\frac{1}{\lambda}}(x-y))\in \D(B_1)$ by Lemma \ref{philem}, $\left\|^i\phi\right\|_{C^{r_\beta}}=\left\|\partial_i\phi\right\|_{C^{r_\beta}}$ is independent of $x,y,\lambda,t$, and $\left\|^i\phi\right\|_{C^{r_\beta}}=\left\|\partial_i\phi\right\|_{C^{r_\beta}}\leq\left\|\phi\right\|_{C^{r_\beta+2}}$. 

Since the condition $\frac{|x-y|}{\lambda}\leq \frac 12$ also implies $|(R_{\frac{1}{\lambda}}(x-y))_i|\leq \frac{|x-y|}{\lambda}$ from the previous identity and (\ref{EPPsin'>MlN}) we obtain that for all $1\leq p<\infty$
$$
\left\|\Pi_{y\beta}(\phi_x^\lambda)-\Pi_{y\beta}(\phi_y^\lambda)\right\|_{\L_p}=\left\|\sum_{i=0}^3\int_0^1\Pi_{y\beta}({^i\phi}_y^\lambda)(R_{\frac{1}{\lambda}}(x-y))_i dt\right\|_{\L_p}\leq
$$
$$
\leq \sum_{i=0}^3\int_0^1\left\|\Pi_{y\beta}({^i\phi}_y^\lambda)\right\|_{\L_p}|(R_{\frac{1}{\lambda}}(x-y))_i| dt\lesssim \left\|\phi\right\|_{C^{r_\beta+2}}\sum_{i=0}^3\int_0^1 \lambda^{|\beta|-1}|x-y|dt\lesssim
$$
\be\lab{1pest}
\lesssim \left\|\phi\right\|_{C^{r_\beta+2}}\lambda^{|\beta|-1}|x-y|\leq \left\|\phi\right\|_{C^{r_\beta+2}}\lambda^{|\beta|-\e'}|x-y|^{\e'},
\ee
where $0<\e'<\e_+(N):=\min\{1,|\gamma|-|\delta|:\gamma,\delta\in \M, |\gamma|>|\delta|, |\gamma|_\prec,|\delta|_\prec <N \}$, the constant in the inequality is independent of $\phi$, $\lambda$, $x$ and $y$, and at the last step we used that $\frac{|x-y|}{\lambda}\leq \left(\frac{|x-y|}{\lambda}\right)^{\e'}$ as $\frac{|x-y|}{\lambda}\leq \frac 12$.

For $\frac{|x-y|}{\lambda}\leq \frac 12$ we can also obtain using (\ref{PiGlN}), (\ref{GestMlN}), (\ref{EPPsin'>MlN}) and the H\"{o}lder inequality that
\be\lab{2pest}
\left\|\Pi_{y\beta}(\phi_x^\lambda)-\Pi_{x\beta}(\phi_x^\lambda)\right\|_{\L_p}\leq\left\|\sum_{\tiny\begin{array}{c}\delta\in \M: \\ |\delta|_\prec < |\beta|_\prec \\ |\delta|< |\beta|\end{array}}\Pi_{x\delta}(\phi_x^\lambda)\circ (\G_{xy})^\delta_\beta\right\|_{\L_p}\lesssim 
\ee
$$
\lesssim \sum_{\tiny\begin{array}{c}\delta\in \M: \\ |\delta|_\prec < |\beta|_\prec \\ |\delta|< |\beta|\end{array}}\left\|\Pi_{x\delta}(\phi_x^\lambda)\right\|_{\L_{2p}}\left\|(\G_{xy})^\delta_\beta\right\|_{\L_{2p}}\lesssim \left\|\phi\right\|_{C^{r_\beta}}\sum_{\tiny\begin{array}{c}\delta\in \M: \\ |\delta|_\prec < |\beta|_\prec \\ |\delta|< |\beta|\end{array}}|x-y|^{|\beta|-|\delta|}\lambda^{|\delta|}=
$$
$$
=\left\|\phi\right\|_{C^{r_\beta}}|x-y|^{\e'}\lambda^{|\beta|-\e'}\sum_{\tiny\begin{array}{c}\delta\in \M: \\ |\delta|_\prec < |\beta|_\prec \\ |\delta|< |\beta|\end{array}}\left(\frac{|x-y|}{\lambda}\right)^{|\beta|-|\delta|-\e'}\lesssim \left\|\phi\right\|_{C^{r_\beta+2}}|x-y|^{\e'}\lambda^{|\beta|-\e'},
$$
where the constant in the inequality is independent of $\phi$, $\lambda$, $x$ and $y$ and at the last step we used that $\left(\frac{|x-y|}{\lambda}\right)^{|\beta|-|\delta|-\e'}\leq \left(\frac{1}{2}\right)^{|\beta|-|\delta|-\e'}$ for $\frac{|x-y|}{\lambda}\leq \frac 12$ as $|\beta|-|\delta|-\e'>0$ by the choice of $\e'$.

Now from (\ref{1pest}) and (\ref{2pest}) we obtain for $|x-y|\leq 2$, and $0<\lambda\leq 8$, $\frac{|x-y|}{\lambda}\leq \frac 12$
\be\lab{3pest}
\left\|\Pi_{y\beta}(\phi_y^\lambda)-\Pi_{x\beta}(\phi_x^\lambda)\right\|_{\L_p}\leq
\ee
$$
\leq \left\|\Pi_{y\beta}(\phi_y^\lambda)-\Pi_{y\beta}(\phi_x^\lambda)\right\|_{\L_p}+\left\|\Pi_{y\beta}(\phi_x^\lambda)-\Pi_{x\beta}(\phi_x^\lambda)\right\|_{\L_p}\lesssim \left\|\phi\right\|_{C^{r_\beta+2}}|x-y|^{\e'}\lambda^{|\beta|-\e'},
$$
where the constant in the inequality is independent of $\phi$, $\lambda$, $x$ and $y$.

For $|x-y|\leq 2$, and $0<\lambda\leq 8$ and $\frac{|x-y|}{\lambda}\geq \frac 12$ we immediately derive from (\ref{EPPsin'>MlN})
\be\lab{4pest}
\left\|\Pi_{y\beta}(\phi_y^\lambda)-\Pi_{x\beta}(\phi_x^\lambda)\right\|_{\L_p}\lesssim \left\|\phi\right\|_{C^{r_\beta}}\lambda^{|\beta|}\lesssim \left\|\phi\right\|_{C^{r_\beta+2}}|x-y|^{\e'}\lambda^{|\beta|-\e'},
\ee 
where the constant in the inequality is independent of $\phi$, $\lambda$, $x$ and $y$, and we used the fact that $\left(\frac{|x-y|}{\lambda}\right)^{\e'}\geq (\frac 12)^{\e'}$ for $\frac{|x-y|}{\lambda}\geq \frac 12$.

Finally from (\ref{3pest}) and (\ref{4pest}) we obtain that for all $x,y\in \R^4$, $|x-y|\leq 2$, $0<\lambda\leq 8$, $1\leq p<\infty$
\be\lab{5pest}
\left\|\Pi_{y\beta}(\phi_y^\lambda)-\Pi_{x\beta}(\phi_x^\lambda)\right\|_{\L_p}\lesssim \left\|\phi\right\|_{C^{r_\beta+2}}|x-y|^{\e'}\lambda^{|\beta|-\e'},
\ee 
where the constant in the inequality is independent of $\phi$,  $\lambda$, $x$ and $y$.

Now for $1\leq p<\infty$ by (\ref{5pest}) and (\ref{Morr1}) with $d+\vartheta=\e'p>2d$, $f(x)=\Pi_{x\beta}(\phi_x^\lambda)$ and for any $v\in \R^4$, $x\in B_1(v)$, $0<\lambda\leq 8$  we have
\be\lab{6pest}
\E(\sup_{x\in B_1(v)}|\Pi_{x\beta}(\phi_x^\lambda)|^p)\lesssim \E|\Pi_{v\beta}(\phi_v^\lambda)|^p+
\ee
$$
+\iint_{B_1(v)\times B_1(v)}\frac{\E|\Pi_{z\beta}(\phi_z^\lambda)-\Pi_{z'\beta}(\phi_{z'}^\lambda)|^p}{|z-z'|^{\e'p}}dzdz'\lesssim
$$
$$
\lesssim  \left\|\phi\right\|_{C^{r_\beta+2}}^p\left(\lambda^{p|\beta|}+\iint_{B_1(v)\times B_1(v)}\frac{|z-z'|^{p\e'}\lambda^{p(|\beta|-\e')}}{|z-z'|^{\e'p}}dzdz\right)\lesssim \left\|\phi\right\|_{C^{r_\beta+2}}^p\lambda^{p(|\beta|-\e')},
$$
where the constant in the inequality is independent of $\phi$, $\lambda$, $x$ and $v$.


Next, observe that by Lemma \ref{philem} for $h\in \R^4$, $|h|\leq \frac 14$ we have ${^h}\varphi\in \D(B_{\frac 12})$ for $\varphi\in \D(B_{\frac 14})$.
Now we use (\ref{Sobt})  with $f(t)=\Pi_{x\beta}(({^h}\varphi)_x^t)$, $a=|\beta|-\e'$, $\bar{\lambda}=8$, $1<p<\infty$ subject to $\e'p>2d$ to estimate using Fubini's theorem
\be\lab{7pest}
\E(\sup_{x\in B_1(v)}\sup_{h\in B_{\frac 14}}\sup_{0<\lambda\leq 8}|\lambda^{-|\beta|+\e'+\frac 1p}\Pi_{x\beta}(({^h}\varphi)_x^\lambda)|^p)\leq 
\ee
$$
\leq\E(\sup_{x\in B_1(v)}\sup_{h\in B_{\frac 14}}\int_0^8(|t^{-a}f(t)|^p+|t^{-a+1}f'(t)|^p)dt)\leq
$$
$$
\leq \int_0^8(\E(\sup_{x\in B_1(v)}\sup_{h\in B_{\frac 14}}|t^{-a}f(t)|^p)+\E(\sup_{x\in B_1(v)}\sup_{h\in B_{\frac 14}}|t^{-a+1}f'(t)|^p))dt.
$$

By Sobolev's embedding $C^0(B_{\frac 14})\supset W_p^1(B_{\frac 14})$, $p>4$ for $F\in W_p^1(B_{\frac 14})$ we have
\be
\sup_{h\in B_{\frac 14}}|F(h)|^p\lesssim \int_{B_{\frac 14}}(|F(h)|^p+\sum_{i=0}^3|\partial_iF(h)|^p)dh.
\ee

Applying this inequality to the first and to the second terms in the right hand side of (\ref{7pest}) with $F=f$ and $F=f'$, respectively, we obtain for $p>4$
\be\lab{7pest'}
\E(\sup_{x\in B_1(v)}\sup_{h\in B_{\frac 14}}\sup_{0<\lambda\leq 8}|\lambda^{-|\beta|+\e'+\frac 1p}\Pi_{x\beta}(({^h}\varphi)_x^\lambda)|^p)\lesssim 
\ee
$$
\lesssim \int_0^8\int_{B_{\frac 14}}(\E(\sup_{x\in B_1(v)}t^{-ap}(|\Pi_{x\beta}(({^h}\varphi)_x^t)|^p+\sum_{i=0}^3|\Pi_{x\beta}(({^h}\partial_i\varphi)_x^t)|^p))+
$$
$$
+\E(\sup_{x\in B_1(v)}t^{(-a+1)p}(|\Pi_{x\beta}(\frac{d}{dt}({^h}\varphi)_x^t)|^p+\sum_{i=0}^3|\Pi_{x\beta}(\frac{d}{dt}({^h}\partial_i\varphi)_x^t)|^p)))dtdh.
$$

Note that by the choice of $\varphi$ and $h$ we have by Lemma \ref{philem} ${^h}\varphi, {^h}\partial_i\varphi\in \D(B_{\frac 12})$, and 
\be\lab{normh}
\left\|{^h}\varphi\right\|_{C^r}=\left\|\varphi\right\|_{C^r}, \left\|{^h}\partial_i\varphi\right\|_{C^r}=\left\|\partial_i\varphi\right\|_{C^r}\leq \left\|\varphi\right\|_{C^{r+2}},~ r\in \N.
\ee
Therefore we can estimate the terms in the right hand side of (\ref{7pest'}) using (\ref{6pest}).

For instance, by (\ref{6pest}) with $\phi={^h}\varphi$ for $0<t\leq 8$ we obtain for $\e'p>2d$
\be\lab{8pest}
\E(\sup_{x\in B_1(v)}|t^{-a}\Pi_{x\beta}(({^h}\varphi)_x^t|^p)=\E(\sup_{x\in B_1(v)}|t^{-(|\beta|-\e')}\Pi_{x\beta}(({^h}\varphi)_x^t|^p)\lesssim 
\ee
$$
\lesssim \left\|\varphi\right\|_{C^{r_\beta+2}}^pt^{-p(|\beta|-\e')}t^{p(|\beta|-\e')}=\left\|\varphi\right\|_{C^{r_\beta+2}}^p,
$$
where the constant in the inequality is independent of $h$, $t$, $v$ and $\varphi$.

Similarly, since by the chain rule
$$
\frac{d}{dt}({^h}\varphi)_x^t=\frac{1}{t}({^h}\hat{\varphi})_x^t,
$$
where
$$
{^h}\hat{\varphi}(z)=-d({^h}\varphi)(z)-2({^h}\partial_0\varphi)(z)z_0-\sum_{i=1}^3({^h}\partial_i\varphi)(z)z_i\in \D(B_{\frac 12}),
$$
we can estimate for $0<t\leq 8$, $\e'p>2d$
\be\lab{9pest}
\E(\sup_{x\in B_1(v)}|t^{-a+1}|\Pi_{x\beta}(\frac{d}{dt}({^h}\varphi)_x^t)|^p)=
\ee
$$
=\E(\sup_{x\in B_1(v)}|t^{-(|\beta|-\e')+1-1}\Pi_{x\beta}(({^h}\hat{\varphi})_x^t)|^p)\lesssim
$$
$$
\lesssim \left\|\varphi\right\|_{C^{r_\beta+4}}^pt^{-p(|\beta|-\e')}t^{p(|\beta|-\e')}=\left\|\varphi\right\|_{C^{r_\beta+4}}^p,
$$
where the constant in the inequality is independent of $t$, $v$ and $\varphi$.

Now using (\ref{8pest}) and (\ref{9pest}) and similar estimates for the remaining terms in the right hand side of (\ref{7pest'}), which can be obtained by replacing $\varphi$ with $\partial_i \varphi$ in (\ref{8pest}) and (\ref{9pest}), we obtain from (\ref{7pest'}), noting the second inequality in (\ref{normh}), 
$$
\E(\sup_{x\in B_1(v)}\sup_{h\in B_{\frac 14}}\sup_{0<\lambda\leq 8}|\lambda^{-|\beta|+\e'+\frac 1p}\Pi_{x\beta}(({^h}\varphi)_x^\lambda)|^p)\lesssim \int_0^8\int_{B_{\frac 14}}\left\|\varphi\right\|_{C^{r_\beta+6}}^p dtdh\lesssim 
$$
\be\lab{10pest}
\lesssim \left\|\varphi\right\|_{C^{r_\beta+6}}^p,
\ee
where $\e'p>2d$, $p>4$, and the constant in the inequality is independent of $v$ and $\varphi$. 

Next, we can obtain the global estimate as follows (compare the second line in the formula below with Definition 5.1 in \cite{BCCH})
$$
\E(\sup_{x\in \R^4}\sup_{h\in B_{\frac 14}}\sup_{0<\lambda\leq 8}|\lambda^{-|\beta|+\e'+\frac 1p}\Pi_{x\beta}(({^h}\varphi)_x^\lambda)\frac{1}{w(x)}|^p)=
$$
$$
=\E(\sup_{\m\in \mathbb{Z}^4}\sup_{x\in B_1^e(\m)}\sup_{h\in B_{\frac 14}}\sup_{0<\lambda\leq 8}|\lambda^{-|\beta|+\e'+\frac 1p}\Pi_{x\beta}(({^h}\varphi)_x^\lambda)|^p\frac{1}{w^p(x)})\leq
$$
$$
\leq C^p\E(\sum_{\m\in \mathbb{Z}^4}\sup_{x\in B_1^e(\m)}\sup_{h\in B_{\frac 14}}\sup_{0<\lambda\leq 8}|\lambda^{-|\beta|+\e'+\frac 1p}\Pi_{x\beta}(({^h}\varphi)_x^\lambda)|^p\frac{1}{w^p(\m)}),
$$
where we used (\ref{bwcd}) with $y=\m$ and the fact that the balls $B_1^e(\m)$, $\m\in \mathbb{Z}^4$ cover $\R^4$ which follows from the inclusions $\Upsilon_\m:=\{x=(x_0,x_1,x_2,x_3)\in \R^4:|x_i-m_i|\leq \frac 12,~i=0,1,2,3\}\subset B_1^e(\m)$, $\m=(m_0,m_1,m_2,m_3)$, and from the obvious identity $\R^4=\bigcup_{\m\in \mathbb{Z}^4}\Upsilon_\m$.

Now we proceed by the monotone convergence theorem which yields
$$
\E(\sup_{x\in \R^4}\sup_{h\in B_{\frac 14}}\sup_{0<\lambda\leq 8}|\lambda^{-|\beta|+\e'+\frac 1p}\Pi_{x\beta}(({^h}\varphi)_x^\lambda)\frac{1}{w(x)}|^p)\leq 
$$
$$
\leq C^p\sum_{\m\in \mathbb{Z}^4}\E(\sup_{x\in B_1^e(\m)}\sup_{h\in B_{\frac 14}}\sup_{0<\lambda\leq 8}|\lambda^{-|\beta|+\e'+\frac 1p}\Pi_{x\beta}(({^h}\varphi)_x^\lambda)|^p)\frac{1}{w^p(\m)},
$$ 

Next, observe that from the definition of the norms $\left\|\cdot\right\|$ and $|\cdot |$ it follows that $B_1^e(\m)\subset B_1(\m)$ for all $\m\in \mathbb{Z}^4$, and hence
\be\lab{10'pest}
\E(\sup_{x\in \R^4}\sup_{h\in B_{\frac 14}}\sup_{0<\lambda\leq 8}|\lambda^{-|\beta|+\e'+\frac 1p}\Pi_{x\beta}(({^h}\varphi)_x^\lambda)\frac{1}{w(x)}|^p)\leq 
\ee
$$
\leq C^p\sum_{\m\in \mathbb{Z}^4}\E(\sup_{x\in B_1(\m)}\sup_{h\in B_{\frac 14}}\sup_{0<\lambda\leq 8}|\lambda^{-|\beta|+\e'+\frac 1p}\Pi_{x\beta}(({^h}\varphi)_x^\lambda)|^p)\frac{1}{w^p(\m)}.
$$ 

Choose $p$ such that $p>\max\{4,\frac{11}{\e_+(N)},\frac{12}{\e_-}\}$ and $\e'=\frac{11}{p}$, so that $\e'p=11>10=2d$, $p>4$ and $\e'=\frac{11}{p}<\e_+(N)$ as required in (\ref{10pest}).
Therefore we can estimate the expectations in the right hand side of (\ref{10'pest}) using (\ref{10pest}),
$$
\E(\sup_{x\in \R^4}\sup_{h\in B_{\frac 14}}\sup_{0<\lambda\leq 8}|\lambda^{-|\beta|+\e'+\frac 1p}\Pi_{x\beta}(({^h}\varphi)_x^\lambda)\frac{1}{w(x)}|^p)\lesssim \sum_{\m\in \mathbb{Z}^4}\left\|\varphi\right\|_{C^{r_\beta+6}}^p\frac{1}{w^p(\m)}\lesssim
$$
\be\lab{10''pest}
\lesssim \left\|\varphi\right\|_{C^{r_\beta+6}}^p,
\ee
where the constant in the inequality is independent of $\varphi$, and at the last step we used (\ref{sumw}) which implies $\sum_{\m\in \mathbb{Z}^4}\frac{1}{w^p(\m)}<\infty$.

Inequality (\ref{10''pest}) obviously yields that almost surely for $x,h\in \R^4$, $|h|\leq \frac 14$, $0<\lambda\leq 8$, 
\be\lab{11pest}
|\Pi_{x\beta}(({^h}\varphi)_x^\lambda)|\lesssim \lambda^{|\beta|-(\e'+\frac 1p)}w(x),
\ee
where the constant in the inequality does not depend on $h$, $\lambda$, and $x$.

Note that by the choice of $p$ we have $\e'+\frac 1p=\frac{12}{p}$. Thus, since $[\beta]\geq 0$ and $\frac{12}{p}< \e_-$, for $0<\lambda\leq 8$ we obtain
$$
\lambda^{|\beta|-(\e'+\frac 1p)}=\lambda^{|\beta|-\frac{12}{p}}\lesssim \lambda^{|\beta|-\e_-}\lesssim \lambda^{|\beta|-\e_-([\beta]+1)},
$$
where the constants in the inequality depends on $\e_-$, $\beta$ and $p$.

By the definition of the homogeneity 
\be\lab{hom-}
|\beta|_-=|\beta|-\e_-([\beta]+1),
\ee 
so we deduce that for $0<\lambda\leq 8$
$$
\lambda^{|\beta|-(\e'+\frac 1p)}\lesssim \lambda^{|\beta|_-},
$$
and hence with our choice of $\e'$ and $p$ inequality (\ref{11pest}) gives (\ref{EPPsin'>MlNloc}). 

By (\ref{1pest})-(\ref{5pest}) the constant in (\ref{6pest}) depends on the model via quantities (\ref{GestMlN}) and (\ref{EPPsin'>MlN}), the dependence being linear for each of them. Therefore the same fact is true for the constant in (\ref{11pest}) and in (\ref{EPPsin'>MlNloc}). This completes the proof.

\epr


\subsection{Pointwise estimates for structure group maps}\lab{ptwsstrgrp}

\setcounter{equation}{0}
\setcounter{theorem}{0}
 
Pointwise estimates for the structure group maps can be obtained using induction over the modified homogeneity and pointwise analogues of Propositions \ref{Gindpr} and \ref{Gestn}.
We keep the notation used in Theorem \ref{mainTes} and in Section \ref{prestrec} and fix $\e_-\in \R$, $0<\e_-<\frac{\e}{100}$, with the help of which the corrected homogeneity is defined.

If $w$ is a weight then for $\beta\in\Ml$ we define 
$$
w_\beta:=\left\{\begin{array}{ll} w & {\rm if}~\beta=0 \\ w^{\beta(\bg)} & {\rm if}~\beta\neq 0\end{array}. \right.
$$

Consider a model of level $N\in \R$ at scale $\bar{\lambda}_0$ satisfying $\bar{\lambda}_{M(N)}\geq 8$ for the regularity structure with model space $T$, with recentered maps $\Pi_{x\beta}$, $\Pi_{x\beta}^-$ and structure group maps $(\G_{xy})_\beta^\gamma$.

We shall obtain pointwise estimates for $(\G_{xy})_\beta^\gamma$ which hold almost surely. 
\bt\lab{GestMNprop}
Let $(\G_{xy})_\mu^\gamma$ be the structure group maps of a model of level $N\in \R$ at scale $\bar{\lambda}_0$ satisfying $\bar{\lambda}_{M(N)}\geq 8$. Then for any weight $w$, for all $\mu,\gamma \in \M$, $|\mu|_\prec,|\gamma|_\prec<N$ the following inequality holds almost surely 
\be\lab{GestMN}
|(\G_{xy})^\gamma_\mu|\lesssim |x-y|^{|\mu|_--|\gamma|_-}w_\mu(x) w_\gamma^{-1}(x)
\ee
for all $x,y\in \R^4$, $|x-y|\leq 2$,  the constant in the inequality does not depend on $x,y$, and it depends on the choice of $w$, and on the model via quantities (\ref{GestMlN}) and (\ref{EPPsin'>MlN}), the last dependence being polynomial, without a constant term.
\et
\br
Note that by (\ref{GbN}) for $\gamma\neq 0$ one has $\mu(\bg)\geq\gamma(\bg)$ in (\ref{GestMN}), and by part (5) in the definition of models in Section \ref{convspace} $(\G_{xy})_\mu^0=\delta_\mu^0{\rm Id}_{W_0}$. Hence $w_\mu(x) w_\gamma^{-1}(x)=w^n$, $n\in \N$.
\er

\bpr
First note that by (\ref{Gxyleq0}) estimates (\ref{GestMN}) hold for $\mu \in \M\setminus \M_{\geq 0}$, $|\mu|_\prec<N$. Indeed, as we observed in Remark \ref{detbound} in this case we can assume that $\gamma\in \M\setminus \M_{\geq 0}$ in (\ref{GestMN}). But in this case $[\mu]=[\gamma]=-1$, and by (\ref{hom-}) $|\mu|_-=|\mu|$, $|\gamma|_-=|\gamma|$, so (\ref{GestMN}) is reduced to (\ref{Gxyleq0}).

Therefore it suffices to consider the case $\mu \in \M_{\geq 0}$. The proof now is by induction over the modified homogeneity of $\mu$.

If $\mu=0$ (\ref{GestMN}) trivially holds by property (4) in the definition of models of level $N$ as it was already observed in Remark \ref{detbound}.

Now the induction step follows from the next two propositions which are analogues of Propositions \ref{Gindpr} and \ref{Gestn}.

Firstly, given that (\ref{GestMN}) is true for all $\mu \in \M$ with $|\mu|_\prec<|\beta|_\prec$ for some $\beta\in \M$, $\beta\neq 0$ we show that (\ref{GestMN}) holds almost surely for $\mu=\beta$ and for all $\gamma\in\M\setminus\M_{pp}$, $|\gamma|_\prec<N$.
\bp\lab{GindprN}
Suppose that for some $\beta\in \M$, $\beta\neq 0$, $|\beta|_\prec<N$ and for all $\gamma\in\M_{pp}$, $\mu \in \M$ with $|\mu|_\prec<|\beta|_\prec$, $|\gamma|_\prec<N$ one has almost surely
\be\lab{GestN}
|(\G_{xy})^\gamma_\mu|\lesssim |x-y|^{|\mu|_--|\gamma|_-}w_\mu(x)w_\gamma^{-1}(x),~x,y\in \R^4,~|x-y|\leq 2, 
\ee
where the constant in the inequality is independent of $x$ and $y$.
Then (\ref{GestN}) holds almost surely for $\mu=\beta$, for all $\gamma\in\M\setminus\M_{pp}$, $|\gamma|_\prec<N$. 
\ep

\bpr
For $\gamma=0$ the statement is void by the first property in part (5) of the definition of models of level $N$ in Section \ref{convspace}.

If $\gamma\neq 0$, $\gamma\in\M\setminus\M_{pp}$ then $\gamma=l\delta_\bg+\sum_{i=1}^r\delta_{{\m}_i}$, $\m_i\in \N^4$, $l\in \N$, $l>0$, and we apply the norm $|\cdot|$ to (\ref{GbN}), and use property (9) in the definition of models of level $N$ and the multiplicative property of the norm to get
\be\lab{gptwsind}
|(\G_{xy})^\gamma_\beta|\lesssim \hspace{-1cm}  \sum_{\tiny\begin{array}{c}\beta_1,\ldots,\beta_r\in \M_{\geq 0}\cup \M_{pp}: \\ \beta_1+\ldots+\beta_r+l\delta_\bg=\beta \\ |\beta_i|_\prec<|\beta|_\prec,i=1,\ldots,r \\ |\beta_i|<|\beta|,i=1,\ldots,r \end{array}} \hspace{-1cm}  |(\G_{xy})_{\beta_1}^{\delta_{\m_1}}|\ldots |(\G_{xy})_{\beta_r}^{\delta_{\m_r}}|.,
\ee
where the constant in the inequality is independent of $x$ and $y$.

By the induction assumption we now obtain that almost surely for all $x,y\in \R^4$, $|x-y|\leq 2$ 
\be\lab{GxyindN}
|(\G_{xy})^\gamma_\beta|\lesssim \hspace{-1cm}  \sum_{\tiny\begin{array}{c}\beta_1,\ldots,\beta_r\in \M_{\geq 0}\cup \M_{pp}: \\ \beta_1+\ldots+\beta_r+l\delta_\bg=\beta \\ |\beta_i|_\prec<|\beta|_\prec,i=1,\ldots,r \\ |\beta_i|<|\beta|,i=1,\ldots,r \end{array}} \hspace{-1cm}|x-y|^{\sum_{i=1}^r(|\beta_i|_--|\delta_{\m_i}|_-)}\prod_{i=1}^r w_{\beta_i}(x), 
\ee
where the constant in the inequality is independent of $x$ and $y$

By the linearity of the function $|\cdot|_--\alpha_-$ we obtain for all terms in the right hand side of the last inequality
\be\lab{bd-N}
\sum_{i=1}^r(|\beta_i|_--|\delta_{\m_i}|_-)=|\beta|_--|\gamma|_-.
\ee

Since $\beta_1+\ldots+\beta_r+l\delta_\bg=\beta$ in the terms in the sum in (\ref{GxyindN}) 
we also have $\beta(\bg)=\sum_{i=1}^r\beta_i(\bg)+l=\sum_{i=1}^r\beta_i(\bg)+\gamma(\bg)$ for each such term. Thus 
$$
\prod_{i=1}^r w_{\beta_i}=w^{\sum_{i=1}^r\beta_i(\bg)}=w_\beta w_\gamma^{-1}.
$$
Substituting this expression and (\ref{bd-N}) into (\ref{GxyindN}) yields the result.
 
\epr

Next, given that (\ref{GestMN}) is true for all $\mu \in \M$ with $|\mu|_\prec<|\beta|_\prec$ for some $\beta\in \M$, $\beta\neq 0$ and for $\mu=\beta$ with arbitrary $\gamma\in\M\setminus\M_{pp}$, $|\gamma|_\prec<N$, we show that (\ref{GestMN}) holds almost surely for $\mu=\beta$ and for all $\gamma\in \M_{pp}$, $|\gamma|_\prec<N$.
\bp\lab{GestnN}
Assume that for some $\beta\in \M_{\geq 0}$, $\beta\neq 0$, $|\beta|_\prec<N$, estimates (\ref{GestN}) hold for $\mu=\beta$, for all $\gamma\in\M\setminus\M_{pp}=\M'$, $|\gamma|_\prec<N$, $x,y\in \R^4$, $|x-y|\leq 2$. 
Then for any $x,y\in \R^4$, $0<|x-y|\leq 2$, $\n\in \N^4$, $|\n|<N$ one has almost surely
\be\lab{GxyestnN}
|(\G_{xy})^{\delta_\n}_\beta|\lesssim |x-y|^{|\beta|_--|\delta_\n|_-}w_\beta(x),
\ee
where the constant in the inequality does not depend on $x,y$.
\ep

\bpr
By properties (3), (4) and  (11) in the definition of models of level $N$  we have
$$
\Pi_{y\beta}(z)=\sum_{\gamma\in \Ml}\Pi_{x\gamma}(z)\circ (\G_{xy})_\beta^\gamma=\sum_{\tiny\begin{array}{c}\gamma\in \M_{\geq 0}\cup\M_{pp}: \\ |\gamma|_\prec<|\beta|_\prec \\ |\gamma|<|\beta| \end{array}}\Pi_{x\gamma}(z)\circ (\G_{xy})_\beta^\gamma+\Pi_{x\beta}(z) =
$$
$$
=\Pi_{x\beta}(z)+\sum_{\tiny\begin{array}{c}\gamma\in \M_{\geq 0}: \\ |\gamma|_\prec<|\beta|_\prec \\ |\gamma|<|\beta| \end{array}}\Pi_{x\gamma}(z)\circ (\G_{xy})_\beta^\gamma+\sum_{\tiny\begin{array}{c}\n\in \N^4: \\ |\n|<|\beta|_\prec \end{array}}\Pi_{x\delta_\n}(z)\circ (\G_{xy})_\beta^{\delta_\n}.
$$
Using the definition $\Pi_{x\delta_\n}(y)=(y-x)^\n{\rm Id}^V_{W_{\delta_\n}}$ we can rewrite the previous identity as follows
$$
\sum_{\tiny\begin{array}{c}\n\in \N^4: \\ |\n|<|\beta|_\prec \end{array}}(z-x)^\n{\rm Id}^V_{W_{\delta_\n}} (\G_{xy})_\beta^{\delta_\n}=\sum_{\tiny\begin{array}{c}\n\in \N^4: \\ |\n|<|\beta|_\prec \end{array}}\Pi_{x\delta_\n}(z)\circ (\G_{xy})_\beta^{\delta_\n}=
$$
$$
=\Pi_{y\beta}(z)-\sum_{\tiny\begin{array}{c}\gamma\in \M_{\geq 0}: \\ |\gamma|_\prec<|\beta|_\prec \\ |\gamma|<|\beta| \end{array}}\Pi_{x\gamma}(z)\circ (\G_{xy})_\beta^\gamma-\Pi_{x\beta}(z).
$$

Now we multiply this identity by the function $\phi^\lambda_x(z)$ and integrate with respect to  $z$ over $\R^4$, where $\phi$ is given in Lemma \ref{phi0lem} with $N$ being the level of the model and arbitrary fixed $\m\in \N^4$, $|\m|<N$. This yields
$$
\lambda^{|\m|}{\rm Id}^V_{W_{\delta_\m}} (\G_{xy})_\beta^{\delta_\m}=\Pi_{y\beta}(\phi^\lambda_x)-\Pi_{x\beta}(\phi^\lambda_x)-\sum_{\tiny\begin{array}{c}\gamma\in \M_{\geq 0}: \\ |\gamma|_\prec<|\beta|_\prec \\|\gamma|<|\beta| \end{array}}\Pi_{x\gamma}(\phi^\lambda_x)\circ (\G_{xy})_\beta^\gamma.
$$

Thus using the multiplicativity property of the norm we obtain
\be\lab{GmN}
\lambda^{|\m|}|(\G_{xy})^{\delta_\m}_\beta|\lesssim |\Pi_{y\beta}(\phi^\lambda_x)|+|\Pi_{x\beta}(\phi^\lambda_x)|
+\sum_{\tiny\begin{array}{c}\gamma\in \M_{\geq 0}: \\ |\gamma|_\prec<|\beta|_\prec \\|\gamma|<|\beta| \end{array}}|\Pi_{x\gamma}(\phi^\lambda_x)| |(\G_{xy})_\beta^\gamma|=
\ee
$$
=|\Pi_{y\beta}({^h}\phi^\lambda_y)|+|\Pi_{x\beta}(\phi^\lambda_x)|
+\sum_{\tiny\begin{array}{c}\gamma\in \M_{\geq 0}: \\ |\gamma|_\prec<|\beta|_\prec \\|\gamma|<|\beta| \end{array}}|\Pi_{x\gamma}(\phi^\lambda_x)| |(\G_{xy})_\beta^\gamma|,
$$
where $h=R_{\frac{1}{\lambda}}(x-y)$.

Now assume that $0<\lambda \leq 8$. If $\frac{|x-y|}\lambda\leq \frac 14$ then we can apply (\ref{EPPsin'>MlNloc}) with $w=w_\beta$, and $h=R_{\frac{1}{\lambda}}(x-y)$ (resp. $h=0$) in the first (resp. in the second) term. We can also use (\ref{EPPsin'>MlNloc}) with $h=0$, $w=w_\gamma$ and (\ref{GestN}) to estimate the factors in the terms of the sum over $\gamma$. This way we obtain from (\ref{GmN}) that almost surely for $x,y\in \R^4$, $|x-y|\leq 2$, $\lambda =4|x-y|$
$$
\lambda^{|\m|}|(\G_{xy})^{\delta_\m}_\beta|\lesssim \lambda^{|\beta|_-}w_\beta(y)+\lambda^{|\beta|_-}w_\beta(x)+
$$
$$
+\sum_{\tiny\begin{array}{c}\gamma\in \M_{\geq 0}: \\ |\gamma|_\prec<|\beta|_\prec \\ |\gamma|<|\beta| \end{array}}|x-y|^{|\beta|_--|\gamma|_-}w_\beta(x) w_\gamma^{-1}(x)\lambda^{|\gamma|_-}w_\gamma(x)\lesssim |x-y|^{|\beta|_-}w_\beta(x),
$$
where at the last step we used (\ref{bwcd}) in the sum.
Finally, (\ref{GxyestnN}) is equivalent to the last estimate.

\epr

The last two propositions justify the induction step which also completes the proof of 
(\ref {GestMN}).

Note that by (\ref{gptwsind}) and (\ref{GmN}) used in the inductive algorithm, the constant in inequality (\ref {GestMN}) depends on the model via quantities (\ref{GestMlN}) and (\ref{EPPsin'>MlN}), the dependence being polynomial, without constant term.
This completes the proof of Theorem \ref{GestMNprop}.

\epr


\subsection{Pointwise estimates for recentered maps}\lab{recestptwsgl}

\setcounter{equation}{0}
\setcounter{theorem}{0}

Consider again a model of level $N\in \R$ at scale $\bar{\lambda}_0$ satisfying $\bar{\lambda}_{M(N)}\geq 8$ for the regularity structure with model space $T$, with recentered maps $\Pi_{x\beta}$, $\Pi_{x\beta}^-$ and structure group maps $(\G_{xy})_\beta^\gamma$.

We shall obtain pointwise estimates for $\Pi_{x\beta}$, $\Pi_{x\beta}^-$ which hold almost surely.
We keep the notation used in Theorem \ref{mainTes} and in Sections \ref{prestrec} and \ref{ptwsstrgrp}, and fix $\e_-\in \R$, $0<\e_-<\frac{\e}{100}$, with the help of which the corrected homogeneity is defined and denote $\alpha_-=\alpha-\e_-$.

\bt\lab{PestNprop}
Let $\Pi_{x\beta}$, $\Pi_{x\beta}^-$ be the recentered maps of a model of level $N\in \R$ at scale $\bar{\lambda}_0$ satisfying $\bar{\lambda}_{M(N)}\geq 8$, $w$ the weight used in Theorem \ref{GestMNprop}. Then for all $\mu \in \M$, $|\mu|_\prec<N$, $\varphi\in \mB^r$, $r\in \N$, $r>2-\alpha_-$, the following inequalities hold almost surely
\be\lab{EphiMN} 
|\Pi_{x\mu}^-(\varphi^{\lambda}_y)|\lesssim \lambda^{\alpha_--2}(\lambda+|x-y|)^{|\mu|_--\alpha_-}w_\mu(x),
\ee
\be\lab{EPPsin'>MN} 
|\Pi_{x\mu}(\varphi_y^\lambda)|\lesssim \lambda^{\alpha_-}(\lambda+|x-y|)^{|\mu|_--\alpha_-}w_\mu(x)
\ee
for all $0<\lambda<\infty$, $x,y\in \R^4$ if $\mu \not\in \M_{\geq 0}$, and if $\mu \in \M_{\geq 0}$ they hold for all  $0<\lambda\leq 1$, $x,y\in \R^4$, $|x-y|\leq 2$.

The constants in inequalities  (\ref{EphiMN}) and (\ref{EPPsin'>MN}) do not depend on $\lambda$, $\varphi$, $x$, $y$, and they depend on $w$ and on the model via quantities (\ref{GestMlN}), (\ref{EphiMlN}) and (\ref{EPPsin'>MlN}), the last dependence being polynomial without a constant term.
\et

\bpr
First note that for $\mu\in \M\setminus \M_{\geq 0}$ one has $[\mu]=-1$, and hence $|\mu|_-=|\mu|$ by (\ref{hom-}). Therefore arguments similar to those used in the proof of (\ref{Piunif}) and the deterministic bounds in Remark \ref{detbound} lead to the following deterministic bounds for $\mu\in \M\setminus \M_{\geq 0}$, $|\mu|_\prec<N$
\be\lab{EphiMlN1det} 
|\Pi_{x\mu}^-(\varphi^{\lambda}_y)|\lesssim \lambda^{\alpha_--2}(\lambda+|x-y|)^{|\mu|_--\alpha_-},
\ee
\be\lab{EPPsin'>MlN1det} 
|\Pi_{x\mu}(\varphi_y^\lambda)|\lesssim \lambda^{\alpha_-}(\lambda+|x-y|)^{|\mu|_--\alpha_-},
\ee
which hold for all $x,y\in \R^4$ and $0<\lambda<\infty$. By the definition of weights these bounds imply (\ref{EphiMN}) and (\ref{EPPsin'>MN}) for all $0<\lambda<\infty$, $x,y\in \R^4$ if $\mu \not\in \M_{\geq 0}$.

Therefore it suffices to establish (\ref{EphiMN}) and (\ref{EPPsin'>MN}) for $\mu\in \M_{\geq 0}$. We shall prove (\ref{EphiMN}). The other inequality is obtained in a similar way. 

Let $\psi=\omega^2\ast \omega$, where $\omega^2$ is $\omega$ rescaled by $2$, and $\omega$ is fixed in (\ref{omegapr*}) with $r\in \N$, $r>2-\alpha_-$.
Observe arguments similar to those used in the proof of (\ref{Piunif}) or Lemma \ref{Pi-int}, with the $\L_p$--norm replaced by $|\cdot |$, and estimates (\ref{GestMN}),
(\ref{EphiMlNloc}) with $\varphi=\psi$, $h=0$ imply that almost surely we have for $x,y\in \R^4$, $|x-y|\leq 2$, $0<\lambda\leq 8$
\be\lab{EphiMNps} 
|\Pi_{x\mu}^-(\psi^{\lambda}_y)|\lesssim \lambda^{\alpha_--2}(\lambda+|x-y|)^{|\mu|_--\alpha_-}w_\mu(x).
\ee
Indeed, by (\ref{Pi-GlN})
\be\lab{pi-estinter}
|\Pi_{x\mu}^-(\psi^{\lambda}_y)|\leq |\Pi_{y\mu}^-(\psi^{\lambda}_y)|+\sum_{\tiny\begin{array}{c}\delta\in \M: \\ |\delta|_\prec < |\mu|_\prec \\ |\delta|< |\mu|\end{array}}|\Pi_{y\delta}^-(\psi^{\lambda}_y)| |(\G_{yx})^\delta_\mu|.
\ee

Now we get almost surely using (\ref{GestMN}) and (\ref{EphiMlNloc}) with $\varphi=\psi$, $h=0$, $x,y\in \R^4$, $|x-y|\leq 2$, $0<\lambda\leq 8$, $w=w_\mu$,
\be\lab{Pi-locint}
|\Pi_{x\mu}^-(\psi^{\lambda}_y)|\lesssim \sum_{\tiny\begin{array}{c}\delta\in \M: \\ |\delta|_\prec \leq |\mu|_\prec \\ |\delta|\leq |\mu|\end{array}}\lambda^{|\delta|_--2}|x-y|^{|\mu|_--|\delta|_-}w_{\delta}(y)w_\mu(x)w_{\delta}^{-1}(y).
\ee

Finally, by (\ref{bwcd})
$$
w_{\delta}(y)w_\mu(x)w_{\delta}^{-1}(y)\lesssim w_{\delta}(x)w_\mu(x)w_{\delta}^{-1}(x)=w_\mu(x),
$$
and hence (\ref{Pi-locint}) takes the form
$$
|\Pi_{x\mu}^-(\psi^{\lambda}_y)|\lesssim \lambda^{\alpha_--2}\sum_{\tiny\begin{array}{c}\delta\in \M: \\ |\delta|_\prec \leq |\mu|_\prec \\ |\delta|\leq |\mu|\end{array}}\lambda^{|\delta|_--\alpha_-}|x-y|^{|\mu|_--|\delta|_-}w_\mu(x)\lesssim 
$$
$$
\lesssim\lambda^{\alpha_--2}(\lambda+|x-y|)^{|\mu|_--\alpha_-}w_\mu(x),
$$
where  at the last step we used (\ref{min}). This proves (\ref{EphiMNps}).

Now we proceed as in the proof of Proposition \ref{Pi-b<} replacing the $\L_p$--norms by $|\cdot |$ and repeating the other arguments verbatim.

Let $\e_n=2^{-(n+1)}$, $n\in \mathbb{Z}$.
We use presentation (\ref{P-appr}) with $0$ replaced by $x$, arbitrary $\varphi\in \mB^r$ and $\e_n=2^{-(n+1)}$, $n\in \N$,
$$
\Pi_{x\mu}^-(\varphi^{\e_{n-1}}_x)=
$$
$$
=\int_{\R^4}\Pi_{x\mu}^-(\psi^{\e_n}_z) (\varphi^{\e_{n-1}}\ast \psi^{\e_{n-1}})(x-z)dz+\sum_{m=n}^{\infty }\int_{\R^4}\Pi_{x\mu}^-(\psi^{\e_m}_z) (\varphi^{\e_{n-1}}\ast \check{\psi}^{\e_m})(x-z)dz.
$$

By (\ref{psisupp}) for $n\in \N$ one has ${\rm supp}(\varphi^{\e_{n-1}}\ast \psi^{\e_{n-1}})\subset B_{\frac{3}{2}\e_{n-1}}\subset B_2$, and by the first property in (\ref{cpsupp}) ${\rm supp}(\varphi^{\e_{n-1}}\ast \check{\psi}^{\e_m})\subset B_{\frac{3}{2}\e_{n-1}}\subset B_2$ for $m\geq n$.
By these observations we have from the previous presentation for $\Pi_{x\mu}^-(\varphi^{\e_{n-1}}_x)$, 
\be\lab{Pi0bn1}
|\Pi_{x\mu}^-(\varphi^{\e_{n-1}}_x)|\leq \int_{B_{\frac{3}{2}\e_{n-1}}(x)}|\Pi_{x\mu}^-(\psi^{\e_n}_z)|dz \left\|\varphi^{\e_{n-1}}\ast \psi^{\e_{n-1}}\right\|_{L_\infty}+
\ee
$$
+\sum_{m=n}^{\infty }\int_{B_{\frac{3}{2}\e_{n-1}}(x)}|\Pi_{x\mu}^-(\psi^{\e_m}_z)|dz \left\|\varphi^{\e_{n-1}}\ast \check{\psi}^{\e_m}\right\|_{L_\infty}.
$$

Using estimates (\ref{est1aux}), (\ref{est2aux})  and (\ref{EphiMNps}) in (\ref{Pi0bn1}) we obtain
\be\lab{Pib-psf1}
|\Pi_{x\mu}^-(\varphi^{\e_{n-1}}_x)|\leq 2^{nd}\int_{B_{\frac{3}{2}\e_{n-1}}(x)}\e_n^{\alpha_--2}(\e_n+|z-x|)^{|\mu|_--\alpha_-}w_\mu(x) dz +
\ee
$$
+\sum_{m=n}^{\infty }2^{n(r+d)-mr}\int_{B_{\frac{3}{2}\e_{n-1}}(x)}\e_m^{\alpha_--2}(\e_m+|z-x|)^{|\mu|_--\alpha_-}w_\mu(x) dz \lesssim
$$
$$
\lesssim (2^{nd-nd}\e_n^{\alpha_--2+|\mu|_--\alpha_-}+\e_n^{\alpha_--2+|\mu|_--\alpha_-}\sum_{m=n}^{\infty}2^{n(r+d)-mr-m(\alpha_--2)+n(\alpha_--2)-nd})w_\mu(x)\lesssim
$$
$$
\lesssim \e_{n-1}^{|\mu|_--2}\sum_{m=n}^{\infty}2^{-(m-n)(r+\alpha_--2)}w_\mu(x)\lesssim \e_{n-1}^{|\mu|_--2}w_\mu(x),
$$ 
where the series in the last expression converges due to the condition $r+\alpha_--2>0$ and 
the constant in the inequality does not depend on $n$, $\varphi$.

Now for any $0<\lambda \leq 1$, let $n\in \N$ be the largest possible such that $\e_{n-1}=2^{-n}\geq \lambda$, so that $\e_n\leq \lambda$, and for $\lambda':=\e_{n-1}^{-1}\lambda$ we have $\frac 12\leq \lambda'\leq 1$. 

Also if $\varphi\in \mB^r$ then $\varphi^{\lambda'}\in C\mB^r$, where $C>0$ is a constant independent of $\lambda'$ and $\varphi$.

Since by the definition $\varphi^\lambda=(\varphi^{\lambda'})^{\e_{n-1}}$ we obtain from the last observation and from (\ref{Pib-psf1}) that
$$
|\Pi_{x\mu}^-(\varphi_x^{\lambda})|\lesssim \e_{n-1}^{|\mu|_--2}w_\mu(x)\leq \lambda^{|\mu|_--2}w_\mu(x),
$$
where the constant in the inequality does not depend on $\lambda$ and $\varphi$.

Finally, the last estimate, which holds for all $\mu\in \M$, $|\mu|_\prec<N$, (\ref{GestMN}) and  arguments similar to those used in the proof of (\ref{EphiMNps}), with $\psi$ replaced by $\varphi$,  imply (\ref{EphiMN}). 

By (\ref{pi-estinter}), (\ref{Pi0bn1}), by Proposition \ref{Pipmptwsprel}, and by Theorem \ref{GestMNprop}  the constant in (\ref{EphiMN}) depends on the model via quantities (\ref{GestMlN}), (\ref{EphiMlN}) and (\ref{EPPsin'>MlN}), the dependence being polynomial without a constant term. This completes the proof.

\epr

Since by Theorem \ref{mainTes} for the smooth models defined in that theorem the quantities (\ref{GestMlN}), (\ref{EphiMlN}) and (\ref{EPPsin'>MlN}) introduced in part (13) of the definition of models of level $N$ are bounded uniformly in the mollification parameter $\rho$, and by Proposition \ref{Pipmptwsprel} and by Theorems \ref{GestMNprop}, \ref{PestNprop} the constants in pointwise estimates (\ref{EphiMN}), (\ref{EPPsin'>MN}) and (\ref{GestMN}) may depend on $\rho$ only via (\ref{GestMlN}), (\ref{EphiMlN}) and (\ref{EPPsin'>MlN}), we obtain the following corollary.

\bc\lab{ptwscor}
The recentered maps and the structure group maps of the smooth model of level $N$  at scale $\bar{\lambda}_0$ with $\bar{\lambda}_{M(N)}\geq 8$ for the regularity structure with model space $T$ defined in Theorem \ref{mainTes} satisfy  estimates (\ref{EphiMN}), (\ref{EPPsin'>MN}) and (\ref{GestMN}) almost surely with arbitrary weight $w$, and the constants in inequalities (\ref{EphiMN}), (\ref{EPPsin'>MN}) and (\ref{GestMN}) in this case do not depend on the mollification parameter $\rho$.
\ec

As a consequence of Theorem \ref{mainmod} and of the definition of models of level $N$ we also obtain another corollary.

\bc
The recentered maps and the structure group maps of the model of level $N$  at scale $\bar{\lambda}_0$ with $\bar{\lambda}_{M(N)}\geq 8$ for the regularity structure with model space $T$ associated to equation (\ref{LangdT}) and defined in Theorem \ref{mainmod} satisfy estimates (\ref{EphiMN}), (\ref{EPPsin'>MN}) and (\ref{GestMN}) almost surely with arbitrary weight $w$.
\ec

The last two corollaries, arguments similar to those presented before Corollary \ref{ptwscor} in conjunction with Theorems \ref{mainTconv}, \ref{mainmod} and with the usual telescoping trick as in the proof of Proposition \ref{indstepprpconv} yield the following statement which is a pointwise counterpart of Theorem \ref{mainmod}.
\bc
Let $\Pi_{x\mu}(\rho)$, $\Pi_{x\mu}^-(\rho)$, $\G_{xy}(\rho)$ be the recentered maps and the structure group maps of the smooth model of level $N$  at scale $\bar{\lambda}_0$ satisfying $\bar{\lambda}_{M(N)}\geq 8$ for the regularity structure with model space $T$ defined in Theorem \ref{mainTes}, and $\Pi_{x\mu}$, $\Pi_{x\mu}^-$, $\G_{xy}$ the recentered maps and the structure group maps of the model of level $N$  at the scale $\bar{\lambda}_0$ for the regularity structure with model space $T$ associated to equation (\ref{LangdT}) and defined in Theorem \ref{mainmod}. Then for any weight $w$ and $r>2-\alpha_-$ we have almost surely
$$
\max_{\tiny\begin{array}{c} \mu,\gamma \in \M: \\ |\mu|_\prec,|\gamma|_\prec<N\end{array}}\sup_{\tiny\begin{array}{c} x,y\in \R^4: \\ 0<|x-y|\leq 2\end{array}}\frac{|(\G_{xy}-\G_{xy}(\rho))^\gamma_\mu|}{|x-y|^{|\mu|_--|\gamma|_-}w_\mu(x) w_\gamma^{-1}(x)}\to 0 ~{\rm when}~\rho\to 0,
$$
$$
\max_{\tiny\begin{array}{c} \mu \in \M: \\ |\mu|_\prec<N\end{array}}\sup_{x\in \R^4}\sup_{0<\lambda\leq 1}\sup_{\varphi\in \mB^r}
\frac{|\Pi_{x\mu}^-(\varphi^{\lambda}_x)-\Pi_{x\mu}^-(\rho)(\varphi^{\lambda}_x)|}{\lambda^{|\mu|_--2}w_\mu(x)}\to 0 ~{\rm when}~\rho\to 0,
$$
$$
\max_{\tiny\begin{array}{c} \mu \in \M: \\ |\mu|_\prec<N\end{array}}\sup_{x\in \R^4}\sup_{0<\lambda\leq 1}\sup_{\varphi\in \mB^r}
\frac{|\Pi_{x\mu}(\varphi^{\lambda}_x)-\Pi_{x\mu}(\rho)(\varphi^{\lambda}_x)|}{\lambda^{|\mu|_-}w_\mu(x)}\to 0 ~{\rm when}~\rho\to 0.
$$
\ec

\end{document}